\documentclass[12pt, eqno]{article}
\usepackage{bbm}
\usepackage{mathrsfs}
\usepackage{amsfonts}
\usepackage{amssymb}
\usepackage{graphicx}
\usepackage[all]{xy}

\usepackage{psfrag}
\usepackage{subfigure}
\usepackage{color}

\usepackage{amssymb,latexsym}
\usepackage{amsmath,latexsym}
\usepackage{amscd}

\newcommand{\K}{{\mathbb K}}
\newcommand{\F}{{\mathbb F}}
\newcommand{\N}{{\mathbb N}}
\newcommand{\C}{{\mathbb C}}
\newcommand{\R}{{\mathbb R}}
\newcommand{\Q}{{\mathbb Q}}
\newcommand{\Z}{{\mathbb Z}}

\newtheorem{theorem}{Theorem}[section]

\newtheorem{corollary}[theorem]{Corollary}

\newtheorem{remark}[theorem]{Remark}

\newtheorem{lemma}[theorem]{Lemma}

\newtheorem{proposition}[theorem]{Proposition}
\newtheorem{claim}[theorem]{Claim}

\begin{document}

\title{Methods of infinite dimensional Morse theory\\ for geodesics on Finsler manifolds}

\author{\\Guangcun Lu
\thanks{\endgraf
 }\\
{\normalsize School of Mathematical Sciences, Beijing Normal University},\\
{\normalsize Laboratory of Mathematics
 and Complex Systems,  Ministry of
  Education},\\
  {\normalsize Beijing 100875, The People's Republic
 of China}\\
{\normalsize (gclu@bnu.edu.cn)}}
\date{October 1,  2014} \maketitle \vspace{-0.1in}


\abstract{We prove the shifting theorems of the critical groups of
critical points and critical orbits for the energy functionals of
Finsler metrics on Hilbert manifolds of $H^1$-curves,  and two
splitting lemmas for the functionals on Banach manifolds of
$C^1$-curves. Two results on critical groups of iterated closed
geodesics are also proved; their corresponding versions on
Riemannian manifolds are based on the usual splitting lemma
by Gromoll and Meyer (1969). Our approach consists in deforming
the square of the Finsler metric in a Lagrangian
which is smooth also on the zero section  and then in using the splitting
lemma for nonsmooth functionals that the author recently developed
in Lu (2011, 0000, 2013). The argument does not involve
 finite-dimensional approximations and any  Palais'
result in Palais (1966). As an application, we extend to Finsler manifolds a result by V.
Bangert and W. Klingenberg (1983) about the existence of infinitely many, geometrically distinct,
 closed geodesics on a compact Riemannian manifold. } \vspace{-0.1in}
\medskip\vspace{6mm}

 \vspace{2mm}

\tableofcontents

\section{Introduction and results}
\setcounter{equation}{0}

Let $M$ be a smooth manifold of dimension $n$. For an integer $k\ge
5$ or $k=\infty$, a $C^k$ {\bf Finsler metric} on $M$  is a
continuous function $F: TM\to\R$ satisfying the following
properties:
\begin{description}
\item[(i)]  It is $C^k$ and positive in $TM\setminus\{0\}$.

\item[(ii)] $F(x,ty)=tF(x,y)$ for every $t>0$ and $(x,y)\in TM$.

\item[(iii)] $L:=F^2$ is fiberwise strongly convex, that is, for any $(x,y)\in TM\setminus\{0\}$
the symmetric bilinear form (the fiberwise Hessian operator)
$$
g^F(x,y): T_xM\times T_xM\to\R,\,(u,
v)\mapsto\frac{1}{2}\frac{\partial^2}{\partial s\partial
t}\left[L(x, y+ su+ tv)\right]\Bigl|_{s=t=0}
$$
is positive definite.
\end{description}
Because of the positive homogeneity it is easily checked that
$g^F(x, \lambda y)=g^F(x,y)$ for any $(x,y)\in TM\setminus\{0\}$ and
$\lambda>0$. The Euler theorem implies
$g^F(x,y)[y,y]=(F(x,y))^2=L(x,y)$. A smooth manifold $M$ endowed
with a $C^k$ Finsler metric is called a $C^k$ {\bf Finsler
manifold}. Note that $L=F^2$ is of class $C^{2-0}$, and of class
$C^2$ if and only if it is  the square of the norm of a Riemannian metric. A Finsler
metric $F$ is said to be {\bf reversible} if $F(x,-v)=F(x,v)$ for
all $(x,v)\in TM$. We also say that a Finsler metric dominates a
Riemannian metric $g$ on $M$ if $F(x,v)\ge C_0\sqrt{g_x(v,v)}$ for
some $C_0>0$ and all $(x,v)\in TM$, (clearly, every Finsler metric
on a compact manifold dominates a Riemannian metric). The length of
a Lipschitz continuous curve $\gamma:[a,b]\to M$ with respect to the
Finsler structure $F$ is defined by $l_F(\gamma)=\int^b_aF(\gamma(t),\dot
\gamma(t))dt$. However  a non-reversible $F$  only induces a
non-symmetric distance and hence leads to the notions of forward and
backward completeness (and geodesic completeness). A differentiable
$\gamma=\gamma(t)$ is said to have {\bf constant speed} if
$F(\gamma(t), \dot\gamma(t))$ is constant along $\gamma$.
According to \cite{BaChSh}, a regular piecewise $C^k$ curve $\gamma:[a, b]\to M$ is
called a Finslerian {\bf geodesic} if it minimizes the length between two sufficiently
close points on the curve.  In terms of
Chern connection $D$ a constant speed Finslerian geodesic is characterized by
the condition $D_{\dot\gamma}\dot\gamma=0$, where $D_{\dot\gamma}{\dot\gamma}$
is the covariant derivative with reference vector  $\dot\gamma$
(see pages 121-124 on \cite{BaChSh}).

Let $(M,F)$ be a $n$-dimensional $C^k$ Finsler manifold with a
complete Riemannian metric $g$, and let
$N\subset M\times M$ be a smooth submanifold. Denote by $I=[0,1]$
the unit interval, and by $W^{1,2}(I, M)$ the space of absolutely
continuous curves $\gamma: I\to M$ such that
$\int^1_0\langle\dot\gamma(t),\dot\gamma(t)\rangle dt<\infty$, where
$\langle u, v\rangle =g_x(u,v)$ for $u,v\in T_xM$. Then $W^{1,2}(I,
M)\subset C^0(I,M)$. Set
$$
\Lambda_N(M):=\{\gamma\in W^{1,2}(I, M)\,|\,(\gamma(0),\gamma(1))\in
N\}.
$$
 It is a Riemannian-Hilbert
 submanifold of $W^{1,2}(I, M)$ of codimension ${\rm
codim}(N)$. If $\gamma\in C^2(I, M)\cap\Lambda_N(M)$ then the
pull-back bundle $\gamma^\ast TM\to I$ is a $C^1$ vector bundle. Let
$W^{1,2}_N(\gamma^\ast TM)$ be the set of absolutely continuous
sections $\xi: I\to \gamma^\ast TM$ such that
$$
(\xi(0),\xi(1))\in T_{(\gamma(0),\gamma(1))}N\quad\hbox{and}\quad
\int^1_0\langle\nabla^g_{\dot{\gamma}}\xi(t),\nabla^g_{\dot{\gamma}}\xi(t)\rangle
dt<\infty.
$$
Here $\nabla^g_\gamma$ is the covariant derivative along $\gamma$
with respect to the Levi-Civita connection of the metric $g$. Then
$T_\gamma\Lambda_N(M)=W^{1,2}_N(\gamma^\ast TM)$ is equipped with
the inner product
\begin{equation}\label{e:1.1}
\langle\xi,\eta\rangle_1=\int^1_0\langle\xi(t),\eta(t)\rangle dt+
\int^1_0\langle\nabla^g_{\dot{\gamma}}\xi(t),\nabla^g_{\dot{\gamma}}\xi(t)\rangle
dt
\end{equation}
and norm $\|\cdot\|_1=\sqrt{\langle\cdot,\cdot\rangle_1}$. For a
general $\gamma\in\Lambda_N(M)\setminus C^2(I, M)$, note that
$\gamma^\ast TM\to I$ is only a bundle of class $H^1$, that is, it
admits a system of local trivializations whose transition maps are
of class $W^{1,2}$. Fortunately, since $W^{1,2}(I,\R^n)$ is a Banach
algebra, one can still define $W^{1,2}$-section of $\gamma^\ast TM$,
and the set $W^{1,2}(\gamma^\ast TM)$ of such sections is also a
well-defined Hilbert space with the inner product given by
(\ref{e:1.1}) (using the $L^2$ covariant derivative along $\gamma$
associated to the Levi-Civita connection of the metric $g$). Hence
we have also $T_\gamma\Lambda_N(M)=W^{1,2}_N(\gamma^\ast TM)$ in
this case. For more details
about these, see Appendix C of
\cite{McSa}. Consider the energy functional on $(M, F)$,
\begin{equation}\label{e:1.2}
{\cal L}:\Lambda_N(M)\to\R,\;\gamma\mapsto\int^1_0L(\gamma(t),
\dot\gamma(t))dt=\int^1_0[F(\gamma(t),\dot\gamma(t))]^2dt.
\end{equation}

\begin{proposition}\label{prop:1.1}
\begin{description}
\item[(i)] The functional ${\cal L}$ is $C^{2-}$, i.e. ${\cal L}$ is
$C^1$ and its differential $d{\cal L}$ is locally Lipschitz.

\item[(ii)] A curve $\gamma\in\Lambda_N(M)$ is a constant (non-zero)
speed $F$-geodesic satisfying the boundary condition
$$
g^F(\gamma(0),\dot\gamma(0))[u,\dot\gamma(0)]=
g^F(\gamma(1),\dot\gamma(1))[v,\dot\gamma(1)]\quad\forall (u,v)\in
T_{(\gamma(0),\gamma(1))}N
$$
if and only if it is a (nonconstant) critical point of ${\cal L}$.

\item[(iii)] Suppose that $(M,F)$ is
forward (resp. backward) complete and that $N$ is  a closed
submanifold of $M\times M$ such that the first projection (resp. the
second projection) of $N$ to $M$ is compact. Then ${\cal L}$
satisfies the Palais-Smale condition on $\Lambda_N(M)$.
\end{description}
\end{proposition}
When $M$ is compact and $N=\triangle_M$ (diagonal) this proposition
was proved in \cite{Me} by Mercuri.   Kozma, Krist\'aly and Varga
\cite{KoKrVa} proved (i), (ii) if $N$ is a product of two
submanifolds $M_1$ and $M_2$ in $M$, and (iii) if $F$  dominates a
complete Riemannian metric on $M$ and $N\subset M\times M$ is a
closed submanifold of $M\times M$ such that $P_1(N)\subset M$ or
$P_2(N)\subset M$ is compact. The above versions were obtained by
Caponio, Javaloyes and Masiello \cite{CaJaMa3} recently.
 A curve $\gamma\in\Lambda_N(M)$ is
called {\bf regular} if $\dot\gamma\ne 0$ a.e. in $[0,1]$. When
$N=\{p,q\}$ for some $p,q\in M$,  Caponio  proved in
\cite[Prop.3.2]{Ca}: if ${\cal L}$ is twice differentiable at a
regular curve $\gamma\in\Lambda_N(M)$ then for a.e. $s\in [0,1]$,
the function
$$
v\in T_{\gamma(s)}M\to F^2(\gamma(s),v)
$$
is a quadratic positive definite form. This suggests that the
results in \cite{MePa, DeSo} cannot be used for the functional
${\cal L}$ in (\ref{e:1.2}).
Actually, the classical method, finite
dimensional approximations developed by
Morse \cite{Mo} for Riemannian geodesics, has been also used for geodesics
of Finsler metrics \cite{BLo, Ma, Ra, Shen}.
 For example,   the
shifting theorem of Gromoll-Meyer \cite{GM1} was obtained using
the finite dimensional approximations setting
in the generalization of a famous theorem of Gromoll-Meyer
\cite{GM2} to Finsler manifolds in \cite{Ma}. However, the finite
dimensional approximation of the loop space $\Lambda
M=\Lambda_{\triangle_M}(M)$ by spaces of broken geodesics carries
only a $\Z_k$ but not an $S^1$ action. It is hard to apply this
classical method to the studies of some geodesic problems.
\vspace{2mm}

\noindent{\bf Notation.}  For a normed vector space $(E, \|\cdot\|)$
and $\delta>0$ we  write ${\bf B}_\delta(E)=\{x\in
E\,|\,\|x\|<\delta\}$ and $\bar{\bf B}_\delta(E)=\{x\in
E\,|\,\|x\|\le\delta\}$ (since the notations of some spaces involved
are complex). For a continuous symmetric bilinear form (or the
associated self-adjoint operator) $B$ on a Hilbert space we write
${\bf H}^-(B)$, ${\bf H}^0(B)$ and ${\bf H}^+(B)$ for the negative
definite, null and positive definite spaces of it.   $\K$ always
denotes an Abelian group.  \vspace{2mm}

In this paper we only consider the following two cases:

\noindent{$\bullet$} $N=M_0\times M_1$,
where $M_0$ and $M_1$ are
two boundaryless and connected submanifolds of $M$ and also  disjoint if $\dim M_0\dim M_1>0$.
 In this case the
boundary condition in Proposition~\ref{prop:1.1}(ii) becomes
\begin{equation}\label{e:1.3}
\left\{\begin{array}{ll}
&g^F(\gamma(0),\dot\gamma(0))[u,\dot\gamma(0)]=0\quad\forall u\in
T_{\gamma(0)}M_0,\\
&g^F(\gamma(1),\dot\gamma(1))[v,\dot\gamma(1)]=0\quad\forall v\in
T_{\gamma(1)}M_1. \end{array}\right.
\end{equation}

\noindent{$\bullet$} $N=\triangle_M$ and hence
$\Lambda_N(M)=W^{1,2}(S^1,M)$, where $S^1=\R/\Z=[0,1]/\{0,1\}$.

\subsection{The case $N=M_0\times M_1$}\label{sec:1.1}

 Suppose that
$\gamma_0\in\Lambda_N(M)$ is an isolated nonconstant critical point
of ${\cal L}$ on $\Lambda_N(M)$.  By Proposition~\ref{prop:1.1}(ii)
$\gamma_0$ is a $C^k$-smooth nonconstant $F$-geodesics with constant
speed $F(\gamma_0(t),\dot\gamma_0(t))\equiv \sqrt{c}>0$. Clearly,
there exist an open neighborhood of $\gamma_0$ in $\Lambda_N(M)$,
${\cal O}(\gamma_0)$, and a compact neighborhood $K$ of
$\gamma_0(I)$ in $M$ such that each curve of ${\cal O}(\gamma_0)$
has an image set contained in $K$. Thus \textsf{we shall assume $M$
to be compact below.}  Moreover, we can  require the Riemannian
metric $g$ on $M$ to have the property that $M_0$ (resp. $M_1$) is
totally geodesic near $\gamma_0(0)$ (resp. $\gamma_0(1)$). Let
$\exp$ denote the exponential map of $g$, and let ${\bf
B}_{2\rho}(T_{\gamma_0}\Lambda_N(M))= \{\xi\in
T_{\gamma_0}\Lambda_N(M)\,|\,\|\xi\|_{1}<2\rho\}$ for $\rho>0$. If
$\rho$ is small enough, by the so-called omega lemma we get a $C^{k-3}$ coordinate chart around $\gamma_0$
on $\Lambda_N(M)$:
\begin{equation}\label{e:1.4}
{\rm EXP}_{\gamma_0}: {\bf
B}_{2\rho}(T_{\gamma_0}\Lambda_N(M))\to\Lambda_N(M)
\end{equation}
given by ${\rm EXP}_{\gamma_0}(\xi)(t)=\exp_{\gamma_0(t)}(\xi(t))$.
  Then ${\cal L}\circ{\rm EXP}_{\gamma_0}$ is
$C^{2-0}$ and has an isolated critical point $0\in
W^{1,2}_N(\gamma_0^\ast TM)$. Consider the Banach manifold
$${\cal
X}=C^1_N(I, M)=\{\gamma\in C^1(I,M)\,|\, (\gamma(0), \gamma(1))\in
N\}.
$$
Then the tangent space $T_{\gamma_0}{\cal X}=C^1_{TN}(\gamma_0^\ast
TM)=\{\xi\in C^1(\gamma_0^\ast TM)\,|\, (\xi(0), \xi(1))\in {TN} \}$
with usual $C^1$-norm. Let ${\cal A}$ be the restriction of the
gradient of ${\cal L}\circ{\rm EXP}_{\gamma_0}$ to ${\bf
B}_{2\rho}(T_{\gamma_0}\Lambda_N(M))\cap T_{\gamma_0}{\cal X}$.
Observe that ${\bf B}_{2\rho}(T_{\gamma_0}\Lambda_N(M))\cap
T_{\gamma_0}{\cal X}$ is an open neighborhood of $0$ in
$T_{\gamma_0}{\cal X}$ and that ${\cal A}$ is a $C^{k-3}$-map from a
small neighborhood of $0\in T_{\gamma_0}{\cal X}$ to
$T_{\gamma_0}{\cal X}$ by (\ref{e:4.28}) and the $C^{k-3}$ smoothness
of $\tilde{A}$ on ${\bf B}_{2r}(X_V)$ claimed below (4.5).
Moreover, ${\cal A}(0)=\nabla{\cal L}(\gamma_0)|_{T_{\gamma_0}{\cal
X}}$ and
$$
\langle d{\cal A}(0)[\xi],\eta\rangle_{1}=d^2{\cal L}|_{\cal
X}(\gamma_0)[\xi,\eta]\quad\forall\xi,\eta\in T_{\gamma_0}{\cal X}.
$$
The key is that the symmetric bilinear form $d^2{\cal L}|_{\cal
X}(\gamma_0)$ can be extended to a symmetric bilinear form on
$T_{\gamma_0}\Lambda_N(M)$, also denoted by $d^2{\cal L}|_{\cal
X}(\gamma_0)$.  The self-adjoint operator associated to the latter
is Fredholm, and has finite dimensional negative definite and null
spaces ${\bf H}^-(d^2{\cal L}|_{\cal X}(\gamma_0))$ and ${\bf
H}^0(d^2{\cal L}|_{\cal X}(\gamma_0))$ by (\ref{e:4.29}), which are
actually contained in $T_{\gamma_0}{\cal X}$.  There exists an
orthogonal decomposition
\begin{equation}\label{e:1.5}
T_{\gamma_0}\Lambda_N(M)={\bf H}^-(d^2{\cal L}|_{\cal
X}(\gamma_0))\oplus {\bf H}^0(d^2{\cal L}|_{\cal
X}(\gamma_0))\oplus{\bf H}^+(d^2{\cal L}|_{\cal X}(\gamma_0)),
\end{equation}
which  induces a (topological) direct sum decomposition of Banach
spaces
$$
T_{\gamma_0}{\cal X}={\bf H}^-(d^2{\cal L}|_{\cal
X}(\gamma_0))\dot{+} {\bf H}^0(d^2{\cal L}|_{\cal
X}(\gamma_0))\dot{+}\bigr({\bf H}^+(d^2{\cal L}|_{\cal
X}(\gamma_0))\cap T_{\gamma_0}{\cal X}\bigl).
$$
Using the implicit function theorem we get $\delta\in (0, 2\rho]$
and a unique $C^{k-3}$-map
$$
 h: {\bf B}_\delta\bigl({\bf H}^0(d^2{\cal L}|_{\cal
X}(\gamma_0))\bigr)\to {\bf H}^-(d^2{\cal L}|_{\cal
X}(\gamma_0))\dot{+}\bigl({\bf H}^+(d^2{\cal L}|_{\cal
X}(\gamma_0))\cap T_{\gamma_0}{\cal X}\bigr)
$$
such that $h(0)=0$, $dh(0)=0$ and $(I-P^0){\cal A}(\xi+
h(\xi))=0\;\forall\xi\in {\bf B}_\delta\bigl({\bf H}^0(d^2{\cal
L}|_{\cal X}(\gamma_0))\bigr)$, where
$P^\star:T_{\gamma_0}\Lambda_N(M)\to {\bf H}^\star(d^2{\cal
L}|_{\cal X}(\gamma_0))$, $\star=-,0,+$, are the orthogonal
projections according to the decomposition (\ref{e:1.5}). Define
${\cal L}^\circ:{\bf B}_\delta\bigl({\bf H}^0(d^2{\cal L}|_{\cal
X}(\gamma_0))\bigr)\to\R$ by
\begin{equation}\label{e:1.6}
{\cal L}^\circ(\xi)={\cal L}\circ{\rm EXP}_{\gamma_0}\bigl(\xi+
h(\xi)\bigr).
\end{equation}
 It is $C^{k-3}$, has an isolated
critical point $0$, and $d^2{\cal L}^\circ(0)=0$.
 Denote by $C_\ast({\cal L},
\gamma_0;\K)$ (resp. $C_{\ast}({\cal L}^\circ, 0;\K)$)  the critical
group of the functional ${\cal L}$ (resp. ${\cal L}^\circ$) at
$\gamma_0$ (resp. $0$) with coefficient group $\K$.
Note that
$$
m^-(\gamma_0):=\dim {\bf H}^-(d^2{\cal L}|_{\cal
X}(\gamma_0))\quad\hbox{and}\quad m^0(\gamma_0):=\dim {\bf
H}^0(d^2{\cal L}|_{\cal X}(\gamma_0))
$$
do not depend on the choice of the Riemannian metric $g$, called
{\bf Morse index} and {\bf nullity} of $\gamma_0$, respectively.
Here is our first key result.

\begin{theorem}\label{th:1.2}
Under the  above notations, there exists a $C^k$ convex with quadratic
growth Lagrangian $L^\ast:TM\to \R$ such that:
\begin{description}
\item[(i)] $L^\ast\le  L$, $L^\ast(x,v)=L(x,v)$ if $L(x,v)\ge 2c/3$,
 and if $F$ is reversible so is $L^\ast$.

\item[(ii)] The corresponding functional ${\cal L}^\ast$ in
(\ref{e:1.26}) is $C^{2-0}$. All functional ${\cal
L}^\tau=(1-\tau){\cal L}+ \tau{\cal L}^\ast$, $\tau\in [0,1]$, have
only a critical point $\gamma_0$ in some neighborhood of
$\gamma_0\in\Lambda_N(M)$, and  satisfy the Palais-Smale condition.
Moreover, $C_\ast({\cal L}^\ast,\gamma_0;\K)=C_\ast({\cal
L},\gamma_0;\K)$.

\item[(iii)] By shrinking the above $\delta>0$
 there exists  an origin-preserving
homeomorphism $\psi$ from ${\bf B}_\delta(T_{\gamma_0}\Lambda_N(M))$
to an open neighborhood of $0$ in $T_{\gamma_0}\Lambda_N(M)$
 such that
$$
{\cal L}^\ast\circ{\rm EXP}_{\gamma_0}\circ\psi(\xi)=\|P^+\xi\|_1^2-
\|P^-\xi\|^2_1 + {\cal L}^\circ(P^0\xi)
$$
for all $\xi\in{\bf B}_\delta(T_{\gamma_0}\Lambda_N(M))$. Moreover,
$\psi\left((P^-+P^0){\bf B}_\delta(T_{\gamma_0}\Lambda_N(M))\right)$ is contained
in $T_{\gamma_0}{\cal X}$, and $\psi$ is also a homeomorphism from
$(P^-+P^0){\bf B}_\delta(T_{\gamma_0}\Lambda_N(M))$ onto
$\psi\left((P^-+P^0){\bf B}_\delta(T_{\gamma_0}\Lambda_N(M))\right)$
even if the topology on the latter is taken as the induced one by $T_{\gamma_0}{\cal X}$.

\item[(iv)]  For any open neighborhood $W$ of $\gamma_0$ in $\Lambda_N(M)$ and a field $\K$,
 write $W_X:=W\cap{\cal X}$ as an open
 subset of ${\cal X}$, then the inclusion
$$
\left(({\cal L}^\ast|_{\cal X})_c\cap W_X, ({\cal L}^\ast|_{\cal
X})_c\cap W_X\setminus\{0\}\right)\hookrightarrow\left({\cal
L}^\ast_c\cap W, {\cal L}^\ast_c\cap W\setminus\{0\}\right)
$$
induces isomorphisms
$$
H_\ast\left(({\cal L}^\ast|_{\cal X})_c\cap W_X, ({\cal
L}^\ast|_{\cal X})_c\cap W_X\setminus\{0\};\K\right)\to
H_\ast\left({\cal L}^\ast_c\cap W, {\cal L}^\ast_c\cap
W\setminus\{0\};\K\right).
$$
\end{description}
\end{theorem}

For ${\cal L}|_{\cal X}$ we have the splitting lemma at $\gamma_0$
as follows.

\begin{theorem}\label{th:1.3}
Under the  above notations,
 there exist  $\epsilon\in (0, \delta)$ and an origin-preserving  homeomorphism $\varphi$ from
${\bf B}_\epsilon(T_{\gamma_0}{\cal X})$ to an open neighborhood of
$0$ in $T_{\gamma_0}{\cal X}$
 such that for any $\xi\in {\bf B}_\epsilon(T_{\gamma_0}{\cal X})$,
$$
{\cal L}\circ{\rm
EXP}_{\gamma_0}\circ\varphi(\xi)=\frac{1}{2}d^2{\cal L}|_{\cal
X}(\gamma_0)[P^+\xi, P^+\xi]-\|P^-\xi\|^2_1 + {\cal
L}^\circ(P^0\xi).
$$
\end{theorem}

  Corollaries~\ref{cor:A.2},~\ref{cor:A.6},  Theorem~\ref{th:1.2}(iii) and
Theorem~\ref{th:1.3} give  the following two versions of the shifting
theorem, respectively.

\begin{theorem}\label{th:1.4}
 For all $q=0,1,\cdots$ it holds that
\begin{eqnarray}\label{e:1.7}
&&C_q({\cal L}, \gamma_0;\K)=C_{q-m^-(\gamma_0)}({\cal L}^\circ,
0;\K)\quad\hbox{and}\nonumber\\
&&C_q({\cal L}|_{\cal X}, \gamma_0;\K)=C_{q-m^-(\gamma_0)}({\cal
L}^\circ, 0;\K).
\end{eqnarray}
\end{theorem}

As a consequence of this result  we have
\begin{equation}\label{e:1.8}
C_\ast({\cal L}|_{\cal X}, \gamma_0;\K)= C_\ast({\cal L},
\gamma_0;\K).
\end{equation}
 Note that Theorem~\ref{th:1.2}(iv) is
stronger than this. Indeed,  if we choose a small open neighborhood
$V$ of $\gamma_0$ in ${\cal X}$ such that the closure of $V$ in
$\Lambda_N(M)$ (resp. ${\cal X}$) is contained in $W$ (resp. the
open subset $W_X$ of ${\cal X}$) and that
$\min_tL(\gamma(t),\dot\gamma(t))\ge 2c/3$ (and hence ${\cal
L}^\ast(\gamma)={\cal L}(\gamma)$) for any $\gamma\in V$, then
$$
H_\ast\left(({\cal L}^\ast|_{\cal X})_c\cap W_X, ({\cal
L}^\ast|_{\cal X})_c\cap W_X\setminus\{0\};\K\right)=
H_\ast\left(({\cal L}|_{\cal X})_c\cap V, ({\cal L}|_{\cal X})_c\cap
V\setminus\{0\};\K\right).
$$
by the excision property of the relative homology groups. This and
Theorem~\ref{th:1.2}(ii),(iv) imply (\ref{e:1.8}). In applications
it is more effective combining these results together.

When $M_0$ and $M_1$ are two  points, an equivalent version
of Theorem~\ref{th:1.3} (as Theorem~\ref{th:4.5}) was also proved by
Caponio-Javaloyes-Masiello \cite[Th.7]{CaJaMa1} in the spirit of
\cite{GM1, Ch93}.  However they only obtained (\ref{e:1.8}) under
the assumption that $\gamma_0$ is a \textsf{nondegenerate} critical
point. See Remark~\ref{rm:4.6} for a detailed discussion. The
infinite dimensional proof method of Theorem~\ref{th:1.3} with
\cite{Ji} was  presented by the author \cite{Lu} (in the
periodic orbit case as in Theorem~\ref{th:1.6}).

\subsection{The case $N=\triangle_M$}\label{sec:1.2}

 Now $\Lambda_N(M)=\Lambda M:=W^{1,2}(S^1,M)=\{\gamma\in W^{1,2}_{loc}(\R,M)\,|\,
 \gamma(t+1)=\gamma(t)\;\forall t\in\R\}$. Hereafter
$S^1:=\R/\Z=\{[s]\,|\,[s]=s+\Z,\, s\in\R\}$.
There exists an equivariant
and also isometric action of $S^1$ on $W^{1,2}(S^1, M)$
and $TW^{1,2}(S^1, M)$:
\begin{eqnarray}
&&[s]\cdot \gamma(t)=\gamma(s+t), \quad \forall [s]\in S^1,
         \; \gamma\in W^{1,2}(S^1, M),\label{e:1.9}\\
 &&[s]\cdot \xi(t)=\xi(s+t),
\quad \forall [s]\in S^1,
         \; \xi\in T_{\gamma}W^{1,2}(S^1, M),\label{e:1.10}
  \end{eqnarray}
which are  continuous, but not differentiable (\cite[Chp.2,
\S2.2]{Kl}). Since ${\cal L}$ is invariant under the $S^1$-action, a
nonconstant curve $\gamma\in \Lambda M$ cannot be an isolated
critical point of the functional ${\cal L}$. Let $\gamma_0\in\Lambda
M$ be a (nonconstant) critical point of ${\cal L}$ with critical
value $c>0$. Under our assumptions $\gamma_0$ is $C^k$-smooth by
Proposition~\ref{prop:3.1}.
  The orbit $S^1\cdot \gamma_0$ is a
$C^{k-1}$-submanifold in $\Lambda M=W^{1,2}(S^1, M)$ by \cite[page
499]{GM2}, and
 hence a $C^{k-1}$-smooth critical submanifold of ${\cal L}$  in
$\Lambda M$.
  \textsf{We assume that $S^1\cdot
\gamma_0$ is an isolated critical orbit below}.

In this subsection, $\mathcal{X}$ denotes the manifold ${\cal X}=C^1(S^1, M)$.
 Let ${\cal
O}=S^1\cdot\gamma_0$ and let $\pi:N{\cal O}\to{\cal O}$ be the
normal bundle of it in $\Lambda M$. It is a $C^{k-2}$ Hilbert vector
bundle over ${\cal O}$, and $XN{\cal O}:=T_{\cal O}{\cal X}\cap
N{\cal O}$ is a $C^{k-2}$ Banach vector bundle over ${\cal O}$. For
$\varepsilon>0$  we denote by
\begin{equation}\label{e:1.11}
\left.\begin{array}{ll}
 & N{\cal O}(\varepsilon)=\{(x,v)\in N{\cal
 O}\,|\,\|v\|_{1}<\varepsilon\},\,\\
& XN{\cal O}(\varepsilon)= \{(x,v)\in XN{\cal
O}\,|\,\|v\|_{C^1}<\varepsilon\}.
\end{array}\right\}
\end{equation}
Clearly, $XN{\cal O}(\varepsilon)\subset N{\cal O}(\varepsilon)$.
Since ${\cal O}$ has a neighborhood in which all elements have supports contained
in a compact neighborhood of the support of $\gamma_0$ in $M$
we would assume directly that $M$ is compact. Then
 for sufficiently small $\varepsilon>0$, the map
\begin{equation}\label{e:1.12}
{\rm EXP}:T\Lambda M(\varepsilon)=\{(x,v)\in T\Lambda
M\,|\,\|v\|_{1}<\varepsilon\}\to \Lambda M
\end{equation}
defined by ${\rm EXP}(x,v)(t)=\exp_{x(t)}v(t)\;\forall t\in\R$ via  the exponential
map $\exp$,  restricts to an $S^1$-equivariant $C^{k-3}$ diffeomorphism  from the
normal disk bundle $N{\cal O}(\varepsilon)$ onto an $S^1$-invariant
open neighborhood of ${\cal O}$ in $\Lambda M$,
\begin{equation}\label{e:1.13}
\digamma:N{\cal O}(\varepsilon)\to{\cal N}({\cal O},\varepsilon).
\end{equation}
 For conveniences we write
\begin{equation}\label{e:1.14}
{\cal F}:={\cal L}\circ\digamma\quad\hbox{and}\quad {\cal
F}^\ast:={\cal L}^\ast\circ\digamma
\end{equation}
 as functionals on $N{\cal O}(\varepsilon)$. They are $C^{2-0}$,
$S^1$-invariant and satisfy the (PS) condition. Denote by ${\cal
F}^X$ and ${\cal F}^{\ast X}$ the restrictions of ${\cal F}$ and
${\cal F}^\ast$ to the open subset $N{\cal O}(\varepsilon)\cap
XN{\cal O}$ of $XN{\cal O}$, respectively.
 For $x\in{\cal O}$ let $N{\cal
O}(\varepsilon)_x$ and $XN{\cal O}(\varepsilon)_x$ be the fibers of
$N{\cal O}(\varepsilon)$ and $XN{\cal O}(\varepsilon)$ at $x\in{\cal
O}$, respectively. Also write ${\cal F}_x$, ${\cal F}^\ast_x$ and
${\cal F}^X_x$, ${\cal F}^{\ast X}_x$ as restrictions of ${\cal F}$,
${\cal F}^\ast$ and ${\cal F}^X$, ${\cal F}^{\ast X}$ to the fibers
at $x$, respectively.
 Denote by $A_x$  the restriction of
the gradient $\nabla{\cal F}_x$ to $N{\cal O}(\varepsilon)_x\cap
XN{\cal O}_x$. Then for $\delta>0$ small enough $A_x$  is a $C^1$
map from $XN{\cal O}(\delta)_x$ to $XN{\cal O}_x$ (and so ${\cal
F}^X_x$ is $C^2$ on $XN{\cal O}(\delta)_x$). Clearly,
\begin{equation}\label{e:1.15}
A_{s\cdot x}(s\cdot v)=s\cdot A_x(v)\quad\forall s\in S^1,\; v\in
N{\cal O}(\varepsilon)_x\cap XN{\cal O}_x.
\end{equation}
  Denote by $B_x$ the symmetric bilinear form
$d^2{\cal F}^X_x(0)$. By (i) above Claim~\ref{cl:5.6} we see that it
can be extended to a symmetric bilinear
form on $N{\cal O}_x$, also denoted by
$B_x$, whose associated self-adjoint operator is Fredholm, and has
finite dimensional negative definite and null spaces ${\bf
H}^-(B_x)$ and ${\bf H}^0(B_x)$. Moreover, ${\bf H}^-(B_x)+ {\bf
H}^0(B_x)$ is contained in $XN{\cal O}_x$. As before there exists an
orthogonal decomposition
\begin{equation}\label{e:1.16}
N{\cal O}_x={\bf H}^-(B_x)\oplus {\bf H}^0(B_x)\oplus{\bf H}^+(B_x).
\end{equation}
  Since $B_{s\cdot x}(s\cdot\xi, s\cdot\eta)=B_x(\xi,\eta)$
for any $s\in S^1$ and $x\in{\cal O}$, (\ref{e:1.16}) leads to a
natural $C^{k-2}$ Hilbert vector bundle orthogonal decomposition
\begin{equation}\label{e:1.17}
N{\cal O}={\bf H}^-(B)\oplus {\bf
H}^0(B)\oplus{\bf H}^+(B)
\end{equation}
with ${\bf H}^\star(B)_x={\bf H}^\star(B_x)$ for $x\in{\cal O}$ and
$\star=+,0,-$, which induces a $C^{k-2}$ Banach vector bundle (topological)
direct sum decomposition
\begin{equation}\label{e:1.18}
 XN{\cal O}={\bf H}^-(B)\dot{+} {\bf
H}^0(B)\dot{+}({\bf H}^+(B)\cap XN{\cal O}).
\end{equation}
We call $m^-({\cal O}):=\dim {\bf H}^-(d^2{\cal L}|_{\cal
X}(\gamma_0))={\rm rank}{\bf H}^-(B)$ and
$$
m^0({\cal
O}):=\dim {\bf
H}^0(d^2{\cal L}|_{\cal X}(\gamma_0))-1={\rm rank}{\bf H}^0(B)
$$
 {\bf Morse index} and {\bf nullity} of
${\cal O}=S^1\cdot\gamma_0$,  respectively.
 When $m^0({\cal O})=0$ the orbit ${\cal O}$ is called {\bf
nondegenerate}. Moreover we have always $0\le m^0(\mathcal{O})\le
2n-1$ because $m^0(\mathcal{O})={\rm rank}{\bf
H}^0({B})=m^0(\gamma_0)-1\le 2n-1$ by (\ref{e:5.2}) and the
inequality above Theorem~3.1 of \cite{Lu0}.

Let ${\bf P}^\star$ be the orthogonal bundle projections from
$N{\cal O}$ onto ${\bf H}^\star(B)$, $\star=+,0,-$, and let ${\bf
H}^0(B)(\epsilon)={\bf H}^0(B)\cap N{\cal O}(\epsilon)$ for
$\epsilon>0$. Note that ${\bf H}^0(B)(\epsilon)\subset XN{\cal O}$
and that  we may choose $\epsilon>0$ so small that ${\bf
H}^0(B)(\epsilon)\subset XN{\cal O}(\delta)$ since ${\rm rank}{\bf
H}^0(B)<\infty$ and  ${\cal O}$ is compact. By the implicit function
theorem,  if $0<\epsilon<\varepsilon$ is sufficiently small for each
$x\in{\cal O}$ there exists a unique $S_x^1$-equivariant $C^1$ map
\begin{equation}\label{e:1.19}
 \mathfrak{h}_x:{\bf H}^0(B)(\epsilon)_x\to {\bf H}^-(B)_x\dot{+}({\bf H}^+(B)_x\cap
XN{\cal O}_x)
\end{equation}
(actually at least $C^2$ because $k\ge 4$),  such that
$\mathfrak{h}_x(0_x)=0_x$, $d\mathfrak{h}_x(0_x)=0_x$ and
\begin{eqnarray}\label{e:1.20}
 ({\bf P}^+_x +{\bf P}^-_x)\circ A_x\bigl(v+ \mathfrak{h}_x(v)\bigr)=0
\quad\forall v\in {\bf H}^0(B)(\epsilon)_x.
\end{eqnarray}
Moreover,  the functional
\begin{eqnarray}\label{e:1.21}
{\cal L}^\circ_\triangle:{\bf H}^0(B)(\epsilon)\ni (x,v)\to{\cal
L}\circ{\rm EXP}_x\bigl(v+
 \mathfrak{h}_x(v)\bigr)\in\R
\end{eqnarray}
is $C^1$, has the isolated critical orbit ${\cal O}$ and also
restricts to a $C^2$ one in each fiber ${\bf H}^0(B)(\epsilon)_x$.
Let us denote by ${\cal L}^\circ_{\triangle x}$
the restriction of ${\cal L}^\circ_{\triangle}$ to ${\bf H}^0(B)(\epsilon)_x$.

\begin{theorem}\label{th:1.5}
Under the  above notations, there exists a $C^k$ convex with quadratic
growth Lagrangian $L^\ast:TM\to \R$ such that:
\begin{description}
\item[(i)] $L^\ast\le  L$, $L^\ast(x,v)=L(x,v)$ if $L(x,v)\ge 2c/3$,
 and if $F$ is reversible so is $L^\ast$.

\item[(ii)] The corresponding functional ${\cal L}^\ast$ in
(\ref{e:1.26}) is $C^{2-0}$ in $\Lambda M$. All functional ${\cal
L}^\tau=(1-\tau){\cal L}+ \tau{\cal L}^\ast$, $\tau\in [0,1]$, have
only a critical orbit ${\cal O}$ in some neighborhood of ${\cal
O}\subset\Lambda M$, and  satisfy the Palais-Smale condition.
Moreover, $C_\ast({\cal L}^\ast,{\cal O};\K)=C_\ast({\cal L},{\cal
O};\K)$.

\item[(iii)] By shrinking the above $\epsilon>0$ (if necessary)
 there exist an $S^1$-invariant open neighborhood $U$
of the zero section of $N{\cal O}$, an $S^1$-equivariant
fiber-preserving,  $C^1$ map $\mathfrak{h}$ given by (\ref{e:1.19})
and (\ref{e:1.20}), and an $S^1$-equivariant fiber-preserving
homeomorphism $\Upsilon:N{\cal O}(\epsilon)\to U$ such that
$$
{\cal F}^\ast\circ\Upsilon(x,u)={\cal L}^\ast\circ{\rm
EXP}\circ\Upsilon(x,u)=\|{\bf P}^+_xu\|^2_1-\|{\bf P}^-_xu\|^2_1+
{\cal L}^\circ_{\triangle x}({\bf P}^0_xu)
$$
for all $(x,u)\in N{\cal O}(\epsilon)$. Moreover,
$\Upsilon\left(({\bf P}^-+{\bf P}^0)N{\cal O}(\epsilon)\right)$ is
contained in $XN{\cal O}$, and $\Upsilon$ is also a homeomorphism
from $({\bf P}^-+{\bf P}^0)N{\cal O}(\epsilon)$ onto
$\Upsilon\left(({\bf P}^-+{\bf P}^0)N{\cal O}(\epsilon)\right)$ even
if the topology on the latter is taken as the induced one by
$XN{\cal O}$. (This implies that $N{\cal O}$ and $XN{\cal O}$ induce
the same topology in $\Upsilon\left(({\bf P}^-+{\bf P}^0)N{\cal
O}(\epsilon)\right)$.)

\item[(iv)]  For any open neighborhood $\mathscr{W}$ of ${\cal O}$ in $\Lambda M$ and a field $\K$,
 write $\mathscr{W}_X=\mathscr{W}\cap{\cal X}$ as  an open
 subset of ${\cal X}$, then the inclusion
$$
\left(({\cal L}^\ast|_{\cal X})_c\cap \mathscr{W}_X, ({\cal
L}^\ast|_{\cal X})_c\cap \mathscr{W}_X\setminus{\cal
O}\right)\hookrightarrow\left({\cal L}^\ast_c\cap \mathscr{W}, {\cal
L}^\ast_c\cap \mathscr{W}\setminus{\cal O}\right)
$$
 induces isomorphisms
$$
H_\ast\left(({\cal L}^\ast|_{\cal X})_c\cap \mathscr{W}_X, ({\cal
L}^\ast|_{\cal X})_c\cap \mathscr{W}_X\setminus{\cal O};\K\right)\to
H_\ast\left({\cal L}^\ast_c\cap \mathscr{W}, {\cal L}^\ast_c\cap
\mathscr{W}\setminus{\cal O};\K\right).
$$
Moreover, the corresponding conclusion can be obtained if
$$
\left(({\cal L}^\ast|_{\cal X})_c\cap \mathscr{W}_X, ({\cal
L}^\ast|_{\cal X})_c\cap \mathscr{W}_X\setminus{\cal
O}\right)\quad\hbox{and}\quad\left({\cal L}^\ast_c\cap \mathscr{W},
{\cal L}^\ast_c\cap \mathscr{W}\setminus{\cal O}\right)
$$
are replaced by $\bigl((\mathring{\cal L}^\ast|_{\cal X})_c\cap
\mathscr{W}_X\cup{\cal O}, (\mathring{\cal L}^\ast|_{\cal X})_c\cap
\mathscr{W}_X\bigr)$ and $\bigl(\mathring{\cal L}^\ast_c\cap
\mathscr{W}\cup{\cal O}, \mathring{\cal L}^\ast_c\cap
\mathscr{W}\bigr)$, respectively, where $\mathring{\cal L}^\ast_c=\{
{\cal L}^\ast<c\}$ and $(\mathring{\cal L}^\ast|_{\cal X})_c=\{{\cal
L}^\ast|_{\cal X}<c\}$.

\item[(v)] If ${\cal L}^{\ast}|_{\mathcal{X}}$ and  ${\cal
L}^\ast$ in (iv) are replaced by ${\cal L}|_{\mathcal{X}}$ and
${\cal L}$, respectively, then the corresponding conclusion also
holds true.
\end{description}
\end{theorem}

Statement (v) in Theorem~\ref{th:1.5} is
based on the following splitting lemma for ${\cal
L}|_{\cal X}$.

\begin{theorem}\label{th:1.6}
Under the notations above, by shrinking the above $\epsilon>0$
 there exist an $S^1$-invariant open neighborhood $V$
of the zero section of $XN{\cal O}$, an $S^1$-equivariant
fiber-preserving, $C^1$ map $\mathfrak{h}$ given by (\ref{e:1.19})
and (\ref{e:1.20}), and an $S^1$-equivariant fiber-preserving
homeomorphism $\Psi:XN{\cal O}(\epsilon)\to V$ such that
$$
{\cal L}\circ{\rm EXP}\circ\Psi(x, v)=\frac{1}{2}d^2{\cal L}|_{\cal
X}(x)[{\bf P}^+_xv, {\bf P}^+_xv]-\|{\bf P}^-_xv\|^2_1+ {\cal
L}^\circ_{\triangle x}({\bf P}^0_xv)
$$
for all $(x,v)\in XN{\cal O}(\epsilon)$.
\end{theorem}

The infinite dimensional proof method of this result using \cite{Ji}
may go back to \cite{Lu}. Clearly,
Theorem~\ref{th:1.5}(v) implies
\begin{eqnarray}\label{e:1.22}
C_\ast({\cal L},{\cal O};\K)=C_\ast({\cal L}|_{\cal X},{\cal O};\K).
\end{eqnarray}
This  also easily follows from Theorem~\ref{th:1.5}(i),(ii) and
(iv). Let ${\bf H}^{0-}(B)={\bf H}^0(B)+{\bf H}^-(B)$ and ${\bf
H}^{0-}(B)(\epsilon)=({\bf H}^0(B)+{\bf H}^-(B))\cap N{\cal
O}(\epsilon)$. Then ${\bf H}^{0-}(B)\subset XN{\cal O}$. Define
$\mathfrak{L}:{\bf H}^{0-}(B)(\epsilon)\to\R$ by $
\mathfrak{L}(x,v)=-\|{\bf P}^-_xv\|^2_1+ {\cal L}^\circ_{\triangle
x}({\bf P}^0_xv)$. By a well known deformation argument we may
derive from Theorem~\ref{th:1.5}(ii)-(iii) and Theorem~\ref{th:1.6}
respectively:
$$
C_\ast({\cal L},{\cal O};\K)=C_\ast(\mathfrak{L},{\cal
O};\K)\quad\hbox{and}\quad  C_\ast({\cal L}|_{\cal X},{\cal
O};\K)=C_\ast(\mathfrak{L},{\cal O};\K).
$$
These give rise to (\ref{e:1.22}) again.

Let $S^1_x\subset S^1$ denote the stabilizer of $x\in\mathcal{O}$.
Since $x$ is nonconstant $S^1_x$ is  a finite cyclic group and the
quotient $S^1/S^1_x\cong S^1\cdot\gamma_0={\cal O}\cong S^1$. (See
\cite[page 499]{GM2}).  Clearly, ${\cal L}^\circ_{\triangle x}$ is
$S^1_x$-invariant. Let $C_\ast({\cal L}^\circ_{\triangle x},
0;\K)^{S^1_x}$ denote the subgroup of all elements in $C_\ast({\cal
L}^\circ_{\triangle x}, 0;\K)$, which are fixed by the induced
action of $S^1_x$ on the homology.  We have the following
generalization of the Gromoll-Meyer shifting theorem for Finsler
manifolds.

\begin{theorem}\label{th:1.7}
Let $\K$ be a field of characteristic $0$ or prime up to order
$|S^1_{\gamma_0}|$ of $S^1_{\gamma_0}$. Then for any $x\in{\cal
O}=S^1\cdot\gamma_0$ and $q=0,1,\cdots$,
\begin{eqnarray*}
&&C_q({\cal L}, {\cal O};\K)\\
&=& \Bigl(H_{m^-({\cal O})}({\bf H}^-(B)_x, {\bf
H}^-(B)_x\setminus\{0_x\};\K)\otimes
C_{q-m^-({\cal O})}({\cal L}^{\circ}_{\triangle x}, 0;\K)\Bigr)^{S^1_x}\nonumber\\
&\oplus& \Bigl(H_{m^-({\cal O})}({\bf H}^-(B)_x, {\bf
H}^-(B)_x\setminus\{0_x\};\K)\otimes C_{q-m^-({\cal O})-1}({\cal
L}^{\circ}_{\triangle x}, 0;\K)\Bigr)^{S^1_x}.
\end{eqnarray*}
provided  $m^-({\cal O})m^0({\cal O})>0$. Moreover,
$$C_q({\cal L},
{\cal O};\K)= \bigl( C_{q-1}({\cal L}^{\circ}_{\triangle x},
0;\K)\bigr)^{S^1_x}\oplus\bigl(C_{q}({\cal L}^{\circ}_{\triangle x},
0;\K)\bigr)^{S^1_x}
$$
if $m^-({\cal O})=0$ and $m^0({\cal O})>0$, and
\begin{eqnarray*}
C_q({\cal L}, {\cal O};\K) &=&H_{q}({\bf H}^-(B), {\bf
H}^-(B)\setminus{\cal O};\K)\nonumber\\
&=&\Bigl(H_{q-1}({\bf H}^-(B)_x, {\bf
H}^-(B)_x\setminus\{0_x\};\K)\Bigr)^{S^1_x}\nonumber\\
&\qquad\oplus& \Bigl(H_{q}({\bf H}^-(B)_x, {\bf
H}^-(B)_x\setminus\{0_x\};\K)\Bigr)^{S^1_x}
\end{eqnarray*}
 if $m^-({\cal O})>0$ and $m^0({\cal O})=0$. Finally,
  $C_q({\cal L}, S^1\cdot\gamma_0;\K)=H_q(S^1;\K)$  for any Abelian group  $\K$
  if   $m^-({\cal O})=m^0({\cal O})=0$.
\end{theorem}

By \cite[Th.I]{Bo} $H_{q}({\bf H}^-(B), {\bf H}^-(B)\setminus{\cal
O};\Z_2)=\Z_2$ for $q=m^-({\cal O}), m^-({\cal O})+1$, and $=0$
otherwise. If $\K=\Z$ and ${\bf H}^-(B)$ is orientable the same is
also true; see  \cite[Cor.2.4.11]{Kl}.

 In the studies of closed geodesics one often define the
critical group
$$
C_\ast({\cal L},
S^1\cdot\gamma_0;\K)=H_\ast(\Lambda(\gamma_0)\cup S^1\cdot\gamma_0,
\Lambda(\gamma_0);\K),
$$
where $\Lambda(\gamma_0)=\{\gamma\in\Lambda M\,|\, {\cal
L}(\gamma)<{\cal L}(\gamma_0)\}$.  Using the excision property of
singular homology and anti-gradient flow it is not hard to prove
that these two kinds of definitions agree (cf.
\cite[Propositions~3.4 and 3.7]{Cor}). Such a version of
Theorem~\ref{th:1.7} was proved in \cite[Prop.3.7]{BLo} by
introducing finite dimensional approximations to $\Lambda M$ as in
 \cite{Ma, Ra, Shen}.

Since one does not know if a generator of the $S^1_x$-action on
${\bf H}^-(B)_x$ reverses orientation or not, no further explicit
version of the formula in Theorem~\ref{th:1.7} can be obtained
through
$$
H_{m^-({\cal O})}({\bf H}^-(B)_x, {\bf
H}^-(B)_x\setminus\{0_x\};\K)=\K.
$$
Recalling that $S^1=\R/\Z=[0,1]/\{0,1\}$, there exists a positive integer
$m=m(\gamma_0)$ such that $1/m$ is the minimal period of
 $\gamma_0$ (since $\gamma_0$ is nonconstant).
 It is equal to the order
of the isotropy group $S^1_{\gamma_0}$, and is called the {\bf
multiplicity} of $\gamma_0$. When $m(\gamma_0)=1$ we say $\gamma_0$
to be {\bf prime}. These can also be described by the $m$-th iterate
operation
\begin{equation}\label{e:1.23}
\varphi_m:\Lambda M\to\Lambda M, \gamma\to\gamma^m
\end{equation}
defined by\footnote{Here $\gamma^m$ is different from $\gamma^m$
appearing in the study of Lagrangian Conley conjecture in \cite{Lo1,
Lu0, Lu1} though we use the same notation.}
$\gamma^m(t)=\gamma(mt)\;\forall t\in \R$ when $\gamma$ is viewed as
a $1$-periodic map $\gamma:\R\to M$. Clearly, there exists a unique
prime curve $\gamma_0^{1/m}\in\Lambda M$ such that
 $(\gamma_0^{1/m})^m$  is
equal to $\gamma_0$. Suppose that $\K$ is a field $\Q$ of rational
numbers and that for each $k\in\N$ the orbit
$S^1\cdot(\gamma_0^{1/m})^k$ is an isolated critical orbit of ${\cal
L}$. We may rewrite the conclusions of Proposition~3.8 in \cite{BLo}
(in our notations) as follows:
\begin{description}
\item[(i)] If $m^0(S^1\cdot\gamma_0)=0$ then
$C_q({\cal L}, S^1\cdot\gamma_0;\K)=\K$ if
$m^-(S^1\cdot\gamma_0)-m^-(S^1\cdot\gamma_0^{1/m})$ is even and
$q\in\{m^-(S^1\cdot\gamma_0), m^-(S^1\cdot\gamma_0)+1\}$, and
$C_q({\cal L}, S^1\cdot\gamma_0;\K)=0$ in other cases.
\item[(ii)] If $m^0(S^1\cdot\gamma_0)>0$ and
$\epsilon(\gamma_0)=(-1)^{m^-(S^1\cdot\gamma_0)-m^-(S^1\cdot\gamma_0^{1/m})}$,
then
\begin{eqnarray*}
C_q({\cal L}, S^1\cdot\gamma_0;\K) &=&
C_{q-m^-(S^1\cdot\gamma_0)-1}({\cal
L}^{\circ}_{\triangle x}, 0;\K)^{S^1_x, \epsilon(\gamma_0)}\nonumber\\
&\oplus&  C_{q-m^-(S^1\cdot\gamma_0)}({\cal L}^{\circ}_{\triangle
x}, 0;\K)^{S^1_x, \epsilon(\gamma_0)}
\end{eqnarray*}
for each $x\in\mathcal{O}$ and $q\in\N\cup\{0\}$, where
 $C_{\ast}({\cal L}^{\circ}_{\triangle x}, 0;\K)^{S^1_x, 1}$ and
$C_{\ast}({\cal L}^{\circ}_{\triangle x}, 0;\K)^{S^1_x, -1}$ are the
eigenspaces of a generator of $S_x^1$ corresponding to $1$ and $-1$,
respectively. Clearly, $C_{\ast}({\cal L}^{\circ}_{\triangle x},
0;\K)^{S^1_x, -1}=0$ if $m$ is odd.
\end{description}

Using Theorem~\ref{th:1.7} we may derive the following result, which
is very important for the proof of Theorem~\ref{th:1.11} in
Section~\ref{sec:8}.

\begin{theorem}\label{th:1.8}
Under the assumptions of Theorem~\ref{th:1.7} suppose that
$m^-(\gamma_0)=0$ and that $C_p({\cal L}, S^1\cdot\gamma_0;\K)\ne 0$
and $C_{p+1}({\cal L}, S^1\cdot\gamma_0;\K)=0$ for some
$p\in\N\cup\{0\}$. Then $p>0$. Furthermore, we have
\begin{description}
\item[(i)] if $p=1$  then each point of $S^1\cdot\gamma_0$
is a local  minimum of  $\mathcal{L}$;

\item[(ii)] if $p\ge 2$ then each point of $S^1\cdot\gamma_0$ is not a local  minimum of
 $\mathcal{L}|_{\cal X}$ (and thus of $\mathcal{L}$).
\end{description}
\end{theorem}

\subsection{Two iteration theorems}\label{sec:1.3}

Our shifting theorems, Theorems~\ref{th:1.4},~\ref{th:1.7}, are
sufficient for computations of critical groups in most of studies about
geodesics on a Finsler manifold. In Riemannian geometry they are
direct consequences of the splitting theorem. So far we have not
obtained the corresponding splitting theorem for the Finsler energy
functional $\mathcal{L}$ on the Hilbert manifold $\Lambda_N(M)$. In
the studies of multiplicity of closed geodesics on a Riemannian
manifold as in \cite{GM2, BKl} etc, one must use the splitting
theorem to deduce a change in the critical groups under the
iteration map $\varphi_m$ (\cite[Th.3]{BKl}). The following is the
corresponding generalization of such results to Finsler manifolds.

\begin{theorem}\label{th:1.9}
For a closed geodesics $\gamma_0$ and some integer $m>1$, suppose
that ${\cal O}=S^1\cdot\gamma_0$ and $\varphi_m({\cal
O})=S^1\cdot\gamma_{0}^m$ are two isolated critical orbits of ${\cal
L}$ in $\Lambda M$ and that $m^-({\cal O})=m^-(\varphi_m({\cal O}))$
and $m^0({\cal O})=m^0(\varphi_m({\cal O}))$. Then $\varphi_m$
induces isomorphisms
\begin{eqnarray*}
(\varphi_m)_\ast: H_\ast\left(\Lambda(\gamma_0)\cup
S^1\cdot\gamma_0, \Lambda(\gamma_0);\K\right)\to
H_\ast\left(\Lambda(\gamma_{0}^m)\cup S^1\cdot\gamma_{0}^m,
\Lambda(\gamma_{0}^m);\K\right).
\end{eqnarray*}
for any field $\K$. (So $\varphi_m$ induces isomorphisms from
$C_\ast({\cal L}, {\cal O};\K)$ to $C_\ast({\cal L}, \varphi_m({\cal
O});\K)$.)
 \end{theorem}

Using Theorem~\ref{th:1.6} it can only be proved that $\varphi_m$
induces isomorphisms
\begin{eqnarray*}
(\varphi_m)_\ast: H_\ast\left(\Lambda(\gamma_0)^X\cup
S^1\cdot\gamma_0, \Lambda(\gamma_0)^X;\K\right)\to
H_\ast\left(\Lambda(\gamma_{0}^m)^X\cup S^1\cdot\gamma_{0}^m,
\Lambda(\gamma_{0}^m)^X;\K\right),
\end{eqnarray*}
where both $\Lambda(\gamma_0)^X:={\cal X}\cap\Lambda(\gamma_0)$ and
$\Lambda(\gamma^m_0)^X:={\cal X}\cap\Lambda(\gamma^m_0)$ are viewed
as topological subspaces of the Banach manifold ${\cal X}$. It is
Theorems~\ref{th:1.5}(v) that make us to derive Theorem~\ref{th:1.9}
from this conclusion. This and (\ref{e:1.22}) cannot lead to such a
result though they are enough in most of applications.

Under the weaker assumption that $m^0({\cal O})=m^0(\varphi_m({\cal
O}))$ some results can also be obtained.
 Since the germ of the map $\mathfrak{h}_x$ at the origin in
 (\ref{e:1.19}) is uniquely determined by ${\mathcal L}$ and the
 metric $g$ we call $\mathscr{N}_x:=\{{\rm EXP}_x\bigl(v+
 \mathfrak{h}_x(v)\bigr)\,|\,v\in {\bf H}^0(B)(\epsilon)_x\}$
a {\bf local characteristic manifold} of ${\cal L}$ at $x$ with
respect to $g$. It is actually a $S^1_x$-invariant, at least $C^{2}$
submanifold of dimension $m^0({\cal O})$ in $H^1(S^1, M)$ and
contains $x$ as its interior point. Moreover $x$ is an isolated
critical point of ${\cal L}|_{\mathscr{N}_x}$ and $d^2({\cal
L}|_{\mathscr{N}_x})(x)=0$.

\begin{theorem}\label{th:1.10}
Let $\K$ be a field. Suppose that $m^0({\cal O})=m^0(\varphi_m({\cal
O}))$. Then
\begin{eqnarray}\label{e:1.24}
\dim C_{q}({\cal L}|_{\mathscr{N}_x}, x;\K)=\dim C_{q}({\cal
L}|_{\mathscr{N}_{x^m}}, x^m;\K)\quad\forall q
\end{eqnarray}
for any $x\in{\cal O}$. If the characteristic of $\K$ is zero or
prime up to order of $S^1_{\gamma_{0}^m}$ then
\begin{eqnarray}\label{e:1.25}
\dim C_{q}({\cal L}|_{\mathscr{N}_x}, x;\K)^{S^1_{ x}}=\dim
C_{q}({\cal L}|_{\mathscr{N}_{x^m}}, x^m;\K)\quad\forall
q\in\{0\}\cup\N.
\end{eqnarray}
\end{theorem}

 (\ref{e:1.24}) may be viewed as the result analogous to
 Theorem~3 of \cite{GM2} on Riemannian
manifolds. With finite-dimensional approximations
Theorem~\ref{th:1.9} and the equivalent forms of
(\ref{e:1.24})-(\ref{e:1.25}) were proved  in \cite[\S 7.1,7.2]{Ra}
and \cite[Th.3.11]{BLo}.

\subsection{An application}\label{sec:1.4}

Using the above theory many results about closed geodesics on
Riemannian manifolds can be generalized to Finsler manifolds in a
straightforward way. For example, repeating the arguments of
\cite{BKl} will lead to similar results. In particular we have the
following generalized version of \cite[Theorem 3]{BKl}.

\begin{theorem}\label{th:1.11}
A connected closed Finsler manifold $(M,F)$ of dimension $n>1$ has
infinitely many geometrically distinct closed geodesics provided
that there  exists a nonconstant closed geodesics $\bar\gamma$ such
that $m^-(\bar\gamma^k)=0$ for all $k\in\N$ and
$H_{\bar{p}}(\Lambda(\bar\gamma)\cup S^1\cdot \bar\gamma,
\Lambda(\bar\gamma);\Q)\ne 0$ for some integer $\bar{p}\ge 2$.
\end{theorem}


The proof of a slightly different version of Theorem 1.11 was outlined by  Rademacher  with
finite-dimensional approximations  \cite[Theorem~7.5]{Ra}. We give a self-contained detailed
proof by improving the arguments in \cite{BKl, Lo1, LoLu, Lu0,Lu1}.
Our proof method is slightly different from \cite{BKl}, and cannot
deal with the case that $\bar{p}=1$ and $\bar\gamma$ is a local
minimum of ${\cal L}$, but not an absolute minimum of ${\cal L}$ in
its free homotopy class. See Remark~\ref{rm:8.8} for comparisons
with the results in \cite{BKl}.

\subsection{Basic idea of the proof}\label{sec:1.5}

Since $W^{1,2}(I, M)\hookrightarrow C^0(I, M)$ is continuous, and
each $\gamma\in W^{1,2}(I, M)$ has compact image, there exists an
open neighborhood ${\cal O}(\gamma)$ of $\gamma$ in $W^{1,2}(I, M)$
such that $\cup\{\alpha(I)\,|\,\alpha\in {\cal O}(\gamma)\}$ is
contained in a compact subset of $M$. Hence as before we may
\textsf{assume that $M$ is compact.}

Recall that  a Lagrangian $L:[0,1]\times TM\to\R$ is
called {\bf convex with quadratic growth} if it satisfies the conditions:
\begin{description}
\item[(L1)] $\exists$ a constant $\ell_0>0$ such that
$\partial_{vv}L(t,x,v)\ge\ell_0 I$,

\item[(L2)] $\exists$ a constant $\ell_1>0$ such that
$\|\partial_{vv}L(t,x,v)\|\le\ell_1$ and
$$\|\partial_{xv}L(t,x,v)\|\le\ell_1(1+ |v|_x),\quad
\|\partial_{xx}L(t,x,v)\|\le \ell_1(1+|v|_x^2)
$$
with respect to some Riemannian metric $\langle\cdot,\cdot\rangle$
(with $|v|^2_x=\langle v,v\rangle_x$).
\end{description}
 Equivalently, there exists a finite atlas on $M$ such that
 under each chart of this atlas the following conditions
 hold for some constants $0<c<C$:
\begin{description}
\item[(L1)] $\sum_{ij}\frac{\partial^2}{\partial v_i\partial
v_j}L(t,x,v)u_iu_j\ge c|{\bf u}|^2\quad\forall {\bf
u}=(u_1,\cdots,u_n)\in\R^n$,

\item[(L2)] $\Bigl| \frac{\partial^2}{\partial x_i\partial
x_j}L(t,x,v)\Bigr|\le C(1+ |v|^2),\quad \Bigl|
\frac{\partial^2}{\partial x_i\partial v_j}L(t,x,v)\Bigr|\le C(1+
|v|),\quad\hbox{and}\\
 \Bigl| \frac{\partial^2}{\partial v_i\partial
v_j}L(t,x,v)\Bigr|\le C$.
\end{description}

\textsf{To show our ideas let us consider Theorem~\ref{th:1.4} for
example.} Given a nontrivial constant speed $F$-geodesic
$\gamma_0:I\to M$ with $(\gamma_0(0), \gamma_0(1))\in N=M_0\times
M_1$, $c=L(\gamma_0,\dot\gamma_0)=[F(\gamma_0,\dot\gamma_0)]^2>0$.
Then we construct a convex with quadratic growth Lagrangian $L^\ast:
TM\to\R$ such that $L^\ast(x,v)= L(x,v)$ if $L(x,v)\ge\frac{2c}{3}$.
Clearly, $\gamma_0$ is a critical point of the functional
\begin{equation}\label{e:1.26}
 {\cal
L}^\ast:\Lambda_N(M)\to\R,\;\gamma\mapsto\int^1_0 L^\ast(\gamma(t),
\dot\gamma(t))dt
\end{equation}
with critical value $c$. Moreover,  $\gamma_0$ is isolated  for
${\cal L}$ if and only if it is so for
 ${\cal L}^\ast$. Define $L^\tau(x,v)=(1-\tau)L(x,v)+
 \tau L^\ast(x,v)$ and
${\cal L}^\tau=(1-\tau){\cal L}+
 \tau {\cal L}^\ast$
for $\tau\in [0,1]$. We shall prove that the family of functionals
$\{{\cal L}^\tau\;|\;\tau\in [0,1]\}$ on $\Lambda_N(M)$
 satisfies the stability property of critical groups
\cite{Ch93, ChGh, Ch05, CiDe, CorHa, MW}, and so
$$
C_q({\cal L}, \gamma_0;{\K})\cong C_q({\cal L}^\ast,
\gamma_0;{\K})\quad\forall q\ge 0.
$$
By Corollary~\ref{cor:A.2} we have a shifting theorem for $C_q({\cal
L}^\ast, \gamma_0;{\K})$ and hence for $C_q({\cal L},
\gamma_0;{\K})$ because ${\cal L}^\ast$ and ${\cal L}$ have the same
characteristic manifold on which they agree.

Our modified Lagrangian $L^\ast$ can be  required to be no more than
$L$. This is very key for the proofs of Theorem~\ref{th:1.5}(v),
Theorem~\ref{th:1.8}(i) and Theorem~\ref{th:1.9}.

\vspace{2mm}

\noindent{\bf Organization of the paper}.  In Section 2 we start
from $L=F^2$ to construct a suitable convex with quadratic growth
Lagrangian $L^\ast$ having the properties outlined above. Then by
considering the corresponding functional family $({\cal
L}^\tau)_{\tau\in [0,1]}$
  with the Lagrangians $L^\tau=(1-\tau)L+
 \tau  L^\ast$ with $\tau\in [0,1]$
we show in Section 3 that the functionals ${\cal L}$ and ${\cal
L}^\ast$ have the same critical groups at $\gamma_0$ and
$S^1\cdot\gamma_0$ in two cases, respectively. The proofs of
Theorems~\ref{th:1.2},~\ref{th:1.3} are given in Section 4, and
those of Theorems~\ref{th:1.5},~\ref{th:1.6},~\ref{th:1.7} and
\ref{th:1.8} are given in Section 5. Section 6 deals with critical
groups of iterated closed geodesics, including the proofs of
Theorems~\ref{th:1.9},~\ref{th:1.10}.  In Section~\ref{sec:7} we
present a  computation method of $S^1$-critical groups. The proof of
Theorem~\ref{th:1.11} is given in Section~\ref{sec:8}. Finally, in
Appendixes A, B we state the splitting lemma and the shifting theorem obtained in
\cite{Lu1,Lu3} and give some related computations, respectively.

Our methods  can also be used to generalize the isometric-invariant
geodesic theory in \cite{GroTa78, Tan82} to Finsler manifolds. They
will be given elsewhere.\vspace{2mm}

\section{The modifications for the energy functionals}\label{sec:2}
\setcounter{equation}{0}

\begin{figure}[!htb]
\centering \psfrag{x}{\footnotesize $x$}\psfrag{y}{\footnotesize
$y$}\psfrag{o}{\footnotesize $O$}\psfrag{varepsilon}{\footnotesize
$\varepsilon$} \psfrag{delta}{\footnotesize
$\delta$}\psfrag{2c}{\footnotesize
$\frac{2c}{3C_1}$}\psfrag{psi}{\footnotesize
$\psi_{\varepsilon,\delta}$}\psfrag{phi}{\footnotesize
$\phi_{\mu,b}+\mu\delta$}
\includegraphics[clip,width=8cm]{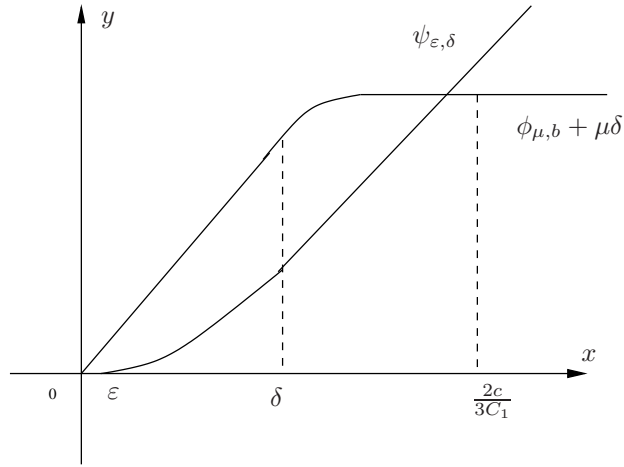}
\caption{The functions $\psi_{\varepsilon,\delta}$ and
$\phi_{\mu,b}+\mu \delta$.}\label{fig:laolu}
\end{figure}

In this section we consider $C^k$ Finsler metrics ($k\ge 2$).
We firstly construct two smooth functions, see
Figure~\ref{fig:laolu}.

\begin{lemma}\label{lem:2.1}
Given positive numbers $c>0$ and $C_1\ge 1$, choose positive
parameters $0<\varepsilon<\delta<\frac{2c}{3C_1}$. Then
\begin{description}
\item[(i)]  there exists
a $C^\infty$ function $\psi_{\varepsilon,\delta}:[0, \infty)\to\R$
such that: $\psi'_{\varepsilon,\delta}>0$ and $\psi$ is convex on
$(\varepsilon, \infty)$, it vanishes in $[0, \varepsilon)$, is equal to
the affine function $\kappa t+ \varrho_0$ on $[\delta, \infty)$,
where $\kappa>0$ and $\varrho_0<0$ are suitable constants;
\item[(ii)] there exists a $C^\infty$ function
$\phi_{\mu,b}:[0, \infty)\to\R$ depending on parameters $\mu>0$ and
$b>0$, such that: $\phi_{\mu,b}$ is nondecreasing and concave (and
hence $\phi_{\mu,b}^{\prime\prime}\le 0$ ), it is equal to the affine
function $\mu t-\mu\delta$ on $[0, \delta]$, is equal to constant
$b>0$ on $[\frac{2c}{3C_1}, \infty)$.
\item[(iii)] Under the above assumptions, $\psi_{\varepsilon,\delta}(t)+
\phi_{\mu,b}(t)-b=\kappa t+ \varrho_0$ for any $t\ge
\frac{2c}{3C_1}$ (and hence for $t\ge\frac{2c}{3}$). Moreover,
$\psi_{\varepsilon,\delta}(t)+ \phi_{\mu,b}(t)-b\ge
-\mu\delta-b\;\forall t\ge 0$, and $\psi_{\varepsilon,\delta}(t)+
\phi_{\mu,b}(t)-b= -\mu\delta-b$ if and only if $t=0$.

\item[(iv)] Under the assumptions (i)-(ii), suppose that the constant $\mu>0$ satisfies
\begin{eqnarray}\label{e:2.1}
\mu+\frac{\varrho_0}{\delta-\varepsilon}>0\quad\hbox{and}\quad\mu\delta+b+\varrho_0>0.
\end{eqnarray}
Then $\psi_{\varepsilon,\delta}(t)+ \phi_{\mu,b}(t)-b\le\kappa t+
\varrho_0\;\forall t\ge\varepsilon$, and
$\psi_{\varepsilon,\delta}(t)+ \phi_{\mu,b}(t)-b\le\kappa t+
\varrho_0\;\forall t\in [0,\varepsilon]$ if $\kappa\ge\mu$.
\end{description}
\end{lemma}

\noindent{\bf Proof of Lemma~\ref{lem:2.1}}. (i)-(iii) are easily
obtained. We only prove (iv).  Observe  that the line connecting
points $(\varepsilon,0)$ to $(\delta, \kappa\delta+ \varrho_0)$
given by
$$
\R\ni
t\mapsto\frac{\kappa\delta+\varrho_0}{\delta-\varepsilon}(t-\varepsilon),
$$
satisfies
\begin{eqnarray}\label{e:2.2} \kappa
t+\varrho_0\le\psi_{\varepsilon,\delta}(t)\le
\frac{\kappa\delta+\varrho_0}{\delta-\varepsilon}(t-\varepsilon)\quad\forall
t\in [\varepsilon,\delta].
\end{eqnarray}
Since $\phi_{\mu,b}(t)-b\le 0\;\forall t\ge 0$ we have
$\psi_{\varepsilon,\delta}(t)+ \phi_{\mu,b}(t)-b\le\kappa t+
\varrho_0\;\forall t\in [\delta, \infty)$ by the definition of
$\psi_{\varepsilon,\delta}$. When $t\in [\varepsilon,\delta]$, since
$\phi_{\mu,b}(t)=\mu t-\mu\delta$ for any $t\le\delta$, by
(\ref{e:2.2}) we only need to prove
\begin{eqnarray}\label{e:2.2.1}
\frac{\kappa\delta+\varrho_0}{\delta-\varepsilon}(t-\varepsilon)+
\mu t-\mu\delta-b\le\kappa t+ \varrho_0\quad\forall t\in
[\varepsilon,\delta].
\end{eqnarray}
Clearly, this is equivalent to
\begin{eqnarray}\label{e:2.2.2}
\left(\frac{\kappa\varepsilon+\varrho_0}{\delta-\varepsilon}+ \mu
\right)t\le
\varrho_0+\mu\delta+b+\frac{\kappa\delta+\varrho_0}{\delta-\varepsilon}\varepsilon
\quad\forall t\in [\varepsilon,\delta]
\end{eqnarray}
because
$\frac{\kappa\delta}{\delta-\varepsilon}-\kappa=\frac{\kappa\varepsilon}{\delta-\varepsilon}$.
The first inequality in (\ref{e:2.1}) implies
$\frac{\kappa\varepsilon+\varrho_0}{\delta-\varepsilon}+\mu\ge
\frac{\kappa\varepsilon}{\delta-\varepsilon}>0$. Hence it suffices
to prove that (\ref{e:2.2.2}), or equivalently,  (\ref{e:2.2.1})
holds for $t=\delta$, that is,
$$
\frac{\kappa\delta+\varrho_0}{\delta-\varepsilon}(\delta-\varepsilon)+
\mu \delta-\mu\delta-b\le\kappa \delta+ \varrho_0.
$$
It is obvious because $b>0$. This proves the first claim in (iv).

For the case $t\in [0, \varepsilon]$ we have
 $\psi_{\varepsilon,\delta}(t)=0$ and
\begin{eqnarray*}
 \phi_{\mu,b}(t)-b
=\mu t-\mu\delta-b \le \kappa t+ \mu\delta+ b+\varrho_0-\mu\delta-b
=\kappa t+ \varrho_0
\end{eqnarray*}
 because the second inequality in (\ref{e:2.1})
implies $(\mu-\kappa)t\le\mu\delta+ b+\varrho_0$ for any
$\mu<\kappa$.  Lemma~\ref{lem:2.1}(iv) is proved.
\hfill$\Box$\vspace{2mm}

\begin{proposition}\label{prop:2.2}
Let $(M, F)$ be a $C^k$ Finsler manifold ($k\ge 2$),  and $L:=F^2$.
Suppose that for some Riemannian metric $g$ on $M$ both
\begin{eqnarray*}
\alpha_g:=\inf_{x\in M,\, |v|_x=1}\inf_{u\ne
0}\frac{g^F(x,v)[u,u]}{g_x(u,u)}\quad\hbox{and}\quad
\beta_g:=\sup_{x\in M,\, |v|_x=1}\sup_{u\ne
0}\frac{g^F(x,v)[u,u]}{g_x(u,u)}
\end{eqnarray*}
are positive numbers, and that for some constant $C_1>0$,
\begin{equation}\label{e:2.3}
|v|^2_x\le L(x, v)\le C_1|v|^2_x\quad\forall (x,v)\in TM.
\end{equation}
Hereafter $|v|_x=\sqrt{g_x(v,v)}$. Then
for a given $c>0$ there exists a $C^k$ Lagrangian $L^\star:TM\to \R$
such that
\begin{description}

\item[(i)] $L^\star(x,v)=\kappa L(x,v)+
\varrho\quad\hbox{if}\;L(x,v)\ge\frac{2c}{3C_1}$,

\item[(ii)] $L^\star$ attains the minimum, and $L^\star(x,v)=\min L^\star\Longleftrightarrow
v=0$,

\item[(iii)]  $L^\star(x,v)\le \kappa L(x,v)+ \varrho_0$ for all
$(x,v)\in TM$,

 \item[(iv)] $\partial_{vv}
L^\star(x,v)[u,u]\ge \min\{2\mu,
\frac{1}{2}\kappa\alpha_g\}|u|_x^2$.
\end{description}
Moreover, if $F$ is reversible, i.e. $F(x,-v)=F(x,v)\;\forall
(x,v)\in TM$, so is $L^\star$.
\end{proposition}

\noindent{\bf Proof}. By the assumptions, for any $(x,v)\in
TM\setminus\{0\}$ and $(x,u)\in TM$ we have
\begin{equation}\label{e:2.4}
\alpha_g |u|_x^2\le g^F(x,v)[u,u]\le\beta_g |u|_x^2.
\end{equation}
Let $\psi_{\varepsilon,\delta}$ and $\phi_{\mu,b}$ be as in
Lemma~\ref{lem:2.1}. Suppose that (\ref{e:2.1}) is satisfied and
that $\kappa\ge\mu$. Consider the function $L^\star:TM\to\R$ defined
by
\begin{eqnarray}\label{e:2.5}
 L^\star(x,v)=\psi_{\varepsilon,\delta}(L(x,v))+
\phi_{\mu,b}(|v|^2_x)-b.
\end{eqnarray}
 Clearly, it is $C^k$ smooth and satisfies the final claim.

In fact, by (\ref{e:2.3}),
$\phi_{\mu,b}(|v|^2_x)\le\phi_{\mu,b}(L(x,v))\;\forall (x,v)$. (i)
 follows from Lemma~\ref{lem:2.1}(iii) immediately, and (ii) comes from
 the fact that both $\psi_{\varepsilon,\delta}$ and $\phi_{\mu,b}$ are nondecreasing.
  Observe that
$$
L^\star(x,v)=\psi_{\varepsilon,\delta}(L(x,v))+
\phi_{\mu,b}(|v|^2_x)-b\le \psi_{\varepsilon,\delta}(L(x,v))+
\phi_{\mu,b}(L(x,v))-b
$$
for any $(x,v)$. We derive (iii) from Lemma~\ref{lem:2.1}(iv).

\underline{Now we are ready to prove (iv)}.  Since
\begin{eqnarray*}
&&\frac{\partial^2}{\partial t\partial
s}\psi_{\varepsilon,\delta}(L(x,v+su+ tu))=\frac{\partial}{\partial
t}\Bigl[\psi'_{\varepsilon,\delta}(L(x,v+su+
tu))\frac{\partial}{\partial s}L(x,v+su+ tu)\Bigr]\\
&&=\psi'_{\varepsilon,\delta}(L(x,v+su+
tu))\frac{\partial^2}{\partial t\partial s}L(x,v+su+ tu)\\
&&+\psi''_{\varepsilon,\delta}(L(x,v+su+
tu))\frac{\partial}{\partial s}L(x,v+su+ tu)\frac{\partial}{\partial
t}L(x,v+su+ tu)
\end{eqnarray*}
for $u\in T_xM$ and $v\in T_xM\setminus\{0\}$, we get
\begin{eqnarray*}
&&\frac{\partial^2}{\partial t\partial
s}\psi_{\varepsilon,\delta}(L(x,v+su+ tu))\Bigm|_{s=0,t=0}\\
&&=\psi'_{\varepsilon,\delta}(L(x,v))\partial_{vv}L(x,v)[u,u]
+\psi''_{\varepsilon,\delta}(L(x,v))\bigl(\partial_vL(x,v)[u]\bigr)^2.
\end{eqnarray*}
Clearly, the left hand side is equal to zero at $v=0$.
 Analogously, one easily  computes
\begin{eqnarray*}
\frac{\partial^2}{\partial t\partial
s}\phi_{\mu,b}(|v+su+tu|^2_x)\Bigm|_{s=0,t=0} =
4\phi''_{\mu,b}(|v|_x^2)\langle v,u\rangle_x^2+
2\phi'_{\mu,b}(|v|_x^2)|u|_x^2.
\end{eqnarray*}
These lead to
\begin{eqnarray}\label{e:2.7}
\partial_{vv}
L^\star(x,v)[u,u]&=&\psi''_{\varepsilon,\delta}(L(x,v))\bigl(\partial_v
L(x,v)[u]\bigr)^2+
\psi'_{\varepsilon,\delta}(L(x,v))\partial_{vv}L(x,v)[u,u]\nonumber\\
&+& 4\phi''_{\mu,b}(|v|_x^2)\langle v,u\rangle_x^2+
2\phi'_{\mu,b}(|v|_x^2)|u|_x^2.
\end{eqnarray}

\noindent{$\bullet$} If $L(x,v)\le\delta$ (and hence
$|v|_x^2\le\delta$ by (\ref{e:2.3})), then
\begin{eqnarray}\label{e:2.8}
\partial_{vv}L^\star(x,v)[u,u]\ge 2\mu|u|_x^2
\end{eqnarray}
because Lemma~\ref{lem:2.1}(i) and Lemma~\ref{lem:2.1}(ii) imply:
$\psi''_{\varepsilon,\delta}(L(x,v))\ge 0$,
$\psi'_{\varepsilon,\delta}(L(x,v))\ge 0$, and
$\phi''_{\mu,b}(|v|_x^2)=0$ and $\phi'_{\mu,b}(|v|_x^2)=\mu$ for
$|v|_x^2\le\delta$.

\noindent{$\bullet$} If $L(x,v)\ge\delta$ (and hence
$|v|_x^2\ge\frac{\delta}{C_1}>\frac{\delta}{3C_1}$ by
(\ref{e:2.3})), then $\psi''_{\varepsilon,\delta}(L(x,v))=0$ and
$$
\psi'_{\varepsilon,\delta}(L(x,v))\partial_{vv}L(x,v)[u,u]=
\kappa\partial_{vv}L(x,v)[u,u]
$$
by Lemma~\ref{lem:2.1}(i). Hence it follows from (\ref{e:2.7}),
Lemma~\ref{lem:2.1}(ii) and (\ref{e:2.4}) that
\begin{eqnarray*}
\partial_{vv}L^\star(x,v)[u,u]&= & \kappa\partial_{vv}L(x,v)[u,u] +
4\phi''_{\mu,b}(|v|_x^2)\langle v,u\rangle_x^2+
2\phi'_{\mu,b}(|v|_x^2)|u|_x^2\\
&\ge & \kappa\partial_{vv}L(x,v)[u,u] +
4\phi''_{\mu,b}(|v|_x^2)\langle
v,u\rangle_x^2\nonumber\\
&\ge & 2\kappa\alpha_g|u|_x^2 +
4\phi''_{\mu,b}(|v|_x^2)|v|_x^2|u|_x^2\nonumber\\
&\ge & 2\kappa\alpha_g|u|_x^2 +
\frac{8c}{3C_1}\phi''_{\mu,b}(|v|_x^2)|u|_x^2\nonumber\\
&=&\left[2\kappa\alpha_g +
\frac{8c}{3C_1}\phi''_{\mu,b}(|v|_x^2)\right]|u|_x^2\nonumber
\end{eqnarray*}
because $\phi''_{\mu,b}\le 0$, and $\phi''_{\mu,b}(|v|_x^2)=0$ for
$|v|_x^2\ge\frac{2c}{3C_1}$. Observe that $\phi''_{\mu,b}(|v|_x^2)$
is bounded for $|v|_x^2\in [\frac{\delta}{3C_1}, \frac{2c}{3C_1}]$.
We may choose $\kappa>0$ so large that
$$
2\kappa\alpha_g + \frac{8c}{3C_1}\phi''_{\mu,b}(|v|_x^2)\ge
\frac{1}{2}\kappa\alpha_g.
$$
This and (\ref{e:2.8}) yield the expected  conclusion.
 \hfill$\Box$\vspace{2mm}

Proposition~\ref{prop:2.2}(iii) and thus the following
Corollary~\ref{cor:2.3}(iii) is only used in  the proof of
Theorem~\ref{th:1.5}(v) at the end of Section~\ref{sec:5.1}.
(i)-(iii) in Lemma~\ref{lem:2.1} and large $\kappa>0$ are sufficient
for other arguments.

If $M$ is compact, for any Riemannian metric $g$ on $M$ both
$\alpha_g$ and $\beta_g$ are positive numbers, and they may be
chosen so that  (\ref{e:2.3}) holds for some constant $C_1>0$. In
this case it is easily proved that $L^\star$ in (\ref{e:2.5}) is a
convex with quadratic growth Lagrangian. Defining
$L^\ast(x,v)=(L^\star(x,v)-\varrho_0)/\kappa$ we get

\begin{corollary}\label{cor:2.3}
Let $(M, F)$ be a compact $C^k$ Finsler manifold ($k\ge 2$), and
$L:=F^2$. Then for a given $c>0$ there exists a  $C^k$ convex with
quadratic growth Lagrangian $L^\ast:TM\to \R$  such that
\begin{description}

\item[(i)] $L^\ast(x,v)= L(x,v)\quad\hbox{if}\;L(x,v)\ge\frac{2c}{3C_1}$,

\item[(ii)] $L^\ast$ attains the minimum, and $L^\ast(x,v)=\min L^\ast\Longleftrightarrow
v=0$,

\item[(iii)]  $L^\ast(x,v)\le  L(x,v)$ for all
$(x,v)\in TM$,

 \item[(iv)] if $F$ is reversible, so is $L^\ast$.
\end{description}
\end{corollary}

 We can also
construct another  $C^k$ convex with quadratic growth Lagrangian
$L^\ast:TM\to \R$ satisfying Corollary~\ref{cor:2.3}(i)(ii)(iv) and
the condition $L^\ast(x,v)\ge  L(x,v)$ for all $(x,v)\in TM$. Such a
$L^\ast$ is not needed in this paper.

Let $L^\ast$ be as in Corollary~\ref{cor:2.3}.
 By the compactness of $M$, (\ref{e:2.3})-(\ref{e:2.4}) and the fact
that $L^\ast$ is convex with quadratic growth we deduce  that there
exist
 $\alpha^\ast_g>0$,
$\beta^\ast_g>0$, $C_j>0$, $j=2,3,4$, such that
\begin{eqnarray}
&&\alpha^\ast_g|u|_x^2\le\partial_{vv}
L^\ast(x,v)[u,u]\le\beta^\ast_g|u|_x^2,\label{e:2.10}\\
&&C_2|v|_x^2-C_4\le L^\ast(x,v)\le C_3(|v|_x^2+1).\label{e:2.11}
\end{eqnarray}
For each $\tau\in [0,1]$ we define $L^\tau:TM\to \R$ by
\begin{equation}\label{e:2.12}
L^\tau(x,v)=(1-\tau)L(x,v)+
 \tau  L^\ast(x,v).
\end{equation}
 Then it is only $C^k$ in $TM\setminus\{0\}$
for $0\le\tau<1$, and by (\ref{e:2.3})-(\ref{e:2.4}) and
(\ref{e:2.10})-(\ref{e:2.11})  one easily proves that in
$TM\setminus\{0\}$,
\begin{eqnarray}
&&\min\{\alpha_g,\alpha^\ast_g\}|u|_x^2\le\partial_{vv}
L^\tau(x,v)[u,u]\le\max\{\beta_g, \beta^\ast_g\}|u|_x^2,\label{e:2.13}\\
&&\min\{C_2,1\}|v|_x^2-C_4\le L^\tau(x,v)\le
(C_1+C_3)(|v|_x^2+1).\label{e:2.14}
\end{eqnarray}
From the mean value theorem it follows that
\begin{eqnarray*}
|\partial_vL^\tau(x,v)\cdot v|&=&|\partial_vL^\tau(x,v)\cdot
v-\partial_vL^\tau(x,0)
\cdot v|\nonumber\\
&=&\left|\int^1_0\partial_{vv}L^\tau(x,sv)[v,v]ds\right|\ge
\min\{\alpha_g,\alpha^\ast_g\}|v|_x^2
\end{eqnarray*}
for any $(x,v)\in TM\setminus\{0\}$, and therefore
\begin{eqnarray}\label{e:2.15}
\|\partial_vL^\tau(x,v)\|\ge\min\{\alpha_g,\alpha^\ast_g\}|v|_x\;\forall
(x,v)\in TM.
\end{eqnarray}

 Consider the Legendre transform
associated with $L^\tau$,
\begin{equation}\label{e:2.16}
{\bf L}^\tau:TM\to T^\ast M,\;(x,v)\mapsto \bigr(x,
\partial_v L^\tau(x,v)\bigl).
\end{equation}

\begin{proposition}\label{prop:2.4}
${\bf L}^\tau$ is a homeomorphism, and a $C^{k-1}$-diffeomorphism for
$\tau=1$. Moreover,  ${\bf L}^\tau$ also restricts to a
$C^{k-1}$-diffeomorphism from $TM\setminus\{0\}$ onto $T^\ast
M\setminus\{0\}$.
\end{proposition}

\noindent{\bf Proof}. The conclusion for $\tau=1$ is standard (see
\cite[Prop.2.1.6]{Fa}).  The case $\tau=0$ had been proved in  \cite[Prop.2.1]{CaJaMa3}. A similar proof yields the case
$0<\tau< 1$. That is, we only need to prove that they are proper
local homeomorphisms by Banach-Mazur theorem (cf.
\cite[Th.5.1.4]{Ber}).

 Since $L^\tau$ is $C^1$, fiberwise strictly convex and superlinear on
$TM$ the map
$$
{\bf L}^\tau_x:T_xM\to T_x^\ast M,\;v\mapsto\partial_vL^\tau(x,v)
$$
is a homeomorphism by Theorem 1.4.6 of \cite{Fa}. Moreover, $L^\tau$
is $C^k$ on $TM\setminus\{0\}$, and
$\partial_{vv}L^\tau(x,v)[u,u]>0$ for any $u,v\in
T_xM\setminus\{0\}$. From the implicit function theorem it follows
that the map ${\bf L}^\tau_x$ is a $C^{k-1}$ diffeomorphism from
$T_xM\setminus\{0_x\}$ to itself.

Now ${\bf L}^\tau:TM\to T^\ast M$ is a continuous bijection.
Consider its inverse
$$
({\bf L}^\tau)^{-1}:TM\to T^\ast M,\;(x,v)\mapsto \bigr(x, ({\bf
L}_x^\tau)^{-1}v\bigl).
$$
 By the positivity of $\partial_{vv}L^\tau$ on $TM\setminus\{0\}$, from the
inverse function theorem we derive that ${\bf L}^\tau$ is locally a
$C^{k-1}$ diffeomorphism on $TM\setminus\{0\}$ and $({\bf L}^\tau)^{-1}$
maps $T^\ast M\setminus\{0\}$ onto $TM\setminus\{0\}$.

 We claim that the continuity of $({\bf L}^\tau)^{-1}$ extends up to the
zero section. Suppose that  $(x_n, w_n)\to (x_0, 0)$ and $({\bf
L}^\tau)^{-1}(x_n, w_n)=(x_n, v_n)$. Then ${\bf
L}^\tau(x_n,v_n)=(x_n, w_n)$ or $w_n=\partial_vL^\tau(x_n,v_n)$. By
(\ref{e:2.15}) we deduce that $|v_n|_{x_n}\to 0$. This leads to the
desired claim.
 Hence
${\bf L}^\tau$ is a homeomorphism from $TM$ onto $T^\ast M$. The
second conclusion is a direct consequence of this fact and the
inverse function theorem. $\Box$\vspace{2mm}

\begin{remark}\label{rm:2.5}
{\rm For a $C^k$ convex with quadratic growth Lagrangian $L^\ast:TM\to\R$
($k\ge 2$), by the convexity of $L^\ast$ we have $\partial_v
L^\ast(x,0)=0\;\forall x$ and therefore (\ref{e:2.15}) holds, it
follows that Proposition~\ref{prop:2.4} is still true.}
\end{remark}

\section{The stability of critical groups}\label{sec:3}
\setcounter{equation}{0}

In this section we assume that $(M, F)$ is a compact $C^k$ Finsler
manifold ($k\ge 5$) and $L:=F^2$. For a $C^k$ convex with quadratic
growth Lagrangian $L^\ast:TM\to \R$ and $\tau\in [0,1]$ we define
$L^\tau(x,v)=(1-\tau)L(x,v)+
 \tau  L^\ast(x,v)$ and
\begin{equation}\label{e:3.1}
{\cal
L}^\tau(\gamma)=\int^1_0L^\tau(\gamma(t),\dot\gamma(t))dt\quad\forall\gamma\in\Lambda_N(M).
\end{equation}
Then ${\cal L}^\tau(\gamma)=(1-\tau){\cal L}(\gamma)+ \tau{\cal
L}^\ast(\gamma)$ for $\tau\in [0,1]$, where the functionals ${\cal
L}$ and ${\cal L}^\ast$ are given by (\ref{e:1.2}) and
(\ref{e:1.26}), respectively. Since  ${\cal L}$ and ${\cal L}^\ast$
on $\Lambda_N(M)$ are $C^{2-0}$ and satisfy the Palais-Smale
condition by Proposition~\ref{prop:1.1} and \cite{AbSc1}, so is each
functional ${\cal L}^\tau$  on $\Lambda_N(M)$.

 The {\bf energy function} of $L^\tau$, $E^\tau:TM\to\R$ is
defined by
\begin{equation}\label{e:3.2}
E^\tau(x,v)= \partial_v L^\tau(x,v)\cdot v-L^\tau(x,v).
\end{equation}
(It is $C^{1-0}$ because  $F^2$ is $C^{2-0}$ on $TM$).

\begin{proposition}\label{prop:3.1}
For the Lagrangian $L^\ast$ in Corollary~\ref{cor:2.3}, if
$\gamma\in\Lambda_N(M)$ is a critical point of the functional ${\cal
L}^\tau$ which is not a constant curve, then $\gamma$ is a $C^k$
regular curve, i.e. $\dot\gamma(t)\ne 0$ for any $t\in I$.
\end{proposition}

\noindent{\bf Proof}. Fix a point $\gamma(\bar t)$ and choose a
coordinate chart around $\gamma(\bar t)$ on $M$, $(V, \chi)$,
$$
\chi:V\to\chi(V)\subset\R^n,\;x\mapsto\chi(x)=(x_1,\cdots,x_n).
$$
Then we get an induced chart on $TM$, $(\pi^{-1}(V), T\chi)$,
$$
T\chi:\pi^{-1}(V)\to \chi(V)\times\R^n,\;(x,v)\mapsto
(x_1,\cdots,x_n;v_1,\cdots,v_n).
$$
Let $I_0$ be a connected component of $\gamma^{-1}(V)$. It has one
of the following three forms: $[0, a)$, $(a,b)$, $(b,1]$. Let
$\zeta(t):=\chi(\gamma(t))$ for $t\in I_0$. Then $\zeta:I_0\to
V\subset\R^n$ is absolutely continuous and
$\dot{\zeta}(t)=d\chi(\gamma(t))(\dot\gamma(t))$ for $t\in I_0$. Set
$$
\tilde L^\tau(x, y):=L^\tau(\chi^{-1}(x), d\chi^{-1}(x)(y))
\quad\forall (x,y)\in\chi(V)\times\R^n.
$$
Since $d{\cal L}^\tau(\gamma)(\xi)=0\;\forall\xi\in T_\gamma\Lambda_N(M)$, we deduce that
$$
\int_{I_0}\Bigl(\partial_x\tilde
L^\tau(\zeta(t),\dot{\zeta}(t))[\xi(t)]+
\partial_y\tilde
L^\tau(\zeta(t),\dot{\zeta}(t))[\dot\xi(t)]\Bigr)dt=0
$$
for any $\xi\in C^\infty_0({\rm Int}(I_0), \R^n)$. Denote by $t_0$
the left end point of $I_0$ and by
$$
H(t)=-\int^t_{t_0}\partial_x\tilde
L^\tau(\zeta(s),\dot{\zeta}(s))ds\quad\forall t\in I_0.
$$
Then $I_0\ni t\mapsto H(t)\in\R^n$ is continuous and
$$
\int_{I_0}\Bigl(H(t)+
\partial_y\tilde
L^\tau(\zeta(t),\dot{\zeta}(t))\Bigr)[\dot\xi(t)]dt=0
$$
for any $\xi\in C^\infty(I_0, \R^n)$ with ${\rm
supp}(\xi)\subset{\rm Int}(I_0)$. It follows that there exists a
constant vector ${\bf u}\in\R^n$ such that $H(t)+
\partial_y\tilde
L^\tau(\zeta(t),\dot{\zeta}(t))={\bf u}$ a.e. on $I_0$. This implies
the map $I_0\ni t\mapsto \partial_y\tilde
L^\tau(\zeta(t),\dot{\zeta}(t))$ is continuous. Moreover, as in
Proposition~\ref{prop:2.4} we can prove that the map
$$
\widetilde{\bf L}^\tau:\chi(V)\times\R^n\to
\chi(V)\times\R^n,\;(x,y)\mapsto \bigr(x, \partial_y\tilde
L^\tau(x,y)\bigl)
$$
is a homeomorphism, and also restricts to a $C^{k-1}$-diffeomorphism
$\widetilde{\bf L}^\tau_0$ from $\chi(V)\times(\R^n\setminus\{0\})$
to $\chi(V)\times(\R^n\setminus\{0\})$.
 So
$$
I_0\ni t\mapsto (\zeta(t),\dot{\zeta}(t))=(\widetilde{\bf
L}^\tau)^{-1}\circ\widetilde{\bf L}^\tau(\zeta(t),\dot{\zeta}(t))=
(\widetilde{\bf L}^\tau)^{-1}\bigl(\zeta(t),\partial_y\tilde
L^\tau(\zeta(t),\dot{\zeta}(t))\bigl)
$$
is continuous. This shows that $\zeta$ is $C^1$. Observe that
$$
\frac{d}{dt}\partial_y\tilde
L^\tau(\zeta(t),\dot{\zeta}(t))=\partial_x\tilde
L^\tau(\zeta(t),\dot{\zeta}(t))\quad\hbox{a.e. on}\;I_0.
$$
We get that the map $I_0\ni t\mapsto\partial_y\tilde
L^\tau(\zeta(t),\dot{\zeta}(t))$ is $C^1$, and hence the composition
$$
\{t\in I_0\,|\,\dot\gamma(t)\ne 0\}\ni t\mapsto\partial_y\tilde
L^\tau(\zeta(t),\dot{\zeta}(t))\stackrel{(\widetilde{\bf
L}^\tau)^{-1}}\longrightarrow\bigl(\zeta(t),\dot{\zeta}(t)\bigr)
$$
is $C^1$.  Summarizing, we
have proved that
{\bf $\gamma$ is $C^1$, and $C^2$ in $\{t\in
I\,|\, \dot\gamma(t)\ne 0\}$.}

Now since $\gamma$ is not constant  the energy
$E^\tau(\gamma,\dot\gamma)= \partial_v
L^\tau(\gamma,\dot\gamma)\cdot \dot\gamma-L^\tau(\gamma,\dot\gamma)$
is constant on every connected component of $\{t\in I\,|\,
\dot\gamma(t)\ne 0\}$.  Corollary~\ref{cor:2.3} implies
\begin{equation}\label{e:3.3}
L^\tau\ge \tau\min L^\ast,\quad L^\tau(x,0)=\tau\min
L^\ast\quad\hbox{and}\quad E^\tau(x, 0)=-\tau \min L^\ast\;\forall
x.
\end{equation}
 For $v\in
T_xM\setminus\{0\}$ and $t>0$, by (\ref{e:3.2}) we have
\begin{eqnarray*}
\frac{d}{dt}E^\tau(x, tv)&=&\frac{d}{dt}\bigl( \partial_v
L^\tau(x, tv)\cdot (tv)\bigl)-\frac{d}{dt}L^\tau(x, tv)\\
&=&\partial_{vv} L^\tau(x, tv)[v, tv]+\partial_v L^\tau(x,
tv)\cdot v-\partial_v L^\tau(x, tv)\cdot v\\
&=&t\partial_{vv} L^\tau(x, tv)[v, v]>0.
\end{eqnarray*}
It follows that $E^\tau(x, v)>E^\tau(x, 0)=-\tau\min L^\ast$ on
$TM\setminus\{0\}$ and hence
$$
E^\tau(x, v)\ge -\tau\min L^\ast\quad\hbox{and}\quad E^\tau(x,
v)=-\tau\min L^\ast\Longleftrightarrow v=0.
$$
Since $E^\tau(\gamma,\dot\gamma)$ is strictly larger than $-\tau\min
L^\ast$ in $I\setminus\{t\in I\,|\, \dot\gamma(t)\ne 0\}$ and $I\ni
t\to E^\tau(\gamma(t),\dot\gamma(t))$ is continuous we must have
 $I=\{t\in
I\,|\, \dot\gamma(t)\ne 0\}$. It easily follows that $\gamma$ is
$C^k$.  $\Box$\vspace{2mm}

\begin{remark}\label{rm:3.2}
{\rm For a $C^k$ convex with quadratic growth Lagrangian $L^\ast:TM\to\R$
($k\ge 2$), if  $L^\ast(x,0)=\min L^\ast\;\forall x$ then
 (\ref{e:2.15}) and therefore
Proposition~\ref{prop:2.4} hold. Moreover we have also (\ref{e:3.3})
and hence Proposition~\ref{prop:3.1} for $k\ge 5$.}
\end{remark}

  Since the weak slope of a $C^1$ functional $f$ on an open subset of a normed space
  is equal to the norm of the differential of $f$ the lower critical
  point (resp. value) of $f$ agrees with the usual one of $f$.
  The following special version of \cite[Th.1.5]{CiDe}
about the stability property of critical groups
  is convenient for us.

\begin{theorem}\label{th:3.3}{\rm (\cite[Theorem~1.5]{CiDe})}
Let $\{f_\tau:\,\tau\in [0, 1]\}$ be a family of $C^1$ functionals
from a Banach space $X$ to $\R$, $U$ an open subset of $X$ and $[0,
1]\ni \tau\mapsto u_\tau\in U$ a continuous map such that $u_0$ and
$u_1$ are critical points, respectively, of $f_0$ and $f_1$. Assume:
\begin{description}
\item[(I)] if $\tau_k\to \tau$ in $[0, 1]$, then $f_{\tau_k}\to f_\tau$
uniformly on $\overline{U}$;

\item[(II)]  for every sequence
$\tau_k\to \tau$ in $[0, 1]$ and $(v_k)$ in $\overline{U}$ with
$f'_{\tau_k}(v_k)\to 0$ and $(f_{\tau_k}(v_k))$ bounded, there
exists a subsequence $(v_{k_j})$ convergent to some $v$ with
$f'_\tau(v)=0$;

\item[(III)] $f'_\tau(u)\ne 0$ for every $\tau\in [0, 1]$ and
$u\in \overline{U}\setminus\{u_\tau\}$

\end{description}
Then $C_q(f_0, u_0;{\K})\cong C_q(f_1, u_1;{\K})$ for every
$q\ge 0$.
\end{theorem}

Actually,   Chang \cite[page 53, Theorem~5.6]{Ch93}, Chang and Ghoussoub
\cite[Theorem~III.4]{ChGh} and Corvellec and Hantoute
\cite[Theorem~5.2]{CorHa}  are sufficient for the proof of our
Theorem~\ref{th:3.8}.

For ${\gamma_0},\tau\in [0, 1]$, since
$L^{\gamma_0}-L^\tau=(\tau-{\gamma_0})L+
 ({\gamma_0}-\tau)L^\ast$ we have
\begin{equation}\label{e:3.4}
|{\cal L}^{\gamma_0}(\gamma)-{\cal L}^\tau(\gamma)|\le
|{\gamma_0}-\tau|\int^1_0\Bigl[\bigl|L(\gamma(t),\dot\gamma(t))\bigr|+
\bigl|L^\ast(\gamma(t),\dot\gamma(t))\bigr|\Bigr]dt
\end{equation}
for all $\gamma\in\Lambda_N(M)$. Note that the condition (L2) in \S1.5 and
the compactness of $M$ imply
\begin{equation}\label{e:3.5}
|L^\ast(x,v)|\le C_2(1+|v|_x^2) \quad\forall (x,v)\in TM
\end{equation}
for some constant $C_2>0$. Let $C_3:=\max\{L(x,v)\,|\, |v|_x=1\}$.
Then $|L(x,v)|\le C_3|v|_x^2$ for all $(x,v)\in TM$. From this and
(\ref{e:3.4})-(\ref{e:3.5}) we immediately get

\begin{claim}\label{cl:3.4}
For any bounded subset $K\subset \Lambda_N(M)$, ${\cal
L}^{\gamma_0}\to{\cal L}^\tau$ uniformly on $K$ as
${\gamma_0}\to\tau$.
\end{claim}

\begin{claim}\label{cl:3.5}
For any $\gamma_0\in\Lambda_N(M)$ there exists a neighborhood ${\cal
U}$ of it such that for every sequence $\tau_k\to \tau$ in $[0, 1]$
and a sequence $(\gamma_k)$ in ${\cal U}$ with $d{\cal
L}^{\tau_k}(\gamma_k)\to 0$ and $({\cal L}^{\tau_k}(\gamma_k))$
bounded, there exists a subsequence $(\gamma_{k_j})$ convergent to
some $\gamma$ with $d{\cal L}^\tau(\gamma)=0$. (Clearly, $\gamma_0$
can be replaced by a compact subset $K\subset\Lambda_N(M)$.)
\end{claim}

\noindent{\bf Proof}. Clearly, this result is of a local nature. By
a well-known localization argument as in \cite{AbSc1, Lu0} (cf.
Section~\ref{sec:4}) the question is reduced to the following case:\\
$\bullet$ $F,  L^\ast:B^n_{2\rho}(0)\times\R^n\to\R$;\\
$\bullet$ $(\gamma_k)\subset W^{1,2}_V([0,1],B^n_{\rho}(0)):=\{ \gamma\in
W^{1,2}([0,1], B^n_{\rho}(0))\,|\, (\gamma(0),\gamma(1))\in
V\}$ is bounded, where $V$ is a linear subspace of $\R^n\times\R^n$,
such that
\begin{equation}\label{e:3.6}
d{\cal L}^{\tau_k}(\gamma_k)\to 0\;\hbox{(as
$k\to\infty$)}\quad\hbox{and}\quad |{\cal L}^{\tau_k}(\gamma_k)|\le
C_4\;\forall k
\end{equation}
for some constant $C_4$ and that both $(\|d{\cal L}(\gamma_k)\|)$ and
$(\|d{\cal L}^\ast(\gamma_k)\|)$ are bounded;\\
$\bullet$ The conclusion is that there exists a subsequence
$(\gamma_{k_j})$ convergent to some $\gamma\in
W^{1,2}_V([0,1],B^n_{2\rho}(0))$ with $d{\cal L}^\tau(\gamma)=0$. In
this time $\gamma$ satisfies
\begin{eqnarray*}
&&\frac{d}{dt}\bigl(\partial_vL^\tau(\gamma(t),\dot\gamma(t))\bigr)-
\partial_xL^\tau(\gamma(t),\dot\gamma(t))=0,\\
&&\partial_vL^\tau(\gamma(0),\dot\gamma(0))\cdot
v_0=\partial_vL^\tau(\gamma(1),\dot\gamma(1))\cdot v_1\quad\forall
(v_0, v_1)\in V.
\end{eqnarray*}

 In order to prove this, note that (\ref{e:3.6})
implies
\begin{eqnarray*}
\|d{\cal L}^{\tau}(\gamma_k)\|&\le& \|d{\cal
L}^{\tau_k}(\gamma_k)\|+
\|d{\cal L}^{\tau}(\gamma_k)-d{\cal L}^{\tau_k}(\gamma_k)\|\\
&\le &\|d{\cal L}^{\tau_k}(\gamma_k)\|+ |\tau_k-\tau|\cdot
\bigl(\|d{\cal L}(\gamma_k)\|+\|d{\cal L}^\ast(\gamma_k)\|\bigr).
\end{eqnarray*}
 We obtain that $\|d{\cal L}^{\tau}(\gamma_k)\|\to 0$ as $k\to\infty$.
Moreover Claim~\ref{cl:3.4} implies that $(|{\cal
L}^{\tau}(\gamma_k)|)$ is bounded too. Combining the proof of
\cite[Th.3.1]{CaJaMa3} and that of \cite[Prop.2.5]{AbSc1} we can
complete the proof of Claim~\ref{cl:3.5}. They are omitted.
$\Box$\vspace{2mm}

\begin{claim}\label{cl:3.6}
For any $C>0$ there exists a $C'>0$, which is independent of
$\tau\in [0,1]$ and $(x,v)\in TM$, such that $E^\tau(x,v)\le
C\Longrightarrow |v|_x\le C'$.
\end{claim}

\noindent{\bf Proof}. Since
$L^\ast(x,v)=(\psi_{\varepsilon,\delta}(L(x,v))+
\phi_{\mu,b}(|v|^2_x)-b-\varrho_0)/\kappa$, we have
\begin{eqnarray*}
E^\ast(x,v)&=&\partial_v  L^\ast(x,v)\cdot v-
L^\ast(x,v)\\
&=&\frac{1}{\kappa}\bigl[\psi'_{\varepsilon,\delta}(L(x,v))\partial_v
L(x,v)\cdot v +2\phi'_{\mu,b}(|v|_x^2)|v|_x^2\\
&&\hspace{20mm}-\psi_{\varepsilon,\delta}(L(x,v))-
\phi_{\mu,b}(|v|^2_x)+b+\varrho_0\bigr].
\end{eqnarray*}
Moreover, $E(x,v)=\partial_v L(x,v)\cdot v- L(x,v)=L(x,v)$. Hence
\begin{eqnarray*}
&& E^\tau(x,v)=(1-\tau)L(x,v)+
\frac{\tau}{\kappa}\bigl[b+\varrho_0+2\psi'_{\varepsilon,\delta}(L(x,v))
L(x,v)\\
&&\hspace{20mm} -\psi_{\varepsilon,\delta}(L(x,v)) +
2\phi'_{\mu,b}(|v|_x^2)|v|_x^2-\phi_{\mu,b}(|v|_x^2)\bigr].
\end{eqnarray*}
Suppose that there exist sequences $(\tau_k)\subset [0,1]$ with
$\tau_k\to\tau$ and $(x_k, v_k)\subset TM$ with $x_k\to x_0$, such
that $E^{\tau_k}(x_k, v_k)\le C\;\forall k$ and
$|v_k|_{x_k}\to\infty$. Since (\ref{e:2.3}) implies
$L(x_k,v_k)\to\infty$ we deduce
\begin{eqnarray*}
&&\psi'_{\varepsilon,\delta}(L(x_k,v_k))=\kappa ,\quad
\psi_{\varepsilon,\delta}(L(x_k,v_k))=\kappa L(x_k,v_k)+
\varrho_0,\\
&&\phi'_{\mu,b}(|v_k|_{x_k}^2)=0,\quad \phi_{\mu,b}(|v_k|_{x_k}^2)=b
\end{eqnarray*}
for $k$ sufficiently large. It follows that
\begin{eqnarray*}
 E^{\tau_k}(x_k,v_k)&=&(1-\tau_k)L(x_k,v_k)+ \frac{\tau_k}{\kappa}\bigl[2\kappa
L(x_k,v_k)- \kappa L(x_k,v_k)]\\
&=&L(x_k,v_k)\to\infty
\end{eqnarray*}
as $k\to\infty$. This contradiction yields the desired claim.
$\Box$\vspace{2mm}

\begin{claim}\label{cl:3.7}
Let $\gamma_0\in\Lambda_N(M)$ be an isolated nonconstant critical
point of ${\cal L}$ on $\Lambda_N(M)$ (and hence  a $C^k$
nonconstant $F$-geodesics with constant speed
$F(\gamma_0(t),\dot\gamma_0(t))\equiv \sqrt{c}>0$). Let
$L^\ast:TM\to\R$ be given by Corollary~\ref{cor:2.3}, and let
$$
L^\tau(x,v)=(1-\tau)L(x,v)+
 \tau L^\ast(x,v)\quad\forall\tau\in [0,1].
 $$
Then there exists a neighborhood of $\gamma_0$ in $\Lambda_N(M)$,
${\cal U}(\gamma_0)$, such that each ${\cal L}^\tau$ has only the
critical point $\gamma_0$ in ${\cal U}(\gamma_0)$.
\end{claim}

\noindent{\bf Proof}. Clearly, if $L(x,v)\ge\frac{2c}{3}$  then
 \begin{equation}\label{e:3.7}
L^\tau(x,v)=L(x,v).
 \end{equation}
 Since $L(\gamma_0(t), \dot\gamma_0(t))=c\;\forall t$, it is easily
checked that $\gamma_0$ is a critical point of the functionals
${\cal L}^\tau$ in (\ref{e:3.1}) on $\Lambda_N(M)$ and
$L^\tau(\gamma_0(t), \dot\gamma_0(t))=c\;\forall t$.

By a contradiction, suppose that there exist sequences
$(\gamma_k)\subset\Lambda_N(M)$, $(\tau_k)\subset [0,1]$ such that
\begin{equation}\label{e:3.8}
\tau_k\to\tau_0,\quad\gamma_k\to\gamma_0,\quad d{\cal
L}^{\tau_k}(\gamma_k)=0\;\forall k.
\end{equation}
Then $\gamma_k$ is nonconstant for each large $k$, and therefore a
$C^k$ regular curve by Proposition~\ref{prop:3.1} (removing finitely
many terms if necessary). Note that
$E^{\tau_0}(\gamma_0(t),\dot{\gamma}_0(t))$ and
 $E^{\tau_k}(\gamma_k(t),\dot{\gamma}_k(t))$ are constants independent
 of $t$. Set $d_0:\equiv E^{\tau_0}(\gamma_0(t),\dot{\gamma}_0(t))$ and
$d_k:\equiv E^{\tau_k}(\gamma_k(t),\dot{\gamma}_k(t))$ for
$k=1,2,\cdots$. Then
\begin{eqnarray*}
d_0&=&\int^1_0 E^{\tau_0}(\gamma_0(t),\dot{\gamma}_0(t))dt\\
&=&\int^1_0 \partial_v
L^{\tau_0}(\gamma_0(t),\dot{\gamma}_0(t))[\dot{\gamma}_0(t)]dt-\int^1_0
L^{\tau_0}(\gamma_0(t),\dot{\gamma}_0(t))dt\\
&=&(1-\tau_0)\left[\int^1_0 \partial_v
L(\gamma_0(t),\dot{\gamma}_0(t))[\dot{\gamma}_0(t)]dt-\int^1_0
L(\gamma_0(t),\dot{\gamma}_0(t))dt\right]\\
&+&\tau_0\left[\int^1_0 \partial_v
L^\ast(\gamma_0(t),\dot{\gamma}_0(t))[\dot{\gamma}_0(t)]dt-\int^1_0
L^\ast(\gamma_0(t),\dot{\gamma}_0(t))dt\right]\\
&=&(1-\tau_0)\int^1_0
L(\gamma_0(t),\dot{\gamma}_0(t))dt\\
&+&\tau_0\left[\int^1_0 \partial_v
L^\ast(\gamma_0(t),\dot{\gamma}_0(t))[\dot{\gamma}_0(t)]dt-\int^1_0
L^\ast(\gamma_0(t),\dot{\gamma}_0(t))dt\right]
\end{eqnarray*}
since $\partial_vL(x,v)[v]=2L(x,v)$. Similarly, we have
\begin{eqnarray*}
d_k&=&\int^1_0 E^{\tau_k}(\gamma_k(t),\dot{\gamma}_k(t))dt\\
&=&\int^1_0 \partial_v
L^{\tau_k}(\gamma_k(t),\dot{\gamma}_k(t))[\dot{\gamma}_k(t)]dt-\int^1_0
L^{\tau_k}(\gamma_k(t),\dot{\gamma}_k(t))dt\\
&=&(1-\tau_k)\int^1_0
L(\gamma_k(t),\dot{\gamma}_k(t))dt\\
&+&\tau_k\left[\int^1_0 \partial_v
L^{\ast}(\gamma_k(t),\dot{\gamma}_k(t))[\dot{\gamma}_k(t)]dt-\int^1_0
L^{\ast}(\gamma_k(t),\dot{\gamma}_k(t))dt\right].
\end{eqnarray*}
Recall that $L^\ast$ is convex with quadratic growth. There exists a
constant $C^\ast>0$ such that $|\partial_vL^\ast(x,v)[v]|\le
C^\ast(1+ |v|^2_x)$ for all $(x,v)\in TM$. From the first two
relations in (\ref{e:3.8}) and a theorem of Krasnosel'skii we deduce
\begin{eqnarray*}
&&\int^1_0 \partial_v
L^{\ast}(\gamma_k(t),\dot{\gamma}_k(t))[\dot{\gamma}_k(t)]dt\to
\int^1_0
\partial_v L^\ast(\gamma_0(t),\dot{\gamma}_0(t))[\dot{\gamma}_0(t)]dt,\\
&&\int^1_0 L^{\ast}(\gamma_k(t),\dot{\gamma}_k(t))dt\to\int^1_0
L^{\ast}(\gamma_0(t),\dot{\gamma}_0(t))dt
\end{eqnarray*}
and hence $d_k\to d_0=c$. Choose $k_0\in\N$ such that $|d_k|<c+1$
for all $k\ge k_0$. Then by Claim~\ref{cl:3.6} we have a constant
$C_5>0$ such that
\begin{equation}\label{e:3.9}
|\dot{\gamma}_k(t)|_{\gamma_k(t)}\le C_5\quad\forall t\in
[0,1]\;\hbox{and}\;k\ge k_0.
\end{equation}

Since $M$ is compact and $L^\ast$ satisfies
  the assumptions (L1)-(L2) in \S1.5,
we may take finite many coordinate charts on $M$,
$$
\varphi_\alpha: V_\alpha\to B^n_{2\rho}(0),\;x\mapsto
(x^\alpha_1,\cdots,x^\alpha_n), \quad \alpha=1,\cdots,m,
$$
and positive constants $C_6>C_7$, such that
$M=\cup^m_{\alpha=1}(\varphi_\alpha)^{-1}(B^n_{\rho}(0))$ and each
$$
L^{\tau}_\alpha: B^n_{2\rho}(0)\times\R^n\to \R,\;(x^\alpha,
v^\alpha)\mapsto L^\tau\bigr(\varphi_\alpha^{-1}(x^\alpha),
d\varphi_\alpha^{-1}(x^\alpha)(v^\alpha)\bigl)
$$
satisfies
\begin{eqnarray}
&&|L^\tau_\alpha(x^\alpha,v^\alpha)|\le C_6(1+ |v^\alpha|^2),\nonumber\\
&&\Bigl| \frac{\partial L^\tau_\alpha}{\partial
x^\alpha_i}(x^\alpha,v^\alpha)\Bigr|\le C_6(1+
|v^\alpha|^2),\quad\Bigl| \frac{\partial L^\tau_\alpha}{\partial
v^\alpha_i}(x^\alpha,v^\alpha)\Bigr|\le C_6(1+ |v^\alpha|),\label{e:3.10}\\
&&\Bigl| \frac{\partial^2 L^\tau_\alpha}{\partial x^\alpha_i\partial
x^\alpha_j}(x^\alpha,v^\alpha)\Bigr|\le C_6(1+ |v^\alpha|^2),\quad
\Bigl| \frac{\partial^2 L^\tau_\alpha}{\partial x^\alpha_i\partial
v^\alpha_j}(x^\alpha,v^\alpha)\Bigr|
\le C_6(1+ |v^\alpha|),\quad\label{e:3.11}\\
&&\Bigl| \frac{\partial^2 L^\tau_\alpha}{\partial v^\alpha_i\partial
v^\alpha_j}(x^\alpha,v^\alpha)\Bigr|\le C_6\quad\hbox{and}\quad
\sum_{ij}\frac{\partial^2 L^\tau_\alpha}{\partial v^\alpha_i\partial
v^\alpha_j}(x^\alpha,v^\alpha)u_iu_j\ge C_7|{\bf u}|^2\label{e:3.12}
\end{eqnarray}
for  $\tau\in [0, 1]$, $(x^\alpha, v^\alpha)\in \bar
B^n_\rho(0)\times\R^n$ and all ${\bf u}=(u_1,\cdots,u_n)\in\R^n$.
Moreover we have also positive constants $C_8, C_9$ and $C_{10}$
such that
\begin{eqnarray}
&&L^\tau(x,v)\ge C_8(|v|_q^2-1)\quad\forall (x,v)\in TM,\label{e:3.13}\\
&&C_9|v^\alpha|\le |v|_x\le C_{10}|v^\alpha|\quad\forall v\in
T_xM\;\hbox{and}\;x\in \varphi_\alpha^{-1}(\bar
B^n_\rho(0)),\label{e:3.14}
\end{eqnarray}
where $v=\sum^n_{i=1}v^\alpha_i\frac{\partial}{\partial
x^\alpha_i}\bigm|_x$.

For each $k\ge k_0$, set
$I_{k,\alpha}:=\gamma_k^{-1}(\varphi_\alpha^{-1}(B^n_{\rho}(0)))$
and $\gamma_{k,\alpha}:=\varphi_\alpha\circ\gamma_k:I_{k,\alpha}\to
B^n_{\rho}(0)$. Then the third condition in (\ref{e:3.8}) implies
that
\begin{eqnarray*}
&&\frac{d}{dt}\partial_{v^\alpha}L^{\tau_k}_\alpha(\gamma_{k,\alpha}(t),\dot{\gamma}_{k,\alpha}(t))
=\partial_{x^\alpha}L^{\tau_k}_\alpha(\gamma_{k,\alpha}(t),\dot{\gamma}_{k,\alpha}(t))
\quad\hbox{or}\nonumber\\
&&\ddot{\gamma}_{k,\alpha}(t)=\left(\partial_{v^\alpha
v^\alpha}L^{\tau_k}_\alpha(\gamma_{k,\alpha}(t),\dot{\gamma}_{k,\alpha}(t))\right)^{-1}\Bigl[
\partial_{x^\alpha}L^{\tau_k}_\alpha(\gamma_{k,\alpha}(t),\dot{\gamma}_{k,\alpha}(t))-\nonumber\\
&&\qquad\qquad-\partial_{x^\alpha
v^\alpha}L^{\tau_k}_\alpha(\gamma_{k,\alpha}(t),\dot{\gamma}_{k,\alpha}(t))[\dot{\gamma}_{k,\alpha}(t)]\Bigr]\;
\forall t\in I_{k,\alpha}.
\end{eqnarray*}
It follows that
\begin{eqnarray*}
|\ddot{\gamma}_{k,\alpha}(t)|^2&=&
\ddot{\gamma}_{k,\alpha}(t)\bigl(\ddot{\gamma}_{k,\alpha}(t)\bigr)^T\\
&=&\left(\partial_{v^\alpha
v^\alpha}L^{\tau_k}_\alpha\right)^{-1}\Bigl[
\partial_{x^\alpha}L^{\tau_k}_\alpha-
\partial_{x^\alpha v^\alpha}L^{\tau_k}_\alpha[\dot{\gamma}_{k,\alpha}(t)]\Bigr]
\bigl(\ddot{\gamma}_{k,\alpha}(t)\bigr)^T\\
&\le&\frac{1}{C_7}\Bigl|
\partial_{x^\alpha}L^{\tau_k}_\alpha-
\partial_{x^\alpha
v^\alpha}L^{\tau_k}_\alpha[\dot{\gamma}_{k,\alpha}(t)]\Bigr|\cdot
|\ddot{\gamma}_{k,\alpha}(t)|
\end{eqnarray*}
because the second inequality of (\ref{e:3.12}) implies that
$|\left(\partial_{v^\alpha
v^\alpha}L^{\tau_k}_\alpha\right)^{-1}v|\le |v|/C_7$. Hence
\begin{eqnarray*}
|\ddot{\gamma}_{k,\alpha}(t)| &\le&\frac{1}{C_7}\Bigl|
\partial_{x^\alpha}L^{\tau_k}_\alpha-
\partial_{x^\alpha
v^\alpha}L^{\tau_k}_\alpha[\dot{\gamma}_{k,\alpha}(t)]\Bigr|\\
&\le& \frac{C_6}{C_7}(1+ |\dot{\gamma}_{k,\alpha}(t)|^2)+
\frac{C_6}{C_7}(1+
|\dot{\gamma}_{k,\alpha}(t)|)|\dot{\gamma}_{k,\alpha}(t)|\\
&\le&  3\frac{C_6}{C_7}(1+ |\dot{\gamma}_{k,\alpha}(t)|^2)\\
&\le&  3\frac{C_6}{C_7}(1+
|\dot{\gamma}_{k}(t)|_{\gamma_k(t)}^2/C_9^2) \le 3\frac{C_6}{C_7}(1+
C_5^2/C_9^2)
\end{eqnarray*}
by (\ref{e:3.10})-(\ref{e:3.12}), (\ref{e:3.14}) and (\ref{e:3.9}).

Note that $I=\cup^m_{\alpha=1}I_{k,\alpha}$ for each $k\ge k_0$. We
may deduce that the sequence $(\gamma_k)$ is bounded in
$C^2([0,1],M)$. Passing to a subsequence we may assume that
$(\gamma_k)$ converges to $\bar\gamma\in C^1([0,1],M)$ in
$C^1([0,1],M)$. Since the sequence $(\gamma_k)$ converges to
$\gamma_0$ in $C^0([0,1],M)$, $\bar\gamma=\gamma_0$. That is,
$(\gamma_k)$ converges to $\gamma_0$ in $C^1([0,1],M)$. It follows
that $\min_{t\in
[0,1]}L(\gamma_k(t),\dot{\gamma}_k(t))>\frac{2c}{3}$ for
sufficiently large $k$. This and (\ref{e:3.7})-(\ref{e:3.8}) lead to
$d{\cal L}(\gamma_k)=0$ for sufficiently large $k$, which
contradicts to the assumption that $\gamma_0\in\Lambda_N(M)$ is an
isolated nonconstant critical point of ${\cal L}$ on $\Lambda_N(M)$.
\hfill$\Box$\vspace{2mm}

Using Theorem~\ref{th:3.3} and Claims~\ref{cl:3.4},~\ref{cl:3.5} and
~\ref{cl:3.7} we immediately get the following key result for the
proof of Theorem~\ref{th:1.2}(ii).

\begin{theorem}\label{th:3.8}
Let $\gamma_0\in\Lambda_N(M)$ be an isolated nonconstant critical
point of ${\cal L}$ on $\Lambda_N(M)$ (and hence  a $C^k$
nonconstant $F$-geodesics with constant speed
$F(\gamma_0(t),\dot\gamma_0(t))\equiv \sqrt{c}>0$). Let
$L^\ast:TM\to\R$ be given by Corollary~\ref{cor:2.3}. Then
$\gamma_0$ is a uniformly isolated critical point of the family of
functionals $\{\mathcal{L}^\tau=(1-\tau)\mathcal{L}+\tau{\cal
L}^\ast\,|\,0\le\tau\le 1\}$ too, and $C_\ast({\cal
L}^\tau,\gamma_0;{\K})$ is independent of $\tau\in [0, 1]$.
\end{theorem}

Next we consider the case $N=\triangle_M$.

\begin{claim}\label{cl:3.9}
Let $S^1\cdot\gamma_0\subset\Lambda M=W^{1,2}(S^1, M)$ be an
isolated nonconstant critical orbit of ${\cal L}$ on $\Lambda M$ (
with $F(\gamma_0(t),\dot\gamma_0(t))\equiv \sqrt{c}>0$). Let
$L^\ast:TM\to\R$ be given by Corollary~\ref{cor:2.3}, and let
$L^\tau(x,v)=(1-\tau)L(x,v)+
 \tau L^\ast(x,v)$ for $\tau\in [0,1]$.
Then there exists a neighborhood of $S^1\cdot\gamma_0$ in $\Lambda
M$, ${\cal U}(S^1\cdot\gamma_0)$, such that the critical set of each
${\cal L}^\tau$ in ${\cal U}(S^1\cdot\gamma_0)$ is the orbit
$S^1\cdot\gamma_0$.
\end{claim}

\noindent{\bf Proof}. From the beginning of the proof of
Claim~\ref{cl:3.7} one easily sees that $S^1\cdot\gamma_0$ is a
critical orbit of each ${\cal L}^\tau$. As in the proof of
Claim~\ref{cl:3.7}, by a contradiction, suppose that there exist
sequences $(\gamma_k)\subset\Lambda M$, $(\tau_k)\subset [0,1]$ and
$s\in S^1$ such that
$$
\tau_k\to\tau_0,\quad\gamma_k\to s\cdot\gamma_0,\quad d{\cal
L}^{\tau_k}(\gamma_k)=0\;\hbox{and}\;\gamma_k\notin
S^1\cdot\gamma_0\;\forall k.
$$
Then it was shown at the end of the proof of Claim~\ref{cl:3.7} that
 $d{\cal L}(\gamma_k)=0$ for sufficiently large $k$.  By the assumption that $S^1\cdot\gamma_0\subset
 \Lambda M$ is an isolated nonconstant critical orbit of ${\cal
L}$ on $\Lambda M$, we obtain $\gamma_k=s_k\cdot\gamma_0$ for some
$s_k\in S^1$ and each sufficiently large $k$. This contradiction
gives our claim. \hfill$\Box$\vspace{2mm}

Theorem~\ref{th:1.5}(ii) follows from Claim~\ref{cl:3.9} and the
following result.

\begin{theorem}\label{th:3.10}
Under the assumptions of Claim~\ref{cl:3.9}, $C_\ast({\cal L}^\tau,
S^1\cdot\gamma_0;{\K})$ is independent of $\tau\in [0, 1]$.
\end{theorem}

\noindent{\bf Proof}.  Since $L^\tau=(1-\tau)L+ \tau L^\ast$ and
$L(s\cdot\gamma_0(t), (s\cdot\gamma_0)'(t))=c=L^\ast(s\cdot\gamma_0(t), (s\cdot\gamma_0)'(t))$ for all $t\in\R$ and
$s\in S^1$,  we have
\begin{equation}\label{e:3.15}
{\cal L}^\tau(s\cdot\gamma_0)=c\quad\forall \tau\in [0,1],\;s\in
S^1.
\end{equation}
 Claim~\ref{cl:3.9} yields a neighborhood of
$S^1\cdot\gamma_0$ in $\Lambda M$, ${\cal U}(S^1\cdot\gamma_0)$,
such that
\begin{equation}\label{e:3.16}
{\cal U}(S^1\cdot\gamma_0)\cap K({\cal
L}^\tau)=S^1\cdot\gamma_0\quad\forall\tau\in [0,1],
\end{equation}
where $K({\cal L}^\tau)$ denotes the critical set of ${\cal
L}^\tau$.

 Let ${\cal
O}=S^1\cdot \gamma_0$. Since $\gamma_0$ is nonconstant, it
 has minimal period $1/m$ for some $m\in\N$, and
${\cal O}$ is a $1$-dimensional $C^{k-1}$ submanifold diffeomorphic to
the circle (\cite[page 499]{GM2}). For every $s\in [0, 1/m]\subset
S^1$ the tangent space $T_{s\cdot \gamma_0}(S^{1}\cdot \gamma_0)$ is
$\R(s\cdot \gamma_0)'$, and the fiber $N{\cal O}_{s\cdot \gamma_0}$
at $s\cdot \gamma_0$ of the normal bundle $N{\cal O}$ of ${\cal O}$
is a subspace of codimension $1$ which is orthogonal to $(s\cdot
\gamma_0)'$ in $T_{s\cdot \gamma_0}\Lambda M$. For the
diffeomorphism $\digamma$ in (\ref{e:1.13}), by shrinking
$\varepsilon>0$ we may require
 that ${\cal N}({\cal O},\varepsilon)\subset {\cal
U}(S^1\cdot\gamma_0)$ (hence ${\cal N}({\cal O},\varepsilon)$
contains no other critical orbit besides ${\cal O}$), and that
$\digamma(\{y\}\times N{\cal O}(\varepsilon)_y\bigl)$ and ${\cal O}$
have a unique intersection point $0_y=y$, where  ${\cal
U}(S^1\cdot\gamma_0)$ is as in Claim~\ref{cl:3.9}.
 For ${\cal L}^\tau$ in (\ref{e:3.15}), define
\begin{equation}\label{e:3.17}
{\cal F}^\tau: N{\cal O}(\varepsilon)\to \R,\;(y,v)\mapsto {\cal
L}^\tau\circ\digamma(y,v).
\end{equation}
It is $S^1$-invariant, $C^{2-0}$, and  satisfies the (PS) condition.
Moreover ${\cal F}^0={\cal F}$ and ${\cal F}^1={\cal F}^\ast$. Since
$C_\ast({\cal L}^\tau, \mathcal{O};{\K})=C_\ast({\cal F}^\tau,
\mathcal{O};{\K})\;\forall\tau\in [0, 1]$, we only need to prove
that
\begin{equation}\label{e:3.18}
C_\ast({\cal F}^\tau, {\cal O};{\K})\quad\hbox{is independent of
$t$}.
\end{equation}
By (\ref{e:3.4}) and $\nabla{\cal L}^{\gamma_0}(\gamma)-\nabla{\cal
L}^\tau(\gamma)=(\tau-{\gamma_0})(\nabla{\cal L}(\gamma)-\nabla{\cal
L}^\ast(\gamma))$ for ${\gamma_0},\tau\in [0, 1]$ and $\gamma\in
\Lambda M$, we derive that  there exists a constant $C_{11}>0$ such
that
\begin{eqnarray}\label{e:3.19}
&&\sup\bigl\{|{\cal F}^{\gamma_0}(y,v)-{\cal
F}^\tau(y,v)|+|\nabla{\cal F}^{\gamma_0}(y,v)-\nabla{\cal
F}^\tau(y,v)|\,:\,(y,v)\in N({\cal O})(\varepsilon)\bigr\}\nonumber\\
&&\le C_{11}|{\gamma_0}-\tau|
\end{eqnarray}
after shrinking $\varepsilon>0$ (if necessary) because $\mathcal{O}$
is compact. Using these (\ref{e:3.18}) easily follows from Chang and
Ghoussoub \cite[Th.III.4]{ChGh}. We may also prove (\ref{e:3.18}) as
follows. For each fixed  $\tau\in [0,1]$, following
\cite[Th.2.3]{Wa} we may construct a Gromoll-Meyer pair of ${\cal
O}$ as a critical submanifold of ${\cal F}^\tau$ on $N({\cal
O})(\varepsilon)$ with respect to $-\nabla{\cal F}^\tau$,
$\bigl(W({\cal O}), W({\cal O})^-\bigr)$, such that
$$
\bigl(W({\cal O})_y, W({\cal O})_y^-\bigr):=\bigl(W({\cal O})\cap
N({\cal O})(\varepsilon)_y,\, W({\cal O})^-\cap
  N({\cal O})(\varepsilon)_y\bigr)
$$
is a  Gromoll-Meyer pair of  ${\cal F}^\tau|_{N({\cal
O})(\varepsilon)_y}$ at its isolated critical point $0=0_y$
satisfying
$$
\bigl(W({\cal O})_{s\cdot y},\, W({\cal O})^-_{s\cdot y}\bigr) =
\bigl(s\cdot W({\cal O})_y,\, s\cdot W({\cal
O})^-_y\bigr)\quad\forall (s, y)\in S^1\times {\cal O}.
$$

Slightly modifying the proof of Lemma~5.2 on the page 52 of
\cite{Ch93} we may show that $\bigl(W({\cal O}), W({\cal
O})^-\bigr)$ is also a
 Gromoll-Meyer pair of ${\cal O}$ as a critical submanifold of ${\cal F}^{\gamma_0}$ on $N({\cal
O})(\varepsilon)$ with respect to certain pseudo-gradient vector
field of ${\cal F}^{\gamma_0}$ if ${\gamma_0}\in [0,1]$ is
sufficiently close to $\tau$
 because of (\ref{e:3.19}). So  we
may get an open neighborhood $J_\tau$  of  $\tau\in [0,1]$ in
$I=[0,1]$ such that $C_\ast({\cal F}^\tau, {\cal
O};{\K})=C_\ast({\cal F}^{\gamma_0}, {\cal O};{\K})$ for any
${\gamma_0}\in J_\tau$. Then (\ref{e:3.18}) follows from this and
the compactness of $[0,1]$. \hfill$\Box$\vspace{2mm}

\section{Proofs of Theorems~\ref{th:1.2},\ref{th:1.3}}\label{sec:4}
\setcounter{equation}{0}

\noindent{\bf Step 1}. {\it Prove the corresponding versions of
Theorems~\ref{th:1.2}, ~\ref{th:1.3} under a new chart around
$\gamma_0$.}

For conveniences of computations we need to consider a coordinate
chart around $\gamma_0$ different from that of (\ref{e:1.4}). Recall
that $M_0$ (resp. $M_1$) is totally geodesic near $\gamma_0(0)$
(resp. $\gamma_0(1)$) with respect to the chosen Riemannian metric
$g$ on $M$.
 Since $\gamma_0$ is of class $C^k$ we may
 take a parallel orthogonal $C^k$ frame field along $\gamma_0$ with respect
to the metric $g$,  $I\ni t\to (e_1(t),\cdots, e_n(t))$. For a small
open ball $B^n_{2\rho}(0)\subset\R^n$ we get a $C^k$ map
\begin{equation}\label{e:4.1}
\phi:I\times B^n_{2\rho}(0)\to M,\;(t,v)\mapsto
\exp_{\gamma_0(t)}\left(\sum^n_{i=1}v_ie_i(t)\right).
\end{equation}
Since there exist linear subspaces $V_i\subset\R^n$, $i=0,1$, such
that $v\in V_i$ if and only if $\sum^n_{k=1}v_ke_k(i)\in
T_{\gamma_0(i)}M_i$, $i=0,1$, by shrinking $\rho>0$ (if necessary)
we get that $ v\in V_i\cap B^n_{2\rho}(0)$ if and only if $
\phi(i,v)\in M_i$, $i=0,1$. Set $V:=V_0\times V_1$ and
\begin{eqnarray*}
&&H_V=W^{1,2}_{V}(I, \R^n):=\{\zeta\in W^{1,2}(I, \R^n)\,|\,
(\zeta(0),\zeta(1))\in V\},\\
&& X_V=C^1_{V}(I, \R^n):=\{\zeta\in C^1(I, \R^n)\,|\,
(\zeta(0),\zeta(1))\in V\}.
\end{eqnarray*}
As usual we use $(\cdot,\cdot)_{W^{1,2}}$ and $\|\cdot\|_{W^{1,2}}$
to denote the inner product and norm in $H_V$. Let ${\bf
B}_{2\rho}(H_V):=\{\zeta\in H_V\,|\,\|\zeta\|_{W^{1,2}}<2\rho\}$.
 Then by the omega lemma the map
\begin{equation}\label{e:4.2}
\Phi:{\bf B}_{2\rho}(H_V) \to\Lambda_N(M)
\end{equation}
defined by $\Phi(\zeta)(t)=\phi(t,\zeta(t))$, gives a $C^{k-3}$
coordinate chart around $\gamma_0$ on $\Lambda_N(M)$ with ${\rm
Im}(\Phi)\subset {\cal O}(\gamma_0)$. Define $\tilde F:I\times
B^n_{2\rho}(0)\times\R^n\to\R$ by
$$
\tilde F(t, x, v)=F\bigr(\phi(t,x), d\phi(t,x)[(1,v)]\bigl).
$$
Then $\tilde F(t, 0, 0)=F(\phi(t,0),
\partial_t\phi(t,0)[1])=F(\gamma_0(t),\dot\gamma_0(t))\equiv\sqrt{c}$.
Moreover the $C^{2-0}$ function $\tilde L:=\tilde F^2$ satisfies
$\tilde L(t, x, v)=L\bigr(\phi(t,x), d\phi(t,x)[(1,v)]\bigl)$, and
is $C^k$ in $(I\times B^n_{2\rho}(0)\times\R^n)\setminus{\cal Z}$,
where
$$
{\cal Z}:=\bigl\{(t,x,v)\in I\times B^n_{2\rho}(0)\times\R^n\,|\,
\partial_x\phi(t,x)[v]=-\partial_t\phi(t,x)\bigr\},
$$
a closed subset in $I\times B^n_{2\rho}(0)\times\R^n$.
 Since $\gamma_0$ is regular, i.e.,
$\partial_t\phi(t,0)=\dot\gamma_0(t)\ne 0$ at each $t\in I$, and
$\partial_x\phi(t,x)$ is injective, we deduce that $(t, 0,
0)\notin{\cal Z}\;\forall t\in I$.  It follows that $I\times
B^n_{2r}(0)\times B^n_{2r}(0)\subset I\times
B^n_{2\rho}(0)\times\R^n\setminus{\cal Z}$ for some $0<r<\rho$. We
also require $r>0$ so small that
\begin{equation}\label{e:4.3}
\tilde{L}(t,x,v)\ge\frac{2c}{3}\quad\forall (t,x,v)\in I\times
B^n_{r}(0)\times B^n_{r}(0).
\end{equation}

 We conclude  that the boundary condition (\ref{e:1.3}) becomes
\begin{equation}\label{e:4.4*}
\partial_v\tilde L(i, 0, 0)[v]=0\qquad\forall v\in V_i,\;i=0,1.
\end{equation}
In fact, for any $X\in T_{\gamma(0)}M_0$ there
exists a unique $v=(v_1,\cdots,v_n)\in V_0$ such that
$X=\sum^n_{k=1}v_ke_k(0)=\partial_v\phi(0,0)[v]$. Since
$\gamma_0(0)=\phi(0,0)$ and $\dot\gamma_0(0)=\partial_t\phi(0,0)$ we
get
\begin{eqnarray*}
0&=&g^F(\gamma_0(0),\dot\gamma_0(0))[X,\dot\gamma_0(0)]\\
&=&\frac{1}{2}\frac{d^2}{ds
d\tau}\Bigl|_{s=0,\tau=0}L(\gamma_0(0),\dot\gamma_0(0)+
sX+\tau\dot\gamma_0(0))\\
&=&\frac{1}{2}\frac{d^2}{ds
d\tau}\Bigl|_{s=0,\tau=0}L\left(\phi(0,0),\partial_t\phi(0,0)+
s\partial_v\phi(0,0)[v]+\tau\partial_t\phi(0,0)\right)\\
&=&\frac{1}{2}\frac{d^2}{ds d\tau}\Bigl|_{s=0,\tau=0}\left[(1+\tau)^2
 L\Big(\phi(0,0),\partial_t\phi(0,0)+
\frac{s}{1+\tau}\partial_v\phi(0,0)[v]\Big)\right]\\
&=&\frac{1}{2}\frac{d^2}{ds d\tau}\Bigl|_{s=0,\tau=0}\left[(1+\tau)^2
 L\Big(\phi(0,0), d\phi(0,0)[1,\frac{s}{1+\tau}v]\Big)\right]\\
&=&\frac{1}{2}\frac{d^2}{ds d\tau}\Bigl|_{s=0,\tau=0}\left[(1+\tau)^2
 \tilde L\Big(0, 0, \frac{s}{1+\tau}v\Big)\right]\\
&=&\frac{1}{2}\frac{d}{ds }\Bigl|_{s=0}\left[2
 \tilde L(0, 0, sv)
 -s\partial_v\tilde L(0, 0, sv)[v]
 \right]=\partial_v\tilde L(0, 0, 0)[v].
\end{eqnarray*}
Similarly, we may prove that
$g^F(\gamma_0(1),\dot\gamma_0(1))[X,\dot\gamma_0(1)]=0\;\forall X\in
T_{\gamma_0(1)}M_1$ if and only if $\partial_v\tilde L(1, 0,
0)[v]=0\;\forall v\in V_1$.

From now on we write $U={\bf B}_{2r}(H_V)$ and
\begin{center}
\textsf{$U_X:={\bf B}_{2r}(H_V)\cap X_V$ as an open subset of
$X_V$}.
\end{center}
 Define
the action functional
\begin{equation}\label{e:4.5}
\tilde{\cal L}: {\bf B}_{2r}(H_V)\to\R,\;\zeta\mapsto \tilde{\cal
L}(\zeta)=\int^1_0\tilde L(t,\zeta(t),\dot\zeta(t))dt,
\end{equation}
that is, $\tilde{\cal L}={\cal L}\circ\Phi$. It is $C^{2-0}$, has
$0\in {\bf B}_{2r}(H_V)$ as a unique
 critical point. Denote by $\tilde{\cal L}_X$ the restriction of
 $\tilde{\cal L}$ on $U_X$, and by $\tilde A$  the restriction of the gradient
 $\nabla\tilde{\cal L}$ to
$U_X$. As in the proof of \cite[Lemma~3.2]{Lu1} or \cite[Lemma~6]{CaJaMa1}
we can prove that $\tilde A(U_X)\subset X_V$. (In fact, this can be seen
from the following (\ref{e:4.13}), (\ref{e:4.14}) and (\ref{e:4.15})
if we replace $\tilde{L}^\ast$ and $\tilde{\cal L}^\ast$ there by
 $\tilde{L}$ and $\tilde{\cal L}$, respectively.)
 Recall that
${\bf B}_{2r}(X_V)=\{\zeta\in X_V\,|\, \|\zeta\|_{C^1}<2r\}\subset
{\bf B}_{2r}(H_V)$. By the choice of $r$ above (\ref{e:4.3}),
$\tilde A:{\bf B}_{2r}(X_V)\to X_V$ is $C^{k-3}$. (This implies that the
restriction of $\tilde{\cal L}$ to $U_X$, denoted by $\tilde{\cal
L}_X$, is $C^{k-2}$ on ${\bf B}_{2r}(X_V)$.) The continuous symmetric
bilinear form $d^2\tilde{\cal L}_X(0)$ on $ X_V$ can be extended
into a symmetric bilinear form on $H_V$ whose associated self-adjoint operator is
Fredholm, and has finite dimensional negative definite and null
spaces $H^-_V$ and $H^0_V$, which are actually contained in $X_V$.
Let $H^+_V$ be the corresponding positive definite space. Then the
orthogonal decomposition $H_V=H^-_V\oplus H^0_V\oplus H^+_V$ induces
  a (topological) direct sum decomposition $X_V=X^-_V\dot{+} X^0_V\dot{+}
X^+_V$, where as sets $X^0_V=H^0_V={\rm Ker}(d\tilde A(0))$,
$X^-_V=H^-_V$ and $X^+_V=X_V\cap H^+_V$. Note that $H_V$ and $X_V$
induce equivalent norms on $H_V^0=X_V^0$. By the implicit function
theorem we get a $\tau\in (0, r]$ and a $C^{k-3}$-map $\tilde h: {\bf
B}_\tau(H_V^0)\to X^-_V\dot{+} X^+_V$ with $\tilde h(0)=0$ and
$d\tilde h(0)=0$ such that
\begin{equation}\label{e:4.6}
\zeta+ \tilde h(\zeta)\in {\bf B}_{2r}(X^0_V)\quad\hbox{and}\quad
(I-P^0_V)\tilde A(\zeta+ \tilde h(\zeta))=0
\end{equation}
for each $\zeta\in {\bf B}_\tau(H^0_V)$, where $P^0_V:H_V\to H^0_V$
is the orthogonal projection. It is not hard to prove that the Morse
index $m^-(\gamma_0)$ and nullity $m^0(\gamma_0)$ of $\gamma_0$ are
equal to  $\dim X^-_V$ and $\dim X^0_V$, respectively.

Let $\tilde L^\ast(t, x, v)=L^\ast\bigr(\phi(t,x),
d\phi(t,x)[(1,v)]\bigl)$. Since $L^\ast(x,v)= L(x,v)$
 if $L(x,v)\ge\frac{2c}{3}$,
$L^\ast(\gamma_0(t), \dot\gamma_0(t))=c\;\forall t$, and
\begin{equation}\label{e:4.7}
\tilde L^\ast(t, x,v)= \tilde L(t, x,v)\quad
 \hbox{if}\quad \tilde L(t,x,v)\ge\frac{2c}{3}.
\end{equation}
In particular,  $\tilde L(t, 0, 0)=c$ implies that $\tilde L^\ast(t,
0,0)\equiv c\;\forall t$. Define the action functional
$$
\tilde{\cal L}^\ast: {\bf B}_{2r}(H_V)\to\R,\;\zeta\mapsto
\int^1_0\tilde L^\ast(t,\zeta(t),\dot\zeta(t))dt,
$$
that is, $\tilde{\cal L}^\ast={\cal L}^\ast\circ\Phi$. It is
$C^{2-0}$, and has $0\in {\bf B}_{2r}(H_V)$ as a unique
 critical point. Let $\tilde{\cal L}^\ast_X$ be the restriction of
$\tilde{\cal L}^\ast$ to $U_X={\bf B}_{2r}(H_V)\cap X_V$.

 We firstly  check that $\tilde{\cal L}^\ast$ satisfies the conditions of
Theorem~\ref{th:A.1}. This can be obtained by almost repeating the
arguments in \cite[\S3]{Lu1}. For the sake of completeness we also
give them.  It is easily computed (cf. \cite[\S3]{Lu1}) that
\begin{eqnarray*}
d\tilde{\cal L}^\ast(\zeta)[\xi]
 \!\!\!&=&\!\!\!\int_0^{1}  \left( \partial_{q} \tilde
L^\ast(t,\zeta(t),\dot{\zeta}(t))\cdot\xi(t) + \partial_{v} \tilde
L^\ast(t,\zeta(t),\dot{\zeta}(t))\cdot\dot{\xi}(t) \right) \, dt
\end{eqnarray*}
for any $\zeta\in {\bf B}_{2r}(H_V)$ and $\xi\in H_V$.  Let us
compute the gradient $\nabla\tilde{\cal L}^\ast(\zeta)$. Define
\begin{eqnarray}\label{e:4.8}
G(\zeta)(t):=\int^t_0 \left[ \partial_v \tilde
L^\ast(s,\zeta(s),\dot{\zeta}(s))-c_0\right]ds\quad\forall t\in
[0,1],
\end{eqnarray}
where $c_0=\int^1_0\partial_v \tilde
L^\ast(s,\zeta(s),\dot{\zeta}(s))ds$. Then
$G(\zeta)(0)=0=G(\zeta)(1)$, and hence $G(\zeta)\in
W_0^{1,2}(I,\R^n)\subset H_V$. Moreover
 \begin{eqnarray*}
&&\quad\int_0^{1}  \left( \partial_q \tilde
L^\ast(t,\zeta(t),\dot{\zeta}(t))\cdot\xi(t) + \partial_v \tilde
L^\ast(t,\zeta(t),\dot{\zeta}(t))\cdot\dot{\xi}(t)
\right)  \, dt\\
&&=(G(\zeta),\xi)_{W^{1,2}}+ c_0\int^1_0\dot\xi(t)dt+\int_0^{1}
\left( \partial_q \tilde
L^\ast(t,\zeta(t),\dot{\zeta}(t))-G(\zeta)(t)\right)\cdot\xi(t)\,dt.
\end{eqnarray*}
By the Riesz theorem one may get a unique $F(\zeta)\in H_V$ such
that
\begin{eqnarray}\label{e:4.9}
c_0\int^1_0\dot\xi(t)dt+\int_0^{1}  \left( \partial_q \tilde
L^\ast(t,\zeta(t),\dot{\zeta}(t))- G(\zeta)(t)\right)\cdot\xi(t)\,dt
=\bigl(F(\zeta),\xi\bigr)_{W^{1,2}}
\end{eqnarray}
for any $\xi\in H_V$. Hence $d\tilde{\cal
L}^\ast(\zeta)[\xi]=(G(\zeta),\xi)_{W^{1,2}}+
(F(\zeta),\xi)_{W^{1,2}}$ and thus
\begin{equation}\label{e:4.10}
\nabla\tilde{\cal L}^\ast(\zeta)= G(\zeta)+ F(\zeta).
\end{equation}
Since
$$
\bigl(F(\zeta),\xi\bigr)_{W^{1,2}}=\int^1_0\left(F(\zeta)(t)\cdot\xi(t)+
\frac{d}{dt}F(\zeta)(t)\cdot\dot\xi(t)\right)dt,
$$
(\ref{e:4.9}) becomes
\begin{eqnarray}\label{e:4.11}
&&\int_0^{1}  \left( \partial_q \tilde
L^\ast(t,\zeta(t),\dot{\zeta}(t))-
G(\zeta)(t)-F(\zeta)(t)\right)\cdot\xi(t)\,dt\nonumber\\
 &&=\int^1_0
\left(\frac{d}{dt}F(\zeta)(t)-c_0\right)\cdot\dot\xi(t)dt\quad\forall\xi\in
H_V.
\end{eqnarray}

\begin{lemma}\label{lem:4.1}
For $f\in L^1(I, \R^n)$  the equation $\ddot{x}(t)- x(t)=f(t)$ has
the general solution of the following form
\begin{eqnarray*}
x(t)= e^t\int^t_0\left[
e^{-2s}\int^s_0e^{\tau}f(\tau)d\tau\right]ds+ c_1e^t+ c_2e^{-t},
\end{eqnarray*}
where $c_i\in\R^n$, $i=1,2$, are constant vectors.
\end{lemma}

Setting $y(t)=\dot x(t)-\int^t_0x(\tau)d\tau$, this lemma can easily be proved by
the standard methods. Let constant vectors $c_1, c_2\in\R^n$ be such
that the function
\begin{equation}\label{e:4.12}
z(t):=
e^t\int^t_0\left[ e^{-2s}\int^s_0e^{\tau}f(\tau)d\tau\right]ds+
c_1e^t+ c_2e^{-t}
\end{equation}
 satisfies $z(0)=F(\zeta)(0)$ and
$z(1)=F(\zeta)(1)-c_0$ with
\begin{eqnarray}\label{e:4.13}
 f(t)&=&- \partial_q \tilde
L^\ast(t,\zeta(t),\dot{\zeta}(t))+
G(\zeta)(t)+ c_0t\nonumber\\
&=& - \partial_q \tilde L^\ast(t,\zeta(t),\dot{\zeta}(t))+ \int^t_0
\partial_v \tilde L^\ast(s,\zeta(s),\dot{\zeta}(s))ds .
\end{eqnarray}
 (We identify each
element of $W^{1,2}(I,\R^n)$ with its unique continuous
representation as usual).
 Then for any $\xi\in
C^1(I,\R^n)$ with $\xi(0)=\xi(1)=0$ it holds that
\begin{eqnarray*}
\int^1_0\dot z(t)\cdot\dot\xi(t)
dt&=&\dot z(t)\cdot\xi(t)\Bigm|^{t=1}_{t=0}-\int^1_0\ddot
z(t)\cdot\xi(t)dt\\
&=&-\int^1_0(z(t)+ f(t))\cdot\xi(t) dt\\
&=&\int^1_0\bigl(\partial_q \tilde
L^\ast(t,\zeta(t),\dot{\zeta}(t))- G(\zeta)(t)-c_0t-
z(t)\bigr)\cdot\xi(t) dt.
\end{eqnarray*}
From this and (\ref{e:4.11}) it follows that
$$
\int^1_0\bigl(F(\zeta)(t)-c_0t-z(t)\bigr)\cdot\xi(t)dt=-\int^1_0
\biggl(\frac{d}{dt}F(\zeta)(t)-c_0-\dot z(t)\biggr)\cdot\dot\xi(t)dt
$$
for any $\xi\in C^1(I,\R^n)$ with $\xi(0)=\xi(1)=0$. Since
$F(\zeta)(t)-c_0t-z(t)$ is equal to  zero at $t=0,1$, by Theorem
8.7 in \cite{Bre} there exists a sequence $(u_k)\in C_0^\infty(\R)$
such that $(u_k|_I)_k$ converges to the function
$F(\zeta)(t)-c_0t-z(t)$ in $W^{1,2}(I,\R^n)$ (and hence in
$C(I,\R^n)$). In particular we have $u_k(i)\to 0$ because
$F(\zeta)(i)-c_0\cdot i-z(i)=0$, $i=0,1$. Define $v_k:I\to\R^n$ by
$v_k(t)=u_k(t)-u_k(0)-t(u_k(0)-u_k(1))$ for each $k\in\N$. Then
$v_k\in C^\infty(I,\R^n)$ and $v_k(0)=v_k(1)=0$ for each $k$, and
the sequence $(v_k)$  converges to $F(\zeta)(t)-c_0t-z(t)$ in
$W^{1,2}(I,\R^n)$. Let $k\to\infty$ in
$$
\int^1_0\bigl(F(\zeta)(t)-c_0t-z(t)\bigr)\cdot v_k(t)dt=-\int^1_0
\biggl(\frac{d}{dt}F(\zeta)(t)-c_0-\dot z(t)\biggr)\cdot\dot
v_k(t)dt,
$$
 we obtain
$$
\int^1_0|F(\zeta)(t)-c_0t-z(t)|^2 dt=-\int^1_0
\biggl|\frac{d}{dt}F(\zeta)(t)-c_0-\dot z(t)\biggr|^2dt
$$
and therefore $F(\zeta)(t)=c_0t+ z(t)\;\forall t\in I$ since both
$F(\zeta)$ and $z$ are continuous on $I$. By (\ref{e:4.12}),
(\ref{e:4.8}) and (\ref{e:4.10}) we arrive at
\begin{eqnarray}
F(\zeta)(t)&=& e^t\int^t_0\left[
e^{-2s}\int^s_0e^{\tau}f(\tau)d\tau\right]ds+ c_1e^t+
c_2e^{-t}+ c_0t,\nonumber\\
\nabla\tilde{\cal L}^\ast(\zeta)(t)&=&e^t\int^t_0\left[
e^{-2s}\int^s_0e^{\tau}f(\tau)d\tau\right]ds  + c_1e^t+
c_2e^{-t}\nonumber\\
&&\qquad +\int^t_0 \partial_v \tilde
L^\ast(s,\zeta(s),\dot{\zeta}(s))ds ,\label{e:4.14}
\end{eqnarray}
where $c_1, c_2\in\R^n$ are suitable constant vectors and $f(t)$ is
given by (\ref{e:4.13}). By (\ref{e:4.14})  the function
$\nabla\tilde{\cal L}^\ast(\zeta)(t)$ is differentiable almost
everywhere, and for a.e. $t\in I$,
\begin{eqnarray}\label{e:4.15}
\frac{d}{dt}\nabla\tilde{\cal L}^\ast(\zeta)(t)&=& e^t\int^t_0\left[
e^{-2s}\int^s_0e^{\tau}f(\tau)d\tau\right]ds
+e^{-t}\int^t_0e^{\tau}f(\tau)d\tau
\nonumber\\
&&\quad + c_1e^t -c_2e^{-t}+ \partial_v \tilde
L^\ast(t,\zeta(t),\dot{\zeta}(t)).
\end{eqnarray}
Let $\tilde A^\ast$ denote the restriction of the gradient
$\nabla\tilde{\cal L}^\ast$ to $U_X$.
 Clearly, (\ref{e:4.14}) and (\ref{e:4.15}) imply that
$\tilde A^\ast(\zeta)\in X_V$ for $\zeta\in U_X$, and that
$U_X\ni\zeta\mapsto \tilde A^\ast(\zeta)\in X_{V}$ is continuous.
From the expression of $d\tilde{\cal L}^\ast(\zeta)$ above
(\ref{e:4.8}) it is easily seen that  $\tilde{\cal L}^\ast_X$
 is at least $C^2$ and
\begin{eqnarray}\label{e:4.16}
 d^2 \tilde{\cal L}^\ast_X
  (\zeta)[\xi,\eta]
   = \int_0^{1} \Bigl(\!\! \!\!\!&&\!\!\!\!\!\partial_{vv}
    \tilde L^\ast(t,\zeta(t),\dot{\zeta}(t))
\bigl[\dot{\xi}(t), \dot{\eta}(t)\bigr] \nonumber\\
&&+ \partial_{qv} \tilde
  L^\ast(t,\zeta(t), \dot{\zeta}(t))
\bigl[\xi(t), \dot{\eta}(t)\bigr]\nonumber \\
&& + \partial_{vq} \tilde
  L^\ast(t,\zeta(t),\dot{\zeta}(t))
\bigl[\dot{\xi}(t), \eta(t)\bigr] \nonumber\\
&&+  \partial_{qq} \tilde L^\ast(t,\zeta(t), \dot{\zeta}(t))
\bigl[\xi(t), \eta(t)\bigr]\Bigr) \, dt
\end{eqnarray}
for any $\zeta\in U_X,\; \xi,\eta\in X_V$. In fact, as in the proof
of \cite[Lemma 3.2]{Lu1} we have also

\begin{lemma}\label{lem:}\label{lem:4.2}
 The map $\tilde
A^\ast:U_X\to X_V$
 is continuously differentiable.
\end{lemma}
From (\ref{e:4.16}) it  easily follows that
\begin{description}
\item[(i)] for any $\zeta\in U_X$ there
exists a constant $C(\zeta)$ such that
$$
\big|d^2 {\tilde{\cal L}}^\ast_{X}
  (\zeta)[\xi,\eta]\big|\le
  C(\zeta)\|\xi\|_{W^{1,2}}\cdot\|\eta\|_{W^{1,2}}
\quad\forall \xi,\eta\in X_V;
$$
\item[(ii)] $\forall\varepsilon>0,\;\exists\;\delta_0>0$, such that
for any $\zeta_1, \zeta_2\in U_X$ with
$\|\zeta_1-\zeta_2\|_{C^1}<\delta_0$,
$$
\big|d^2 {\tilde{\cal L}}^\ast_{X}
  (\zeta_1)[\xi,\eta]- d^2 {\tilde{\cal L}}^\ast_{X}
  (\zeta_2)[\xi,\eta]\big|\le
  \varepsilon \|\xi\|_{W^{1,2}}\cdot\|\eta\|_{W^{1,2}}
\quad\forall \xi,\eta\in X_V.
$$
\end{description}
 (i) shows that the right side of (\ref{e:4.16}) is also a bounded symmetric
bilinear form on $H_V$. As in \cite[\S3]{Lu1} we have a map $\tilde
B^\ast: U_X\to L_s(H_V)$, which is uniformly continuous, such that
\begin{equation}\label{e:4.17}
\bigl(d\tilde A^\ast(\zeta)[\xi], \eta\bigr)_{W^{1,2}}=d^2
{\tilde{\cal L}}^\ast_{X}
  (\zeta)[\xi,\eta]=\bigl(\tilde B^\ast(\zeta)\xi,
  \eta\bigr)_{W^{1,2}}
\end{equation}
for  any $\zeta\in U_X$ and $\xi, \eta\in X_V$. Namely (\ref{e:A.2})
is satisfied. Almost repeating the arguments in \cite[\S3]{Lu1} we
may check that the map $\tilde B^\ast$ satisfy the conditions {\bf
(B1)} and {\bf (B2)} in Appendix~\ref{app:A}. Summarizing the above
arguments, $(\tilde{\cal L}, \tilde A^\ast, \tilde B^\ast)$
satisfies the conditions of Theorem~\ref{th:A.1} (resp.
Theorem~\ref{th:A.5}) around the critical point $0\in H_V$ (resp.
$0\in X_V$). By (\ref{e:4.3}) and (\ref{e:4.7}) for any $\zeta\in
{\bf B}_{r}(X_V)$ we have
\begin{equation}\label{e:4.18}
\tilde L^\ast(t, \zeta(t), \dot{\zeta}(t))= \tilde L(t,
\zeta(t),\dot{\zeta}(t))\quad\forall t\in I.
\end{equation}
This implies that
\begin{equation}\label{e:4.19}
\tilde A^\ast(\zeta)= \tilde A(\zeta)\quad\forall\zeta\in {\bf
B}_{r}(X_V),\quad\hbox{and}\quad\tilde B^\ast(0)= \tilde B(0),
\end{equation}
where the map $\tilde B: {\bf B}_{2r}(X_V)\to L_s(H_V)$ is
determined by the equation
$$
d^2{\tilde{\cal L}}_{X}(\zeta)[\xi,\eta]=\bigl(\tilde B(\zeta)\xi,
  \eta\bigr)_{W^{1,2}}\quad\forall \zeta\in {\bf
B}_{2r}(X_V),\; \xi, \eta\in X_V.
$$
({\bf Note}: the domain of $\tilde B$ is different from the one of
$\tilde{B}^\ast$!) Shrink $\tau$ in (\ref{e:4.6}) so small that
\begin{equation}\label{e:4.20}
\zeta+ \tilde h(\zeta)\in {\bf B}_r(X_V)\quad\forall \zeta\in {\bf
B}_\tau(H^0_V).
\end{equation}
Then this and (\ref{e:4.19}) imply  that the $C^{k-1}$-map $\tilde h$ in
(\ref{e:4.6}) satisfies
$$
(I-P^0_V)\tilde A^\ast(\zeta+ \tilde h(\zeta))=0\;\forall\zeta\in
{\bf B}_\tau(H^0_V).
$$
Let $\tilde{\cal L}^{\ast\circ}, \tilde{\cal L}^\circ:{\bf
B}_\tau(H_V^0)\to\R$ be defined by
 \begin{equation}\label{e:4.21}
\tilde{\cal L}^{\ast\circ}(\zeta)=\tilde{\cal L}^\ast\bigl(\zeta+
\tilde h(\zeta)\bigr)\quad\hbox{and}\quad  \tilde{\cal
L}^\circ(\zeta)=\tilde{\cal L}\bigl(\zeta+ \tilde h(\zeta)\bigr),
\end{equation}
respectively. Then (\ref{e:4.20}) and (\ref{e:4.18}) lead to
\begin{equation}\label{e:4.22}
\tilde{\cal L}^{\ast\circ}(\zeta)=\tilde{\cal L}^\ast\bigl(\zeta+
\tilde h(\zeta)\bigr)=\tilde{\cal L}\bigl(\zeta+ \tilde
h(\zeta)\bigr)=\tilde{\cal L}^{\circ}(\zeta)
\end{equation}
for all $\zeta\in {\bf B}_\tau(H^0_V)$. By  Theorem~\ref{th:A.1} we
obtain the following splitting theorem.
\begin{theorem}\label{th:4.3}
Under the notations above, there exists a ball ${\bf
B}_\eta(H_V)\subset{\bf B}_\tau(H_V)$,  an origin-preserving local
homeomorphism $\tilde\psi$ from ${\bf B}_\eta(H_V)$ to an open
neighborhood of $0\in H_V$  such that
$$
\tilde{\cal
L}^\ast\circ\tilde\psi(\zeta)=\|P_V^+\zeta\|^2_{W^{1,2}}- \|
P^-_V\zeta\|^2_{W^{1,2}} + \tilde{\cal
L}^{\ast\circ}(P^0_V\zeta)\quad\forall\zeta\in{\bf B}_\eta(H_V).
$$
\end{theorem}

Now Corollary~\ref{cor:A.2} gives rise to
 \begin{equation}\label{e:4.23}
C_q({\cal L}^\ast, \gamma_0;\K)=C_q(\tilde{\cal L}^\ast,
0;\K)=C_{q-m^-(\gamma_0)}(\tilde{\cal L}^{\ast\circ},
0;\K)\;\quad\forall q=0,1,\cdots,
\end{equation}
and $ C_{\ast}(\tilde{\cal L}^{\ast\circ},
0;\K)=C_{\ast}(\tilde{\cal L}^{\circ}, 0;\K)$ by (\ref{e:4.22}).
From (\ref{e:4.23}) and Theorem~\ref{th:3.8} we arrive at  the
following shifting theorem.

\begin{theorem}\label{th:4.4}
 $C_q({\cal L}, \gamma_0;\K)=C_q(\tilde{\cal L},
0;\K)=C_{q-m^-(\gamma_0)}(\tilde{\cal L}^\circ, 0;\K)\;\forall
q=0,1,\cdots$.
\end{theorem}

By  Theorem~\ref{th:A.5} there exists a ball ${\bf
B}_\mu(X_V)\subset {\bf B}_\tau(X_V)$, an origin-preserving local
homeomorphism $\tilde\varphi$ from ${\bf B}_\mu(X_V)$ to an open
neighborhood of $0$ in $X_V$ with $\tilde\varphi({\bf
B}_\mu(X_V))\subset {\bf B}_r(X_V)$ such that
\begin{equation}\label{e:4.24}
\tilde{\cal L}^\ast_X\circ\tilde\varphi(\zeta)=\frac{1}{2}(\tilde
B^\ast(0)\zeta^\bot, \zeta^\bot)_{W^{1,2}}+ \tilde{\cal
L}^\ast(\tilde h(\zeta^0)+\zeta^0)
\end{equation}
for any $\zeta\in {\bf B}_\mu(X_V)$, where $\zeta^0=P^0_V(\zeta)$
and $\zeta^\bot=\zeta-\zeta^0$. Since $\mu\le\tau$, from
(\ref{e:4.7}) and (\ref{e:4.18})-(\ref{e:4.19}) we derive that
(\ref{e:4.24}) becomes:
\begin{equation}\label{e:4.25}
\tilde{\cal L}_X\circ\tilde\varphi(\zeta)=\frac{1}{2}(\tilde
B(0)\zeta^\bot, \zeta^\bot)_{W^{1,2}}+ \tilde{\cal
L}^\circ(\zeta^0)\quad\forall\zeta\in {\bf B}_\mu(X_V).
\end{equation}
Since $X_V^\star=H_V^\star$, $\star=0,-$, as in the arguments below
Theorem~\ref{th:A.5} the following splitting theorem may be derived
from (\ref{e:4.25}) by changing $\mu>0$ and $\tilde\varphi$
suitably.

\begin{theorem}\label{th:4.5}
Under the notations above, there exists a ball ${\bf
B}_\mu(X_V)\subset {\bf B}_\tau(X_V)$,  an origin-preserving local
homeomorphism $\tilde\varphi$ from ${\bf B}_\mu(X_V)$ to an open
neighborhood of $0$ in $X_V$  such that
\begin{equation}\label{e:4.26}
\tilde{\cal L}_X\circ\tilde\varphi(\zeta)=\frac{1}{2}(\tilde
B(0)P^+_V\zeta, P^+_V\zeta)_{W^{1,2}}-\|P^-_V\zeta\|^2_{W^{1,2}}+
\tilde{\cal L}^\circ(P_V^0\zeta)
\end{equation}
for any $\zeta\in {\bf B}_\mu(X_V)$, where  $\tilde{\cal L}^\circ$
is as in (\ref{e:4.21}).
\end{theorem}

\begin{remark}\label{rm:4.6}
{\rm  It is easily shown  that $H_V=W^{1,2}_{V}(I,
\R^n)=W^{1,2}_{0}(I, \R^n)$ and $X_V=C^1_{V}(I, \R^n)=C^1_{0}(I,
\R^n)$ if $M_0$ and $M_1$ are two disjoint points. For this
case along the proof lines of \cite[Ch.I, Th.5.1]{Ch93},
Caponio-Javaloyes-Masiello proved (\ref{e:4.26}) in
\cite[Th.7]{CaJaMa1}  and hence the shifting theorem
$C_\ast(\tilde{\cal L}_X, 0;\K)=C_{\ast-m^-(\gamma_0)}(\tilde{\cal
L}^\circ, 0;\K)$. When the critical point $0$ is nondegenerate, they
also claimed that (29) of \cite{CaJaMa1}, or equivalently
$C_\ast(\tilde{\cal L}, 0;\K)=C_\ast(\tilde{\cal L}_X, 0;\K)$, can
be obtained with Palais' theorems 16 and 17 in \cite{Pa} as in
\cite{Ch94}. However, a detailed proof of such a claim is not
trivial and was recently given in \cite{CaJaMa2} by combining Chang's
ideas of \cite{Ch83} with the technique of \cite{AbSc1}, and the
nondegenerate assumption of the critical point $0$ was used in an
essential way. }
\end{remark}


\noindent{\bf Step 2}. {\it Complete the proofs  of
Theorem~\ref{th:1.2},~\ref{th:1.3}.}

Note that the differential at $0$ of the chart $\Phi$ in
(\ref{e:4.2}),
$$
d\Phi(0):H_V\to T_{\gamma_0}\Lambda_N(M)=W^{1,2}(\gamma_0^\ast
TM),\;\zeta\mapsto\sum^n_{i=1}\zeta_ie_i
$$
 is a Hilbert space isomorphism and that for any $\zeta\in {\bf B}_{2r}(H_V)$,
$$
{\rm
EXP}_{\gamma_0}(d\Phi(0)[\zeta])(t)=\exp_{\gamma_0(t)}\left((d\Phi(0)[\zeta])(t)\right)=\Phi(\zeta)(t),
$$
 i.e., ${\rm EXP}_{\gamma_0}\circ
d\Phi(0)=\Phi$ on ${\bf B}_{2r}(H_V)$.  Since $\tilde{\cal L}={\cal
L}\circ\Phi$ on ${\bf B}_{2r}(H_V)$ by (\ref{e:4.5}), we get
\begin{equation}\label{e:4.27}
{\cal L}\circ{\rm EXP}_{\gamma_0}\circ d\Phi(0)={\cal
L}\circ\Phi=\tilde{\cal L}\quad\hbox{on}\; {\bf B}_{2r}(H_V)
\end{equation}
and thus $\nabla({\cal L}\circ{\rm
EXP}_{\gamma_0})(d\Phi(0)[\zeta])=d\Phi(0)[\nabla\tilde{\cal
L}(\zeta)]\;\forall\zeta\in {\bf B}_{2r}(H_{V})$. It follows that
\begin{eqnarray}
&&{\cal A}(d\Phi(0)[\zeta])=d\Phi(0)[\tilde
A(\zeta)]\quad\forall\zeta\in {\bf
B}_{2r}(H_{V})\cap X_{V},\label{e:4.28}\\
&&d{\cal A}(0)\circ d\Phi(0)=d\Phi(0)\circ d\tilde
A(0)\label{e:4.29}
\end{eqnarray}
because the Hilbert space isomorphism $d\Phi(0):H_{V}\to
T_{\gamma_0}\Lambda_N(M)$  induces a Banach space isomorphism from
$X_{V}$ to $T_{\gamma_0}{\cal X}=T_{\gamma_0}C^1_N(I,M)$.
(\ref{e:4.29}) implies that $d\Phi(0)(H_V^\star)= {\bf
H}^\star(d^2{\cal L}|_{\cal X}(\gamma_0))$ for $\star=-,0,+$, and so
\begin{eqnarray}\label{e:4.30}
d\Phi(0)\circ P_V^\star=P^\star\circ d\Phi(0),\quad \star=-,0,+.
\end{eqnarray}
Shrinking $\delta>0$ above (\ref{e:1.6}) so that $\delta<\tau$, for
$\xi\in {\bf B}_\delta\bigl({\bf H}^0(d^2{\cal L}|_{\cal
X}(\gamma_0))\bigr)$ we derive
\begin{eqnarray*}
0&=&d\Phi(0)\circ (I-P^0_V)\tilde A(d\Phi(0)^{-1}\xi+ \tilde
h(d\Phi(0)^{-1}\xi))\\
&=& (I-P^0)\circ d\Phi(0)\circ\tilde A\circ d\Phi(0)^{-1}\bigl(\xi+
d\Phi(0)\circ\tilde
h(d\Phi(0)^{-1}\xi)\bigr)\\
&=& (I-P^0)\circ {\cal A}\bigl(\xi+ d\Phi(0)\circ\tilde
h(d\Phi(0)^{-1}\xi)\bigr)
\end{eqnarray*}
from (\ref{e:4.30}), (\ref{e:4.28}) and (\ref{e:4.6}). But we know that $(I-P^0){\cal A}(\xi+ h(\xi))=0$ for any $\xi\in {\bf
B}_\delta\bigl({\bf H}^0(d^2{\cal L}|_{\cal X}(\gamma_0))\bigr)$
by the formula above (\ref{e:1.6}). By the uniqueness of $h$ there,
$h(\xi)=d\Phi(0)\circ\tilde h(d\Phi(0)^{-1}\xi)$ and hence
\begin{eqnarray}\label{e:4.31}
{\cal L}^\circ(d\Phi(0)[\zeta])&=&{\cal L}\circ{\rm
EXP}_{\gamma_0}\bigl(d\Phi(0)[\zeta]+ h(d\Phi(0)[\zeta])\bigr)\nonumber\\
&=&{\cal L}\circ{\rm EXP}_{\gamma_0}\circ d\Phi(0)\bigl(\zeta+
d\Phi(0)^{-1}\circ h(d\Phi(0)[\zeta])\bigr)\nonumber\\
&=&{\cal L}\circ{\rm EXP}_{\gamma_0}\circ d\Phi(0)\bigl(\zeta+
\tilde
h(\zeta)\bigr)\nonumber\\
&=&\tilde{\cal L}\bigl(\zeta+ \tilde h(\zeta)\bigr)=\tilde{\cal
L}^\circ(\zeta)\quad\forall\zeta\in {\bf B}_\delta(H_V^0)
\end{eqnarray}
by (\ref{e:1.6}), (\ref{e:4.27}) and the definition of $\tilde{\cal
L}^\circ$ in (\ref{e:4.21}). Hence $C_\ast({\cal L}^\circ,
0;\K)=C_\ast(\tilde{\cal L}^\circ,0;\K)$. (That is,
 this and Theorem~\ref{th:4.4}  give the first equality in Theorem~\ref{th:1.4} too.)

\underline{In order to prove Theorem~\ref{th:1.2}}, by
Claims~\ref{cl:3.5},~\ref{cl:3.7} and Theorem~\ref{th:3.8} it
suffices to prove (iii)-(iv). Shrink $\delta>0$ above (\ref{e:1.6})
so that $\delta\le\eta$,  where $\eta$ is the radius
of the ball in Theorem~\ref{th:4.3}. Note that $d\Phi(0)\left({\bf
B}_\delta(H_V)\right)={\bf B}_\delta(T_{\gamma_0}\Lambda_N(M))$ and
(\ref{e:4.29}) implies that
\begin{eqnarray}\label{e:4.32}
(P_V^\star\zeta, P_V^\star\zeta)_{W^{1,2}}&=&\langle d\Phi(0)\circ
P_V^\star\zeta, d\Phi(0)\circ
P_V^\star\zeta\rangle_{1}\nonumber\\
&=&\langle P^\star\circ d\Phi(0)\zeta, P^\star\circ
d\Phi(0)\zeta\rangle_{1}
\end{eqnarray}
for $\zeta\in {\bf B}_\delta(H_V)\subset {\bf B}_\eta(H_V)$ and
$\star=+,-$. Moreover, as in (\ref{e:4.27}) and (\ref{e:4.31}) we have
$$
{\cal L}^\ast\circ{\rm EXP}_{\gamma_0}\circ d\Phi(0)={\cal
L}^\ast\circ\Phi=\tilde{\cal L}^\ast\quad\hbox{on}\; {\bf
B}_{2r}(H_V).
$$
Define $\psi: {\bf B}_\delta(T_{\gamma_0}\Lambda_N(M))\to
T_{\gamma_0}\Lambda_N(M)$ by $\psi=d\Phi(0)\circ\tilde\psi\circ
d\Phi(0)^{-1}$. For $\xi\in {\bf
B}_\delta(T_{\gamma_0}\Lambda_N(M))$ and $\zeta=d\Phi(0)^{-1}\xi$ we
get
\begin{eqnarray*}
&&{\cal L}^\ast\circ{\rm EXP}_{\gamma_0}\circ\psi(\xi)={\cal
L}^\ast\circ{\rm EXP}_{\gamma_0}\circ d\Phi(0)\circ\tilde\psi\circ
d\Phi(0)^{-1}\xi=\tilde{\cal
L}^\ast\circ\tilde\psi(\zeta),\\
&&{\cal L}^\circ(P^0\xi)={\cal L}^\circ(P^0\circ
d\Phi(0)\zeta)={\cal L}^\circ(d\Phi(0)\circ P^0_V\zeta)=\tilde{\cal
L}^\circ(P^0_V\zeta)
\end{eqnarray*}
by (\ref{e:4.31}). These, (\ref{e:4.32}), (\ref{e:4.22}) and Theorem~\ref{th:4.3}
give Theorem~\ref{th:1.2}(iii).

As to Theorem~\ref{th:1.2}(iv),  for any open neighborhood $W$ of
$0$ in  $U={\bf B}_{2r}(H_V)$ and a field $\K$, writing  $W_X=W\cap
X_V$ as  an open subset of $X_V$ and using Theorem~\ref{th:A.9} we
deduce that the inclusion
$$
\big(\tilde{\cal L}^\ast_c\cap W_X, \tilde{\cal L}^\ast_c\cap W_X
\setminus\{0\}\big)\hookrightarrow\big(\tilde{\cal L}^\ast_c\cap W,
\tilde{\cal L}^\ast_c\cap W\setminus\{0\}\big)
$$
 induces isomorphisms
$$H_\ast\big(\tilde{\cal L}^\ast_c\cap W_X,
\tilde{\cal L}^\ast_c\cap W_X\setminus\{0\};\K\big)\cong
H_\ast\big(\tilde{\cal L}^\ast_c\cap W, \tilde{\cal L}^\ast_c\cap
W\setminus\{0\};\K\big).
$$
The expected conclusion follows from this immediately.

\underline{Finally, let us prove Theorem~\ref{th:1.3}}. Since
$d\Phi(0)$ is a Banach isomorphism from $X_V$ to $T_{\gamma_0}{\cal
X}=C^1_{TN}(\gamma_0^\ast TM)$ we may choose $\epsilon>0$ such that
$$
d\Phi(0)^{-1}\left({\bf B}_\epsilon(T_{\gamma_0}{\cal
X})\right)\subset  {\bf B}_\mu(X_V).
$$
By (\ref{e:4.29}), $d{\cal A}(0)\circ d\Phi(0)=d\Phi(0)\circ d\tilde
A(0)= d\Phi(0)\circ \tilde B(0)$. So for  $\xi\in {\bf
B}_\epsilon(T_{\gamma_0}{\cal X})$ and $\zeta=[d\Phi(0)]^{-1}\xi$ we
have
\begin{eqnarray*}
 d^2{\cal L}|_{\cal X}(\gamma_0)[\xi^+,
\xi^+]&=&d^2{\cal L}|_{\cal X}(\gamma_0)[P^+\xi,
P^+\xi]\\
&=&\langle d{\cal A}(0)[P^+\xi], P^+\xi\rangle_{1}\\
&=&\langle d{\cal A}(0)\circ d\Phi(0)\circ P^+_V\zeta, d\Phi(0)\circ
P^+_V\zeta\rangle_{1}\\
&=&\langle d\Phi(0)\circ\tilde B(0)\zeta^+, d\Phi(0)
\zeta^+\rangle_{1}=(\tilde B(0)\zeta^+,  \zeta^+)_{W^{1,2}},\\
\|P^-\xi\|^2_{1}&=&\langle P^-\xi,
P^-\xi\rangle_{1}=\langle P^-\circ d\Phi(0)\zeta, P^-\circ d\Phi(0)\zeta\rangle_{1}\\
&=&\langle d\Phi(0)\circ P^-_V\zeta, d\Phi(0)\circ P^-_V\zeta\rangle_{1}\\
&=&(P^-_V\zeta,  P^-_V\zeta)_{W^{1,2}}=\|P^-_V\zeta\|^2_{W^{1,2}}
\end{eqnarray*}
by (\ref{e:4.30}). Define $\varphi: {\bf
B}_\epsilon(T_{\gamma_0}{\cal X})\to T_{\gamma_0}{\cal X}$ by
$\varphi=d\Phi(0)\circ\tilde\varphi\circ d\Phi(0)^{-1}$. As above
these and Theorem~\ref{th:4.5} yield Theorem~\ref{th:1.3}.
 \hfill$\Box$\vspace{2mm}

\section{Proofs of Theorems~\ref{th:1.5},~\ref{th:1.6},~\ref{th:1.7}, ~\ref{th:1.8}}\label{sec:5}
\setcounter{equation}{0}

\subsection{Proofs of Theorems~\ref{th:1.5},~\ref{th:1.6},~\ref{th:1.7}}\label{sec:5.1}

By  Claim~\ref{cl:3.9} and Theorem~\ref{th:3.10},
Theorem~\ref{th:1.5}(i)-(ii) are clear.

When ${\cal L}$ is replaced by ${\cal L}^\ast$ in the statement of
Theorem~\ref{th:1.2}, we have the corresponding maps $A^\ast_x$ and
$B^\ast_x$ for $x\in{\cal O}$. Shrinking $\varepsilon>0$ if
necessary, by Corollary~\ref{cor:2.3} we have
\begin{equation}\label{e:5.1}
{A}^\ast_{x}(v)= {A}_{x}(v)\quad\forall (x,v)\in XN{\cal
O}(\varepsilon)_x
\end{equation}
  and hence ${B}_x^\ast= {B}_x$
 for any $x\in{\cal O}$.
The latter implies that
\begin{equation}\label{e:5.2}
{\bf H}^\star({B}^\ast)={\bf H}^\star({B}),\quad \star=+,0,-.
\end{equation}
It follows from this, (\ref{e:5.1}) and (\ref{e:1.20}) that the
$C^1$ map ${\mathfrak{h}}_x$ in (\ref{e:1.19}) satisfies
\begin{eqnarray*}
 ({\bf P}^+_x + {\bf P}^-_x)\circ {A}^\ast_x\bigl(v+ {\mathfrak{h}}_x(v)\bigr)=0
\quad\forall v\in {\bf H}^0(B)(\epsilon)_x
\end{eqnarray*}
(by shrinking $\epsilon\in (0,\varepsilon)$ if necessary). Define
the functional
\begin{eqnarray}\label{e:5.3}
{\cal L}^{\ast\circ}_\triangle:{\bf H}^0({B}^\ast)(\epsilon)\ni
(x,v)\to{\cal L}^\ast\circ{\rm EXP}_x\bigl(v+
 {\mathfrak{h}}_x(v)\bigr)\in\R.
\end{eqnarray}
By Corollary~\ref{cor:2.3} and (\ref{e:1.21}) we have
\begin{eqnarray}\label{e:5.4}
{\cal L}_{\triangle}^{\ast\circ}(x,v)={\cal L}_{\triangle
x}^{\ast\circ}(v)&=&{\cal L}^\ast\circ{\rm EXP}_x\bigl(v+
 {\mathfrak{h}}_x(v)\bigr)\nonumber\\
 &=&{\cal L}\circ{\rm EXP}_x\bigl(v+
 {\mathfrak{h}}_x(v)\bigr)={\cal
L}^{\circ}_{\triangle x}(v)={\cal L}^{\circ}_{\triangle}(x, v)
\end{eqnarray}
 for any $v\in {\bf H}^0({B})(\epsilon)_x={\bf
 H}^0({B}^\ast)(\epsilon)_x$.
   Hence ${\cal L}_{\triangle}^{\ast\circ}$ is $C^1$ and has the isolated critical orbit ${\cal O}$.
 Moreover, ${\cal L}^{\ast\circ}_{\triangle x}$ is $C^2$, $S_x^1$-invariant and
 \begin{equation}\label{e:5.5}
C_{q}({\cal L}^\circ_{\triangle x}, 0;\K)=C_{q}({\cal
L}^{\ast\circ}_{\triangle x}, 0;\K)\quad\forall x\in{\cal
O},\;q\in\N\cup\{0\}.
\end{equation}

As below Lemma~\ref{lem:4.2}, for every $x\in {\cal O}$ there exists
a continuous map ${\bf B}^\ast_x:{N}{\cal O}(\epsilon)_x\cap
X{N}{\cal O}_x\to L_s({N}{\cal O}_x)$ with respect to the topology
of $X{N}{\cal O}_x$ such that
$$
\langle d{A}^\ast_{x}(\zeta)[\xi], \eta\rangle_{1}=d^2 {\cal
F}^{\ast X}_{x}
  (\zeta)[\xi,\eta]=\langle{\bf B}_x^\ast(\zeta)\xi,
  \eta\rangle_{1}
$$
for  any $\zeta\in{N}{\cal O}(\epsilon)_x\cap X{N}{\cal O}_x$ and
$\xi, \eta\in {N}{\cal O}_x$. Clearly, ${\bf
B}_x^\ast(0)={B}_x^\ast$.  The proof of the following proposition
will be postponed  to Section~\ref{sec:5.3}.

\begin{proposition}\label{prop:5.1}
 $({N}{\cal O}_{\gamma_{0}}, X{N}{\cal O}_{\gamma_{0}},
 {\cal F}^\ast_{\gamma_{0}},
  {A}^\ast_{\gamma_{0}}, {\bf B}^\ast_{\gamma_{0}})$ satisfies
the conditions of Theorem~\ref{th:A.1} (and hence
Theorem~\ref{th:A.5}) around the critical point $0=0_{\gamma_{0}}$.
\end{proposition}

Clearly, Proposition~\ref{prop:5.1} implies that $(N{\cal O}_{x},
X{N}{\cal O}_{x}, {\cal F}^\ast_x,
  {A}^\ast_{x}, {\bf B}^\ast_{x})$ satisfies
the conditions of Theorems~\ref{th:A.1},~\ref{th:A.5} around $0=0_x$
for each $x\in{\cal O}$. (Indeed, let $s\in S^1$ such that
$x=s\cdot\gamma_{0}$. Using the Hilbert isomorphism ${N}{\cal
O}_{\gamma_{0}}\to{N}{\cal O}_{x}$ and
Theorems~\ref{th:A.4},~\ref{th:A.8} one easily proves them.)

By Theorem~\ref{th:A.1} and (\ref{e:5.1}) we obtain a
$S_x^1$-invariant open neighborhood ${U}_{x}$ of $0_x$ in ${N}{\cal
O}(\varepsilon)_x$ and a $S_x^1$-equivariant origin-preserving
homeomorphism
\begin{equation}\label{e:5.6}
{\Upsilon}_{x}: {N}{\cal O}(\epsilon)_x\to {U}_{x}
\end{equation}
(shrinking $\epsilon\in (0,\varepsilon)$ if necessary),   such that
\begin{equation}\label{e:5.7}
{\cal F}^\ast_x\circ{\Upsilon}_{x}(u)=\|{\bf P}_{x}^+u\|^2_1-\|{\bf
P}_{x}^-u\|^2_1+ {\cal F}^\ast_x({\bf P}_{x}^0u+
{\mathfrak{h}}_x({\bf P}_{x}^0u))
\end{equation}
for all $u\in {N}{\cal O}(\epsilon)_x$. $\Upsilon_x$ also maps
$({\bf P}^-_x+{\bf P}^0_x){N}{\cal O}(\epsilon)_x$ into $X{N}{\cal
O}_x$ and is a homeomorphism from $({\bf P}^-_x+{\bf P}^0_x){N}{\cal
O}(\epsilon)_x$ to $\Upsilon_x\left(({\bf P}^-_x+{\bf
P}^0_x){N}{\cal O}(\epsilon)_x\right)$
 even if the topology on the latter is taken as the induced one
 by $XN{\cal O}_x$.
Moreover, the $S_x^1$-invariant functional
\begin{equation}\label{e:5.8}
{\cal F}_{x}^{\ast\circ}:{\bf H}^0({B}^\ast)(\epsilon)_x\to\R,\;
z\mapsto {\cal F}^\ast_x(z+ {\mathfrak{h}}_x(z)),
\end{equation}
is $C^2$, has the isolated critical point $0_x\in {\bf
H}^0({B}^\ast)_x$ and $d^2{\cal F}_{x}^{\ast\circ}(0_x)=0$. Observe
that
\begin{eqnarray}\label{e:5.9}
{\cal L}^{*\circ}_{\triangle x}(v)&=&{\cal L}^\ast\circ{\rm
EXP}_x\bigl(v+
 {\mathfrak{h}}_x(v)\bigr)={\cal F}^\ast_x\bigl(v+
 {\mathfrak{h}}_x(v)\bigr)\nonumber\\
 &=&{\cal
F}_{x}^{\ast\circ}(v)\quad\forall v\in{\bf
H}^0({B}^\ast)(\epsilon)_x
\end{eqnarray}
because of (\ref{e:5.3}) and the definition of ${\cal F}^\ast$ at
the beginning of this section.

Let ${U}=\cup_{x\in{\cal O}}U_x$. It is an $S^1$-invariant tubular
open neighborhood of the zero section of ${N}{\cal O}$, and in fact
a fiber bundle over ${\cal O}$.
 Define maps $\Upsilon: {N}{\cal O}(\epsilon)\to {U}$ and
$$
{\mathfrak{h}}:{\bf H}^0({B}^\ast)(\epsilon)\to {\bf
H}^-({B}^\ast)\oplus({\bf H}^+({B}^\ast)\cap X{N}{\cal O})
$$
by $\Upsilon|_{{N}{\cal O}(\epsilon)_x}=\Upsilon_x$ and
${\mathfrak{h}}|_{{\bf H}^0({B}^\ast)(\epsilon)_x}={\mathfrak{h}}_x$
for any $x\in{\cal O}$, respectively. Here both ${\mathfrak{h}}_x$
and ${\Upsilon}_{x}$ are given by (\ref{e:1.19}) and (\ref{e:5.6}),
respectively. Since
\begin{eqnarray*}
&&{A}^\ast_{s\cdot x}(s\cdot v)= s\cdot {A}^\ast_x(v)\quad\forall
s\in S^1,\; v\in
N{\cal O}(\varepsilon)_x\cap XN{\cal O}_x,\\
&&{B}^\ast_{s\cdot x}(s\cdot\xi,
s\cdot\eta)={B}^\ast_x(\xi,\eta)\quad\forall s\in S^1,\;x\in{\cal
O},\;\xi,\eta\in N{\cal O}_x,
\end{eqnarray*}
as in the proofs of Theorem~7.3 and Corollary~7.1 in \cite[page
72]{Ch93} it follows from these, (\ref{e:1.19})-(\ref{e:1.20}) and
(\ref{e:5.4}) that ${\mathfrak{h}}$ is an $S^1$-equivariant
fiber-preserving  $C^1$ map  and that $\Upsilon$ is an $S^1$-equivariant fiber-preserving homeomorphism from ${N}{\cal
O}(\epsilon)$ onto ${U}$. \textsf{So  the first claim in
Theorem~\ref{th:1.5}(iii) is proved.} The second may be derived from
either the above arguments or the corresponding conclusion in
Theorem~A.2 of \cite{Lu2, Lu3}.

By completely similar arguments \textsf{we may use
Theorem~\ref{th:A.5} to derive Theorem~\ref{th:1.6}.}

\underline{Now we prove Theorem~\ref{th:1.5}(iv)-(v) and
Theorem~\ref{th:1.7}.}

\textsf{Firstly, we prove the latter.} If $m^0({\cal O})=m^-({\cal
O})=0$, by Theorem~\ref{th:1.5}(iii) we have ${\cal
F}^\ast\circ{\Upsilon}(u)=c+\|u\|^2_1$ for all $u\in {N}{\cal
O}(\epsilon)$, and hence
\begin{eqnarray*}
&&C_\ast({\cal L}, {\cal O};\K)=C_\ast({\cal F}^\ast, {\cal
O};\K)=C_\ast({\cal F}^\ast\circ{\Upsilon}, {\cal
O};\K)\\
&=&H_\ast(\{{\cal F}^\ast\circ{\Upsilon}<c\}\cup {\cal O}, \{{\cal
F}^\ast\circ{\Upsilon}<c\};\K)=H_\ast ({\cal O};\K) =H_\ast (S^1;\K)
\end{eqnarray*}
because  $\gamma_0$ is nonconstant.
 If either $m^0({\cal
O})=0$ and $m^-({\cal O})>0$ or $m^0({\cal O})>0$ and $m^-({\cal
O})=0$ one easily sees the desired conclusions from the following
proof in the case $m^0({\cal O})>0$ and $m^-({\cal O})>0$. For the
proof of the final case it suffices to prove:

\begin{claim}\label{cl:5.2}
Let $\K$ be a field of characteristic $0$ or prime up to order
$|S^1_{\gamma_{0}}|$ of $S^1_{\gamma_{0}}$. For any $x\in{\cal O}$
and $q=0,1,\cdots$, it holds that
\begin{eqnarray*}
&&C_q({\cal F}^\ast, {\cal O};\K)\\
&=& \Bigl(H_{m^-({\cal O})}({\bf H}^-({B})_x, {\bf
H}^-({B})_x\setminus\{0_x\};\K)\otimes C_{q-m^-({\cal O})}({\cal
F}^{\ast\circ}_{x},
0;\K)\Bigr)^{S^1_x}\\
&\oplus& \Bigl(H_{m^-({\cal O})}({\bf H}^-({B})_x, {\bf
H}^-({B})_x\setminus\{0_x\};\K)\otimes C_{q-m^-({\cal O})-1}({\cal
F}^{\ast\circ}_{x}, 0;\K)\Bigr)^{S^1_x}.
\end{eqnarray*}
\end{claim}

\noindent{\bf Proof}.
 Since the map $S^1\times {N}{\cal O}_x\to {N}{\cal O}$ with $(s,v)\mapsto s\cdot v$
is a normal covering with group of covering transformations $S^1_x$
(\cite[page 500]{GM2}), so is the map $S^1\times {W}_x\to {W}$ with
$(s,v)\mapsto s\cdot v$ for any subset ${W}\subset {N}{\cal O}$ with
properties
\begin{equation}\label{e:5.10}
{W}_{s\cdot x}=s\cdot {W}_x\quad\forall x\in{\cal O}, s\in S^1,
\end{equation}
where ${W}_x={W}\cap {N}{\cal O}_x$, Then ${W}=
(S^1\times{W}_x)/S^1_x$, where $S_x^1$ acts on $S^1\times{W}_x$ by
covering transformations as described above. Note that
 \begin{eqnarray*}
 {W}:= {N}{\cal
O}(\epsilon),\quad{W}\cap \{{\cal F}^\ast\le c\}
\quad\hbox{and}\quad ({W}\setminus{\cal O})\cap \{{\cal F}^\ast\le
c\}
\end{eqnarray*}
 satisfy (\ref{e:5.10}). It follows that
\begin{eqnarray}\label{e:5.11}
&&C_q({\cal F}^\ast, {\cal O};\K)=H_q\left({W}\cap \{{\cal
F}^\ast\le
c\}, ({W}\setminus{\cal O})\cap \{{\cal F}^\ast\le c\};\K\right)\nonumber\\
&=&H_q\Bigl((S^1\times ({W}_x\cap \{{\cal F}_x^\ast\le c\}))/S_x^1,
(S^1\times (({W}_x\setminus\{0_x\})\cap \{{\cal F}_x^\ast\le
c\}))/S_x^1;\K\Bigr)\nonumber\\
&=&H_q\bigl((S^1\times\triangle_x)/S_x^1,
(S^1\times\triangle^\prime_x)/S_x^1;\K\bigr),
\end{eqnarray}
where $\triangle_x= {N}{\cal O}(\epsilon)_x\cap\{{\cal
F}^\ast_x\circ\Upsilon_x\le c\}$ and $\triangle^\prime_x=  ({N}{\cal
O}(\epsilon)_x \setminus\{0_x\})\cap\{{\cal
F}^\ast_x\circ\Upsilon_x\le c\}$. The final equality in
(\ref{e:5.11}) comes from the fact that $\Upsilon_x^{-1}$  is a
$S_x^1$-equivariant homeomorphism from
\begin{eqnarray*}
\bigl(S^1\times ({W}_x\cap \{{\cal F}_x^\ast\le c\}), S^1\times
(({W}_x\setminus\{0_x\})\cap \{{\cal F}_x^\ast\le c\})\bigr)
\end{eqnarray*}
to $(S^1\times\triangle_x, S^1\times\triangle^\prime_x)$ by
(\ref{e:5.7}). Let ${N}{\cal O}(\epsilon)^{-0}_x:={N}{\cal
O}(\epsilon)_x\cap ({\bf H}^-({B})_x^-\oplus {\bf H}^0({B})_x)$,
which is a finite dimensional $C^3$-smooth manifold contained in
$X{N}{\cal O}$. Define
\begin{equation}\label{e:5.12}
({\cal F}_x^\ast\circ\Upsilon_x)^{-0}:{N}{\cal
O}(\epsilon)^{-0}_x\to\R
\end{equation}
 by
$({\cal F}_x^\ast\circ\Upsilon_x)^{-0}(v^-+v^0)= -\|v^-\|^2_1+ {\cal
F}^\ast_x(v^0+ {\mathfrak{h}}_x(v^0))= -\|v^-\|^2_1+ {\cal
F}^{\ast\circ}_x(v^0)$. It is $C^2$ because of (\ref{e:5.8}).
Observe that $(S^1\times\triangle_x, S^1\times\triangle^\prime_x)$
can be retracted $S^1_x$-equivariantly into $(S^1\times
\triangle_x^{-0}, S^1\times \triangle_x^{\prime-0})$, where
$\triangle_x^{-0}={N}{\cal O}(\epsilon)^{-0}_x\cap\{({\cal
F}^\ast_x\circ\Upsilon_x)^{-0}\le c\}$,
\begin{eqnarray*}
\triangle_x^{\prime-0}=({N}{\cal O}(\epsilon)^{-0}_x
\setminus\{0_x\})\cap\{({\cal F}^\ast_x\circ\Upsilon_x)^{-0}\le c\}.
\end{eqnarray*}
From this and (\ref{e:5.11}) we derive that
\begin{eqnarray}\label{e:5.13}
C_q({\cal F}^\ast, {\cal
O};\K)=H_q\bigl((S^1\times\triangle_x^{-0})/S_x^1,
(S^1\times\triangle^{\prime-0}_x)/S_x^1;\K\bigr).
\end{eqnarray}
 Using either Satz~6.6 on the page 57
of \cite{Ra} or the proof of  Lemma~3.6 in \cite{BLo} (with
Theorem~7.2 on the page 142 of \cite{Br}) we know that the transfer
is an isomorphism
$$
H_\ast\bigl((S^1\times\triangle_x^{-0})/S_x^1,
(S^1\times\triangle^{\prime-0}_x)/S_x^1;\K\bigr)=H_\ast\bigl(S^1\times\triangle_x^{-0},
S^1\times\triangle^{\prime-0}_x;\K\bigr)^{S_x^1}.
$$
This, (\ref{e:5.13}) and the K\"unneth formula lead to
\begin{eqnarray}\label{e:5.14}
C_q({\cal F}^\ast, {\cal O};\K)=H_{q-1}\bigl(\triangle_x^{-0},
\triangle^{\prime-0}_x;\K\bigr)^{S_x^1}\oplus
H_q\bigl(\triangle_x^{-0}, \triangle^{\prime-0}_x;\K\bigr)^{S_x^1}.
\end{eqnarray}
As in the proof of Theorem~5.5  on the page 51 of \cite{Ch93} we may
obtain
\begin{eqnarray*}
&&H_q\bigl(\triangle_x^{-0},
\triangle^{\prime-0}_x;\K\bigr)\\
&=&\hspace{-2mm}H_{m^-({\cal O})}\bigl({\bf H}^-({B})(\epsilon)_x,
\partial {\bf H}^-({B})(\epsilon)_x;\K\bigr)\otimes
C_{q-m^-({\cal O})}\bigl({\cal F}_x^{\ast\circ}, 0_x;\K\bigr)\\
\hspace{-2mm}&=&\hspace{-2mm}H_{m^-({\cal O})}\bigl({\bf
H}^-({B})_x, {\bf H}^-({B})_x\setminus\{0_x\};\K\bigr)\otimes
C_{q-m^-({\cal O})}\bigl({\cal F}_x^{\ast\circ}, 0_x;\K\bigr)
\end{eqnarray*}
for all $q=0,1,\cdots$. (See the proof of Lemma~4.6 in \cite{Lu0}).
This and (\ref{e:5.14}) give Claim~\ref{cl:5.2} immediately.
\hfill$\Box$\vspace{2mm}

Now by Theorem~\ref{th:1.5}(ii) for $q\in\N\cup\{0\}$ we have
$$
C_q({\cal L}, {\cal O};{\K})=C_q({\cal L}^\ast, {\cal
O};{\K})=C_q({\cal F}^\ast, {\cal O};{\K}).
$$
Moreover, (\ref{e:5.4}) and (\ref{e:5.9}) imply ${\cal
L}^{\circ}_{\triangle x}(v)={\cal L}_{\triangle
x}^{\ast\circ}(v)={\cal F}_{x}^{\ast\circ}(v)$ for any $x\in{\cal
O}$ and $v\in {\bf H}^0({B})(\epsilon)_x={\bf
H}^0({B}^\ast)(\epsilon)_x$.
 \textsf{Theorem~\ref{th:1.7} follows from these and Claim~\ref{cl:5.2}
  immediately.}\vspace{2mm}

\noindent{\bf Another proof of Claim~\ref{cl:5.2}}.  By
Theorem~\ref{th:1.6}, for any $q\in\N\cup\{0\}$,
\begin{eqnarray}\label{e:5.15}
C_q({\cal F}^X, {\cal O};\K)=C_q({\cal F}^X\circ{\Psi}, {\cal O};\K)
=C_q(({\cal F}^X\circ{\Psi})^{-0}, {\cal O};\K),
\end{eqnarray}
where $({\cal F}^X\circ{\Psi})^{-0}:{N}{\cal O}^{-0}(\rho)={N}{\cal
O}(\rho)\cap\bigl({\bf H}^-({B}){\oplus}{\bf H}^0({B})\bigr)\to\R$
is defined by
\begin{equation}\label{e:5.16}
({\cal F}^X\circ{\Psi})^{-0}(x,v)=-\|{\bf P}^-_xv\|^2_1+ {\cal
L}^\circ_{\triangle x}({\bf P}^0_xv)
\end{equation}
for a small $\rho\in (0,\epsilon)$, and the second equality in
(\ref{e:5.15}) is obtained by the standard deformation method as
done above (3.6) of \cite{Lu2}. Taking ${W}$ to be $\{({\cal
F}^X\circ{\Psi})^{-0}\le c\}$ or $\{({\cal F}^X\circ{\Psi})^{-0}\le
c\}\setminus\{{\cal O}\}$ we have ${W}=(S^1\times{W}_x)/S_x^1$  and
hence
\begin{eqnarray*}
&C_q(({\cal F}^X\circ{\Psi})^{-0}, {\cal
O};\K)\\
&\hspace{-3mm}=H_q\left(\bigl(S^1\times\{({\cal
F}^X\circ{\Psi})^{-0}_x\le c\}\bigr)/S_x^1, \bigl(S^1\times(\{({\cal
F}^X\circ{\Psi})^{-0}_x\le
c\}\setminus\{0_x\})\bigr)/S_x^1;\K\right).
\end{eqnarray*}
Almost repeating the arguments below (\ref{e:5.13}) in the proof of
Claim~\ref{cl:5.2}  we arrive at

\begin{claim}\label{cl:5.3}
Let $\K$ be a field of characteristic $0$ or prime up to order
$|S^1_{\gamma_{0}}|$ of $S^1_{\gamma_{0}}$. For any $x\in{\cal O}$
and $q=0,1,\cdots$,
 it holds that
\begin{eqnarray*}
&&C_q({\cal F}^X, {\cal O};\K)\\
&=& \Bigl(H_{m^-({\cal O})}({\bf H}^-({B})_x, {\bf
H}^-({B})_x\setminus\{0_x\};\K)\otimes C_{q-m^-({\cal O})}({\cal
L}^\circ_{\triangle x},
0;\K)\Bigr)^{S^1_x}\\
&\oplus& \Bigl(H_{m^-({\cal O})}({\bf H}^-({B})_x, {\bf
H}^-({B})_x\setminus\{0_x\};\K)\otimes C_{q-m^-({\cal O})-1}({\cal
L}^\circ_{\triangle x}, 0;\K)\Bigr)^{S^1_x}.
\end{eqnarray*}
\end{claim}

As in the reasoning of (\ref{e:5.15}), using
Theorem~\ref{th:1.5}(iii) we derive
$$
C_q({\cal F}^\ast, {\cal O};\K)=C_q({\cal F}^\ast\circ{\Upsilon},
{\cal O};\K) =C_q(({\cal F}^\ast\circ{\Upsilon})^{-0}, {\cal O};\K)
 $$
for each $q\in\N\cup\{0\}$, where $({\cal
F}^\ast\circ{\Upsilon})^{-0}$ is defined by (\ref{e:5.12}). By
(\ref{e:5.9}) and (\ref{e:5.4}), ${\cal F}^{\ast\circ}(x,v)={\cal
L}_{\triangle}^{\ast\circ}(x,v)={\cal L}^{\circ}_{\triangle }(x,v)$
for $(x,v)\in{N}{\cal O}^{-0}(\rho)$ if $\rho>0$ is small enough. So
\begin{eqnarray*}
 ({\cal
F}^\ast\circ{\Upsilon})^{-0}(x,v^-+v^0) =-\|v^-\|^2_1+ {\cal
L}^\circ_{\triangle x}(v^0) =({\cal F}^X\circ{\Psi})^{-0}(x,v^-+
v^0)
\end{eqnarray*}
for $(x,v^-+v^0)\in{N}{\cal O}^{-0}(\rho)$ by (\ref{e:5.16}). It
follows that
$$
C_q(({\cal F}^\ast\circ{\Upsilon})^{-0}, {\cal O};\K)=C_q(({\cal
F}^X\circ{\Psi})^{-0}, {\cal O};\K)
$$
 and thus
$C_q({\cal F}^\ast, {\cal O};\K)=C_q({\cal F}^X, {\cal
O};\K)\;\forall q\in\N\cup\{0\}$ by (\ref{e:5.15}). This and
Claim~\ref{cl:5.3} lead to
Claim~\ref{cl:5.2}.\hfill$\Box$\vspace{2mm}

\noindent\underline{Proof of Theorem~\ref{th:1.5}(iv).} For each
$q\in\N\cup\{0\}$ we have $C_q({\cal F}^{\ast}, {\cal
O};\K)=C_q({\cal F}^X, {\cal O};\K)$ by
Claims~\ref{cl:5.2},~\ref{cl:5.3}, and $C_q({\cal F}^{\ast X}, {\cal
O};\K)=C_q({\cal F}^X, {\cal O};\K)$ by  Theorem~\ref{th:1.5}(i) and
the excision property of relative homology groups. Hence $C_q({\cal
F}^\ast, {\cal O};\K)=C_q({\cal F}^{\ast X}, {\cal O};\K)$ for any
$q\in\N\cup\{0\}$. This can also be obtained as a direct consequence
of shifting theorems for ${\cal F}^\ast$ and ${\cal F}^{\ast X}$ at
${\cal O}$.

\begin{claim}\label{cl:5.4}
For a field $\K$,  the group $C_q({\cal F}^{\ast X}, {\cal
O};\K)=C_q({\cal F}^\ast, {\cal O};\K)$ is a finite dimensional
vector space over $\K$ for each $q\in\N\cup\{0\}$.
\end{claim}

\noindent{\bf Proof}. As in the proof of \cite[(7.15)]{BiFrPi} we use Proposition~A.1 in \cite{BiFrPi}
to derive that
\begin{eqnarray*}
\dim C_q({\cal F}^\ast, {\cal O};\K)\le 2\left(\dim C_{q-m^-({\cal
O})}({\cal F}^{\ast\circ}_{x}, 0;\K)+ \dim C_{q-m^-({\cal
O})-1}({\cal F}^{\ast\circ}_{x}, 0;\K) \right)
\end{eqnarray*}
for any $x\in{\cal O}$ and $q\in\N\cup\{0\}$. By
\cite[Remark~4.6]{Lu3} or \cite[Remark~2.24]{Lu2} we know that $\dim
C_q({\cal F}^{\ast\circ}_{x},0;\K)<\infty\;\forall q$ and $\dim
C_q({\cal F}^{\ast\circ}_{x},0;\K)=0$ for almost all $q$. Of course,
this claim can also be proved with the method therein.
\hfill$\Box$\vspace{2mm}

We write $f_c=\{f\le c\}$ for $f={\cal
F}, {\cal F}^X$ and ${\cal F}^\ast$, ${\cal F}^{\ast X}$.

\begin{claim}\label{cl:5.5}
  For any open neighborhood
${W}$ of ${\cal O}$ in $N{\cal O}(\varepsilon)$,
 write ${W}_X={W}\cap XN{\cal O}$ as an open
 subset of $XN{\cal O}$, then  the inclusion
$$
\left({\cal F}^{\ast X}_c\cap {W}_X, {\cal F}^{\ast X}_c\cap
{W}_X\setminus{\cal O}\right)\hookrightarrow\left({\cal
F}^\ast_c\cap {W}, {\cal F}^\ast_c\cap {W}\setminus{\cal O}\right),
$$
 induces surjective
homomorphisms
$$ H_\ast\left({\cal F}^{\ast X}_c\cap {W}_X,
{\cal F}^{\ast X}_c\cap {W}_X\setminus{\cal O};\K\right) \to
H_\ast\left({\cal F}^\ast_c\cap {W}, {\cal F}^\ast_c\cap
{W}\setminus{\cal O};\K\right)
$$ for any Abelian group $\K$.
\end{claim}

This may be obtained by Corollary~3.3 of \cite{Lu2}. For clearness
we present its proof with the proof method  of \cite[Cor.2.5]{Lu3}
since it shows that we need not assume the normal bundle of
${\cal O}$ to be trivial in Theorem~3.17 of \cite{Lu2}.\vspace{2mm}

\noindent{\bf Proof of Claim~\ref{cl:5.5}}. If the closure of ${W}$
in $N{\cal O}$ is contained in the interior of $N{\cal
O}(\varepsilon)$, it is easy to show that the closure of ${W}_X$ in
$XN{\cal O}$ is also contained in the interior of $N{\cal
O}(\varepsilon)\cap XN{\cal O}$ when the latter is equipped with
topology of $XN{\cal O}$. By the excision of relative homology
groups we only need to prove Claim~\ref{cl:5.5} for some open
neighborhood ${W}$ of ${\cal O}$ in $N{\cal O}(\varepsilon)$. Note
that $\Upsilon(N{\cal O}(\epsilon))=U\subset N{\cal O}(\varepsilon)$
by (\ref{e:5.6}).

Since ${\bf H}^-({B})(\nu)\oplus {\bf H}^0({B})(\nu)\oplus{\bf
H}^+({B})(\nu)\subset{N}{\cal O}(\epsilon)$ for $\nu\in
(0,\epsilon/3)$, we may take
$$
{W}=\Upsilon\bigl({\bf H}^-({B})(\nu)\oplus {\bf
H}^0({B})(\nu)\oplus{\bf H}^+({B})(\nu)\bigr),\quad
{V}=\Upsilon\bigl({\bf H}^-({B})(\nu)\oplus {\bf
H}^0({B})(\nu)\bigr).
$$
 Consider the deformation
 $\eta: {W}\times [0, 1]\to {W}$ given by
 $$
\eta\bigl(\Upsilon(x, u^-+ u^0+ u^+), t\bigr)=\Upsilon(x, u^-+ u^0+
tu^+).
 $$
It gives a deformation retract from $\bigl({\cal F}^\ast_c\cap {W},
{\cal F}^\ast_c\cap {W}\setminus{\cal O}\bigr)$ onto $\bigl({\cal
F}^\ast_c\cap {V}, {\cal F}^\ast_c\cap {V}\setminus{\cal O}\bigr)$.
Hence  the inclusion $I:\bigl({\cal F}^\ast_c\cap {V}, {\cal
F}^\ast_c\cap {V}\setminus{\cal O}\bigr)\hookrightarrow \bigl({\cal
F}^\ast_c\cap {W}, {\cal F}^\ast_c\cap {W}\setminus{\cal O}\bigr)$
induces isomorphisms between their relative singular groups.
 This means that  for each nontrivial $\alpha\in H_q\bigl({\cal F}^\ast_c\cap {W}, {\cal F}^\ast_c\cap
{W}\setminus{\cal O};\K\bigr)$ we can choose a relative singular
cycle representative of it, $c=\sum_jg_j{\gamma_0}_j$, such that
$$
|c|:=\cup_j{\gamma_0}_j(\triangle^q)\subset {\cal F}^\ast_c\cap
{V}\quad\hbox{and}\quad |\partial c|\subset {\cal F}^\ast_c\cap
{V}\setminus{\cal O}.
$$
By the second conclusion of Theorem~\ref{th:1.5}(iii) we deduce that
$c$ is also a relative singular cycle in $\bigl({\cal F}^{\ast
X}_c\cap {V}, {\cal F}^{\ast X}_c\cap {V}\setminus{\cal O}\bigr)$,
 denoted by $c_X$ for clearness. Precisely, if $\imath^V$
denotes the identity map from $\bigl({\cal F}^{\ast X}_c\cap {V},
{\cal F}^{\ast X}_c\cap {V}\setminus{\cal O}\bigr)$ to $\bigl({\cal
F}^{\ast}_c\cap {V}, {\cal F}^{\ast}_c\cap {V}\setminus{\cal
O}\bigr)$ then $\imath^V(c_X)=c$. Let $\jmath^W$ denote the
inclusion map from $\bigl({\cal F}^{\ast X}_c\cap {W}, {\cal
F}^{\ast X}_c\cap {W}\setminus{\cal O}\bigr)$ to $\bigl({\cal
F}^{\ast}_c\cap {W}, {\cal F}^{\ast}_c\cap {W}\setminus{\cal
O}\bigr)$. Then $I$, $\imath^V$,
 $\jmath^W$ and the inclusion
$$
I^X:\bigl({\cal F}^{\ast X}_c\cap {V}, {\cal F}^{\ast X}_c\cap
{V}\setminus{\cal O}\bigr)\hookrightarrow \bigl({\cal F}^{\ast
X}_c\cap {W}, {\cal F}^{\ast X}_c\cap {W}\setminus{\cal O}\bigr)
$$
 satisfy $I\circ\imath^V=\jmath^W\circ
I^X$. Since $I_\ast([c])=\alpha$ we obtain
$$
\alpha=I_\ast\circ(\imath^V)_\ast[c_X]=(\jmath^W)_\ast\circ
(I^X)_\ast[c_X]=(\jmath^W)_\ast\bigl((I^X)_\ast[c_X]\bigr).
$$
Claim~\ref{cl:5.5} is proved. \hfill$\Box$\vspace{2mm}

Since any surjective (or injective) homomorphism between vector
spaces of same finite dimension is an isomorphism, for a field $\K$
those surjective homomorphisms in Claim~\ref{cl:5.5} are
isomorphisms  by Claim~\ref{cl:5.4}. Taking
$\mathscr{W}=\digamma(W)$ the first claim of
Theorem~\ref{th:1.5}(iv) immediately follows from the commutative
diagram
$$
\begin{CD}
\left({\cal F}^{\ast X}_c\cap {W}_X, {\cal F}^{\ast X}_c\cap
{W}_X\setminus{\cal O}\right) @>{\rm Inclusion}
>> \left({\cal
F}^\ast_c\cap {W}, {\cal F}^\ast_c\cap {W}\setminus{\cal O}\right)
 \\
@V {\rm EXP} VV @VV {\rm EXP} V \\
\left(({\cal L}^\ast|_{\cal X})_c\cap \mathscr{W}_X, ({\cal
L}^\ast|_{\cal X})_c\cap \mathscr{W}_X\setminus{\cal O}\right)@>{\rm
Inclusion}
>>\left({\cal
L}^\ast_c\cap \mathscr{W}, {\cal L}^\ast_c\cap
\mathscr{W}\setminus{\cal O}\right)
\end{CD}
$$
Similarly, the second claim  of Theorem~\ref{th:1.5}(iv) can be
proved as that of Claim~\ref{cl:5.5}. \hfill$\Box$\vspace{2mm}

\noindent\underline{Proof of Theorem~\ref{th:1.5}(v).} As above we
shall prove an equivalent version of it on the bundle $XN{\cal O}$.
Take a small open neighborhood $\mathcal{V}$ of ${\cal O}$ in
$XN{\cal O}$ such that its closure is contained in ${\rm
Int}({W}_X)$ and that ${\cal F}^{\ast X}$ and ${\cal F}^{X}$ are
same in $\mathcal{V}$. Then
\begin{equation}\label{e:5.17}
\bigl({\cal F}^{\ast X}_c\cap \mathcal{V}, {\cal F}^{\ast X}_c\cap
\mathcal{V}\setminus{\cal O}\bigr)= \bigl({\cal F}^{X}_c\cap
\mathcal{V}, {\cal F}^{X}_c\cap \mathcal{V}\setminus{\cal O}\bigr).
\end{equation}
 By the excision of
relative homology groups we deduce that  the inclusion
$$
\bigl({\cal F}^{\ast X}_c\cap \mathcal{V}, {\cal F}^{\ast X}_c\cap
\mathcal{V}\setminus{\cal O}\bigr)\hookrightarrow \bigl({\cal
F}^{\ast X}_c\cap {W}_X, {\cal F}^{\ast X}_c\cap {W}_X\setminus{\cal
O}\bigr)
$$
induces isomorphisms between their  homology groups. This and the
above equivalent version of Theorem~\ref{th:1.5}(iv) imply that the
same claim holds true for
 the inclusion
 \begin{equation}\label{e:5.18}
I^{vw}:\bigl({\cal F}^{\ast X}_c\cap \mathcal{V}, {\cal F}^{\ast
X}_c\cap \mathcal{V}\setminus{\cal O}\bigr)
\hookrightarrow\bigl({\cal F}^\ast_c\cap {W}, {\cal F}^\ast_c\cap
{W}\setminus{\cal O}\bigr).
\end{equation}
By Corollary~\ref{cor:2.3}, $L^\ast\le L$ and hence
$\mathcal{L}^\ast\le\mathcal{L}$ and ${\cal F}^\ast\le {\cal F}$.
The latter implies
 \begin{equation}\label{e:5.19}
 \bigl({\cal
F}_c\cap {W}, {\cal F}_c\cap {W}\setminus{\cal
O}\bigr)\subset\bigl({\cal F}^\ast_c\cap {W}, {\cal F}^\ast_c\cap
{W}\setminus{\cal O}\bigr).
\end{equation}
Moreover, by (\ref{e:5.17}) we have also the inclusion
\begin{equation}\label{e:5.20}
J^{vw}:\bigl({\cal F}^{X}_c\cap \mathcal{V}, {\cal F}^{X}_c\cap
\mathcal{V}\setminus{\cal O}\bigr)\hookrightarrow \bigl({\cal
F}_c\cap {W}, {\cal F}_c\cap {W}\setminus{\cal O}\bigr).
\end{equation}
Hence $I^{vw}=\mathfrak{I}\circ J^{vw}$ and thus
$I^{vw}_\ast=\mathfrak{I}_\ast\circ J^{vw}_\ast$, where
$\mathfrak{I}$ is the inclusion in (\ref{e:5.19}). Since the
homology groups of pairs in (\ref{e:5.18})-(\ref{e:5.20}) are
isomorphic vector spaces of finite dimensions by Claim~\ref{cl:5.4},
and $I^{vw}_\ast$ are isomorphisms by the above equivalent version
of Theorem~\ref{th:1.5}(iv), both $h_\ast$ and $J^{vw}_\ast$ must be
isomorphisms as well.

As above the excision leads to that the inclusion
$$
 \bigl({\cal F}^{X}_c\cap
\mathcal{V}, {\cal F}^{ X}_c\cap \mathcal{V}\setminus{\cal
O}\bigr)\hookrightarrow \bigl({\cal F}^{X}_c\cap {W}_X, {\cal
F}^{X}_c\cap {W}_X\setminus{\cal O}\bigr)
$$
induces isomorphisms
between their homology groups. Composing the inclusion with
\begin{equation}\label{e:5.21}
\bigl({\cal F}^{X}_c\cap {W}_X, {\cal F}^{X}_c\cap
{W}_X\setminus{\cal O}\bigr)\hookrightarrow \bigl({\cal F}_c\cap
{W}, {\cal F}_c\cap {W}\setminus{\cal O}\bigr)
\end{equation}
 we get
$J^{vw}$. It follows that the inclusion in (\ref{e:5.21}) induces
isomorphisms between their homology groups. The desired conclusion
is proved.\hfill$\Box$\vspace{2mm}

\subsection{Proof of Theorem~\ref{th:1.8}}\label{sec:5.2}

 \textsf{Step 1. Proving $p>0$}. By a
contradiction argument we assume $p=0$. Write ${\cal
O}=S^1\cdot\gamma_0$ as before.

 \noindent{\bf Case 1}.
$m^0({\cal O})=0$. Since $m^-(\gamma_0)=0$, by Theorem~\ref{th:1.7}
we have $C_q(\mathcal{L},\mathcal{O};\K)=H_q(S^1;\K)$
 and hence $C_q({\cal L},{\cal O};\K)=0$ for $q\notin\{0,1\}$. By
the assumption $p$ must be larger than zero.

\noindent{\bf Case 2}. $m^0({\cal O})>0$. Since $m^-(\gamma_0)=0$ we
derive from Theorem~\ref{th:1.7} that
\begin{eqnarray*}
&& C_0({\cal L}, {\cal O};\K)= \bigl(C_{0}({\cal
L}^{\circ}_{\triangle x}, 0;\K)\bigr)^{S^1_x}\quad\forall x\in {\cal
O},\\
&&C_1({\cal L}, {\cal O};\K)= \bigl( C_{0}({\cal
L}^{\circ}_{\triangle x}, 0;\K)\bigr)^{S^1_x}\oplus\bigl(C_{1}({\cal
L}^{\circ}_{\triangle x}, 0;\K)\bigr)^{S^1_x}\quad\forall x\in {\cal
O}.
\end{eqnarray*}
Hence $p$ must be more than zero too.

\noindent\textsf{Step 2. Proving (i).} If $m^0({\cal O})=0$, then
${\cal L}^\ast\circ\digamma\circ\Upsilon(u)=\|u\|_1^2$ for any $u\in
N{\cal O}(\epsilon)$ by Theorem~\ref{th:1.5}(iii). Moreover ${\cal
L}^\ast\le{\cal L}$ and ${\cal L}^\ast={\cal L}$ on ${\cal O}$. The
claim follows directly.

If $m^0({\cal O})>0$, by Theorem~\ref{th:1.7} for any $x\in {\cal
O}$ we have
\begin{eqnarray*}
&&0\ne C_1({\cal L}, {\cal O};\K)= \bigl( C_{0}({\cal
L}^{\circ}_{\triangle x}, 0;\K)\bigr)^{S^1_x}\oplus\bigl(C_{1}({\cal
L}^{\circ}_{\triangle x}, 0;\K)\bigr)^{S^1_x},\\
&&0=C_{2}({\cal L}, {\cal O};\K)= \bigl( C_{1}({\cal
L}^{\circ}_{\triangle x}, 0;\K)\bigr)^{S^1_x}\oplus\bigl(C_{2}({\cal
L}^{\circ}_{\triangle x}, 0;\K)\bigr)^{S^1_x}.
\end{eqnarray*}
It follows that $C_{0}({\cal L}^{\circ}_{\triangle x}, 0;\K)\ne 0$.
By Theorem~4.6 on the page 43 of \cite{Ch93} $0=0_x$ must be a local
minimum of ${\cal L}^{\circ}_{\triangle x}$. This and
Theorem~\ref{th:1.5}(iii) imply that $x$ is a local minimum of
${\cal L}^\ast$ and hence of ${\cal L}$ because ${\cal
L}^\ast\le{\cal L}$ and ${\cal L}^\ast={\cal L}$ on ${\cal O}$ as
above.

\textsf{Step 3. Proving (ii)}. If $m^0({\cal O})=0$, by the proof of
Case 1 in Step 1 we have
$$
C_{p}({\cal L},{\cal O};\K)=H_{p}(S^1;\K)=0\quad\hbox{since}\;p\ge
2.
$$
This case cannot occur. Hence it must hold that $m^0({\cal O})>0$.
By Theorem~\ref{th:1.7}
\begin{eqnarray*}
&&0\ne C_p({\cal L}, {\cal O};\K)= \bigl( C_{p-1}({\cal
L}^{\circ}_{\triangle x}, 0;\K)\bigr)^{S^1_x}\oplus\bigl(C_{p}({\cal
L}^{\circ}_{\triangle x}, 0;\K)\bigr)^{S^1_x},\\
&&0=C_{p+1}({\cal L}, {\cal O};\K)= \bigl( C_{p}({\cal
L}^{\circ}_{\triangle x},
0;\K)\bigr)^{S^1_x}\oplus\bigl(C_{p+1}({\cal L}^{\circ}_{\triangle
x}, 0;\K)\bigr)^{S^1_x}
\end{eqnarray*}
for any $x\in {\cal O}$. It follows that $C_{p-1}({\cal
L}^{\circ}_{\triangle x}, 0;\K)\ne 0$. By Example~1 on the page 33
of \cite{Ch93} $0_x$ may not be a local minimum of ${\cal
L}^{\circ}_{\triangle x}$. This and Theorem~\ref{th:1.6} imply that
$x$ is not a local minimum of ${\cal L}|_{\cal X}$ and thus of
${\cal L}$.
 \hfill$\Box$\vspace{2mm}

\subsection{Proof of Proposition~\ref{prop:5.1}}\label{sec:5.3}

 Since  $\gamma_0$ is a $C^k$-map to $M$ with
$k\ge 5$, starting with a unit orthogonal
 frame at $T_{\gamma(0)}M$ and using the parallel
transport along $\gamma_0$ with respect to the Levi-Civita
connection of the Riemannian metric $g$ we get a unit orthogonal
parallel $C^k$ frame field $\R\to \gamma_0^\ast TM,\;t\mapsto
(e_1(t),\cdots, e_n(t))$.
 Note that  there exists a unique orthogonal matrix
 $E_{\gamma_0}$ such that
$(e_1(1),\cdots, e_n(1))=(e_1(0),\cdots, e_n(0))E_{\gamma_0}$.
(\textsf{All vectors in $\R^n$ will be understood as row vectors.}) By the elementary matrix theory there
exists an orthogonal matrix $\Xi$ such that
$$
\Xi^{-1}E_{\gamma_0}\Xi={\rm
diag}(S_1,\cdots,S_{\sigma})\in\R^{n\times n},
$$
where each $S_j$ is either $1$, or $-1$, or
${\scriptscriptstyle\left(\begin{array}{cc}
\cos\theta_j& \sin\theta_j\\
-\sin\theta_j& \cos\theta_j\end{array}\right)}$, $0<\theta_j<\pi$,
and their orders satisfy: ${\rm ord}(S_1)\ge\cdots\ge{\rm
ord}(S_{\sigma})$.
 So replacing $(e_1,\cdots,e_n)$ by $(e_1,\cdots,e_n)\Xi$ we may
assume
$$
E_{\gamma_0}={\rm diag}(S_1,\cdots,S_{\sigma})\in\R^{n\times n}.
$$
Since $\gamma_0(t+1)=\gamma_0(t)\;\forall t\in\R$ and $\R\ni
t\mapsto (e_1(t),\cdots, e_n(t))E_{\gamma_0}$ is also a unit
orthogonal parallel $C^k$ frame field along $\gamma_0$ it is easily
proved that
$$
(e_1(t+1),\cdots, e_n(t+1))=(e_1(t),\cdots,
e_n(t))E_{\gamma_0}\quad\forall t\in\R.
$$
Such a frame field is called  $E_{\gamma_0}$-{\bf 1-invariant}.

A curve $\xi=(\xi_1,\cdots,\xi_n):\R\to\R^n$ is called
$E_{{\gamma_0}}$-{\bf 1-invariant} if
$\xi(t+1)^T=E_{{\gamma_0}}\xi(t)^T$ for all $t\in\R$. Here
$\xi(t)^T$ denotes  the transpose of the matrix $\xi(t)$ as usual.
Let $X_{{\gamma_0}}$ be the Banach space of all
$E_{{\gamma_0}}$-1-invariant $C^1$ curves from $\R$ to $\R^n$
according to the usual $C^1$-norm. Denote by $H_{{\gamma_0}}$  the
completion of $X_{{\gamma_0}}$ with respect to the norm
$\|\cdot\|_{1,2}$ induced by the inner product
\begin{eqnarray}\label{e:5.22}
\langle\xi,\eta\rangle_{1,2}=\int_0^1[(\xi(t),\eta(t))_{\R^n}+
(\dot\xi(t),\dot\eta(t))_{\R^n}]dt.
\end{eqnarray}
 Since
$\sum^n_{j=1}\xi_j(t+1){e}_j(t+1)=\sum^n_{j=1}\xi_j(t){e}_j(t)$ for
all $t\in\R$, we obtain a Hilbert space isomorphism
 $$
\mathfrak{I}_{{\gamma_0}}: ({H}_{\gamma_0},
\langle\cdot,\cdot\rangle_{1,2})\to (T_{\gamma_{0}}{\Lambda M},
\langle\cdot,\cdot\rangle_1),\quad\xi\mapsto
\mathfrak{I}_{{\gamma_0}}(\xi)=\sum^n_{j=1}\xi_j{e}_j.
$$
 With the exponential map  $\exp$ of $g$  we define the $C^k$ map
$$
\phi_{{\gamma_0}}:\R\times B^n_{2\rho}(0)\to
M,\;(t,x)\mapsto\exp_{\gamma_{0}(t)}\left(\sum^n_{i=1}x_i
e_i(t)\right)
$$
 for some small open ball
$B^n_{2\rho}(0)\subset\R^n$. Then
$\phi_{{\gamma_0}}(t+1,x)=\phi_{{\gamma_0}}(t, (E_{{\gamma_0}}
x^T)^T)$ and
\begin{eqnarray*}
d\phi_{{\gamma_0}}(t+1,x)[(1,v)]=d\phi_{{\gamma_0}}(t,
(E_{{\gamma_0}} x^T)^T)[(1, (E_{{\gamma_0}} v^T)^T)]
\end{eqnarray*}
for $(t,x)\in \R\times B^n_{2\rho}(0)$. This yields  a $C^{k-3}$
coordinate chart around $\gamma_{0}$ on ${\Lambda M}$,
$$
\Phi_{{\gamma_0}}:{\bf B}_{2\rho}({H}_{{\gamma_0}}):=\{\xi\in
 H_{{\gamma_0}}\,|\,\|\xi\|_{1,2}<2\rho\} \to{\Lambda M}
$$
given by $\Phi_{{\gamma_0}}(\xi)(t)=\phi_{{\gamma_0}}(t,\xi(t))$ for
$\xi\in {\bf B}_{2\rho}({H}_{{\gamma_0}})$. Note that
$d\Phi_{{\gamma_0}}(0)={\mathfrak{I}}_{{\gamma_0}}$.

Let us compute the expression of ${\cal L}^\ast$ in this chart.
Define
\begin{eqnarray}\label{e:5.23}
 L^\ast_{{\gamma_0}}(t, x, v)=L^\ast\bigr(\phi_{{\gamma_0}}(t,x),
d\phi_{{\gamma_0}}(t,x)[(1,v)]\bigl).
\end{eqnarray}
Then for any $t\in\R$ we have $L^\ast_{{\gamma_0}}(t,0,0)\equiv  c$,
\begin{eqnarray*}
L^\ast_{{\gamma_0}}(t+1, x,
v)&=&L^\ast\bigr(\phi_{{\gamma_0}}(t+1,x),
d\phi_{{\gamma_0}}(t+1,x)[(1,v)]\bigl)\\
&=&L^\ast\bigr(\phi_{{\gamma_0}}(t, (E_{{\gamma_0}} x^T)^T),
d\phi_{{\gamma_0}}(t, (E_{{\gamma_0}} x^T)^T)[(1, (E_{{\gamma_0}} v^T)^T)]\bigl)\\
&=&L^\ast_{{\gamma_0}}(t, (E_{{\gamma_0}} x^T)^T, (E_{{\gamma_0}}
v^T)^T)\quad\hbox{and}\\
\partial_x L^\ast_{{\gamma_0}}(t+1, x, v)&=&\partial_x L^\ast_{{\gamma_0}}(t,
(E_{{\gamma_0}} x^T)^T,
(E_{{\gamma_0}} v^T)^T)E_{{\gamma_0}},\\
\partial_v L^\ast_{{\gamma_0}}(t+1, x, v)&=&\partial_v L^\ast_{{\gamma_0}}(t,
(E_{{\gamma_0}} x^T)^T, (E_{{\gamma_0}} v^T)^T)E_{{\gamma_0}}.
\end{eqnarray*}
It follows from these that
\begin{eqnarray}
 L^\ast_{{\gamma_0}}(t+1,\xi(t+1),\dot\xi(t+1))&=&
L^\ast_{{\gamma_0}}(t,
(E_{{\gamma_0}} \xi(t+1)^T)^T, (E_{{\gamma_0}} \dot\xi(t+1)^T)^T)\nonumber\\
&=& L^\ast_{{\gamma_0}}(t,\xi(t),\dot\xi(t))\quad\forall t\in\R,\label{e:5.24}\\
\partial_x L^\ast_{{\gamma_0}}(t+1,\xi(t+1),\dot\xi(t+1))&=&\partial_x
L^\ast_{{\gamma_0}}(t,\xi(t),\dot\xi(t))E_{{\gamma_0}}\quad\forall t\in\R,\label{e:5.25}\\
\partial_v L^\ast_{{\gamma_0}}(t+1,\xi(t+1),\dot\xi(t+1))&=&\partial_v
L^\ast_{{\gamma_0}}(t,\xi(t),\dot\xi(t))E_{{\gamma_0}}\quad\forall
t\in\R\label{e:5.26}
\end{eqnarray}
 for any $\xi\in {\bf
B}_{2\rho}({H}_{{\gamma_0}})$ since
$\xi(t+1)^T=E_{{\gamma_0}}\xi(t)^T$ and
$\dot\xi(t+1)^T=E_{{\gamma_0}}\dot\xi(t)^T$. Define  the action
functional ${\cal L}^\ast_{{\gamma_0}}: {\bf
B}_{2\rho}({H}_{{\gamma_0}})\to\R$ by
$$
{\cal L}^\ast_{{\gamma_0}}(\xi)=\int^1_0
L^\ast_{{\gamma_0}}(t,\xi(t),\dot\xi(t))dt.
$$
 Then
${\cal L}^\ast_{{\gamma_0}}={\cal L}^\ast\circ\Phi_{{\gamma_0}}$ on
${\bf B}_{2\rho}({H}_{{\gamma_0}})$. We claim that
 it  is $C^{2-0}$. Since  the gradient  $\nabla{\cal
 L}^\ast_{{\gamma_0}}$ of ${\cal
 L}^\ast_{{\gamma_0}}$ at  $\xi\in {\bf B}_{2\rho}({H}_{{\gamma_0}})$ is given
 by (\ref{e:B.34}),  we derive
\begin{eqnarray}\label{e:5.27}
\frac{d}{dt}\nabla{\cal
L}^\ast_{{\gamma_0}}(\xi)(t)&=&\frac{e^{t}}{2}\int^\infty_te^{-s}\left(
\partial_x
L^\ast_{{\gamma_0}}\bigl(s,\xi(s),\dot{\xi}(s)\bigr)-\mathfrak{R}^{\xi}(s) \right)\, ds\nonumber\\
&-&\frac{e^{-t}}{2}\int^t_{-\infty}e^{s}\left(\partial_x
L^\ast_{{\gamma_0}}\bigl(s,\xi(s),\dot{\xi}(s)\bigr)-\mathfrak{R}^{\xi}(s)
\right)\, ds \nonumber\\
&+&\partial_v L^\ast_{{\gamma_0}}\bigl(t,\xi(t),\dot{\xi}(t)\bigr)-
\left(\int^{1}_0\partial_v
L^\ast_{{\gamma_0}}\bigl(s,\xi(s),\dot{\xi}(s)\bigr)ds\right)\times\nonumber\\
&&\times \biggl[\frac{d}{dt}\left(\oplus_{l\le
p}\frac{\sin\theta_l}{2-2\cos\theta_l}{\scriptscriptstyle\left(\begin{array}{cc}
0& -1\\
1&
0\end{array}\right)}-\frac{1}{2}I_{2p}\right)\nonumber\\
&&\qquad\oplus\, {\rm diag}(\dot{a}_{p+1}(t), \cdots,
\dot{a}_\sigma(t))\biggr]
\end{eqnarray}
provided $2={\rm ord}(S_p)>{\rm ord}(S_{p+1})$ for some
$p\in\{0,\cdots,\sigma\}$, where $\mathfrak{R}^{\xi}$ is given by
(\ref{e:B.35}).
 Let  $\mathbb{A}^\ast_{{\gamma_0}}$ be the restriction of the gradient
 $\nabla{\cal L}^\ast_{{\gamma_0}}$ to
${\bf B}_{2\rho}({H}_{{\gamma_0}})\cap X_{{\gamma_0}}$. Then
$\mathbb{A}^\ast_{{\gamma_0}}({\bf B}_{2\rho}({H}_{{\gamma_0}})\cap
X_{{\gamma_0}})\subset X_{{\gamma_0}}$.  As in the proof of
\cite[Lemma 3.2]{Lu1} it follows from (\ref{e:B.34}), (\ref{e:B.35})
and (\ref{e:5.27})  that $\mathbb{A}^\ast_{{\gamma_0}}$ is $C^1$ as
a map from the open subset ${\bf B}_{2\rho}({H}_{{\gamma_0}})\cap
X_{{\gamma_0}}$ of $X_{{\gamma_0}}$ to $X_{{\gamma_0}}$ (and hence
the restriction ${\cal L}^{\ast X}_{{\gamma_0}}$ of ${\cal
L}^\ast_{{\gamma_0}}$ to ${\bf B}_{2\rho}({H}_{{\gamma_0}})\cap
X_{{\gamma_0}}$ is $C^2$). Moreover
\begin{eqnarray*}
 d^2 {\cal L}^{\ast X}_{{\gamma_0}}
  (\zeta)[\xi,\eta]
   = \int_0^{1} \Bigl(\!\! \!\!\!&&\!\!\!\!\!\partial_{vv}
     L^\ast_{{\gamma_0}}\bigl(t,\zeta(t),\dot{\zeta}(t)\bigr)
\bigl[\dot{\xi}(t), \dot{\eta}(t)\bigr] \nonumber\\
&&+ \partial_{qv}
  L^\ast_{{\gamma_0}}\bigl(t,\gamma(t), \dot{\gamma}(t)\bigr)
\bigl[\xi(t), \dot{\eta}(t)\bigr]\nonumber \\
&& + \partial_{vq}
  L^\ast_{{\gamma_0}}\bigl(t,\zeta(t),\dot{\zeta}(t)\bigr)
\bigl[\dot{\xi}(t), \eta(t)\bigr] \nonumber\\
&&+  \partial_{qq} L^\ast_{{\gamma_0}}\bigl(t,\zeta(t),
\dot{\zeta}(t)\bigr) \bigl[\xi(t), \eta(t)\bigr]\Bigr) \, dt
\end{eqnarray*}
for any $\zeta\in {\bf B}_{2\rho}({H}_{{\gamma_0}})\cap
X_{{\gamma_0}}$ and $\xi,\eta\in X_{{\gamma_0}}$,  it is easily
checked that the corresponding properties to (i) and (ii) below
Lemma~\ref{lem:4.2} hold for ${\cal L}^{\ast X}_{{\gamma_0}}$, that
is,
\begin{description}
\item[(i)] for any $\zeta\in {\bf B}_{2\rho}({H}_{{\gamma_0}})\cap X_{{\gamma_0}}$ there
exists a constant $C(\zeta)$ such that
$$
\big|d^2 {\cal L}^{\ast X}_{\gamma_0} (\zeta)[\xi,\eta]\big|\le
  C(\zeta)\|\xi\|_{1,2}\cdot\|\eta\|_{1,2}
\quad\forall \xi,\eta\in X_{\gamma_0};
$$
\item[(ii)] $\forall\varepsilon>0,\;\exists\;\delta_0>0$, such that
for all $\zeta_1, \zeta_2\in {\bf B}_{2\rho}({H}_{\gamma_0})\cap
X_{{\gamma_0}}$ with $\|\zeta_1-\zeta_2\|_{C^1}<\delta_0$,
 $$
\big|d^2 {\cal L}^{\ast X}_{{\gamma_0}} (\zeta_1)[\xi,\eta]- d^2
{\cal L}^{\ast X}_{{\gamma_0}}(\zeta_2)[\xi,\eta]\big|\le
  \varepsilon \|\xi\|_{1,2}\cdot\|\eta\|_{1,2}
\quad\forall \xi,\eta\in X_{\gamma_0}.
$$
\end{description}
It follows that there exists a  continuous map
$\mathbb{B}^\ast_{{\gamma_0}}:
 {\bf B}_{2\rho}({H}_{{\gamma_0}})\cap X_{{\gamma_0}}\to
{\cal L}_s({H}_{{\gamma_0}})$ with respect to the induced topology
on ${\bf B}_{2\rho}({H}_{{\gamma_0}})\cap X_{{\gamma_0}}$ by
$X_{{\gamma_0}}$, such that
$$
\langle d\mathbb{A}^\ast_{{\gamma_0}}(\zeta)[\xi],
\eta\rangle_{1,2}=d^2 {\cal L}^{\ast X}_{{\gamma_0}}
  (\zeta)[\xi,\eta]=\langle\mathbb{B}^\ast_{{\gamma_0}}(\zeta)\xi,
  \eta\rangle_{1,2}\quad\forall \xi, \eta\in X_{{\gamma_0}}.
$$
 By these and similar arguments to those
of \cite{Lu1} we may obtain

\begin{claim}\label{cl:5.6}
 Around the critical point $0\in
{H}_{{\gamma_0}}$, $({H}_{{\gamma_0}}, X_{{\gamma_0}}, {\cal
L}^\ast_{{\gamma_0}}, \mathbb{A}^\ast_{{\gamma_0}},
\mathbb{B}^\ast_{{\gamma_0}})$ satisfies the conditions of
Theorem~\ref{th:A.1} except that the critical point $0$ is not
isolated.
\end{claim}

Observe that ${\cal L}^\ast_{\gamma_0}$ has an one-dimensional
critical manifold $S:= \Phi^{-1}_{{\gamma_0}}\bigl({\cal O}\cap{\rm
Im}(\Phi_{\gamma_0})\bigr)$, and that $T_{\gamma_0}{\cal
O}=\dot\gamma_0\R\subset W^{1,2}(\gamma_0^\ast TM)$. Since
$d\Phi_{\gamma_0}(0)=\mathfrak{I}_{{\gamma_0}}$
  is an isomorphism there exists a
unique $\zeta_0\in H_{\gamma_0}$ satisfying
${\mathfrak{I}}_{{\gamma_0}}(\zeta_0)=\dot\gamma_{0}$,
  that is, for any $t\in \R$,
$$
\dot\gamma_0(t)=\frac{d}{ds}\Bigl|_{s=0}\Phi_{\gamma_0}(s\zeta_0)(t)=
\frac{d}{ds}\Bigl|_{s=0}\exp_{\gamma_0(t)}\left(\sum^n_{k=1}s\zeta_{0k}(t)e_k(t)\right)
=\sum^n_{k=1}\zeta_{0k}(t)e_k(t).
$$
Hence $g(\dot\gamma_0(t), e_j(t))=\zeta_{0j}(t)$ for any $t\in \R$
and $j=1,\cdots,n$. Clearly, $T_0S=\zeta_0\R$. The normal space of
$S$ at $0\in S$ is  the orthogonal complementary of $\zeta_0\R$ in
the Hilbert space $H_{\gamma_0}$, denoted by $H_{{\gamma_0},0}$.
Note that $\zeta_0\in X_{\gamma_0}$ actually. Let
${X}_{{\gamma_0},0}:= H_{{\gamma_0},0}\cap X_{{\gamma_0}}$.  Denote
by ${\cal L}^\ast_{{\gamma_0},0}$ the restriction of ${\cal
L}^\ast_{{\gamma_0}}$ to ${\bf B}_{2\rho}({H}_{{\gamma_0},0})$. Then
for $\xi\in {\bf B}_{2\rho}(H_{{\gamma_0},0})$ we have
$$
\nabla{\cal L}^\ast_{{\gamma_0},0}(\xi)=\nabla{\cal
L}^\ast_{{\gamma_0}}(\xi)-\frac{\langle \nabla{\cal
L}^\ast_{{\gamma_0}}(\xi),  \zeta_0\rangle_{1,2}}{\|
\zeta_0\|_{1,2}^2} \zeta_0.
$$
 It follows that the
restriction of $\nabla{\cal L}^\ast_{{\gamma_0},0}$ to ${\bf
B}_{2\rho}( H_{{\gamma_0},0})\cap X_{{\gamma_0},0}$, denoted by
$\mathbb{A}^\ast_{{\gamma_0},0}$, is a $C^1$-map into
$X_{{\gamma_0},0}$ (with respect to the $C^1$-topology on ${\bf
B}_{2\rho}( H_{{\gamma_0},0})\cap X_{{\gamma_0},0}$), and
\begin{equation}\label{e:5.28}
\mathbb{A}^\ast_{{\gamma_0},0}(\xi)=\mathbb{A}^\ast_{{\gamma_0}}(\xi)-\frac{\langle
 \mathbb{A}^\ast_{{\gamma_0}}(\xi),
\zeta_0\rangle_{1,2}}{\| \zeta_0\|_{1,2}^2} \zeta_0\quad\forall
\xi\in {\bf B}_{2\rho}(H_{{\gamma_0},0})\cap X_{{\gamma_0},0}.
\end{equation}
 Let ${\cal L}^{\ast X}_{{\gamma_0},0}$ be the restriction of ${\cal L}^{\ast
X}_{{\gamma_0}}$ to ${\bf
B}_{2\rho}({H}_{{\gamma_0},0})\cap{X}_{{\gamma_0},0}$, and for
$\xi\in{\bf B}_{2\rho}({H}_{{\gamma_0},0})\cap{X}_{{\gamma_0},0}$
let
 $\mathbb{B}^\ast_{{\gamma_0},0}(\xi)$  be the extension of the continuous
symmetric bilinear form $d^2({\cal L}^{\ast X}_{{\gamma_0},0})(\xi)$
on ${H}_{{\gamma_0},0}$ .   Then
$$
\mathbb{B}^\ast_{{\gamma_0},0}(\xi)\eta=\mathbb{B}^\ast_{{\gamma_0}}(\xi)\eta-\frac{\langle
\mathbb{B}^\ast_{{\gamma_0}}(\xi)\eta,
\zeta_0\rangle_{1,2}}{\|\zeta_0\|_{1,2}^2}\zeta_0\quad\forall\eta\in
{H}_{{\gamma_0},0}.
$$
Clearly,  $\mathbb{B}^\ast_{{\gamma_0},0}: {\bf
B}_{2\rho}({H}_{{\gamma_0},0})\cap {X}_{{\gamma_0},0}\to {\cal
L}_s({H}_{{\gamma_0},0})$ is  continuous with respect to the
$C^1$-topology on ${\bf B}_{2\rho}({H}_{{\gamma_0},0})\cap
{X}_{{\gamma_0},0}$, and for  any $\xi\in {\bf
B}_{2\rho}({H}_{{\gamma_0},0})\cap {X}_{{\gamma_0},0}$ and $\zeta,
\eta\in {X}_{{\gamma_0},0}$,
$$
\langle d\mathbb{A}^\ast_{{\gamma_0},0}(\xi)[\zeta],
\eta\rangle_{1,2}=d^2({\cal L}^{\ast
X}_{{\gamma_0},0})(\xi)[\zeta,\eta]=\langle\mathbb{B}^\ast_{{\gamma_0},0}(\xi)\zeta,
  \eta\rangle_{1,2}.
$$
 By Claim~\ref{cl:5.6} and
Theorem~\ref{th:A.3} we arrive at

\begin{claim}\label{cl:5.7}
 $({H}_{{\gamma_0},0}, {X}_{{\gamma_0},0}, {\cal
L}^\ast_{{\gamma_0},0}, \mathbb{A}^\ast_{{\gamma_0},0},
\mathbb{B}^\ast_{{\gamma_0},0})$
 satisfies the conditions of Theorem~\ref{th:A.1}
around the critical point $0\in{H}_{{\gamma_0},0}$.
\end{claim}

\begin{remark}\label{rm:5.8}
{\rm From the arguments above Claim~\ref{cl:5.6} one easily sees
that  the conditions of Theorem~\ref{th:A.5} are satisfied for
$({H}_{{\gamma_0}}, X_{{\gamma_0}}, {\cal L}^\ast_{{\gamma_0}},
\mathbb{A}^\ast_{{\gamma_0}}, \mathbb{B}^\ast_{{\gamma_0}})$ around
the critical point $0$. By Theorem~\ref{th:A.7}
$({H}_{{\gamma_0},0}, {X}_{{\gamma_0},0}, {\cal
L}^\ast_{{\gamma_0},0}, \mathbb{A}^\ast_{{\gamma_0},0},
\mathbb{B}^\ast_{{\gamma_0},0})$ satisfies the conditions of
Theorem~\ref{th:A.5} around the critical point $0$ too. }
\end{remark}

Let $\digamma_{\gamma_0}$ be the restriction of the map $\digamma$
in (\ref{e:1.13}) to $N{\cal O}(\varepsilon)_{\gamma_{0}}$. Let
$\delta=\min\{\varepsilon, 2\rho\}$. For $\xi\in {\bf
B}_{\delta}({H}_{{\gamma_0}})$, by the arguments above
(\ref{e:5.23}) we get
\begin{eqnarray}\label{e:5.29}
\Phi_{{\gamma_0}}(\xi)(t)&=&\phi_{{\gamma_0}}(t,\xi(t))=\exp_{\gamma_{0}(t)}\left(\sum^n_{i=1}\xi_i(t)
e_i(t)\right)\nonumber\\
&=&\exp_{\gamma_{0}(t)}\left(\mathfrak{I}_{{\gamma_0}}(\xi)(t)\right)=
\exp_{\gamma_{0}(t)} \left((d\Phi_{{\gamma_0}}(0)[\xi])(t)\right)
\nonumber\\
&=&\digamma_{\gamma_{0}}(d\Phi_{{\gamma_0}}(0)[\xi])(t).
\end{eqnarray}
 That is,
$\digamma_{\gamma_{0}}\circ d\Phi_{{\gamma_0}}(0)=\Phi_{{\gamma_0}}$
on ${\bf B}_{\delta}({H}_{{\gamma_0}})$. Recall that ${\cal
L}^\ast_{{\gamma_0}}={\cal L}^\ast\circ\Phi_{{\gamma_0}}$ on ${\bf
B}_{2\rho}({H}_{{\gamma_0}})$ and ${\cal F}^\ast_{\gamma_{0}}={\cal
L}^\ast\circ\digamma_{\gamma_{0}}$ by the definition below
(\ref{e:1.14}). From these and (\ref{e:5.29}) we derive ${\cal
F}^\ast_{\gamma_{0}}\circ d\Phi_{{\gamma_0}}(0)={\cal
L}^\ast\circ\digamma_{\gamma_{0}}\circ d\Phi_{{\gamma_0}}(0)={\cal
L}^\ast\circ\Phi_{{\gamma_0}}={\cal L}^\ast_{{\gamma_0}}$ on ${\bf
B}_{\delta}({H}_{{\gamma_0}})$, and hence
\begin{equation}\label{e:5.30}
{\cal F}^\ast_{\gamma_{0}}\circ d\Phi_{{\gamma_0}}(0)={\cal
L}^\ast_{{\gamma_0},0}\quad\hbox{on}\quad {\bf
B}_{\delta}({H}_{{\gamma_0},0}).
\end{equation}
Since  ${\mathfrak{I}}_{{\gamma_0}}=d\Phi_{{\gamma_0}}(0)$ restricts
to a Hilbert space isomorphism ${\mathfrak{I}}_{{\gamma_0},0}$ from
$ H_{{\gamma_0},0}$ to $N{\cal O}_{\gamma_{0}}$ and a Banach space
isomorphism ${\mathfrak{I}}^X_{{\gamma_0},0}$ from
$X_{{\gamma_0},0}$ to $X N{\cal O}_{\gamma_{0}}$,  (\ref{e:5.30})
leads to
\begin{equation}\label{e:5.31}
{A}^\ast_{\gamma_{0}}={\mathfrak{I}}^X_{{\gamma_0},0}\circ\mathbb{A}^\ast_{{\gamma_0},0}\circ
({\mathfrak{I}}^X_{{\gamma_0},0})^{-1}
\end{equation}
and ${\bf
B}^\ast_{\gamma_{0}}(\xi)={\mathfrak{I}}_{{\gamma_0},0}\circ\mathbb{B}^\ast_{{\gamma_0},0}\left(
({\mathfrak{I}}^X_{{\gamma_0},0})^{-1}\xi\right)\circ
({\mathfrak{I}}_{{\gamma_0},0})^{-1}$ for $\xi\in {\bf
B}_{\delta}({H}_{{\gamma_0},0})\cap {X}_{{\gamma_0},0}$. By
Claim~\ref{cl:5.7} and Theorem~\ref{th:A.4}, $({N}{\cal
O}_{\gamma_{0}}, X{N}{\cal O}_{\gamma_{0}},
 {\cal F}^\ast_{\gamma_{0}},
  {A}^\ast_{\gamma_{0}}, {B}^\ast_{\gamma_{0}})$ satisfies
the conditions of Theorem~\ref{th:A.1} (and hence
Theorem~\ref{th:A.5}) around the critical point
$0\equiv0_{\gamma_{0}}$. Proposition~\ref{prop:5.1} is proved.
 \hfill$\Box$\vspace{2mm}

\begin{remark}\label{rm:5.9}
{\rm   With $E_\sigma:={\rm diag}(\sigma,1,\cdots,1)\in\R^{n\times
n}$, where $\sigma=\sigma(\gamma_0)=1$ if $\gamma_0^\ast TM\to S^1$
is trivial, and $-1$ otherwise, it was claimed in \cite[\S5.3]{Lu4}
that by parallel transport we may obtain a unit orthogonal parallel
frame field along $\gamma_0$, $\R\ni t\mapsto \gamma_0^\ast
TM,\;t\mapsto (e_1(t),\cdots, e_n(t))$, satisfying
$(e_1(t+1),\cdots, e_n(t+1))=(e_1(t),\cdots,
e_n(t))E_\sigma\;\forall t\in\R$. This is not true in general. In
fact, one can only get a unit orthogonal
 frame field which may not be parallel. So the inner product
in (\ref{e:5.22}) must be replaced by
$$
\langle\xi,\eta\rangle_{1,2,\ast}=\int_0^1[(\xi(t),\eta(t))_{\R^n}+
(\dot\xi(t),\dot\eta(t))_{\R^n}]dt+
\sum^{n}_{i,j=1}\int^1_0g(\nabla_{\dot\gamma_0}e_i,
\nabla_{\dot\gamma_0}e_j)\xi_i\eta_jdt
$$
so that the map $\mathfrak{I}_{{\gamma_0}}$ below (\ref{e:5.22}) is
 a Hilbert space isomorphism from
 $({H}_{\gamma_0},\langle\cdot,\cdot\rangle_{1,2,\ast})$ to $(T_{\gamma_{0}}{\Lambda
M}, \langle\cdot,\cdot\rangle_1)$. We prefer the above method to
this correction.}
\end{remark}

\section{Proofs of Theorems~\ref{th:1.9},~\ref{th:1.10}}\label{sec:6}
\setcounter{equation}{0}

\subsection{Proof of Theorem~\ref{th:1.9}}\label{sec:6.1}

Clearly,  Theorem~\ref{th:1.9} is equivalent to the following
theorem.

\begin{theorem}\label{th:6.1}
 With
$c={\cal L}|_{\cal O}$ there exist $S^1$-invariant neighborhoods
${\cal U}_i$ of $\varphi_i({\cal O})$, $i=1,m$, $\varphi_m({\cal
U}_1)\subset{\cal U}_m$, such that $\varphi_m$ induces isomorphisms
\begin{eqnarray*}
&&(\varphi_m)_\ast: H_\ast\bigl({\cal L}_c\cap{\cal U}_1, {\cal
L}_c\cap({\cal U}_1\setminus\{\cal O\});\K\bigr)\to\\
&&\hspace{30mm}H_\ast\bigl({\cal L}_{m^2c}\cap{\cal U}_m, {\cal
L}_{m^2c}\cap({\cal U}_m\setminus\{\varphi_m(\cal O)\});\K\bigr).
\end{eqnarray*}
 \end{theorem}

Write ${\cal U}_i^X={\cal U}_i\cap {\mathcal{X}}$, $i=1,m$,  as open
subsets of ${\mathcal{X}}=C^1(S^1, M)$. Then  it is
Theorem~\ref{th:1.5}(v) that shows the inclusions
\begin{eqnarray*}
&&\left(({\cal L}|_{\mathcal{X}})_c\cap {\cal U}_1^X, ({\cal
L}|_{\mathcal{X}})_c\cap {\cal U}_1^X\setminus{\cal
O}\right)\hookrightarrow\left({\cal L}_c\cap {\cal U}_1, {\cal
L}_c\cap {\cal U}\setminus{\cal O}\right)\quad{and}\\
&&\bigl(({\cal L}|_{\mathcal{X}})_{m^2c}\cap {\cal U}_m^X, ({\cal
L}|_{\mathcal{X}})_{m^2c}\cap {\cal U}_m^X\setminus\varphi_m({\cal
O})\bigr)\hookrightarrow\bigl({\cal L}_{m^2c}\cap {\cal U}_m, {\cal
L}_{m^2c}\cap {\cal U}_m\setminus\varphi_m({\cal O})\bigr),
\end{eqnarray*}
induce isomorphisms among their  homology groups. Moreover these two
inclusions commute with $\varphi_m$. Hence Theorem~\ref{th:6.1} is
equivalent to

\begin{claim}\label{cl:6.2}
 $\varphi_m$ induces isomorphisms
\begin{eqnarray*}
&&(\varphi_m)_\ast: H_\ast\bigl(({\cal L}|_{\mathcal{X}})_c\cap
{\cal U}_1^X, ({\cal
L}|_{\mathcal{X}})_c\cap {\cal U}_1^X\setminus{\cal O};\K\bigr)\to\\
&&\hspace{30mm}H_\ast\bigl(({\cal L}|_{\mathcal{X}})_{m^2c}\cap
{\cal U}_m^X, ({\cal L}|_{\mathcal{X}})_{m^2c}\cap {\cal
U}_m^X\setminus\varphi_m({\cal O});\K\bigr).
\end{eqnarray*}
\end{claim}

We shall use the  splitting lemma Theorem~\ref{th:1.6} around ${\cal
O}$ for  ${\cal L}|_{\mathcal{X}}$ to prove it. As in the case of
Riemannian geometry in \cite{GM2},
 we need to introduce  the equivalent
Riemannian-Hilbert structure on $T{\Lambda M}$:
\begin{equation}\label{e:6.1}
(\xi,\eta )_m=m^2\int^1_0\langle\xi(t),\eta(t)\rangle dt+
\int^1_0\langle\nabla^g_{\dot{\gamma}}\xi(t),\nabla^g_{\dot{\gamma}}\eta(t)\rangle
dt.
\end{equation}
Denote by $\|\cdot\|_m$ the norm of it. Then the $m$-th iteration
map $\varphi_m: (\Lambda M, \langle\cdot,\cdot\rangle_1)\to (\Lambda
M, (\cdot,\cdot)_m)$ is an isometry up to the factor $m^2$, and
satisfies
\begin{equation}\label{e:6.2}
\left.\begin{array}{ll} &\varphi_m({\cal O})=S^1\cdot x^m, \quad
{\cal L}\circ\varphi_m=m^2{\cal L},\\
&\varphi_m([ms]\cdot\alpha)=[s]\cdot\varphi_m(\alpha)\quad\forall
([s],\alpha)\in S^1\times\Lambda M.
\end{array}\right\}
\end{equation}
Clearly, it induces a bundle embedding $\tilde\varphi_m: T\Lambda
M\to T{\Lambda M}$ given by
\begin{equation}\label{e:6.3}
(x,v)\mapsto (x^m, v^m),\qquad v^m(t)=v(mt)\;\hbox{for $t\in S^1$}.
\end{equation}
 Denote by $\hat N\varphi_m({\cal O})$ the normal bundle of
$\varphi_m({\cal O})$ with respect to the metric (\ref{e:6.1}). Its
fiber $\hat N\varphi_m({\cal O})_{x^m}$ at $x^m$ is the orthogonal
complementary of $\dot x^m\R$ in $(T_{x^m}{\Lambda M},
(\cdot,\cdot)_m)$. (By shrinking $\varepsilon>0$ if necessary) the
functional
\begin{equation}\label{e:6.4}
\hat{\cal F}: \hat N\varphi_m({\cal O})(\sqrt{m}\varepsilon)\to
\R,\;(y,v)\mapsto {\cal L}\circ{\rm EXP}(y,v)
\end{equation}
is well-defined, $S^1$-invariant, of class $C^{2-0}$, and satisfies
the (PS) condition. Consider the Banach vector subbundle of
$T_{\varphi_m({\cal O})}{\cal X}$,
  $X\hat{N}\varphi_m({\cal O}):=T_{\varphi_m({\cal O})}{\cal X}\cap
\hat{N}\varphi_m({\cal O})$, and denote by $\hat{\cal F}^X$ the
restriction of $\hat{\cal F}$ on the open subset $\hat
N\varphi_m({\cal O})(\sqrt{m}\varepsilon)\cap
X\hat{N}\varphi_m({\cal O})$ of $X\hat{N}\varphi_m({\cal O})$.
Unlike in the case of Riemannian geometry we need to prove a
splitting theorem for $\hat{\cal F}^X$ around $\varphi_m({\cal O})$
corresponding to Theorem~\ref{th:1.6} though the ideas are same. To
this end, for $r>0$ let
\begin{eqnarray*}
 && \hat{N}\varphi_m({\cal O})(r)=\{(y,v)\in \hat{N}\varphi_m({\cal O})\,|\,\|v\|_{m}<r\},\\
&& X\hat{N}\varphi_m({\cal O})(r)= \{(y,v)\in
X\hat{N}\varphi_m({\cal O})\,|\,\|v\|_{C^1}<r\},
\end{eqnarray*}
and denote by $\hat{N}\varphi_m({\cal O})(r)_y$ and
$X\hat{N}\varphi_m({\cal O})(r)_y$
 their fibers  at $y\in\varphi_m({\cal O})$, respectively.
  Let $\hat{\cal F}_y$ and $\hat{\cal F}^X_y$
denote the restrictions of the functionals $\hat{\cal F}$ and
$\hat{\cal F}^X$ to $\hat N\varphi_m({\cal
O})(\sqrt{m}\varepsilon)_y$ and $\hat N\varphi_m({\cal
O})(\sqrt{m}\varepsilon)_y\cap X\hat{N}\varphi_m({\cal O})_y$,
respectively.
 Denote by $\hat{\nabla}\hat{\cal F}_y$  the gradient of $\hat{\cal F}_y$ with
respect to the inner product given by (\ref{e:6.1}) on $\hat
N\varphi_m({\cal O})_y$, and by
 $\hat A_y$ be the restriction of
$\hat\nabla\hat{\cal F}_y$ to $\hat{N}\varphi_m({\cal
O})(\sqrt{m}\varepsilon)_y\cap X\hat{N}\varphi_m({\cal O})_y$. It is
clear that
\begin{equation}\label{e:6.5}
\hat A_{s\cdot y}(s\cdot v)=s\cdot \hat A_y(v)\quad\forall s\in
S^1,\; (y,v)\in \hat{N}\varphi_m({\cal O})(\sqrt{m}\varepsilon)\cap
X\hat{N}\varphi_m({\cal O}).
\end{equation}

\begin{claim}\label{cl:6.3}
 $\hat{\cal F}^X_y$ satisfies the
conditions of Theorem~\ref{th:A.5} on $X\hat{N}\varphi_m({\cal
O})(\delta)_y$ if $\delta\in (0, \varepsilon)$ is small enough.
\end{claim}

We postpone its proof to the end of this subsection. In particular,
it means: (i) the map $\hat A_y$  is $C^1$-smooth from
$X\hat{N}\varphi_m({\cal O})(\delta)_y$ to $X\hat{N}\varphi_m({\cal
O})_y$ (and so $\hat{\cal F}^X_y$ is $C^2$ in
$X\hat{N}\varphi_m({\cal O})(\delta)_y$); (ii) the symmetric
bilinear form $d^2\hat{\cal F}^X_y(0)$ has a continuous extension
$\hat{B}_y$ on $\hat{N}\varphi_m({\cal O})_y$; (iii) ${\bf
H}^-(\hat{B}_y)+ {\bf H}^0(\hat{B}_y)\subset X\hat{N}\varphi_m({\cal
O})_y$ is of finite dimension and there exist an orthogonal
decomposition
\begin{eqnarray}
\hat{N}\varphi_m({\cal O})_y={\bf H}^-(\hat{B}_y)\hat\oplus {\bf
H}^0(\hat{B}_y)\hat\oplus{\bf H}^+(\hat{B}_y)\label{e:6.6}
 \end{eqnarray}
with respect to the metric in (\ref{e:6.1}) and an induced  Banach
space direct sum decomposition $X\hat{N}\varphi_m({\cal O})_y={\bf
H}^-(\hat{B}_y){\dot{+} } {\bf H}^0(\hat{B}_y){\dot{+} }({\bf
H}^+(\hat{B}_y)\cap X\hat{N}\varphi_m({\cal O})_y)$.

For $y=x^m$, since both $\dim{\bf H}^-(B_y)$ and $\dim{\bf
H}^-(\hat{B}_y)$ (resp. $\dim{\bf H}^0(B_y)+1$ and $\dim{\bf
H}^0(\hat{B}_y)+1$)
 are equal to the Morse index (resp. nullity) of the symmetric bilinear form
$d^2({\cal L}\circ{\rm EXP}|_{T_y\mathcal{X}})(0_y)$, we obtain
$m^-(\varphi_m({\cal O}))={\rm dim}{\bf H}^-(B_y)={\rm dim}{\bf
H}^-(\hat{B}_y)$ and $m^0(\varphi_m({\cal O}))={\rm dim}{\bf
H}^0(B_y)={\rm dim}{\bf H}^0(\hat{B}_y)$.
 Let $\hat{\bf P}^\star_y$ be the orthogonal projections from
$\hat{N}\varphi_m({\cal O})_y$ onto ${\bf H}^\star(\hat{B}_y)$ in
(\ref{e:6.6}), $\star=+,0,-$. Since $\dim{\bf H}^0(\hat{B}_y)$ is
finite we may shrink $\epsilon\in (0,\varepsilon)$ in (\ref{e:1.19})
so small that
$${\bf H}^0(\hat{B})(\sqrt{m}\epsilon)_y:={\bf
H}^0(\hat{B}_y)\cap \hat{N}\varphi_m({\cal
O})(\sqrt{m}\epsilon)_y\subset X\hat{N}\varphi_m({\cal
O})(\delta)_y
$$
 and use  the implicit function theorem to get a
unique
 $C^1$ map
\begin{equation}\label{e:6.7}
 \hat{\mathfrak{h}}_y:{\bf H}^0(\hat{B})(\sqrt{m}\epsilon)_y\to
 {\bf H}^-(\hat{B}_y){\dot{+}}({\bf H}^+(\hat{B}_y)\cap
X\hat{N}\varphi_m({\cal O})_y)
\end{equation}
satisfying $\hat{\mathfrak{h}}_y(0)=0$, $d\hat{\mathfrak{h}}_y(0)=0$
and
\begin{eqnarray}\label{e:6.8}
 (\hat{\bf P}^+_y +\hat{\bf P}^-_y)\circ \hat{A}_y\bigl(v+ \hat{\mathfrak{h}}_y(v)\bigr)=0
\quad\forall v\in {\bf H}^0(\hat{B})(\sqrt{m}\epsilon)_y.
\end{eqnarray}
By (\ref{e:6.5}) the map $\hat{\mathfrak{h}}_y$ is also
$S_y^1$-equivariant. Define the functional
\begin{eqnarray}\label{e:6.9}
\hat{\cal L}^\circ_{\triangle y}:{\bf
H}^0(\hat{B})(\sqrt{m}\epsilon)_y\ni v\to{\cal L}\circ{\rm
EXP}_y\bigl(v+
 \hat{\mathfrak{h}}_y(v)\bigr)\in\R.
\end{eqnarray}
It is  $C^{2}$ and has the isolated critical point $0$. As in the
proof of Theorem~\ref{th:1.6} we may derive from these and
Theorem~\ref{th:A.5}:

\begin{proposition}\label{prop:6.4}
 By shrinking the above $\epsilon>0$ (if necessary)
 there exist an $S^1$-invariant open neighborhood $\widehat{V}$
of the zero section of $X\hat{N}\varphi_m({\cal O})$, an $S^1$-equivariant fiber-preserving,  $C^1$ map $\hat{\mathfrak{h}}$
given by (\ref{e:6.7}) and (\ref{e:6.8}), and an $S^1$-equivariant
fiber-preserving homeomorphism
$\widehat{\Psi}:X\hat{N}\varphi_m({\cal O})(\sqrt{m}\epsilon)\to
\widehat{V}$ such that
$$
\hat{\cal
F}^X\circ\widehat{\Psi}(x^m,v)=\frac{1}{2}\hat{B}_{x^m}(\hat{\bf
P}^+_{x^m}v, \hat{\bf P}^+_{x^m}v)-\|\hat{\bf P}^-_{x^m}v\|^2_m+
\hat{\cal L}^\circ_{\triangle x^m}(\hat{\bf P}^0_{x^m}v)
$$
for all $(x^m,v)\in X\hat{N}\varphi_m({\cal O})(\sqrt{m}\epsilon)$.
\end{proposition}

Let $\hat\beta_{x^m}(v)=\frac{1}{2}\hat{B}_{x^m}(\hat{\bf
P}^+_{x^m}v, \hat{\bf P}^+_{x^m}v)-\|\hat{\bf P}^-_{x^m}v\|^2_m$ and
$\hat\alpha_{x^m}(v)=\hat{\cal L}^\circ_{\triangle x^m}(\hat{\bf
P}^0_{x^m}v)$. Then
\begin{eqnarray}
\hat{\cal
F}^X_{x^m}\circ\widehat\Psi_{x^m}(v)&=&\hat\beta_{x^m}(v)+\hat\alpha_{x^m}(v)\quad\forall
v\in X\hat{N}\varphi_m({\cal O})(\sqrt{m}\epsilon)_{x^m},\label{e:6.10}\\
{\cal F}^X_{x}\circ\Psi_{x}(v)&=&\frac{1}{2}B_{x}({\bf P}^+_{x}v,
{\bf P}^+_{x}v)-\|{\bf P}^-_{x}v\|^2_1+
{\cal L}^\circ_{\triangle x}({\bf P}^0_{x}v)\nonumber\\
&\equiv&\beta_{x}(v)+ \alpha_{x}(v)\quad\forall v\in XN{\cal
O}(\epsilon)_{x}\label{e:6.11}
\end{eqnarray}
 by Theorem~\ref{th:1.6}.
Note that $\tilde\varphi_m(N{\cal O}(r)_x)\subset
\hat{N}\varphi_m({\cal O})(\sqrt{m}r)_{x^m}$ and
$\tilde\varphi_m(XN{\cal O}(r)_x)\subset X\hat{N}\varphi_m({\cal
O})(r)_{x^m}$ for any $x\in{\cal O}$ and $r>0$. Clearly,
(\ref{e:6.2})-(\ref{e:6.4}) imply
\begin{eqnarray}\label{e:6.12}
\hat{\cal F}_{x^m}(\tilde\varphi_m(v))=m^2{\cal
F}_{x}(v)\quad\forall (x,v)\in N{\cal O}(\varepsilon).
 \end{eqnarray}
We conclude
\begin{eqnarray}\label{e:6.13}
\hat\nabla\hat{\cal
F}_{x^m}(\tilde\varphi_m(v))=\tilde{\varphi}_{m}(\nabla{\cal
F}_{x}(v))\quad \forall (x,v)\in N{\cal O}(\varepsilon).
\end{eqnarray}
In fact, since
$(\tilde\varphi_m(\xi),\tilde\varphi_m(\eta))_m=m^2\langle\xi,\eta\rangle_1$
for any $\xi,\eta\in N{\cal O}_x$, by (\ref{e:6.12}) we only need to
prove that $\hat\nabla\hat{\cal F}_{x^m}(\tilde\varphi_m(v))\in{\rm
Im}(\tilde{\varphi}_{m})$ for all $v\in N{\cal O}(\varepsilon)_x$.
Consider the isometric action  on $\hat{N}\varphi_m({\cal
O})(\sqrt{m}\varepsilon)_{x^m}$ of group $\Z_m=\{e^{2\pi
ip/m}\,|\,p=0,\cdots,m-1\}$ given by
$$
(e^{2\pi ip/m}\cdot u)(t)=u(t+\frac{p}{m})\quad\forall
t\in\R,\;p=0,\cdots,m-1.
$$
Then $\tilde{\varphi}_{m}(N{\cal O}(\varepsilon)_x)$ is the fixed
point set of the action. Moreover it is not hard to see that the
functional $\hat{\cal F}_{x^m}$ is invariant under this
$\Z_m$-action. Hence
$$e^{2\pi ip/m}\cdot \hat\nabla\hat{\cal
F}_{x^m}(\tilde\varphi_m(v))=\hat\nabla\hat{\cal
F}_{x^m}(\tilde\varphi_m(v)),\quad p=0,\cdots,m-1,
$$
that is, $\hat\nabla\hat{\cal F}_{x^m}(\tilde\varphi_m(v))\in{\rm
Im}(\tilde{\varphi}_{m})$. (See \cite[page
23]{Pa1}.)

Clearly, $e^{2\pi ip/m}\cdot \hat A_{x^m}(\tilde\varphi_m(v))=\hat A_{x^m}(\tilde\varphi_m(v))\;
\forall v\in N{\cal O}(\varepsilon)_x\cap XN{\cal O}_x$, $p=0,\cdots,m-1$.
From these and the density of $X\hat{N}\varphi_m({\cal
O})_{x^m}$ in $\hat{N}\varphi_m({\cal O})_{x^m}$ we derive
  \begin{eqnarray*}
  (\hat{B}_{x^m}e^{2\pi ip/m}\cdot u, e^{2\pi ip/m}\cdot v)_m
 =(\hat{B}_{x^m}u, v)_m\quad\forall u,v\in \hat{N}\varphi_m({\cal O})_{x^m}
 \end{eqnarray*}
 and thus $\hat{B}_{x^m}e^{2\pi ip/m}\cdot u=e^{2\pi ip/m}\cdot(\hat{B}_{x^m} u)\;\forall u\in \hat{N}\varphi_m({\cal O})_{x^m}$, $p=0,1,\cdots,m$. The latter equalities imply
$\hat{B}_{x^m}(\tilde\varphi_m(N{\cal O}_x))\subset \tilde\varphi_m(N{\cal O}_x)$.
Moreover,  (\ref{e:6.13}) implies
\begin{eqnarray}\label{e:6.14}
 \hat A_{x^m}(\tilde\varphi_m(v))=\tilde{\varphi}_{m}(
A_{x}(v))\quad \forall (x,v)\in  N{\cal
O}(\varepsilon)\cap XN{\cal O}.
\end{eqnarray}
It easily follows from this that
$$
(\hat{B}_{x^m}\tilde\varphi_m(v), \tilde\varphi_m(u))_m=m^2\langle B_xv,u\rangle_1\quad\forall u,v\in N{\cal O}_x.
$$
 Using this and $\hat{B}_{x^m}(\tilde\varphi_m(N{\cal O}_x))\subset \tilde\varphi_m(N{\cal O}_x)$ we can obtain
\begin{eqnarray}\label{e:6.15}
 \hat B_{x^m}\circ\tilde\varphi_m=\tilde{\varphi}_{m}\circ
B_{x}\quad\forall x\in{\cal O},
\end{eqnarray}
which implies that
 $\tilde\varphi_m({\bf H}^\star(B_{x}))\subset{\bf
H}^\star(\hat{B}_{x^m})$, $\star=-,0,+$, i.e., $\tilde\varphi_m$
preserves the decompositions in (\ref{e:1.16}) and (\ref{e:6.6}).
For $v\in {\bf H}^0(B)(\epsilon)_{x}$ we have $A_{x}(v+
\mathfrak{h}_{x}(v))\in {\bf H}^0(B_{x})$ by (\ref{e:1.20}), and so
${\bf H}^0(\hat{B}_{x^m})\ni \tilde\varphi_{m}\left(A_{x}(v+
\mathfrak{h}_{x}(v))\right)=\hat{A}_{x^m}(\tilde\varphi_m(v)+
\tilde\varphi_{m}(\mathfrak{h}_{x}(v)))$ by (\ref{e:6.14}) after
shrinking $\epsilon>0$ so that $u+ \mathfrak{h}_{x}(u)\in XN{\cal
O}(\delta)_{x}\;\forall u\in {\bf H}^0(B)(\epsilon)_{x}$. Hence
$$
(\hat{\bf P}_{x^m}^-+\hat{\bf P}_{x^m}^+)
\hat{A}_{x^m}(\tilde\varphi_m(v)+
\tilde\varphi_{m}(\mathfrak{h}_{x}(v))) =0\quad\forall v\in {\bf
H}^0(B)(\epsilon)_{x}.
$$
By the implicit function theorem this and (\ref{e:6.8}) imply that
\begin{eqnarray}\label{e:6.16}
\hat{\mathfrak{h}}_{x^m}(\tilde\varphi_m(v))=
\tilde\varphi_m\left({\mathfrak{h}}_{x}(v)\right) \quad\forall
v\in{\bf H}^0(B)(\epsilon)_{x}.
 \end{eqnarray}
It follows from this and (\ref{e:6.9})-(\ref{e:6.12}) that for all
$v\in{\bf H}^0(B)(\epsilon)_{x}$,
\begin{eqnarray}\label{e:6.17}
\hat{\cal L}^\circ_{\triangle x^m}(v^m)&=&\hat{\cal
F}_{x^m}\bigl(v^m+
 \hat{\mathfrak{h}}_{x^m}(v^m)\bigr)
=\hat{\cal F}_{x^m}\bigl(v^m+
 \tilde\varphi_m({\mathfrak{h}}_{x}(v))\bigr)\nonumber\\
&=&m^2{\cal F}_{x}\bigl(v+ \mathfrak{h}_{x}(v)\bigr)=m^2{\cal
L}^\circ_{\triangle x}(v).
\end{eqnarray}
This and (\ref{e:6.15})-(\ref{e:6.16}) lead to
\begin{eqnarray}
&&\hat{\alpha}_{x^m}(\tilde\varphi_m(v))= m^2\alpha_{x}(v)
\quad\forall v\in {\bf H}^0(B)(\epsilon)_{x},\label{e:6.18}\\
&&\hat{\beta}_{x^m}(\tilde\varphi_m(v))= m^2\beta_{x}(v)
\quad\forall v\in ({\bf P}^+_{x}+{\bf P}^-_{x})(N{\cal
O}_{x}).\label{e:6.19}
\end{eqnarray}
Carefully checking the proof of Theorem~2.5 in \cite{Ji} it follows
from (\ref{e:6.12})-(\ref{e:6.15}) that
\begin{eqnarray}
\widehat{\Psi}_{x^m}(\tilde\varphi_m(v))=
\tilde\varphi_m\left(\Psi_{x}(v) \right)\quad\forall v\in XN{\cal
O}(\epsilon)_{x}.\label{e:6.20}
\end{eqnarray}
Here $\widehat{\Psi}_{x^m}$ and $\Psi_{x}$ are as in (\ref{e:6.11})
and (\ref{e:6.10}).

 Take $0<\epsilon_1\ll\epsilon$ such that
\begin{eqnarray*}
&&R_{x}:={\bf H}^0(B)(\epsilon_1)_{x}\dot{+} {\bf
H}^-(B)(\epsilon_1)_{x}\dot{+} ({\bf H}^+(B)_{x}\cap XN{\cal
O}(\epsilon_1)_{x})
\subset XN{\cal O}(\epsilon)_{x}\\
&&\hbox{(which implies $\tilde\varphi_m(R_{x})\subset X\hat{N}{\cal
O}(\epsilon)_{x^m}$ and}\\
&&\tilde\varphi_m({\bf H}^+(B)_{x}\cap XN{\cal
O}(\epsilon_1)_{x})\subset{\bf H}^+(\hat{B})_{x^m}\cap
X\hat{N}\tilde\varphi_m({\cal
O})(\epsilon_1)_{x^m}\;\hbox{), and}\\
&&\hat{R}_{x^m}:=\tilde\varphi_m({\bf
H}^0(B)(\epsilon_1)_{x})\dot{+} \tilde\varphi_m({\bf
H}^-(B)(\epsilon_1)_{x})\dot{+} ({\bf
H}^+(\hat{B})_{x^m}\cap X\hat{N}\tilde\varphi_m({\cal O})(\epsilon_1)_{x^m})\nonumber\\
&&\hspace{10mm}\subset X\hat{N}\tilde\varphi_m({\cal
O})(\sqrt{m}\epsilon)_{x^m}.
\end{eqnarray*}
Set $R=S^1\cdot R_{x}$ and $\hat{R}=S^1\cdot \hat{R}_{x^m}$. Let $V$
and $\widehat{V}$ be as in Theorem~\ref{th:1.6} and
Proposition~\ref{prop:6.4}, respectively. By (\ref{e:6.20}) we have
the commutative diagrams
$$
\begin{CD}
R @>\Psi >> \Psi(R)\subset V
 \\
@V \tilde\varphi_m VV @VV \tilde\varphi_m V \\
\hat{R}@>\widehat\Psi
>>\widehat\Psi(\hat{R})\subset \widehat{V}
\end{CD}
$$
and thus {\tiny
$$
\begin{CD}
H_\ast\bigl(({\cal F}^X\circ\Psi)_{c}\cap R, ({\cal
F}^X\circ\Psi)_{c}\cap (R\setminus{\cal O});\K\bigr) @>\Psi_\ast
>> H_\ast\bigl({\cal F}^X_c\cap\Psi(R),
{\cal F}^X_c\cap(\Psi(R)\setminus{\cal O});\K\bigr)
 \\
@V (\tilde\varphi_m)_\ast VV @VV (\tilde\varphi_m)_\ast V \\
H_\ast\bigl((\hat{\cal F}^X\circ\widehat\Psi)_{m^2c}\cap \hat{R},
(\hat{\cal F}^X\circ\widehat\Psi)_{m^2c}\cap
(\hat{R}\setminus\varphi_m({\cal O}));\K\bigr)@>(\widehat\Psi)_\ast
>>H_\ast\bigl(\hat{\cal F}^X_{m^2c}\cap\widehat\Psi(\hat{R}),
\hat{\cal
F}^X_{m^2c}\cap(\widehat\Psi(\hat{R})\setminus\varphi_m({\cal
O}));\K\bigr).
\end{CD}
$$}
It is the assumptions $m^-({\cal O})=m^-(\varphi_m({\cal O}))$ and
$m^0({\cal O})=m^0(\varphi_m({\cal O}))$ that imply  $\hat{R}$ is an
open neighborhood of $\varphi_m({\cal O})$ in
$X\hat{N}\varphi_m({\cal O})(\sqrt{m}\epsilon)$. Take
${V}_1=\Psi(R)$ and $\hat{V}_m=\widehat\Psi(\hat{R})$ and note that
${\cal F}^X\circ\Psi=\beta+\alpha$, $\hat{\cal
F}^X\circ\widehat\Psi=\hat\beta+\hat\alpha$ by (\ref{e:6.10}) and
(\ref{e:6.11}), and that $\Psi_\ast$ and $\widehat\Psi_\ast$  are
isomorphisms. We deduce

\begin{proposition}\label{prop:6.5}
$\tilde\varphi_m$ induces isomorphisms
\begin{eqnarray}
&&(\tilde\varphi_m)_\ast: H_\ast\bigl({\cal F}^{X}_c\cap {V}_1,
{\cal F}^{X}_c\cap {V}_1\setminus{\cal
O};\K\bigr)\to\nonumber\\
&&\hspace{30mm}H_\ast\bigl(\hat{\cal F}^{X}_{m^2c}\cap \hat{V}_m,
\hat{\cal F}^{X}_{m^2c}\cap \hat{V}_m\setminus\varphi_m({\cal
O});\K\bigr).\label{e:6.21}
\end{eqnarray}
if and only if $\tilde\varphi_m$ induces  isomorphisms
\begin{eqnarray}\label{e:6.22}
&& H_\ast\bigl((\beta+ \alpha)_{c}\cap R, (\beta+ \alpha)_{c}\cap
(R\setminus{\cal O});\K\bigr)\to\nonumber\\
&&\hspace{20mm} H_\ast\bigl((\hat\beta+ \hat\alpha)_{m^2c}\cap
\hat{R}, (\hat\beta+ \hat\alpha)_{m^2c}\cap
(\hat{R}\setminus\varphi_m({\cal O}));\K\bigr).
\end{eqnarray}
\end{proposition}

Since the deformation  retracts
\begin{eqnarray*}
&&{\bf H}^0(B)\dot{+} {\bf H}^-(B)\dot{+} \bigl({\bf H}^+(B)\cap
XN{\cal O}\bigr) \times [0, 1]\to{\bf H}^0(B)\oplus {\bf
H}^-(B)\dot{+} \bigl({\bf
H}^+(B)\cap XN{\cal O}\bigr)\\
&&\hspace{20mm}(x, v^0+ v^-+ v^+)\mapsto (x, v^0+ v^-+ tv^+),\\
&&{\bf H}^0(\hat{B})\dot{+} {\bf H}^-(\hat{B})\dot{+}
\bigl({\bf H}^+(\hat{B})\cap X\hat{N}\varphi_m({\cal O})\bigr)
\times [0,
1]\\
&&\hspace{50mm}\to{\bf H}^0(\hat{B})\dot{+} {\bf
H}^-(\hat{B})\dot{+} \bigl({\bf
H}^+(\hat{B})\cap X\hat{N}\varphi_m({\cal O})\bigr)\\
&&\hspace{20mm}(x, v^0+ v^-+ v^+)\mapsto (x, v^0+ v^-+ tv^+)
\end{eqnarray*}
commute with $\tilde\varphi_m$, (\ref{e:6.22}) are isomorphisms if
and only if $\tilde\varphi_m$ induces  isomorphisms from
$H_\ast\bigl((\beta+ \alpha)_{c}\cap \Box, (\beta+ \alpha)_{c}\cap
(\Box\setminus{\cal O});\K\bigr)$ to
\begin{eqnarray*}
H_\ast\bigl((\hat\beta+ \hat\alpha)_{m^2c}\cap
\tilde\varphi_m(\Box),  (\hat\beta+ \hat\alpha)_{m^2c}\cap
(\tilde\varphi_m(\Box)\setminus\varphi_m({\cal O}));\K\bigr),
\end{eqnarray*}
where $\Box={\bf H}^0(B)(\epsilon_1)\oplus {\bf
H}^-(B)(\epsilon_1)$. The latter is clear since
$\tilde\varphi_m:\Box\to \tilde\varphi_m(\Box)$ is a linear
diffeomorphism and $(\hat\beta+
\hat\alpha)(\tilde\varphi_m(\xi))=m^2(\beta+ \alpha)(\xi)\;\forall
\xi\in\Box$ by (\ref{e:6.18}) and (\ref{e:6.19}). It is not hard to
see that Claim~\ref{cl:6.2} is equivalent to that (\ref{e:6.21})
are isomorphisms. Hence the proof of Theorem~\ref{th:6.1} is
completed once Claim~\ref{cl:6.3} is proved.
\hfill$\Box$\vspace{2mm}

\noindent{\bf Proof of Claim~\ref{cl:6.3}}. Corresponding to the
maps $\hat{\cal F}_y$, $\hat{\cal F}^X_y$, $\hat A_y$ and $\hat B_y$
between (\ref{e:6.4}) and (\ref{e:6.6}) we have maps $\hat{\cal
F}^\ast_y$, $\hat{\cal F}^{\ast X}_y$, $\hat A^\ast_y$ and $\hat
B^\ast_y$ if $L=F^2$ is replaced by the $L^\ast$. Let $\hat{\bf
B}_x^\ast$ be defined by $d^2 \hat{\cal F}^{\ast X}_{y}
  (\zeta)[\xi,\eta]=\bigl(\hat{\bf B}_y^\ast(\zeta)\xi,
  \eta\bigr)_{m}$. Then $\hat{\bf
B}_x^\ast(0)=\hat B_y$. As in (\ref{e:5.1})-(\ref{e:5.2}) it holds
that $\hat{A}^\ast_{y}(v)= \hat{A}_{y}(v)\;\forall (y,v)\in
X\hat{N}\varphi_m({\cal O})(\sqrt{m}\varepsilon)_y$ (by shrinking
$\varepsilon>0$ if necessary). In particular, this implies
$\hat{B}_y^\ast= \hat{B}_y\;\forall y\in\varphi_m({\cal O})$, and
hence ${\bf H}^\star(\hat{B}^\ast)={\bf
H}^\star(\hat{B}),\;\star=+,0,-$. It suffices to prove the
corresponding conclusions for $\hat{\cal F}^{\ast X}_y$ on
$X\hat{N}\varphi_m({\cal O})(\delta)_y$ with $y=\gamma_0^m$ for
$\delta>0$ small enough. The proof is almost as same as that of
Proposition~\ref{prop:5.1}. We only outline it. Following the
notations in Section~\ref{sec:5.3} we may obtain  a
$E_{{\gamma_0}^m}$-1-invariant unit orthogonal parallel frame field
along $\gamma_{0}^m$, $(\hat e_1(t),\cdots, \hat
e_n(t)):=(e_1(mt),\cdots, e_n(mt))$, that is, it satisfies $(\hat
e_1(t+1),\cdots, \hat e_n(t+1))=(\hat e_1(t),\cdots, \hat
e_n(t))E_{\gamma_0^m}\;\forall t\in\R$. (It is easily checked that
$E_{\gamma_0^m}=(E_{\gamma_0})^m$). A curve
$\xi=(\xi_1,\cdots,\xi_n):\R\to\R^n$ is called $E_{\gamma_0^m}$-{\bf
1-invariant} if $\xi(t+1)^T=E_{\gamma_0^m}\xi(t)^T\;\forall t\in\R$.
 Let $X_{\gamma_0^m}$ be the Banach space of all
$E_{\gamma_0^m}$-1-invariant  $C^1$ curves from $\R$ to $\R^n$
according to the usual $C^1$-norm.
 Let
 $\hat{H}_{\gamma_0^m}$ be
the Hilbert space of all $E_{\gamma_0^m}$-1-invariant $W^{1,2}_{\rm
loc}$ curves from $\R$ to $\R^n$  with inner product
\begin{eqnarray}\label{e:6.23}
\langle\xi,\eta\rangle_{1,2,m}=\int_0^1[m^2(\xi(t),\eta(t))_{\R^n}+
(\dot\xi(t),\dot\eta(t))_{\R^n}]dt.
\end{eqnarray}
It contains $X_{\gamma_0^m}$ as a dense subset. Write
$\|\cdot\|_{1,2,m}$ the induced norm on $\hat{H}_{\gamma_0^m}$.
There exists a Hilbert space isomorphism
 $\mathfrak{I}_{\gamma_0^m}: (\hat{H}_{\gamma_0^m},
\langle\cdot,\cdot\rangle_{1,2,m})\to (T_{\gamma_{0}^m}{\Lambda M},
(\cdot,\cdot)_m)$
 given by $\mathfrak{I}_{\gamma_0^m}(\xi)=\sum^n_{j=1}\xi_j\hat{e}_j$.
Under the isomorphisms $\mathfrak{I}_{\gamma_0}$ below
(\ref{e:5.22}) and ${\mathfrak{I}}_{\gamma_0^m}$ the iteration
$\tilde\varphi_m: T_{\gamma_0}{\Lambda M}\to
T_{\gamma_{0}^m}{\Lambda M}$ in (\ref{e:6.3}) corresponds to the map
\begin{equation}\label{e:6.24}
\tilde\varphi_{m{\gamma_0}}:H_{\gamma_0}\to \hat{H}_{\gamma_0^m}
\end{equation}
given by
$\tilde\varphi_{m{\gamma_0}}(\xi)(t)=\xi^m(t)=\xi(mt)\;\forall
t\in\R$.
 So $\tilde\varphi_m\circ{\mathfrak{I}}_{\gamma_0^m}=
\tilde\varphi_{m{\gamma_0}}\circ\mathfrak{I}_{\gamma_0}$.
 With the exponential map  $\exp$ of $g$  we define the $C^k$ map
$$
\phi_{\gamma_0^m}:\R\times B^n_{2\rho}(0)\to
M,\;(t,u)\mapsto\exp_{\gamma_{0}^m(t)}\left(\sum^n_{i=1}u_i\hat
e_i(t)\right)
$$
 for some small open ball
$B^n_{2\rho}(0)\subset\R^n$. It  satisfies
\begin{eqnarray*}
&&\phi_{\gamma_0^m}(t,z)=\phi_{\gamma_0}(mt,z),\quad d\phi_{\gamma_0^m}(t,z)[(1,u)]=d\phi_{\gamma_0}(mt,z)[(m,u)],\\
&&\phi_{\gamma_0^m}(t+1,z)=\phi_{\gamma_0^m}(t, (E_{\gamma_0^m} z^T)^T)\quad\hbox{and}\\
&&d\phi_{\gamma_0^m}(t+1,z)[(1,u)]=d\phi_{\gamma_0^m}(t,
(E_{\gamma_0^m} z^T)^T)[(1, (E_{\gamma_0^m} u^T)^T)]
\end{eqnarray*}
for $(t,z,u)\in \R\times B^n_{2\rho}(0)\times\R^n$ (by shrinking
$\rho>0$ if necessary), where $\phi_{\gamma_0}$ is as in
Section~\ref{sec:5.3}. This yields a $C^{k-3}$ coordinate chart around
$\gamma_{0}^m$ on ${\Lambda M}$,
$$
\Phi_{\gamma_0^m}:{\bf B}_{2\rho}(\hat{H}_{\gamma_0^m}):=\{\xi\in
 H_{\gamma_0^m}\,|\,\|\xi\|_{1,2,m}<2\rho\} \to{\Lambda M}
$$
given by $\Phi_{\gamma_0^m}(\xi)(t)=\phi_{\gamma_0^m}(t,\xi(t))$ for
$\xi\in {\bf B}_{2\rho}(\hat{H}_{\gamma_0^m})$. Clearly
$d\Phi_{\gamma_0^m}(0)={\mathfrak{I}}_{\gamma_0^m}$. Define
$$
 L^\ast_{\gamma_0^m}(t, z, u)=L^\ast\bigr(\phi_{\gamma_0^m}(t,z),
d\phi_{\gamma_0^m}(t,z)[(1,u)]\bigl)\quad\forall (t,z,u)\in \R\times
B^n_{2\rho}(0)\times\R^n.
$$
Then  $L^\ast_{\gamma_0^m}(t,0,0)\equiv m^2 c\;\forall t\in\R$ and
as in (\ref{e:5.24})-(\ref{e:5.26}) we have
\begin{eqnarray*}
 L^\ast_{\gamma_0^m}(t+1,\xi(t+1),\dot\xi(t+1))
&=& L^\ast_{\gamma_0^m}(t,\xi(t),\dot\xi(t))\quad\forall t\in\R,\\
\partial_x L^\ast_{\gamma_0^m}(t+1,\xi(t+1),\dot\xi(t+1))&=&\partial_x
L^\ast_{\gamma_0^m}(t,\xi(t),\dot\xi(t))E_{\gamma_0^m}\quad\forall t\in\R,\\
\partial_v L^\ast_{\gamma_0^m}(t+1,\xi(t+1),\dot\xi(t+1))&=&\partial_v
L^\ast_{\gamma_0^m}(t,\xi(t),\dot\xi(t))E_{\gamma_0^m}\quad\forall
t\in\R
\end{eqnarray*}
 for any $\xi\in {\bf
B}_{2\rho}(\hat{H}_{\gamma_0^m})$ since
$\xi(t+1)^T=E_{\gamma_0^m}\xi(t)^T$ and
$\dot\xi(t+1)^T=E_{\gamma_0^m}\dot\xi(t)^T$. Define  the action
functional ${\cal L}^\ast_{\gamma_0^m}: {\bf
B}_{2\rho}(\hat{H}_{\gamma_0^m})\to\R$ by
$$
{\cal L}^\ast_{\gamma_0^m}(\xi)=\int^1_0
L^\ast_{\gamma_0^m}(t,\xi(t),\dot\xi(t))dt.
$$
Then ${\cal L}^\ast_{\gamma_0^m}={\cal
L}^\ast\circ\Phi_{\gamma_0^m}$ on ${\bf
B}_{2\rho}(\hat{H}_{\gamma_0^m})$. As in Section~\ref{sec:5.3}
 using (\ref{e:B.15})-(\ref{e:B.16}) and (\ref{e:B.30})-(\ref{e:B.33})   we may prove that
 it  is $C^{2-0}$ and that the restriction $\hat{\mathbb{A}}^\ast_{\gamma_0^m}$ of the gradient of ${\cal
 L}^\ast_{\gamma_0^m}$ on $\hat{H}_{\gamma_0^m}$  to
${\bf B}_{2\rho}(\hat{H}_{\gamma_0^m})\cap X_{\gamma_0^m}$ is a
$C^1$ map from the open subset ${\bf
B}_{2\rho}(\hat{H}_{\gamma_0^m})\cap X_{\gamma_0^m}$ of
$X_{\gamma_0^m}$ to $X_{\gamma_0^m}$
 (and hence the restriction of ${\cal L}^\ast_{\gamma_0^m}$ to ${\bf
B}_{2\rho}(\hat{H}_{\gamma_0^m})\cap X_{\gamma_0^m}$, denoted by
${\cal L}^{\ast X}_{\gamma_0^m}$, is $C^2$). Since
\begin{eqnarray*}
 d^2 {\cal L}^{\ast X}_{\gamma_0^m}
  (\zeta)[\xi,\eta]
   = \int_0^{1} \Bigl(\!\! \!\!\!&&\!\!\!\!\!\partial_{vv}
     L^\ast_{\gamma_0^m}\bigl(t,\zeta(t),\dot{\zeta}(t)\bigr)
\bigl[\dot{\xi}(t), \dot{\eta}(t)\bigr] \nonumber\\
&&+ \partial_{qv}
  L^\ast_{\gamma_0^m}\bigl(t,\zeta(t), \dot{\zeta}(t)\bigr)
\bigl[\xi(t), \dot{\eta}(t)\bigr]\nonumber \\
&& + \partial_{vq}
  L^\ast_{\gamma_0^m}\bigl(t,\zeta(t),\dot{\zeta}(t)\bigr)
\bigl[\dot{\xi}(t), \eta(t)\bigr] \nonumber\\
&&+  \partial_{qq} L^\ast_{\gamma_0^m}\bigl(t,\zeta(t),
\dot{\zeta}(t)\bigr) \bigl[\xi(t), \eta(t)\bigr]\Bigr) \, dt
\end{eqnarray*}
for any $\zeta\in {\bf B}_{2\rho}(\hat{H}_{\gamma_0^m})\cap
X_{\gamma_0^m}$ and $\xi,\eta\in X_{\gamma_0^m}$,  it is easily
checked that the corresponding properties to (i) and (ii) below
Lemma~\ref{lem:4.2} hold for ${\cal L}^{\ast X}_{\gamma_0^m}$, and
so
 there exists a  continuous map $\hat{\mathbb{B}}^\ast_{\gamma_0^m}:
 {\bf B}_{2\rho}(\hat{H}_{\gamma_0^m})\cap X_{\gamma_0^m}\to
{\cal L}_s(\hat{H}_{\gamma_0^m})$ with respect to the induced
topology on ${\bf B}_{2\rho}(\hat{H}_{\gamma_0^m})\cap
X_{\gamma_0^m}$ by $X_{\gamma_0^m}$, such that
$$
\langle d\hat{\mathbb{A}}^\ast_{\gamma_0^m}(\zeta)[\xi],
\eta\rangle_{1,2,m}=d^2 {\cal L}^{\ast X}_{\gamma_0^m}
  (\zeta)[\xi,\eta]=\langle\hat{\mathbb{B}}^\ast_{\gamma_0^m}(\zeta)\xi,
  \eta\rangle_{1,2,m}\quad\forall \xi, \eta\in X_{\gamma_0^m}.
$$
From these we may obtain the following similar result to
Claim~\ref{cl:5.6}:  \textsf{Around the critical point $0\in
\hat{H}_{\gamma_0^m}$, $(\hat{H}_{\gamma_0^m}, X_{\gamma_0^m}, {\cal
L}^\ast_{\gamma_0^m}, \hat{\mathbb{A}}^\ast_{\gamma_0^m},
\hat{\mathbb{B}}^\ast_{\gamma_0^m})$ satisfy the conditions of
Theorem~\ref{th:A.1} except that the critical point $0$ is not
isolated.} Suitably changing the arguments below Claim~\ref{cl:5.6}
we may prove that $(\hat{N}\varphi_m({\cal O})_{y},
X\hat{N}\varphi_m({\cal O})_{y},
 \hat{\cal F}^\ast_{y},
  \hat{A}^\ast_{y}, \hat{\bf B}^\ast_{y})$ with $y=\gamma^m_{0}$ satisfies
the conditions of Theorem~\ref{th:A.1}  around the critical point
$0=0_{y}$, which implies that the conditions of Theorem~\ref{th:A.5}
are satisfied for $\hat{\cal F}^{\ast X}_y$ (and hence $\hat{\cal
F}^{X}_y$)  on $X\hat{N}\varphi_m({\cal O})(\delta)_y$.
\hfill$\Box$\vspace{2mm}

\subsection{Proof of Theorem~\ref{th:1.10}}\label{sec:6.2}

When ${\cal O}$ is replaced by $\varphi_m({\cal O})$, for
$y\in\varphi_m({\cal O})$ let $A_y, B_y, \mathfrak{h}_y, {\bf
P}^+_y, {\bf P}^-_y,  {\cal L}^\circ_{\triangle y}$ be the
corresponding maps with those in (\ref{e:1.15})-(\ref{e:1.21}).
Write $x=\gamma_0$ in the following. Then for all $q\in\{0\}\cup\N$
we have
\begin{eqnarray*}
&&\dim C_{q}({\cal L}|_{\mathscr{N}_x}, x;\K)=\dim C_{q}({\cal
L}^\circ_{\triangle x}, 0;\K),\\
&&\dim C_{q}({\cal L}|_{\mathscr{N}_x}, x;\K)^{S_x^1}=\dim
C_{q}({\cal L}^\circ_{\triangle x}, 0;\K)^{S_x^1}\\
&&\dim C_{q}({\cal L}|_{\mathscr{N}_{x^m}}, x^m;\K)=\dim C_{q}({\cal
L}^\circ_{\triangle x^m}, 0;\K),\\
&&\dim C_{q}({\cal L}|_{\mathscr{N}_{x^m}}, x^m;\K)^{S_{x^m}^1}=\dim
C_{q}({\cal L}^\circ_{\triangle x^m}, 0;\K)^{S_{x^m}^1}.
\end{eqnarray*}
So the proofs of (\ref{e:1.24}) and (\ref{e:1.25}) are reduced to
prove that for any $q\in\{0\}\cup\N$,
\begin{eqnarray}
&&\dim C_{q}({\cal L}^\circ_{\triangle x}, 0;\K)=\dim C_{q}({\cal
L}^\circ_{\triangle x^m}, 0;\K),\label{e:6.25}\\
&&\dim C_{q}({\cal L}^\circ_{\triangle x}, 0;\K)^{S_x^1}=\dim
C_{q}({\cal L}^\circ_{\triangle x^m}, 0;\K)^{S_x^1}.\label{e:6.26}
\end{eqnarray}

\begin{claim}\label{cl:6.6}
$\dim C_{\ast}({\cal L}^{\circ}_{\triangle x}, 0;\K)=\dim
C_{\ast}(\hat{\cal L}^{\circ}_{\triangle x^m},0;\K)$ for $\hat{\cal
L}^{\circ}_{\triangle x^m}$ in (\ref{e:6.9}).
\end{claim}

\noindent{\bf Proof}. Since $m^0({\cal O})=m^0(\varphi_m({\cal
O}))$, $\tilde\varphi_m$ restricts to a linear homeomorphism
$\tilde\varphi^x_m: {\bf H}^0(B)_{x}\to{\bf H}^0(\hat{B})_{x^m}$. By
(\ref{e:6.16}) it also restricts to homeomorphism
\begin{equation}\label{e:6.27}
({\cal L}_{\triangle x}^\circ)_c\cap {\bf
H}^0(B)({\epsilon/\sqrt{m}})_x\to (\hat{\cal L}_{\triangle
x^m}^\circ)_{m^2c}\cap {\bf H}^0(\hat{B})(\epsilon)_{x^m},
\end{equation}
which  induces isomorphisms
\begin{eqnarray*}
&&\qquad\quad H_\ast\left(({\cal L}_{\triangle x}^\circ)_c\cap {\bf
H}^0(B)({\epsilon/\sqrt{m}})_x, ({\cal L}_{\triangle x}^\circ)_c\cap
{\bf H}^0(B)({\epsilon/\sqrt{m}})_x\setminus\{0\};\K\right)\nonumber\\
&&\to H_\ast\left((\hat{\cal L}_{\triangle x^m}^\circ)_{m^2c}\cap
{\bf H}^0(\hat{B})(\epsilon)_{x^m}, (\hat{\cal L}_{\triangle
x^m}^\circ)_{m^2c}\cap {\bf
H}^0(\hat{B})(\epsilon)_{x^m}\setminus\{0\};\K\right).
\end{eqnarray*}

\begin{claim}\label{cl:6.7}
$\dim C_{\ast}({\cal L}^{\circ}_{\triangle x}, 0;\K)^{S^1_{ x}}=\dim
C_{\ast}(\hat{\cal L}^{\circ}_{\triangle x^m}, 0;\K)^{S^1_{x^m}}$.
\end{claim}

\noindent{\bf Proof}.   Recall that
$S^1=\R/\Z=\{[s]\,|\,[s]=s+\Z,\;s\in\R\}$. Since $S^1_{x}$ is a
finite cyclic subgroup of $S^1$, we may assume $S^1_{x}=\{[j/l]\,|\,
j=0,\cdots,l-1\}$ for some $l\in\N$. Then
$S^1_{x^m}=\{[\frac{j}{ml}]\,|\, j=0,\cdots,ml-1\}$. Clearly,
$\tilde\varphi^x_{m}([1/l]\cdot\xi)=[\frac{1}{ml}]\cdot\tilde\varphi^x_{m}(\xi)$.
Observe   that ${\mathfrak{h}}_{x}$ (resp.
$\hat{\mathfrak{h}}_{x^m}$) is $S^1_{x}$ (resp. $S^1_{x^m}$)
equivariant. From  (\ref{e:6.16}) and (\ref{e:6.17}) we derive that
the homeomorphism in (\ref{e:6.27}) induces a homeomorphism
$$
({\cal L}_{\triangle x}^\circ)_c\cap {\bf
H}^0(B)({\epsilon/\sqrt{m}})_x/S^1_x\to(\hat{\cal L}_{\triangle
x^m}^\circ)_{m^2c}\cap {\bf H}^0(\hat{B})(\epsilon)_{x^m}/S^1_{x^m}
$$
and therefore  isomorphisms
\begin{eqnarray*}
&&\qquad H_\ast\left(({\cal L}_{\triangle x}^\circ)_c\cap {\bf
H}^0(B)({\epsilon/\sqrt{m}})_x/S^1_x, \bigr(({\cal L}_{\triangle
x}^\circ)_c\cap
{\bf H}^0(B)({\epsilon/\sqrt{m}})_x\setminus\{0\}\bigl)/S^1_x;\K\right)\nonumber\\
&&\to H_\ast\left((\hat{\cal L}_{\triangle x^m}^\circ)_{m^2c}\cap
{\bf H}^0(\hat{B})(\epsilon)_{x^m}/S^1_{x^m}, \bigl((\hat{\cal
L}_{\triangle x^m}^\circ)_{m^2c}\cap {\bf
H}^0(\hat{B})(\epsilon)_{x^m}\setminus\{0\}\bigr)/S^1_{x^m};\K\right).
\end{eqnarray*}
As in the arguments  below (\ref{e:5.13}), this may induce
isomorphisms
\begin{eqnarray*}
&&\qquad H_\ast(({\cal L}_{\triangle x}^\circ)_c\cap {\bf
H}^0(B)({\epsilon/\sqrt{m}})_x, ({\cal L}_{\triangle x}^\circ)_c\cap
{\bf H}^0(B)({\epsilon/\sqrt{m}})_x\setminus\{0\};\K)^{S^1_x}\nonumber\\
&&\to H_\ast((\hat{\cal L}_{\triangle x^m}^\circ)_{m^2c}\cap {\bf
H}^0(\hat{B})(\epsilon)_{x^m}, (\hat{\cal L}_{\triangle
x^m}^\circ)_{m^2c}\cap {\bf
H}^0(\hat{B})(\epsilon)_{x^m}\setminus\{0\};\K)^{S^1_{x^m}}
\end{eqnarray*}
if  the characteristic of $\K$ is zero or prime up to orders of
$S^1_{\gamma_0}$ and $S^1_{\gamma_{0}^m}$.
 \hfill$\Box$\vspace{1mm}

By Claims~\ref{cl:6.6},~\ref{cl:6.7} for proofs of (\ref{e:6.25})
and (\ref{e:6.26}) it suffices to show
\begin{eqnarray}\label{e:6.28}
&&\dim C_{\ast}(\hat{\cal L}^{\circ}_{\triangle x^m},0;\K)=\dim
C_{\ast}({\cal L}^{\circ}_{\triangle x^m},0;\K),\\
&&\dim C_{\ast}(\hat{\cal L}^{\circ}_{\triangle x^m},
0;\K)^{S^1_{x^m}}=\dim C_{\ast}({\cal L}^{\circ}_{\triangle x^m},
0;\K)^{S^1_{x^m}}\label{e:6.29}
\end{eqnarray}
if the characteristic of $\K$ is zero or prime up to orders of
$S^1_{\gamma_0}$ and $S^1_{\gamma_{0}^m}$.

\noindent{\bf Proof of of  (\ref{e:6.28})}. By Claim~\ref{cl:6.3}
and Proposition~\ref{prop:5.1}
 $\hat{\cal F}^X_{x^m}$ (resp. ${\cal F}^X_{x^m}$)
satisfies  Theorem~\ref{th:A.5} around $0$ on $X\hat
N\tilde\varphi_m({\cal O})(\delta)_{x^m}$ (resp.
$N\tilde\varphi_m({\cal O})(\delta)_{x^m}$), respectively.
 Let $m^-=m^-({\cal
O})=m^-(\varphi_m({\cal O}))$. Applying Corollary~\ref{cor:A.6} to
${\cal F}^X_{x^m}$ and $\hat{\cal F}^X_{x^m}$ gives
\begin{eqnarray*}
C_{q-m^-}({\cal L}^\circ_{\triangle x^m},0;\K)= C_q({\cal
F}^X_{x^m},0;\K)\quad\hbox{and}\quad
 C_{q-m^-}(\hat{\cal L}^\circ_{\triangle x^m},0;\K)=C_q(\hat{\cal
F}^X_{x^m},0;\K)
\end{eqnarray*}
for all $q\in\N\cup\{0\}$.
  So it suffices to prove that
\begin{equation}\label{e:6.30}
C_\ast({\cal F}^X_{x^m},0;\K)=C_\ast(\hat{\cal F}^X_{x^m},0;\K).
\end{equation}
Since $C_\ast({\cal F}^X_{x^m},0;\K)=C_{\ast}({\cal L}|_{\bf
S},x^m;\K)$ and $C_\ast(\hat{\cal F}^X_{x^m},0;\K)=C_{\ast}({\cal
L}|_{\hat{\bf S}},x^m;\K)$, where
\begin{eqnarray*}
{\bf S}={\rm EXP}_{x^m}\bigl(XN(\varphi_m({\cal
O}))(\delta)_{x^m}\bigr)\quad\hbox{and}\quad\hat{\bf S}={\rm
EXP}_{x^m}\bigl(X\hat N(\varphi_m({\cal O}))(\delta)_{x^m}\bigr)
\end{eqnarray*}
are two at least $C^{k-3}$-smooth submanifolds of codimension one in
${\cal X}$ and both are transversely intersecting at $x^m$ with
$S^1\cdot x^m=\varphi_m({\cal O})$, we only need to prove that
\begin{equation}\label{e:6.31}
C_{\ast}({\cal L}|_{\bf S}, x^m;\K)=C_{\ast}({\cal L}|_{\hat{\bf
S}}, x^m;\K).
\end{equation}
Note that the $S^1$-action on ${\cal X}$ is $C^1$. We get a $C^1$
map
$$
\mathfrak{T}:S^1\times{\bf S}\to{\cal X},\; ([s], x)\mapsto [s]\cdot
x.
$$
Since ${\bf S}$ is transversal to $\varphi_m({\cal O})$ at $x^m$,
the differential $d\mathfrak{T}([0], x^m)$ is an isomorphism. So
there exists a neighborhood ${\cal N}(x^m)$ of $x^m$ in ${\bf S}$
and $0<r\ll 1/2$ such that $\mathfrak{T}$ is a diffeomorphism from
$\{[s]\,|\, s\in (-r,r)\}\times {\cal N}(x^m)$ onto a neighborhood
of $x^m$ in ${\cal X}$. Because $\hat{\bf S}$ is transversal to
$\varphi_m({\cal O})$ at $x^m$ as well $\mathfrak{T}^{-1}(\hat{\bf
S})$ is a submanifold in $\{[s]\,|\, s\in (-r,r)\}\times {\cal
N}(x^m)$ which is transversal to $\{[s]\,|\, s\in (-r,r)\}\times
\{x^m\}$. It follows that there exists a $C^1$-map
$\mathfrak{S}:{\cal N}(x^m)\to \{[s]\,|\, s\in (-r,r)\}$ such that
the graph ${\rm Gr}(\mathfrak{S})$ is a neighborhood of $x^m\equiv
([0], x^m)$ in $\mathfrak{T}^{-1}(\hat{\bf S})$
 by shrinking $\{[s]\,|\, s\in
(-r,r)\}\times {\cal N}(x^m)$ if necessary. This implies that the
composition
$$
\Theta:{\cal N}(x^m)\ni x\mapsto (\mathfrak{S}(x),x)\in {\rm
Gr}(\mathfrak{S})\stackrel{\mathfrak{T}}{\longrightarrow}{\cal X}
$$
is a $C^1$ diffeomorphism from ${\cal N}(x^m)$ onto a neighborhood
${\cal W}(x^m)$ of $x^m$ in $\hat{\bf S}$. Clearly,
$\Theta(x^m)=x^m$ and ${\cal L}|_{\hat{\bf S}}(\Theta(x))={\cal
L}([\mathfrak{S}(x)]\cdot x)= {\cal L}(x)={\cal L}|_{{\bf S}}(x)$
for any $x\in {\cal N}(x^m)$. Hence $\Theta$ induces an isomorphism
from $C_{\ast}({\cal L}|_{\bf S},x^m;\K)$ to $C_{\ast}({\cal
L}|_{\hat{\bf S}}, x^m;\K)$. This proves (\ref{e:6.31}).
\vspace{1mm}

\noindent{\bf Proof of of  (\ref{e:6.29})}. Let $d=m^2c$. Recall
that $S^1_{x^m}=\{[\frac{j}{ml}]\,|\, j=0,\cdots,ml-1\}$ for some
$l\in\N$. It is easily checked that $S^1_{x^m}$ acts on
$XN\tilde\varphi_m({\cal O})(\delta)_{x^m}$ (resp. $X\hat
N\tilde\varphi_m({\cal O})(\delta)_{x^m}$) and therefore on ${\bf
S}$ (resp. $\hat{\bf S}$). We shrink the above $r>0$ so that
$r<\frac{1}{4ml}$. Take a $S^1_{x^m}$-invariant open neighborhood
${\cal N}(x^m)^*$ of $x^m$ in ${\bf S}$ such that ${\cal
N}(x^m)^*\subset {\cal N}(x^m)$. Given a $p\in {\cal N}(x^m)^*$, by
the construction of $\Theta$ we have
$\Theta(p)=[\mathfrak{S}(p)]\cdot p$ and
\begin{eqnarray*}
\Theta([\frac{1}{ml}]\cdot p)&=&[\mathfrak{S}([\frac{1}{ml}]\cdot
p)]\cdot([\frac{1}{ml}]\cdot p)=[\mathfrak{S}([\frac{1}{ml}]\cdot p)+ \frac{1}{ml}]\cdot p\\
&=&[\mathfrak{S}([\frac{1}{ml}]\cdot p)-\mathfrak{S}(p)+
\frac{1}{ml}]\cdot \Theta(p).
\end{eqnarray*}
 Note that
$\mathfrak{S}([\frac{1}{ml}]\cdot p)$ and $-\mathfrak{S}(p)$ belong
to $(-r,r)$, and that  $0<r<\frac{1}{4ml}$. We deduce that
 $\mathfrak{S}([\frac{1}{ml}]\cdot p)-\mathfrak{S}(p)+ \frac{1}{ml}$
sits in $(-\frac{1}{ml}, \frac{2}{ml})$. Moreover, it is also one of
the numbers $\frac{j}{ml}$, $j=0,1,\cdots,ml-1$. Hence
$\mathfrak{S}([\frac{1}{ml}]\cdot p)-\mathfrak{S}(p)+
\frac{1}{ml}=\frac{1}{ml}$ or $0$. If $p\ne x^m$ the second case
cannot occur because $\Theta$ is a diffeomorphism. This shows that
$\Theta([\frac{1}{ml}]\cdot p) =[\frac{1}{ml}]\cdot\Theta(p)$ for
any $p\in {\cal N}(x^m)^*$. Namely, $\Theta$ is an equivariant
diffeomorphism from ${\cal N}(x^m)^*$ to some $S^1_{x^m}$-invariant
open neighborhood of $x^m$ in $\hat{\bf S}$. It follows that there
exist $S^1_{x^m}$-invariant open neighborhoods of $x^m$ in ${\bf S}$
and $\hat{\bf S}$, $U$ and $\hat{U}$, such that
\begin{eqnarray*}
&&H_\ast\left(({\cal L}|_{{\bf S}})_d\cap U/S^1_{x^m};({\cal
L}|_{\hat{\bf S}})_d\cap
(U\setminus\{x^m\})/S^1_{x^m};\K\right)\nonumber\\
&=&H_\ast\left(({\cal L}|_{\hat{\bf S}})_d\cap
\hat{U}/S^1_{x^m};({\cal L}|_{\hat{\bf S}})_d\cap
(\hat{U}\setminus\{x^m\})/S^1_{x^m};\K\right).
\end{eqnarray*}
This implies that  there exist $S^1_{x^m}$-invariant open
neighborhoods  $V$ and $\hat{V}$ of $0$ in $XN\tilde\varphi_m({\cal
O})(\delta)_{x^m}$ and $X\hat N\tilde\varphi_m({\cal
O})(\delta)_{x^m}$ respectively, such that
\begin{eqnarray}\label{e:6.32}
&&H_\ast\left(({\cal F}^X_{x^m})_d\cap V/S^1_{x^m};({\cal
F}^X_{x^m})_d\cap
(V\setminus\{x^m\})/S^1_{x^m};\K\right)\nonumber\\
&=&H_\ast\left((\hat{\cal F}^X_{x^m})_d\cap
\hat{V}/S^1_{x^m};(\hat{\cal F}^X_{x^m})_d\cap
(\hat{V}\setminus\{x^m\})/S^1_{x^m};\K\right).
\end{eqnarray}
 For the field $\K$ of characteristic
zero or prime up to order of $S^1_{\gamma_0^m}$, as in the proof of the
equality above (\ref{e:5.14}) using Theorem~\ref{th:A.5} (for ${\cal
F}^X_{x^m}$ and $\hat{\cal F}^X_{x^m}$) we can deduce that
\begin{eqnarray}
&&H_q\left(({\cal F}^X_{x^m})_d\cap V/S^1_{x^m};({\cal
F}^X_{x^m})_d\cap
(V\setminus\{x^m\})/S^1_{x^m};\K\right)\label{e:6.33}\\
&=& \Bigl(H_{m^-}({\bf H}^-({B})_{x^m}, {\bf
H}^-({B})_{x^m}\setminus\{0\};\K)\otimes C_{q-m^-}({\cal
L}^\circ_{\triangle x^m},
0;\K)\Bigr)^{S^1_{x^m}},\nonumber\\
&&H_\ast\left((\hat{\cal F}^X_{x^m})_d\cap
\hat{V}/S^1_{x^m};(\hat{\cal F}^X_{x^m})_d\cap
(\hat{V}\setminus\{x^m\})/S^1_{x^m};\K\right)\label{e:6.34}\\
&=&\Bigl(H_{m^-}({\bf H}^-(\hat{B})_{x^m}, {\bf
H}^-(\hat{B})_{x^m}\setminus\{0\};\K)\otimes C_{q-m^-}(\hat{\cal
L}^\circ_{\triangle x^m}, 0;\K)\Bigr)^{S^1_{x^m}}\nonumber
\end{eqnarray}
for each $q=0,1,\cdots$.  Consider the following homotopy of
equivalent Riemannian-Hilbert structures on $T{\Lambda M}$ from
$\langle\cdot,\cdot\rangle_1$ to $(\cdot,\cdot)_m$,
\begin{equation}\label{e:6.35}
[0,1]\ni\tau\mapsto\langle\!\langle\cdot,\cdot\rangle\!\rangle_{\tau}
=(1-\tau)\langle\cdot,\cdot\rangle_1+\tau(\cdot,\cdot)_m,
\end{equation}
Let $\mathfrak{B}^{(\tau)}\in{\cal L}(T_{x^m}{\Lambda M})$ be
determined by $d^2({\cal L}\circ{\rm
EXP}_{x^m})(0)[\xi,\eta]=\langle\!\langle\mathfrak{B}^{(\tau)}\xi,\eta\rangle\!\rangle_{\tau}$.
It is Fredholm and self-adjoint with respect to the inner product
$\langle\!\langle\cdot,\cdot\rangle\!\rangle_{\tau}$. Observe that
$\dim{\bf H}^-(\mathfrak{B}^{(\tau)})$ does not depend on $\tau\in
[0, 1]$ and that ${\bf H}^-(\mathfrak{B}^{(0)})={\bf
H}^-({B})_{x^m}$ and  ${\bf H}^-(\mathfrak{B}^{(1)})={\bf
H}^-(\hat{B})_{x^m}$. Moreover, each ${\bf
H}^-(\mathfrak{B}^{(\tau)})$ is $S^1_{x^m}$-invariant. It is easy to
find continuous paths $e_i:[0, 1]\to T_{x^m}{\Lambda M}$,
$i=1,\cdots,m^-$, such that $e_1(\tau),\cdots, e_{m^-}(\tau)$ form a
basis in ${\bf H}^-(\mathfrak{B}^{(\tau)})$ for each $\tau\in
[0,1]$. Let $[s]\in S^1_{x^m}$ be a generator. Then the determinant
of the transformation matrix from $e_1(\tau),\cdots, e_{m^-}(\tau)$
to $[s]\cdot e_1(\tau),\cdots, [s]\cdot e_{m^-}(\tau)$ is a nonzero
continuous function of $\tau\in [0, 1]$.  Hence  a generator of the
$S^1_{x^m}$-action on ${\bf H}^-(B)_{x^m}$ reverses orientation if
and only if its action on ${\bf H}^-(\hat{B})_{x^m}$ reverses
orientation. This and (\ref{e:6.32})-(\ref{e:6.34}) imply
(\ref{e:6.29}). \hfill$\Box$\vspace{2mm}

\begin{remark}\label{rm:6.8}
{\rm There is another possible method to prove (\ref{e:6.28}) and
(\ref{e:6.29}). Let $N\varphi_m({\cal O})^\tau$ denote the normal
bundle of $\varphi_m({\cal O})$ with respect to the metric
$\langle\!\langle\cdot,\cdot\rangle\!\rangle_{\tau}$ in
(\ref{e:6.35}). Note that there exists a natural bundle isomorphism
$\mathcal {I}^\tau: N\varphi_m({\cal O})\to N\varphi_m({\cal
O})^\tau$ whose restriction $\mathcal {I}^\tau_{x^m}$ on fiber
$N\varphi_m({\cal O})_{x^m}$ at $x^m$ is given by
$$
 \xi-\frac{\langle\xi,
\dot{x}^m\rangle_1}{\langle\dot{x}^m,
\dot{x}^m\rangle_1}\dot{x}^m\to\xi-\frac{\langle\!\langle\xi,
\dot{x}^m\rangle\!\rangle_\tau}{\langle\!\langle\dot{x}^m,
\dot{x}^m\rangle\!\rangle_\tau}\dot{x}^m
$$
for any $\xi\in T_{x^m}{\Lambda M}$. Clearly, $\mathcal
{I}^\tau_{x^m}$ is only a Banach space isomorphism from
${N}\varphi_m({\cal O})_{x^m}$ onto $N\varphi_m({\cal
O})^\tau_{x^m}$. For a small $\epsilon>0$ we have a smooth homotopy
from ${\cal F}_{x^m}$ to $\hat{\cal F}_{x^m}\circ(\mathcal
{I}^1_{x^m}|_{{N}\varphi_m({\cal O})(\epsilon)_{x^m}})$, ${\cal
L}\circ{\rm EXP}_{x^m}\circ(\mathcal
{I}^\tau_{x^m}|_{{N}\varphi_m({\cal O})(\epsilon)_{x^m}})$ with
$0\le\tau\le 1$. Theorem~\ref{th:3.3} implies $C_\ast({\cal
F}_{x^m},0;\K)=C_\ast(\hat{\cal F}_{x^m}\circ(\mathcal
{I}^1_{x^m}|_{{N}\varphi_m({\cal O})(\epsilon)_{x^m}}),0;\K)=
C_\ast(\hat{\cal F}_{x^m},0;\K)$. This and a corresponding
conclusion with (\ref{e:1.8}) yield (\ref{e:6.30}). A similar
argument with the stability of critical groups  for continuous
functionals on metric spaces can give (\ref{e:6.29}) as in
section~\ref{sec:7}. }
\end{remark}

\section{Computations of $S^1$-critical groups}\label{sec:7}
\setcounter{equation}{0}

Let $\Gamma$ be a subgroup of $S^1$ and let $A\subset\Lambda M$ be a
$\Gamma$-invariant subset. We denote by $A/\Gamma$ the quotient
space of $A$ with respect to the action of $\Gamma$.
  Rademacher \cite[\S6.1]{Ra} defined
$S^1$-{\it critical group} of $\gamma_0$ by
$$
\overline{C}_\ast(\mathcal{L},
\gamma_0;\K)=H_\ast\left((\Lambda(\gamma_0)\cup
S^1\cdot\gamma_0)/S^1, \Lambda(\gamma_0)/S^1;\K\right),
$$
which is important in studies of closed geodesics. It was proved in
\cite[(3.10)]{BLo} that the  groups $\overline{C}_\ast(\mathcal{L},
\gamma_0;\K)$ and $C_\ast(\mathcal{L},
\gamma_0;\K):=H_\ast\left(\Lambda(\gamma_0)\cup \{\gamma_0\},
\Lambda(\gamma_0);\K\right)$ have relations
\begin{equation}\label{e:7.1}
\overline{C}_\ast(\mathcal{L}, \gamma_0;\K)=C_\ast(\mathcal{L},
\gamma_0;\K)^{S^1_{\gamma_0}}.
\end{equation}
Let us outline a proof of it with our previous methods. For
simplicity we write $\gamma_0$ as $\gamma$. Define a family of inner
products on $T_{\gamma}\Lambda M$ by
$$
\langle\xi,\eta\rangle_{\tau}=\int^1_0\langle\xi(t),\eta(t)\rangle
dt+ \tau
\int^1_0\langle\nabla^g_{\dot{\gamma}}\xi(t),\nabla^g_{\dot{\gamma}}\eta(t)\rangle
dt,\quad \tau\in [0,1].
$$
For a small $\epsilon>0$ let $\Omega_{2\epsilon}=\{\xi\in{\bf
B}_{2\epsilon}(T_{\gamma}\Lambda
M)\,|\,\langle\dot\gamma,\xi\rangle_0=0\}$. By Lemma~2.2.8 of
\cite{Kl} the set ${\Gamma_0}_{2\epsilon}:={\rm
EXP}(\Omega_{2\epsilon})=\{\exp_{\gamma}(\xi)\,|\,\xi\in
\Omega_{2\epsilon}\}$ is a slice of the $S^1$-action on $\Lambda M$.
This means:(i)
$s\cdot{\Gamma_0}_{2\epsilon}={\Gamma_0}_{2\epsilon}\;\forall s\in
S^1_{\gamma}$, (ii) if
$(s\cdot{\Gamma_0}_{2\epsilon})\cap{\Gamma_0}_{2\epsilon}\ne\emptyset$
for $s\in S^1$ then $s\in S^1_{\gamma}$, (iii) $\exists\,0<\nu\ll 1$
such that
\begin{equation}\label{e:7.2}
\mathfrak{H}:(-\nu,\nu)\times{\Gamma_0}_{2\epsilon}\ni (s,x)\mapsto
s\cdot x\in\Lambda M
\end{equation}
 is a homeomorphism onto an open neighborhood of $\gamma$ in
$\Lambda M$. The latter implies that
${\Gamma_0}_{2\epsilon}/S^1_{\gamma}\ni S^1_{\gamma}\cdot x\to
S^1\cdot x\in \Lambda M/S^1$ is a homeomorphism onto
$(S^1\cdot{\Gamma_0}_{2\epsilon})/S^1$ which is an open neighborhood
of the equivalence class of $\gamma$ in $\Lambda M/S^1$. Hence
\begin{eqnarray}\label{e:7.3}
\overline{C}_\ast(\mathcal{L},
\gamma;\K)&=&H_\ast\bigl((\Lambda(\gamma)\cap{\Gamma_0}_{2\epsilon}\cup\{\gamma\})/S^1_{\gamma},
\Lambda(\gamma)\cap{\Gamma_0}_{2\epsilon}/S^1_{\gamma};\K\bigr)\nonumber\\
&=&H_\ast\bigl(\mathcal{L}_c\cap\overline{{\Gamma_0}_{\epsilon}}/S^1_{\gamma},
\mathcal{L}_c\cap\overline{{\Gamma_0}_{\epsilon}}\setminus\{\gamma\}/S^1_{\gamma}
;\K\bigr),
\end{eqnarray}
where $\overline{{\Gamma_0}_{\epsilon}}={\rm
EXP}(\overline{\Omega_{\epsilon}})$ is the closure of
${\Gamma_0}_{\epsilon}$. Since
$\overline{{\Gamma_0}_{\epsilon}}/S^1_{\gamma}$ is a complete metric
space, using the stability of critical groups of isolated critical
point for continuous functionals on metric spaces in
\cite[Th.5.2]{CorHa} or \cite[Th.1.5]{CiDe} we may prove that
\begin{eqnarray}\label{e:7.4}
&&H_\ast\bigl(\mathcal{L}_c\cap\overline{{\Gamma_0}_{\epsilon}}/S^1_{\gamma},
\mathcal{L}_c\cap\overline{{\Gamma_0}_{\epsilon}}\setminus\{\gamma\}/S^1_{\gamma}
;\K\bigr)\nonumber\\
&=&H_\ast\bigl((\mathcal{L}^\ast)_c\cap\overline{{\Gamma_0}_{\epsilon}}/S^1_{\gamma},
(\mathcal{L}^\ast)_c\cap\overline{{\Gamma_0}_{\epsilon}}\setminus\{\gamma\}/S^1_{\gamma}
;\K\bigr),
\end{eqnarray}
see \cite{LuW} for details.  Let $2\epsilon$ be less than
$\varepsilon$ in (\ref{e:1.12}). Using the chart in (\ref{e:1.12})
we write $\mathcal{E}^\ast:=\mathcal{L}^\ast\circ{\rm EXP}$. Then
$\mathcal{E}^\ast|_{\overline{\Omega_{\epsilon}}}=\mathcal{L}^\ast|_{\overline{{\Gamma_0}_{\epsilon}}
}\circ{\rm EXP}|_{\overline{\Omega_{\epsilon}}}$ and
\begin{eqnarray}\label{e:7.5}
&&H_\ast\bigl((\mathcal{L}^\ast)_c\cap\overline{{\Gamma_0}_{\epsilon}}/S^1_{\gamma},
(\mathcal{L}^\ast)_c\cap\overline{{\Gamma_0}_{\epsilon}}\setminus\{\gamma\}/S^1_{\gamma}
;\K\bigr)\nonumber\\
&=&H_\ast\bigl((\mathcal{E}^\ast|_{\overline{\Omega_{\epsilon}}})_c/S^1_{\gamma},
(\mathcal{E}^\ast|_{\overline{\Omega_{\epsilon}}})_c\setminus\{\gamma\}/S^1_{\gamma}
;\K\bigr).
\end{eqnarray}
Note that $\gamma$ is at least $C^2$ in our case. Denote by
$\langle\dot\gamma\rangle^\bot_\tau$ the orthogonal complementary of
$\dot\gamma\R$ with respect to the inner product
$\langle\cdot,\cdot\rangle_{\tau}$. Then
$\langle\dot\gamma\rangle^\bot_1=N\mathcal{O}_{\gamma}$ and
$$
\langle\dot\gamma\rangle^\bot_0= \{\xi\in T_{\gamma}\Lambda
M\,|\,\langle\dot\gamma,\xi\rangle_0=0\}\quad\hbox{and hence}\quad
\Omega_{2\epsilon}={\bf B}_{2\epsilon}(T_{\gamma}\Lambda
M)\cap\langle\dot\gamma\rangle^\bot_0.
$$
For each $\tau\in [0, 1]$ there exists an obvious
$S^1_{\gamma}$-equivariant linear diffeomorphism
$$
h_\tau:\langle\dot\gamma\rangle^\bot_0\to\langle\dot\gamma\rangle^\bot_\tau,\;
\xi\mapsto\xi-\langle\xi,\dot\gamma\rangle_\tau\dot\gamma.
$$
And $|h_\tau(\xi)|_1^2\le
(1+|\dot\gamma|_1^2)^2|\xi|_1^2\;\forall\xi\in
\langle\dot\gamma\rangle^\bot_0$. So we can take $0<\rho\ll\epsilon$
such that $h_\tau(\xi)\in{\bf B}_{2\epsilon}(T_{\gamma}\Lambda M)$
for all $\xi\in {\bf B}_{2\rho}(T_{\gamma}\Lambda M)$. Let
$\Omega^\tau_{2\epsilon}={\bf B}_{2\epsilon}(T_{\gamma}\Lambda
M)\cap\langle\dot\gamma\rangle^\bot_\tau$ and consider the pullback
$h^\ast_\tau(\mathcal{E}^\ast|_{\Omega^\tau_{2\epsilon}})$ for
$\tau\in [0,1]$. They are well-defined on
$\Omega_{2\rho}=\Omega^0_{2\rho}$ and have an isolated critical
point $\gamma$. As above using the stability of critical groups for
continuous functionals on metric spaces in \cite[Th.5.2]{CorHa} or
\cite[Th.1.5]{CiDe} we derive
\begin{eqnarray}\label{e:7.6}
&&H_\ast\bigl((\mathcal{E}^\ast|_{\overline{\Omega_{\epsilon}}})_c/S^1_{\gamma},
(\mathcal{E}^\ast|_{\overline{\Omega_{\epsilon}}})_c\setminus\{\gamma\}/S^1_{\gamma}
;\K\bigr)\nonumber\\
&=&H_\ast\bigl((h_0^\ast\mathcal{E}^\ast|_{\overline{\Omega_{\rho}}})_c/S^1_{\gamma},
(h_0^\ast\mathcal{E}^\ast|_{\overline{\Omega_{\rho}}})_c\setminus\{\gamma\}/S^1_{\gamma}
;\K\bigr)\nonumber\\
&=&H_\ast\bigl((h_1^\ast\mathcal{E}^\ast|_{\overline{\Omega_{\rho}}})_c/S^1_{\gamma},
(h_1^\ast\mathcal{E}^\ast|_{\overline{\Omega_{\rho}}})_c\setminus\{\gamma\}/S^1_{\gamma}
;\K\bigr)\nonumber\\
&=&H_\ast\bigl((\mathcal{E}^\ast|_{N\mathcal{O}(\epsilon)_{\gamma}})_c/S^1_{\gamma},
(\mathcal{E}^\ast|_{N\mathcal{O}(\epsilon)_{\gamma}})_c\setminus\{\gamma\}/S^1_{\gamma}
;\K\bigr)\nonumber\\
&=&H_\ast\bigl(\mathcal{F}^\ast_c/S^1_{\gamma},
\mathcal{F}^\ast_c\setminus\{\gamma\}/S^1_{\gamma} ;\K\bigr),
\end{eqnarray}
where the excision property of relative homology groups are used in
the first, second and fourth equality and we may require that
$\epsilon$ is less than $\varepsilon$ in (\ref{e:3.17}). By
(\ref{e:7.3})-(\ref{e:7.6}) we get
\begin{eqnarray*}
\overline{C}_\ast(\mathcal{L},
\gamma;\K)&=&H_\ast\bigl((\Lambda(\gamma)\cap{\Gamma_0}_{2\epsilon}\cup\{\gamma\})/S^1_{\gamma},
\Lambda(\gamma)\cap{\Gamma_0}_{2\epsilon}/S^1_{\gamma};\K\bigr)\nonumber\\
&=&H_\ast\bigl(\mathcal{F}^\ast_c/S^1_{\gamma},
\mathcal{F}^\ast_c\setminus\{\gamma\}/S^1_{\gamma} ;\K\bigr).
\end{eqnarray*}
As in the proofs (\ref{e:5.13}) and (\ref{e:5.14})  using this and
Theorem~\ref{th:1.5}(iii) we can derive

\begin{proposition}\label{prop:7.1}
If $\K$ is a field of characteristic $0$ or prime up to order
$|S^1_{\gamma}|$ of $S^1_{\gamma}$,
 \begin{eqnarray*}
&&\overline{C}_q(\mathcal{L},
\gamma;\K)\\
&=& \Bigl(H_{m^-({\cal O})}({\bf H}^-({B})_{\gamma}, {\bf
H}^-({B})_{\gamma}\setminus\{0\};\K)\otimes C_{q-m^-({\cal
O})}({\cal L}^\circ_{\triangle\gamma}, 0;\K)\Bigr)^{S^1_{\gamma}}.
\end{eqnarray*}
\end{proposition}

Let $\mathfrak{H}$ denote the homeomorphism in (\ref{e:7.2}). Then
\begin{eqnarray}\label{e:7.7}
 C_q(\mathcal{L}, \gamma;\K)&=&H_q\bigl(\Lambda(\gamma)\cup \{\gamma\},
\Lambda(\gamma);\K\bigr)\nonumber\\
&=&H_q\bigl({\rm Im}(\mathfrak{H})\cap\Lambda(\gamma)\cup
\{\gamma\},
{\rm Im}(\mathfrak{H})\cap\Lambda(\gamma);\K\bigr)\nonumber\\
&=&H_q\bigl((-\nu,\nu)\times({\Gamma_0}_{2\epsilon}\cap\Lambda(\gamma)\cup
\{\gamma\}),
(-\nu,\nu)\times({\Gamma_0}_{2\epsilon}\cap\Lambda(\gamma));\K\bigr)\nonumber\\
&=&H_q\bigl({\Gamma_0}_{2\epsilon}\cap\Lambda(\gamma)\cup
\{\gamma\},
{\Gamma_0}_{2\epsilon}\cap\Lambda(\gamma);\K\bigr)\nonumber\\
&=&H_q\bigl({\Gamma_0}_{2\epsilon}\cap\mathcal{L}_c,{\Gamma_0}_{2\epsilon}\cap\mathcal{L}_c\setminus
\{\gamma\};\K\bigr).
\end{eqnarray}
As before using Theorem~\ref{th:3.3}, the excision and
Corollary~\ref{cor:A.6} we can prove
\begin{eqnarray*}
(\ref{e:7.7})&=&H_q\bigl(\mathcal{L}_c\cap\overline{{\Gamma_0}_{\epsilon}},
\mathcal{L}_c\cap\overline{{\Gamma_0}_{\epsilon}}\setminus\{\gamma\}
;\K\bigr)\nonumber\\
&=&H_q\bigl(\mathcal{L}^\ast_c\cap\overline{{\Gamma_0}_{\epsilon}},
\mathcal{L}^\ast_c\cap\overline{{\Gamma_0}_{\epsilon}}\setminus\{\gamma\}
;\K\bigr)\nonumber\\
&=&H_q\bigl((\mathcal{E}^\ast|_{\overline{\Omega_{\epsilon}}})_c,
(\mathcal{E}^\ast|_{\overline{\Omega_{\epsilon}}})_c\setminus\{\gamma\}
;\K\bigr)\nonumber\\
&=&H_q\bigl((\mathcal{E}^\ast|_{N\mathcal{O}(\epsilon)_{\gamma}})_c,
(\mathcal{E}^\ast|_{N\mathcal{O}(\epsilon)_{\gamma}})_c\setminus\{\gamma\};\K\bigr)\nonumber\\
&=&H_q\bigl(\mathcal{F}^\ast_c,\mathcal{F}^\ast_c\setminus\{\gamma\}
;\K\bigr)\nonumber\\
&=&H_{m^-({\cal O})}({\bf H}^-({B})_{\gamma}, {\bf
H}^-({B})_{\gamma}\setminus\{0\};\K)\otimes C_{q-m^-({\cal
O})}({\cal L}^\circ_{\triangle\gamma}, 0;\K).
\end{eqnarray*}
Hence we arrive at
\begin{eqnarray*}
&& C_q(\mathcal{L}, \gamma;\K)\nonumber\\
&=&H_{m^-({\cal O})}({\bf H}^-({B})_{\gamma}, {\bf
H}^-({B})_{\gamma}\setminus\{0\};\K)\otimes C_{q-m^-({\cal
O})}({\cal L}^\circ_{\triangle\gamma}, 0;\K).
\end{eqnarray*}
Then (\ref{e:7.1}) follows from this and Proposition~\ref{prop:7.1}.

\section{Proof of Theorem~\ref{th:1.11}}\label{sec:8}
\setcounter{equation}{0}

Though as done in \cite{Lu4} the proof can be completed following
the arguments in Refs. \cite{BKl, Lo1, LoLu, Lu0,Lu1} we
shall present a simplified proof in the sense that  complicated
quantitative  estimates about Bangert homotopies in the original
methods of \cite{BKl,Lo1} can be removed with more ``soft''
topological arguments. The same methods can be used to simplify the
arguments in \cite{Lu0,Lu1}.

 Firstly,  the key Lemmas~1,2 in
\cite{GM2} also hold true for closed geodesics on $(M, F)$ (cf.
\cite[\S4.2, \S 7.1]{Ra}), that is,

\begin{lemma}\label{lem:8.1}
For a closed geodesic $\gamma$ on $(M,F)$ one has:
\begin{description}
\item[(a)] Either $m^-(\gamma^k)=0$ for all $k$ or there exist
numbers $a>0$, $b>0$ such that $m^-(\gamma^{k+l})-m^-(\gamma^k)\ge
la-b$ for all $k,l$;
 \item[(b)] There are positive integers $k_1,\cdots, k_s$ and sequence $n_{ji}\in\N$,
$i>0$, $j=1,\cdots, s$, such that the numbers $n_{ji}k_j$ are
mutually distinct, $n_{j1}=1$, $\{n_{ji}k_j\}=\N$, and
$m^0(\gamma^{n_{ji}k_j})=m^0(\gamma^{k_j})$.
\end{description}
\end{lemma}

For every integer $q\ge 0$ let $\triangle^q$ be the standard
$q$-simplex, i.e., $\triangle^q=\langle e_0,\cdots,e_q\rangle$,
where $e_0=0$ and $e_1,\cdots,e_q$ is the standard basis of $\R^q$.
Our methods depend on the following generalization of \cite[Lemma
1]{BKl}. (Actually the second part of it is exactly the content of
\cite[Lemma 1]{BKl}.)

\begin{lemma}\label{lem:8.2}
Let $(X, A)$ be a pair of topological spaces and
$\alpha=\sum^\ell_{i=0}n_i{\gamma_0}_i$ a singular relative integral
$p$-cycle of $(X,A)$. Let ${\Gamma_0}(\alpha)$ denote the set of
singular simplices of $\alpha$ together with all their faces.
Suppose to every ${\gamma_0}\in{\Gamma_0}$,
${\gamma_0}:\triangle^q\to X$, $0\le q\le p$, there is assigned a
homotopy $P({\gamma_0}):\triangle^q\times [0, 1]\to X$ such that
\begin{description}
\item[(i)] $P({\gamma_0})(z,0)={\gamma_0}(z)$ for $z\in\triangle^q$,
\item[(ii)] $P({\gamma_0})(z,t)={\gamma_0}(z)$ if ${\gamma_0}(\triangle^q)\subset
A$ with $q=p-1$,
\item[(iii)] $P({\gamma_0})\circ (e^i_q\times id)=P({\gamma_0}\circ e^i_q)$
for $0\le i\le q$.
\end{description}
Then the homology classes $[\alpha]$ and $[\bar\alpha]$ are same in
$H_p(X,A;\Z)$, where $\bar\alpha=\sum^\ell_{i=0}n_i\bar{\gamma_0}_i$
with $\bar{\gamma_0}_i=P({\gamma_0}_i)_1=P({\gamma_0}_i)(\cdot,1)$.
In particular, $[\alpha]=0$ provided that the above assigned
homotopy every $P({\gamma_0})$ is also required to satisfy
\begin{description}
\item[(ii')] $P({\gamma_0})(z,t)={\gamma_0}(z)$ if ${\gamma_0}(\triangle^q)\subset
A$ (including the above {\bf (ii)}),
\item[(iv)] $P({\gamma_0})(\triangle^q\times\{1\})\subset A$.
\end{description}
\end{lemma}

\noindent{\bf Proof}.\quad Note that the cycles $\alpha$ and
$\bar\alpha$ can be considered (continuous) maps from a
$p$-dimensional finite oriented simplicial complex $K$ with boundary
to $X$ such that $\alpha(\partial K)\subset A$, $\bar\alpha(\partial
K)\subset A$ and that $\alpha_\ast$ and $\bar\alpha_\ast$ map the
fundamental class $[K,\partial K]\in H_{p}(K,
\partial K;\Z)$ to $[\alpha]\in H_p(X,A;\Z)$ and $[\bar\alpha]\in H_p(X,A;\Z)$
respectively. Moreover, all homotopies $P({\gamma_0})$,
${\gamma_0}\in{\Gamma_0}(\alpha)$, give a homotopy between
$\alpha:(K,\partial K)\to (X,A)$ and $\bar\alpha:(K,\partial K)\to
(X,A)$ relative to $\partial K$. Hence $\alpha_\ast[K,\partial
K]=\bar\alpha_\ast[K,\partial K]$ in $H_p(X,A;\Z)$.
\hfill$\Box$\vspace{2mm}

 For $k\in\N\cup\{0\}$ recall that a pair $(X,A)$ consisting of a
topological space $X$ and a subspace $A$  is said to {\bf
$k$-connected} if $\pi_0(A)\to\pi_0(X)$ is surjective and
$\pi_j(X,A,a)=0$ for $j\in\{1,\cdots,k\}$ and each $a\in A$, in
particular the $0$-connectedness of the pair $(X,A)$ is equivalent
to the condition that every point of $X$ is joined by a path to some
point of $A$ (cf. for example \cite[p. 143]{Di}). Let
$Cl(\triangle^q)$ be the closure complex of $\triangle^q$, i.e., the
standard simplicial complex consisting of all faces of
$\triangle^q$, and  let $Cl(\triangle^q)_k$ be the $p$-skeleton of
the complex $Cl(\triangle^q)$ (and so
$Cl(\triangle^q)_k=Cl(\triangle^q)$ for $k\ge q$). We need the
following standard result in Eilenberg and Blakers homology groups.
The reader is referred to one of Refs. \cite[Th.9.5.1]{Di},
\cite[Chap.7, \S4, Th.8]{Spa} and \cite[Chap.4, Th.5.1]{Whi} for the
proof.

\begin{lemma}\label{lem:8.3}
Suppose that a pair $(X, A)$ of topological spaces is $k$-connected.
Then   every simplex  ${\gamma_0}:\triangle^q\to X$, ($0\le q\le k$
or $q>k$), can be assigned a homotopy
$P({\gamma_0}):\triangle^q\times [0, 1]\to X$ such that
\begin{description}
\item[(i)] $P({\gamma_0})_0={\gamma_0}$, i.e., $P({\gamma_0})(z,0)={\gamma_0}(z)\;\forall z\in\Delta^q$,
\item[(ii)] ${\gamma_0}(\triangle^q)\subset
A\Longrightarrow P({\gamma_0})_t={\gamma_0}\;\forall t\in [0, 1]$,
i.e., $P({\gamma_0})(z,t)={\gamma_0}(z)\;\forall
(z,t)\in\Delta^q\times [0,1]$,
\item[(iii)] $P({\gamma_0})\bigl(Cl(\triangle^q)_k\times\{1\}\bigr)\subset A$,
\item[(iv)] $P({\gamma_0})\circ (e^i_q\times id)=P({\gamma_0}\circ e^i_q)$
for $0\le i\le q$, where $e^i_q$ is the $i$-th face of
$\triangle^q$.
\end{description}
\end{lemma}

Actually, we may use the proof method of this lemma to prove:

\begin{lemma}\label{lem:8.4}
For a pair $(X, A)$ of topological spaces, a  simplex
${\gamma_0}:\triangle^p\to X$, $p>0$,  and an integer $0\le k<p$,
suppose that each $q$-simplex ${\gamma_0}\in Cl(\triangle^p)_k$ has
been
 assigned a homotopy $P({\gamma_0}):\triangle^q\times [0, 1]\to X$ satisfying
the conditions (i)-(ii), (iv) in Lemma~\ref{lem:8.3} and
$P({\gamma_0})\bigl(\triangle^q\times\{1\}\bigr)\subset A$.  Then
each $r$-simplex ${\gamma_0}\in Cl(\triangle^p)\setminus
Cl(\triangle^p)_k$ may be  assigned a homotopy
$P({\gamma_0}):\triangle^r\times [0, 1]\to X$  such that
\begin{description}
\item[(i)] $P({\gamma_0})_0={\gamma_0}$,
\item[(ii)] ${\gamma_0}(\triangle^r)\subset
A\Longrightarrow P({\gamma_0})_t={\gamma_0}\;\forall t\in [0, 1]$,
\item[(iii)] $P({\gamma_0})\bigl(Cl(\triangle^r)_k\times\{1\}\bigr)\subset A$,
\item[(iv)] $P({\gamma_0})\circ (e^i_r\times id)=P({\gamma_0}\circ e^i_r)$
for $0\le i\le r$, where $e^i_r$ is the $i$-th face of
$\triangle^r$.
\end{description}
\end{lemma}

\noindent{\bf Proof}.\quad  Let ${\gamma_0}\in
Cl(\triangle^p)\setminus Cl(\triangle^p)_k$ be an arbitrary
$(k+1)$-simplex. If ${\gamma_0}(\triangle^{k+1})\subset A$ we define
$P({\gamma_0}):\triangle^{k+1}\times [0,1]\to X$ by
$P({\gamma_0})(z,t)={\gamma_0}(z)$ for all $(z,t)\in
\triangle^{k+1}\times [0,1]$. If
${\gamma_0}(\triangle^{k+1})\not\subset A$, the homotopies
$P({\gamma_0}\circ e^i_{k+1})$, $0\le i\le k+1$, induce a continuous
map $f:\partial\triangle^{k+1}\times [0, 1]\to X$ such that
$f\circ(e^i_{k+1}\times id)=P({\gamma_0}\circ e^i_{k+1})$. Since
$\partial\triangle^{k+1}\subset\triangle^{k+1}$ is a cofibration,
there exists a retraction $r:\triangle^{k+1}\times [0,1]\to
(\triangle^{k+1}\times \{0\})\cup (\partial\triangle^{k+1}\times
[0,1])$. By the homotopy extension property we obtain a continuous
map $P({\gamma_0}):\triangle^{k+1}\times [0, 1]\to X$ such that
$P({\gamma_0})(\cdot,0)={\gamma_0}(\cdot)$ and
$P({\gamma_0})(z,t)=f(z,t)$ for all $(z,t)\in
\partial\triangle^{k+1}\times
[0,1]$. Clearly, $P({\gamma_0})$ satisfies (i)-(ii) and (iv). For
$\mu\in Cl(\triangle^{k+1})_k$, since
$\mu\subset\partial\triangle^{k+1}$ we have $\mu\subset {\rm
Im}(e^i_{r+1})$ for some $0\le i\le k+1$, and hence
$$
P({\gamma_0})(\mu\times\{1\})=f(\mu\times\{1\})=P({\gamma_0}\circ
e^i_{k+1})(\mu\times\{1\})\subset A
$$
because (iii) holds for $r=k$. If $k+1<p$ a simple induction
argument leads to the desired proof. \hfill$\Box$\vspace{2mm}

Replacing \cite[Th.1]{BKl} we need the following analogue of it,
whose corresponding version in the content of the Lagrangian Conley
conjecture was firstly proved by Y. Long \cite[Prop. 5.1]{Lo1} for
the finite energy homology on $M=T^n$, and then by the author
\cite[Prop. 5.6]{Lu0} for the  singular homology on closed
manifolds.

\begin{proposition}\label{prop:8.5}
 Let $\Lambda=\Lambda M$ and $\mathring{\cal L}_{d}=\{\gamma\in\Lambda\,|\,{\cal
L}(\gamma)<d\}$ for $d>0$. Then for a  $C^1$-smooth $q$-simplex
$\eta:(\triangle^q,\partial\triangle^q)\to (\Lambda M,\mathring{\cal
L}_{d})$, there exists an integer $m(\eta)>0$ such that for every
integer $m\geq m(\eta)$, the $q$-simplex
$$
\eta^m\equiv \varphi_m(\eta):(\triangle^q,\partial\triangle^q)
            \to (\Lambda M,\mathring{\cal L}_{m^2d})
            $$
is  homotopic to a  singular $q$-simplex
$\eta_m:(\triangle^q,\partial\triangle^q)\to
 (\mathring{\cal L}_{m^2d}, \mathring{\cal L}_{m^2d})
$
  with $\eta^m = \eta_m$ on $\partial\triangle^q$
 and the homotopy fixes $\eta^m|_{\partial\triangle^q}$.
\end{proposition}

This can be proved by almost repeating the proof of
\cite[Prop.5.6]{Lu0}.  We shall outline its proof at the end of this
section in order to show that the $C^1$-smoothness of $\eta$ bring
conveniences and simplifications.

Recall that $\Lambda(\gamma):=\{\beta\in\Lambda\,|\,{\cal
L}(\beta)<{\cal L}(\gamma)\}$ for $\gamma\in\Lambda$. Then
$\Lambda(\gamma)=\mathring{\cal L}_{d}$ if ${\cal L}(\gamma)=d$.

Now we begin with proof of Theorem~\ref{th:1.11}. By indirect
arguments we make\vspace{2mm}

\noindent{\bf Assumption F}: There is only a finite number of
  distinct closed geodesics.

Then each critical orbit of ${\cal L}$ is isolated. Since $0\le
m^0(S^1\cdot\gamma)\le 2n-1$ for each closed geodesic $\gamma$, from
Theorem~\ref{th:1.7} it follows that $H_{q}(\Lambda(\gamma)\cup
S^1\cdot \gamma, \Lambda(\gamma);\Q)=C_q(\mathcal{L},
S^1\cdot\gamma;\Q)\ne 0$ and $m^-(\gamma)=0$ imply $q\in[0, 2n-1]$.
Moreover, by the assumption there  exist a  closed geodesics
$\bar\gamma$ and an integer $\bar{p}\ge 2$ such that
$m^-(\bar{\gamma}^k)=0\;\forall k\in\mathbb{N}$ and
$H_{\bar{p}}(\Lambda(\bar\gamma)\cup S^1\cdot \bar\gamma,
\Lambda(\bar\gamma);\Q)\ne 0$.  Hence we can find an integer
$p\ge\bar{p}$ and a closed geodesic $\gamma_0$ such that
\begin{eqnarray}\label{e:8.1}
\left.\begin{array}{ll} &m^-(\gamma_0^k)=0\;\forall k\in\mathbb{N},\quad
H_p(\Lambda(\gamma_0)\cup
S^1\cdot\gamma_0, \Lambda(\gamma_0);\Q)\ne 0,\\
&H_{q}(\Lambda(\gamma)\cup S^1\cdot \gamma,
\Lambda(\gamma);\Q)=0\quad\forall q>p\\
&\hbox{for each closed geodesic $\gamma$ with $m^-({\gamma}^k)=0\;\forall k\in\mathbb{N}$.}
\end{array}\right\}
\end{eqnarray}
 By Lemma~\ref{lem:8.1}(a)  we can find $A>0$ such that
 every closed geodesics $\gamma$ with $\mathcal{L}(\gamma)>A$
 either satisfies $m^-(\gamma)>{p}+1$ or $m^-(\gamma^k)=0\;\forall k\in\mathbb{N}$.
From this, (\ref{e:8.1}) and Theorem~1.7 it follows that
\begin{equation}\label{e:8.2}
H_{p+1}(\Lambda(\gamma)\cup S^1\cdot \gamma, \Lambda(\gamma);\Q)=0
\end{equation}
for every closed geodesic $\gamma$ with $\mathcal{L}(\gamma)>A$. We
conclude
\begin{equation}\label{e:8.3}
H_{p+1}(\Lambda, \Lambda(\gamma^m_0)\cup S^1\cdot
\gamma^m_0;\,\Q)=0\quad\hbox{if}\;{\cal L}(\gamma^m_0)>A.
\end{equation}
Otherwise, take a nonzero class $\mathfrak{B}\in H_{p+1}(\Lambda,
\Lambda(\gamma_0^m)\cup S^1\cdot \gamma_0^m;\,\Q)$ and a singular
cycle representative $Z$ of it in $(\Lambda, \Lambda(\gamma_0^m)\cup
S^1\cdot \gamma_0^m)$.
 Since the support of $Z$ is compact we can choose  a large regular value
 $b>{\cal L}(\gamma_0^m)=m^2{\cal L}(\gamma_0)$
 so that ${\rm supp}(Z)\subset\mathring{\cal L}_b$. Clearly, $Z$  cannot be
homologous to zero in $(\mathring{\cal L}_b,\,
\Lambda(\gamma_0^m)\cup S^1\cdot \gamma_0^m)$ (otherwise it is
homologous to zero in $(\Lambda, \Lambda(\gamma_0^m)\cup S^1\cdot
\gamma_0^m)$.) Hence $H_{p+1}(\mathring{\cal L}_b,
\Lambda(\gamma_0^m)\cup S^1\cdot \gamma_0^m;\,\Q)\ne 0$. As showed
in the proof of \cite[Theorem 3]{BKl}, under Assumption F the
standard Morse theoretic arguments yield  a closed geodesic
$\gamma'$ such that $b>{\cal L}(\gamma')\ge {\cal L}(\gamma_0^m)>A$
and
$$
H_{p+1}(\Lambda(\gamma')\cup S^1\cdot\gamma',\,
\Lambda(\gamma');\Q)\ne 0.
$$
This contradicts to (\ref{e:8.2}). (\ref{e:8.3}) is proved.

Applying Lemma~\ref{lem:8.1}(b) to $\gamma_0$ we get integers
$k_i\ge 2$, $i=1,\cdots, s$, such that
\begin{equation}\label{e:8.4}
m^0(\gamma_0)=m^0(\gamma^k_0) \quad\forall k\in
\N\setminus\cup^s_{i=1}k_i\N.
\end{equation}

For $m\in\N$ let $\iota_m: (\Lambda(\gamma_0^m)\cup S^1\cdot
\gamma_0^m,\;\Lambda(\gamma_0^m))\to
(\Lambda,\;\Lambda(\gamma_0^m))$ denote the inclusion map. As stated
at the end of the proof of \cite[Theorem 3]{BKl}, using
(\ref{e:8.3}) and the exact sequence for the triple $(\Lambda,\,
\Lambda(\gamma_0^m)\cup S^1\cdot \gamma_0^m,\,
\Lambda(\gamma_0^m))$,
\begin{eqnarray*}
{\scriptstyle \to H_{p+1}(\Lambda(\gamma_0^m)\cup S^1\cdot
\gamma_0^m,\, \Lambda(\gamma_0^m);\Q)\to H_{p+1}(\Lambda,
\Lambda(\gamma_0^m);\Q)\to
H_{p+1}(\Lambda, \Lambda(\gamma_0^m)\cup S^1\cdot \gamma_0^m;\Q)}\\
{\scriptstyle\to H_{p}(\Lambda(\gamma_0^m)\cup S^1\cdot
\gamma_0^m,\,
\Lambda(\gamma_0^m);\Q)\stackrel{\iota_{m\ast}}{\longrightarrow}
H_{p}(\Lambda, \Lambda(\gamma_0^m);\Q)\to H_{p}(\Lambda,
\Lambda(\gamma_0^m)\cup S^1\cdot \gamma_0^m;\Q)\to}
\end{eqnarray*}
we obtain that for each integer $m$ with ${\cal L}(\gamma_0^m)>A$,
the induced homomorphism
 $$
\iota_{m\ast}: H_p(\Lambda(\gamma_0^m)\cup S^1\cdot
\gamma_0^m,\;\Lambda(\gamma_0^m);\,\Q)\to
H_p(\Lambda,\;\Lambda(\gamma_0^m);\,\Q)
 $$
 is injective. This, (\ref{e:8.4}) and Theorem~\ref{th:1.9} lead to the claim on
 page 385 of \cite{BKl}.

\begin{claim}\label{cl:8.6}
   For any $m\in \N\setminus\cup^s_{i=1}k_i\N$ with ${\cal L}(\gamma_0^m)>A$, the composition
$$
H_p(\Lambda(\gamma_0)\cup S^1\cdot\gamma_0,
\Lambda(\gamma_0);\Q)\xrightarrow{\varphi_{m\ast}}
H_p(\Lambda(\gamma_0^m)\cup S^1\cdot \gamma_0^m,
\Lambda(\gamma_0^m);\Q)\xrightarrow{\iota_{m\ast}} H_p(\Lambda,
\Lambda(\gamma_0^m);\Q)
$$
is injective. Here $\varphi_{m\ast}$ is the  homomorphism induced by
$\varphi_m: (\Lambda(\gamma_0)\cup S^1\cdot\gamma_0,
\Lambda(\gamma_0))\to (\Lambda, \Lambda(\gamma_0^m))$ between their
singular homology groups.
\end{claim}

Hence as in the proof of \cite[Theorem 3]{BKl} we shall be able to
complete the proof of Theorem~\ref{th:1.11} provided that we may
prove that the homomorphism $\iota_{m\ast}\circ\varphi_{m\ast}$
maps some nonzero class $\mathfrak{C}\in H_p(\Lambda(\gamma_0)\cup
S^1\cdot\gamma_0, \Lambda(\gamma_0);\Q)$ to the zero in
$H_p(\Lambda, \Lambda(\gamma_0^m);\Q)$ for some
$m\in\N\setminus\cup^s_{i=1}k_i\N$ with ${\cal L}(\gamma_0^m)>A$.
 However, for such a nonzero class
$\mathfrak{C}\in H_p(\Lambda(\gamma_0)\cup S^1\cdot\gamma_0,
\Lambda(\gamma_0);\Q)$ it is impossible to find a  $C^1$-smooth
relative cycle representative. Fortunately,  we have the following
commutative diagram
$$
\xymatrix{
   (\Lambda(\gamma_0)\cup S^1\cdot\gamma_0,
\Lambda(\gamma_0)) \ar[d]_{\varphi_m} \ar[r]^-{\iota}
                &(\Lambda,
\Lambda(\gamma_0)) \ar[d]^{\varphi_m}  \\
  (\Lambda(\gamma^m_0)\cup S^1\cdot\gamma^m_0,
\Lambda(\gamma^m_0)) \ar[r]^-{\iota_m}
                & (\Lambda,
\Lambda(\gamma_0^m)) }
$$
where $\iota=\iota_1$ is the inclusion. And the class
$\iota_\ast(\mathfrak{C})\in H_p(\Lambda, \Lambda(\gamma_0);\Q)$ has
always a smooth relative cycle representative $\varrho$ in
$(\Lambda, \Lambda(\gamma_0))$ because $(\Lambda,
\Lambda(\gamma_0))$ is a pair of smooth Hilbert manifolds, which can
be realized by either finite-dimensional approximations of $\Lambda$
or the axioms for homology theories of Eilenberg and Steenrod as
done in \cite{Lo1}.

Since $\gamma_0$ is not a local minimum of the functional
 $\mathcal{L}$, $\Lambda(\gamma_0)\ne\emptyset$ and thus
we can choose a \textsf{path connected} neighborhood ${\cal U}$ of
$S^1\cdot\gamma_0$ such that each point of ${\cal U}$ can be
connected to a point of ${\cal U}\cap\Lambda(\gamma_0)$ by a smooth
path in ${\cal U}$, that is, the pair $({\cal U}, {\cal
U}\cap\Lambda(\gamma_0))$ is $0$-connected. (Indeed, let ${\cal
O}=S^1\cdot\gamma_0$ and we may take ${\cal U}={\cal N}({\cal
O},\varepsilon)$ as in (\ref{e:1.13}). Then for $x\in{\cal O}$,
$0_x\in N{\cal O}(\varepsilon)_x$ is not a local minimum of the
functional $\mathcal{F}_x$ below (\ref{e:1.14}). The desired claim
follows from this and the convexity of $N{\cal O}(\varepsilon)_x$
immediately.)

We can also require that the range (or carrier)
 of $\varrho$ is contained in the
given neighborhood ${\cal U}$ of $S^1\cdot\gamma_0$. Denote by
${\Gamma_0}(\varrho)$ the set of all $p$-simplices of $\varrho$
 together with all their faces
contained in $\varrho$,  and by
${\Gamma_0}_{j}(\varrho)=\{\mu\in{\Gamma_0}(\varrho)\,|\,\dim\mu\le
j\}$ for $0\le j\le p$. (Recall $p\ge\bar{p}\ge 2$). Let
$K=\{k_1,\cdots,k_s\}$, and so
$\N\setminus\cup^s_{i=1}k_i\N=\N\setminus\N K$. Let $m_0$ be the
smallest integer such that ${\cal L}(\gamma_0^m)>A$.

 According to the original methods in \cite{BKl, Lo1},
  it is to prove:
 \textsf{There exists an integer
  $m\in \N\setminus K\N$ with $m>m_0$ such that for every
$\mu\in{\Gamma_0}(\varrho)$ with $\mu:\triangle^r\to \Lambda$ and
$0\leq r\leq p$, there exists a homotopy
$P(\varphi_m(\mu)):\triangle^r\times [0,1]\to \Lambda$ such that the
properties (i) to (iv) in \cite[Lemma~1]{BKl} (i.e., the second part
of Lemma~\ref{lem:8.2}) hold for $(X,A)=\bigl(\Lambda,
\Lambda(\gamma_0^m)\bigr)$.} (As in \cite[Prop. 5.2]{Lo1} this was
completed in \cite{Lu4}  using Proposition~\ref{prop:8.5},
\cite[Lemma~1]{BKl} and \cite[Lemma~A.4]{Lu0}). It follows from the
second part of Lemma~\ref{lem:8.2} that
$\varphi_{m\ast}(\iota_\ast(\mathfrak{C}))\in H_p(\Lambda,
\Lambda(\gamma_0^m);\Q)$ vanishes. This contradicts to
Claim~\ref{cl:8.6}, and therefore  the proof of
Theorem~\ref{th:1.11} is completed.

These methods actually depend on more things  than what
Proposition~\ref{prop:8.5} can give. It is very trouble to check
that Proposition~\ref{prop:8.5} may apply to the homotopy extension
constructed in each induction step (because some quantitative
estimates about Bangert homotopies or its variants are needed). We
shall bypass them. The main novelty is to prove

\begin{proposition}\label{prop:8.7}
Let
${\Gamma_0}_{p}(\varrho)=\{{\gamma_0}_1,\cdots,{\gamma_0}_\ell\}$
and $\varrho=\sum^{\ell}_{i=1}a_i{\gamma_0}_i$, where $a_i\in\Q$,
$i=1,\cdots,\ell$. Let $k\in\N$ such that $ka_i\in\Z$,
$i=1,\cdots,\ell$. Then there exist positive integers $m_1<\cdots
<m_{p-1}$ (in the case $p>2$) and cycles ${\varrho}_{j,m_j,k}$,
$j=1,\cdots,p-1$ such that each ${\varrho}_{j,m_j,k}$ is smooth on
each simplex of it, satisfies
\begin{equation}\label{e:8.5}
\cup^{j}_{i=0}{\Gamma_0}_i({\varrho}_{j,m_j,k})\subset
\Lambda(\gamma^{m_j}_0)
\end{equation}
and also represents the class
$k(\varphi_{m_{j}})_\ast\circ\iota_\ast(\mathfrak{C})\in
H_p(\Lambda, \Lambda(\gamma^{m_j}_0);\Q)$. Furthermore there exists
an integer $\mathfrak{M}>0$ such that for every integer
$m>\mathfrak{M}$ and for every
$\mu\in{\Gamma_0}({\varrho}_{p-1,m_{p-1},k})$ with
$\mu:\triangle^r\to \Lambda$ and $0\leq r\leq p$, there exists a
homotopy $P(\varphi_m(\mu)):\triangle^r\times [0,1]\to \Lambda$ such
that the properties (i), (ii') and (iii)-(iv) in Lemma~\ref{lem:8.2}
hold for $(X,A)=\bigl(\Lambda, \Lambda(\gamma_0^m)\bigr)$.
Consequently, $\varphi_{m\ast}(\iota_\ast(\mathfrak{C}))\in
H_p(\Lambda, \Lambda(\gamma_0^m);\Q)$ vanishes for every integer
$m>\mathfrak{M}$.
\end{proposition}

Once this is proved. Take an integer $m>\mathfrak{M}+ m_0$ such that
$m\in\N\setminus K\N$. The expected contradiction is obtained, and
so the proof of Theorem~\ref{th:1.11} is completed. It remains to
prove Propositions~\ref{prop:8.5},~\ref{prop:8.7}.

\noindent{\bf Proof of Proposition~\ref{prop:8.7}.}  {\it Step 1}.
\textsf{Find another cycle representative $\tilde{\varrho}_k$ of the
class $k\iota_\ast(\mathfrak{C})\in H_p(\Lambda,
\Lambda(\gamma_0);\Z)$ in $(\Lambda, \Lambda(\gamma_0))$ such that
${\Gamma_0}_0(\tilde{\varrho}_k)\subset \Lambda(\gamma_0)$.}

 Because the  carrier  of $\varrho$ is contained in the given
neighborhood ${\cal U}$, and the pair $({\cal U}, {\cal
U}\cap\Lambda(\gamma_0))$ is $0$-connected, by Lemma~\ref{lem:8.3}
 we may assign to every simplex ${\gamma_0}\in
{\Gamma_0}_q(k\varrho)$ with $q\in\{0,1,\cdots,p\}$  a homotopy
$P({\gamma_0}):\triangle^q\times [0, 1]\to \mathcal{U}$ such that
$$
\left.\begin{array}{ll}
 &P({\gamma_0})_0={\gamma_0}, \\
&{\gamma_0}(\triangle^q)\subset \Lambda(\gamma_0)\Longrightarrow
P({\gamma_0})_t={\gamma_0}\;\forall t\in [0, 1],\\
&P({\gamma_0})\circ (e^i_q\times id)=P({\gamma_0}\circ
e^i_q)\;\forall 0\le i\le q.
\end{array}\right\}
$$
It follows from Lemma~\ref{lem:8.2} that
$\tilde{\varrho}_k:=\sum^{\ell}_{i=1}ka_iP({\gamma_0}_i)_1$ is also
a cycle representative of the class $k\iota_\ast(\mathfrak{C})\in
H_p(\Lambda, \Lambda(\gamma_0);\Z)$ in $(\Lambda,
\Lambda(\gamma_0))$ with the  carrier in ${\cal U}$, and
\begin{equation}\label{e:8.6}
{\Gamma_0}_0(\tilde{\varrho}_k)\subset
\Lambda(\gamma_0)\cap\mathcal{U}.
\end{equation}

{\it Step 2}. \textsf{There exists an integer
$\mathfrak{M}(\tilde{\varrho}_k)>m_0$ such that for every integer
$m\geq \mathfrak{M}(\tilde{\varrho}_k)$ the class
$k(\varphi_m)_\ast\circ \iota_\ast(\mathfrak{C}) \in H_p(\Lambda,
\Lambda(\gamma^m_0);\Z)$ has a cycle representative
$\tilde{\varrho}_{1,m,k}$ such that
$\cup^1_{i=0}{\Gamma_0}_i(\tilde{\varrho}_{1,m,k})\subset
\Lambda(\gamma^m_0)$.}

 For each $\mu\in {\Gamma_0}_0(\tilde{\varrho}_k)$ we
define $P(\mu):\triangle^0\times[0,1]\to \Lambda$ by $P(\mu)({\bf
0},t)\equiv\mu$.

By Proposition~\ref{prop:8.5} there exists an integer
$\mathfrak{M}(\tilde{\varrho}_k)>m_0$ such that for every integer
$m\geq \mathfrak{M}(\tilde{\varrho}_k)$ and every $1$-simplex
$\eta\in{\Gamma_0}_1(\tilde{\varrho}_k)$, we have a homotopy
$P(\eta^m):\triangle^1\times [0, 1]\to \Lambda$
 from  the $1$-simplex
$\eta^m\equiv \varphi_m(\eta):(\triangle^1,\partial\triangle^1)
            \to \bigl(\Lambda,
\Lambda(\gamma_0^m)\bigr)$
 to a  singular $1$-simplex
$\eta_m:(\triangle^1,\partial\triangle^1)\to
 \left(\Lambda(\gamma_0^m), \Lambda(\gamma_0^m)\right)$
 to satisfy $P(\eta^m)_t|_{\partial\triangle^1}=\eta^m|_{\partial\triangle^1}\;\forall t\in [0,1]$.
By Lemma~\ref{lem:8.4} we assign  to every simplex ${\gamma_0}\in
{\Gamma_0}_q(\tilde\varrho_k)$, $1<q\le p$,   a homotopy
$P({\gamma_0}^m):\triangle^q\times [0, 1]\to \Lambda$ such that
$$
\left.\begin{array}{ll}
 &P({\gamma_0})_0={\gamma_0}, \\
&{\gamma_0}(\triangle^q)\subset \Lambda(\gamma_0^m)\Longrightarrow
P({\gamma_0})_t={\gamma_0}\;\forall t\in [0, 1],\\
&P({\gamma_0})\circ (e^i_q\times id)=P({\gamma_0}\circ
e^i_q)\;\forall 0\le i\le q.
\end{array}\right\}
$$
It follows from (\ref{e:8.6}), these and Lemma~\ref{lem:8.2} that
$$
\tilde{\varrho}_{1,m,k}:=\sum^{\ell}_{i=1}ka_iP\left(\varphi_m(P({\gamma_0}_i)_1)\right)
$$
is  a cycle representative of the class
$(\varphi_m)_{\ast}([\tilde{\varrho}_k])=k(\varphi_m)_\ast\circ
\iota_\ast(\mathfrak{C}) \in H_p(\Lambda, \Lambda(\gamma^m_0);\Z)$,
and that
\begin{equation}\label{e:8.7}
\cup^1_{i=0}{\Gamma_0}_i(\tilde{\varrho}_{1,m,k})\subset
\Lambda(\gamma^m_0).
\end{equation}
Unfortunately, the simplices of
${\Gamma_0}_1(\tilde{\varrho}_{1,m,k})$ are not necessarily smooth.
These obstruct direct applications of Proposition~\ref{prop:8.5}. It
is next step that helps us to overcome this difficulty.

{\it Step 3}. \textsf{For every integer $m\geq
\mathfrak{M}(\tilde{\varrho}_k)$ the class $k(\varphi_m)_\ast\circ
\iota_\ast(\mathfrak{C}) \in H_p(\Lambda, \Lambda(\gamma^m_0);\Z)$
has a cycle representative  ${\varrho}_{1,m,k}$ that is smooth on
each simplex, such that
$\cup^1_{i=0}{\Gamma_0}_i({\varrho}_{1,m,k})\subset
\Lambda(\gamma^m_0)$.}

 As in the proof of Lemma~\ref{lem:8.2} we consider the integral cycle
 $k\tilde{\varrho}_{1,m,k}$  a (continuous) map from a $p$-dimensional
finite oriented simplicial complex $K$ with boundary, which may be
assumed to be an Euclid complex in some $\R^N$, to $\Lambda$ such
that $\tilde{\varrho}_{1,m,k}(\partial K)\subset
\Lambda(\gamma_0^m)$. And so $(\tilde{\varrho}_{1,m,k})_\ast$ maps
the fundamental class $[K,\partial K]\in H_{p}(K, \partial K;\Z)$ to
$k(\varphi_m)_\ast\circ \iota_\ast(\mathfrak{C})$.

Now let us approximate $\tilde{\varrho}_{1,m,k}$ by a map
${\varrho}_{1,m,k}:K\to\Lambda$ that is smooth on each simplex.
Because $\Lambda(\gamma^{m}_0)$ is open, by (\ref{e:8.7}) we may
require
$$
\cup^1_{i=0}{\Gamma_0}_i({\varrho}_{1,m,k})\subset
\Lambda(\gamma^{m}_0).
$$
Hence ${\varrho}_{1,m,k}$ determines a new cycle representing  the
class $k(\varphi_{m})_\ast\circ \iota_\ast(\mathfrak{C})$, also
denoted by ${\varrho}_{1,m,k}$.

{\it Step 4}. Fix an integer $m_1\geq
\mathfrak{M}(\tilde{\varrho}_k)$. If $p>2$ we may continue to repeat
the above procedures and obtain integers
$m_{p-1}>m_{p-2}>\cdots>m_1$ and cycles ${\varrho}_{j,m_j,k}$,
$j=2,\cdots,p-1$, such that each ${\varrho}_{j,m_j,k}$ is smooth on
each simplex of it, satisfies (\ref{e:8.5}) and also represents the
class $k(\varphi_{m_{j}})_\ast\circ\iota_\ast(\mathfrak{C})\in
H_p(\Lambda, \Lambda(\gamma^{m_j}_0);\Z)$.

{\it Step 5}. For the cycle ${\varrho}_{p-1,m_{p-1},k}$, by
Proposition~\ref{prop:8.5} we get an integer
$\mathfrak{M}({\varrho}_{p-1,m_{p-1},k})$ such that for every
integer $m\geq \mathfrak{M}({\varrho}_{p-1,m_{p-1},k})$ and every
$p$-simplex $\eta\in{\Gamma_0}_{p}({\varrho}_{p-1,m_{p-1},k})$, we
have a homotopy $P(\eta^m):\triangle^{p}\times [0, 1]\to\Lambda$
 from  the $p$-simplex
$$
\eta^m\equiv \varphi_m(\eta):(\triangle^{p},\partial\triangle^{p})
            \to \left(\Lambda, \Lambda(\gamma_0^m)\right)
            $$
 to a  singular $p$-simplex
$\eta_m:(\triangle^{p},\partial\triangle^{p})\to
 \bigl(\Lambda(\gamma_0^m), \Lambda(\gamma_0^m)\bigr)$
 satisfying
$P(\eta^m)_t|_{\partial\triangle^{p}}=\eta^m|_{\partial\triangle^{p}}\;\forall
t\in [0,1]$. By (\ref{e:8.5}), for every $r$-simplex $\eta\in
\cup^{p-1}_{i=0}{\Gamma_0}_i({\varrho}_{p-1,m_{p-1},k})$ and for
every integer $m\ge\mathfrak{M}({\varrho}_{p-1,m_{p-1},k})$, we
define a homotopy $P(\eta^m):\triangle^{r}\times [0, 1]\to \Lambda$
by $P(\eta^m)(z,t)=\eta^m(z)\;\forall (z,t)\in\triangle^{r}\times
[0, 1]$. These homotopies satisfy the second part of
Lemma~\ref{lem:8.2} (i.e., \cite[Lemma~1]{BKl}) and hence
$$
k(\varphi_{m'})_\ast\circ\iota_\ast(\mathfrak{C})=
k(\varphi_m)_\ast\circ(\varphi_{m_{p-1}})_\ast\circ\iota_\ast(\mathfrak{C})=
k(\varphi_m)_\ast[{\varrho}_{p-1,m_{p-1},k}] =0
$$
in $H_{p}\bigl(\Lambda, \Lambda(\gamma_0^{m'});\Z\bigr)$ with
$m'=m+m_{p-1}$. Define
$\mathfrak{M}:=\mathfrak{M}({\varrho}_{p-1,m_{p-1},k})+ m_{p-1}$,
the desired conclusion holds for it. Proposition~\ref{prop:8.7} is
proved. \hfill$\Box$\vspace{2mm}

\noindent{\bf Proof of Proposition~\ref{prop:8.5}}. Since
$\triangle^q$ and $\partial\triangle^q$ are compact,
\begin{eqnarray*}
\mathfrak{K}_0(\eta):=\max\left\{{\cal
L}(\eta(x))\;|\;x\in\partial\triangle^q \right\}\quad\hbox{and}\quad
\mathfrak{K}_1(\eta):=\max\left\{{\cal
L}(\eta(x))\;|\;x\in\triangle^q\right\}
\end{eqnarray*}
are always finite. Following \cite{Lo1}, for $t,s\in [0,1]$ let
$e(t,s)=(t,\cdots,t,s)\in\R^q\times [0,1]$ and thus the barycenter
of $\triangle^q\times\{s\}\subset\R^q\times\{s\}$ is
$\hat{e}(s):=e(1/(q+1),s)$. Set
$$
\triangle^q(s)=e((1-s)/(q+1),s)+
(s\triangle^q)\times\{0\}\subset\triangle^q\times\{s\}
$$
for $s\in [0,1]$. (Clearly, $\triangle^q(1)=\triangle^q\times\{1\}$
and $\triangle^q(0)=\hat{e}(0)$.) Denote by $L(s)$ the straight line
passing through $e(0,s)$ and $\hat e(s)$ successively in
$\R^q\times\{s\}$, that is, $ L(s)=\{e(t,s)\,|\,t\in\R\}$.
 Then we have an orthogonal subspace decomposition
$\R^q\times\{s\}=V_{q-1,s}\times L(s)$, and each
$w\in\triangle^q\times\{s\}$ may be uniquely written as
$w=(v,\tau)\in V_{q-1,s}\times L(s)$.  For such a
$w=(v,\tau)\in\triangle^q(s)$ denote by $l(v,s)$ the intersection
segment of $\triangle^q(s)$ with the straight line passing through
$w$ and parallel to $L(s)$. Actually $l(v,s)=\{(v, e(t,s))\,|\,
a(v,s)\le t\le b(v,s)\}$, where $a(v,s)$ and $b(v,s)$ depend
piecewise smoothly on $(v,s)$.

For $w=(v,\tau)\in\triangle^q(s)$  define a $v$-parameterized curve
$$
^\sharp\eta_v:l(v,s)\to\Lambda,\;\tau\mapsto\eta(v,\tau)\quad\hbox{for}\;
(v,\tau)=w\in[V_{q-1,s}\times L(s)]\cap\triangle^q(s).
$$
It is  $C^1$, and the corresponding initial value curves $\eta^{\rm
ini}_v:l(v,s)\to M $ is also  $C^1$ as the evaluation
$\Lambda\ni\gamma\mapsto\gamma(0)\in M$ is smooth.
 As in \cite{Lo1} using $\eta^{\rm ini}_v$ we define
$\eta_m:\triangle^q\to\Lambda$ by $\eta_m(w)=(^\sharp\eta_v)_m(t)$,
where $(v,t)$ is the unique $V_{q-1,1}\times L(1)$ decomposition of
$(w,1)\in\triangle^q(1)=\triangle^q\times\{1\}$. It is continuous.
Writing $(w,s)\in\triangle^q\times [0,1]$ as $(v,t,s)\in
[V_{q-1,s}\times L(s)]\times [0,1]$ we define a homotopy
$H:\triangle^q\times [0,1]\to\Lambda$ by
$$
H(w,s)=\left\{\begin{array}{ll} \eta^m(w)
&\hbox{if}\;(w,s)\in(\triangle^q\times\{s\})\setminus\triangle^q(s),\\
(^\sharp\eta_v)_m(t) &\hbox{if}\;(w,s)=(v,t,s)\in\triangle^q(s).
\end{array}\right.
$$
It is continuous and satisfies $H(w,0)=\eta^m$, $H(w,1)=\eta_m$ and
$H(w,s)=\eta^m(w)$ for any $(w,s)\in\partial\triangle^q\times
[0,1]$. Now for
$(w,s)\in(\triangle^q\times\{s\})\setminus\triangle^q(s)$ it holds
that
$$
{\cal L}(H(w,s))=m^2{\cal L}(\eta(w))\le m^2\mathfrak{K}_1(\eta),
$$
and for $(w,s)=(v,t,s)\in\triangle^q(s)$ it follows from
\cite[(1)]{BKl} that
\begin{eqnarray*}
{\cal L}(H(w,s))&=&{\cal L}((^\sharp\eta_v)_m(t))\\
&\le&(m+
2|l(v,s)|)\bigl((m-1)\kappa_0(^\sharp\eta_v)+\kappa_1(^\sharp\eta_v)+2\kappa_2(^\sharp\eta_v)\bigr)\\
&\le&(m+ 4)\bigl((m-1)\mathfrak{K}_0(\eta)+
\mathfrak{K}_1(\eta)+\int^{b(v,s)}_{a(v,s)}\Bigl\|\frac{d}{dt}
\eta(v,t)(0)\Bigr\|^2_1dt\bigr)\\
&\le&(m+ 4)\bigl((m-1)\mathfrak{K}_0(\eta)+
\mathfrak{K}_1(\eta)+\mathscr{C}(\eta)\bigr).
\end{eqnarray*}
Here $\mathscr{C}(\tilde\eta):=\max\{|\nabla_x(\eta(x)(0))|^2\;
 |\;x\in\triangle^q\}<+\infty$ because
  $\triangle^q\ni x\to\eta(x)(0)\in M$ is  $C^1$.
 Proposition~\ref{prop:8.5}(iii)-(v) follows
immediately.\hfill$\Box$\vspace{2mm}

By the arguments below Assumption F it is not hard to see that we
actually prove the following result: if the closed geodesic
$\bar\gamma$ in Theorem~\ref{th:1.11} is big in the sense that
$H_{q}(\Lambda(\gamma)\cup S^1\cdot \gamma,
\Lambda(\gamma);\Q)=0\;\forall q>\bar{p}$ for each other closed
geodesic $\gamma$ with $m^-({\gamma}^k)=0\;\forall k\in\mathbb{N}$, then there exist
infinitely many geometrically distinct closed geodesics on $(M,F)$
such that each of them is free homotopic to one of
$\{\bar\gamma^k\,|\, k\in\N\}$.

\begin{remark}\label{rm:8.8}
{\rm Our proof method is slightly different from \cite{BKl}. It is
key for us that $\gamma_0$ is not a local minimum of the  functional
 $\mathcal{L}$, which comes from (\ref{e:8.1}) and thus our
 assumption that there  exist a  closed geodesics $\bar\gamma$ and an integer
$\bar{p}\ge 2$ such that $m^-(\bar{\gamma}^k)=0\;\forall k\in\mathbb{N}$ and
$H_{\bar{p}}(\Lambda(\bar\gamma)\cup S^1\cdot \bar\gamma,
\Lambda(\bar\gamma);\Q)\ne 0$. According to \cite[Theorem 3]{BKl}
one should make the following weaker: \vspace{1mm}

\noindent{\bf Assumption}: There  exists a  closed geodesics $\gamma$
with $m^-(\gamma^k)=0\;\forall k\in\mathbb{N}$ such that
$H_{q}(\Lambda(\gamma)\cup S^1\cdot \gamma, \Lambda(\gamma);\Q)\ne
0$ for some $q\in\N\cup\{0\}$ and that $\gamma$ is not an absolute
minimum of $\mathcal{L}$ in its free homotopy class.\vspace{1mm}

Let us explain how such an assumption leads to the following
sentence at the beginning of the proof of \cite[Theorem 3, page
385]{BKl} (in our notations and with the $\Z$-coefficient groups
being replaced by $\Q$-coefficient ones): \vspace{1mm}

\noindent{\bf Sentence} (\cite[lines 5-6 on page 385]{BKl}):
\textsf{Choosing a different $\gamma$, if necessary, we can find
$p\in\N$ such that
 $H_p(\Lambda(\gamma)\cup S^1\cdot\gamma,\; \Lambda(\gamma);\Q)\ne 0$ and
 $H_q(\Lambda(\beta)\cup S^1\cdot\beta,\; \Lambda(\beta);\Q)= 0$ for
every $q>p$ and every closed geodesic $\beta$ with
$m^-(\beta^k)=0\;\forall k\in\mathbb{N}$.}\vspace{1mm}

In fact, by  Theorem~\ref{th:1.8}  the following case cannot occur:
$$
H_0(\Lambda(\gamma)\cup S^1\cdot\gamma,\; \Lambda(\gamma);\Q)\ne
0\quad\hbox{and}\quad H_k(\Lambda(\gamma)\cup S^1\cdot\gamma,\;
\Lambda(\gamma);\Q)=0\;\forall k\in\N.
$$
Hence we may choose $s\in\N$ such that $H_s(\Lambda(\gamma)\cup
S^1\cdot\gamma,\; \Lambda(\gamma);\Q)\ne 0$ and
$$
H_q(\Lambda(\gamma)\cup S^1\cdot\gamma,\;
\Lambda(\gamma);\Q)=0\;\forall q>s
$$
because the shifting theorem implies that $H_q(\Lambda(\alpha)\cup
S^1\cdot\alpha,\; \Lambda(\alpha);\Q)=0$ for all $q\notin [0, 2n-1]$
for every closed geodesic $\alpha$ with $m^-(\alpha)=0$.\\
$\bullet$ If $H_q(\Lambda(\beta)\cup S^1\cdot\beta,\;
\Lambda(\beta);\Q)= 0$ for every $q>s$ and every closed geodesic
$\beta$ with $m^-(\beta^k)=0\;\forall k\in\mathbb{N}$, then we may take $p=s$ in the
above sentence. In this case  the closed geodesic $\gamma$ is not a
local minimum of $\mathcal{L}$ if $s>1$ by Theorem~\ref{th:1.8}(ii).
If $s=1$ then $\gamma$ is a local minimum of  $\mathcal{L}$ by
Theorem~\ref{th:1.8}(i), but is not an absolute minimum of
$\mathcal{L}$ in the free homotopy class of $\gamma$ by
the above Assumption.\\
$\bullet$ If there must exist an integer $l>s$ and a closed geodesic
$\bar\gamma$ with $m^-(\bar\gamma^k)=0\;\forall k\in\mathbb{N}$ such that
$H_l(\Lambda(\bar\gamma)\cup S^1\cdot\bar\gamma,\;
\Lambda(\bar\gamma);\Q)\ne 0$, it is easily seen that such a pair
$(l,\bar\gamma)$ can be chosen so that $H_q(\Lambda(\beta)\cup
S^1\cdot\beta,\; \Lambda(\beta);\Q)=0$ for every integer $q>l$ and
for every closed geodesic $\beta$ with $m^-(\beta^k)=0\;\forall k\in\mathbb{N}$. In
this case $\gamma$ and $p$ in the above sentence  can be chosen as
$\bar\gamma$ and $l$, respectively. Since $p=l\ge 2$
Theorem~\ref{th:1.8}(ii) shows that $\gamma=\bar\gamma$ is not a
local  minimum of $\mathcal{L}$.

Summarizing the above arguments we equivalently expressed the above
Assumption as follows:\\
$\bullet$ either $H_1(\Lambda(\gamma)\cup S^1\cdot\gamma,\;
\Lambda(\gamma);\Q)\ne 0$ and  $H_q(\Lambda(\beta)\cup
S^1\cdot\beta,\; \Lambda(\beta);\Q)= 0$ for every $q>1$ and every
closed geodesic $\beta$ with $m^-(\beta^k)=0\;\forall k\in\mathbb{N}$, and $\gamma$ is
not an absolute minimum of $\mathcal{L}$ in the free homotopy class
of $\gamma$;\\
$\bullet$ or $\exists\, p\ge 2$ and a closed geodesic $\gamma'$ such
that $H_p(\Lambda(\gamma')\cup S^1\cdot\gamma',\;
\Lambda(\gamma');\Q)\ne 0$  and that $H_q(\Lambda(\beta)\cup
S^1\cdot\beta,\; \Lambda(\beta);\Q)= 0$ for every $q>p$ and every
closed geodesic $\beta$ with $m^-(\beta^k)=0\;\forall k\in\mathbb{N}$. (Hence
$\gamma'$ not a local  minimum of $\mathcal{L}$.)

Hence the above Assumption may lead to  the above sentence.
Comparing with the method by Bangert and Klingenberg \cite{BKl} ours
can only deal with the latter case. }
\end{remark}

\vspace{2mm}
 \noindent{\bf Acknowledgements}.  I would like to thank
Professors Rong Ruan and Jiequan Li for producing the beautiful
figure. The partial results of this paper were reported on the
International Conference in Geometry, Analysis and Partial
Differential Equations held at Jiaxing University on December 17-21,
2011; I thank Professor Guozhen Lu for invitation and hospitality. I
also thank Professors Jean-Noel Corvellec and Erasmo Caponio for
sending me their related papers respectively. I am very grateful to
Dr. Marco Mazzucchelli for pointing out my misunderstanding for
the precise meaning of a sentence in \cite{BKl} (Sentence 2 in
Section~9 of the original version  of the paper).
 I am deeply grateful to  the anonymous referees for kindly checking this paper,
  giving many useful comments and correcting
many errors in mathematics, English and typos.
 This work was partially supported
by the NNSF   10971014 and 11271044 of China,  RFDPHEC (No. 200800270003)
 and the Fundamental Research Funds for the Central Universities (No. 2012CXQT09).

\appendix

\section{Appendix:\quad The splitting theorems in
\cite{Lu1,Lu2, Lu3, Ji}} \label{app:A}\setcounter{equation}{0}

 Let $H$ be a Hilbert space with inner product
$(\cdot,\cdot)_H$ and the induced norm $\|\cdot\|$, and let $X$ be a
Banach space with norm $\|\cdot\|_X$, such that
\begin{description}
\item[(S)]  $X\subset H$ is dense in $H$ and
 $\|x\|\le \|x\|_X\;\forall x\in X$.
\end{description}
For an open neighborhood $U$ of the origin $0\in H$,   write
$U_X=U\cap X$ as  an open neighborhood of $0$ in $X$. Let  ${\cal
L}\in C^1(U,\mathbb{R})$ have  $0$ as an isolated critical point.
Suppose that there exist maps $A\in C^1(U_X, X)$ and $B\in C(U_X,
L_s(H))$ such that
\begin{eqnarray}
&&{\cal L}'(x)(u)=(A(x), u)_H\quad\forall x\in U_X\;\hbox{and}\;
u\in X,\label{e:A.1}\\
&&(A'(x)u, v)_H=(B(x)u, v)_H\;\forall x\in U_X\;\hbox{and}\; u, v\in
X. \label{e:A.2}
\end{eqnarray}
(These imply: (a) ${\cal L}|_{U_X}\in C^2(U_X, \R)$,  (c) $d^2{\cal
L}|_{U_X}(x)[u,v]=(B(x)u,v)_H$ for any $x\in U_X$ and $u, v\in X$,
(c) $B(x)(X)\subset X\;\forall x\in U_X$). Furthermore we  assume
$B$ to satisfy  the following properties:
\begin{description}
\item[(B1)]  If $u\in H$ such that $B(0)(u)=v$ for some
$v\in X$, then $u\in X$. Moreover,  all eigenvectors of the operator
$B(0)$ that correspond to negative eigenvalues belong to $X$.

\item[(B2)] The map $B:U_X\to
L_s(H)$  has a decomposition $B(x)=P(x)+ Q(x)$ for each $x\in U_X$,
where $P(x):H\to H$ is a positive definitive linear operator and
$Q(x):H\to H$ is a compact linear  operator with the following
properties:
\begin{description}
\item[(i)] For any sequence $\{x_k\}\subset
U_X$ with $\|x_k\|\to 0$ it holds that $\|P(x_k)u-P(0)u\|\to 0$ for
any $u\in H$;

\item[(ii)] The  map $Q:U\cap X\to
L(H)$ is continuous at $0$ with respect to the topology induced from
$H$ on $U\cap X$;

\item[(iii)] There exist positive constants $\eta_0>0$ and  $C_0>0$ such that
$$
(P(x)u, u)\ge C_0\|u\|^2\quad\forall u\in H,\;\forall x\in
X\;\hbox{with}\;\|x\|<\eta_0.
$$
\end{description}
\end{description}

  Let $H^-$, $H^0$ and $H^+$ be the negative definite, null and
  positive definite spaces of $B(0)$.  Then  $H=H^-\oplus
H^0\oplus H^+$. By (B1) and (B2)  both $H^0$ and $H^-$ are
finite-dimensional subspaces contained in $X$.
 Denote by $P^\ast$ the orthogonal
projections onto $H^\ast$, $\ast=+, -, 0$, and by $X^\ast=X\cap
H^\ast=P^\ast(X),\;\ast=+, -$. Then $X^+$ is dense in $H^+$, and
$(I-P^0)|_X=(P^++P^-)|_X: (X, \|\cdot\|_X)\to (X^++X^-, \|\cdot\|)$
is also continuous because all norms are equivalent on a linear
space of finite dimension.  These give the topological direct sum
decomposition: $X=H^-\dot{+} H^0\dot{+} X^+$.  $m^0=\dim H^0$ and
$m^-=\dim H^-$  are called the {\bf nullity} and the {\bf Morse
index} of critical point $0$ of ${\cal L}$, respectively. The
critical point $0$ is called {\bf nondegenerate} if $m^0=0$. The
following is, except for (iv), Theorem~1.1 of \cite{Lu1} , which is a
special version of Theorem~2.1 in \cite{Lu3}.

\begin{theorem}\label{th:A.1}
Under the above assumptions (S) and (B1)-(B2), there exist a
positive number $\epsilon\in\R$,  a  $C^1$ map $h:{\bf
B}_\epsilon(H^0)\to X^++X^-$ satisfying $h(0)=0$ and
$$
 (I-P^0)A(z+ h(z))=0\quad\forall z\in {\bf B}_\epsilon(H^0),
$$
an open neighborhood $W$ of $0$ in $H$ and an origin-preserving
homeomorphism
$$
\Phi: {\bf B}_\epsilon(H^0)\times \left({\bf B}_\epsilon(H^+) + {\bf
B}_\epsilon(H^-)\right)\to W
$$
of form $\Phi(z, u^++ u^-)=z+ h(z)+\phi_z(u^++ u^-)$ with
$\phi_z(u^++ u^-)\in H^\pm$  such that
$$
{\cal L}\circ\Phi(z, u^++ u^-)=\|u^+\|^2-\|u^-\|^2+ {\cal L}(z+
h(z))
$$
for all $(z, u^+ + u^-)\in {\bf B}_\epsilon(H^0)\times \left({\bf
B}_\epsilon(H^+) + {\bf B}_\epsilon(H^-)\right)$, and that
$$
\Phi\left({\bf B}_\epsilon(H^0)\times \left({\bf
B}_\epsilon(H^+)\cap X + {\bf B}_\epsilon(H^-)\right)\right)\subset
X.
$$
Moreover,  $\Phi$, $h$ and the function ${\bf
B}_\epsilon(H^0)\ni z\mapsto {\cal L}^\circ(z):={\cal L}(z+ h(z))$
also satisfy:
\begin{description}
\item[(i)]  For  $z\in {\bf B}_\epsilon(H^0)$, $\Phi(z, 0)=z+ h(z)$,
$\phi_z(u^++ u^-)\in H^-$ if and only if $u^+=0$;

\item[(ii)]  $h'(z)=-[(I-P^0)A'(z+
h(z))|_{X^\pm}]^{-1}(I-P^0)A'(z+ h(z))|_{H^0} \quad\forall z\in {\bf
B}_\epsilon(H^0)$;

\item[(iii)] ${\cal L}^\circ$ is $C^{2}$, has an isolated critical point $0$,   $d^2{\cal L}^\circ(0)=0$ and
 $$
d{\cal L}^\circ(z_0)[z]=(A(z_0+ h(z_0)), z)_H\quad\forall z_0\in
{\bf B}_\epsilon(H^0),\; z\in H^0.
 $$

\item[(iv)] $\Phi$ is also a homeomorphism from ${\bf B}_\epsilon(H^0)\times {\bf B}_\epsilon(H^-)$ to
$\Phi\left({\bf B}_\epsilon(H^0)\times  {\bf
B}_\epsilon(H^-)\right)$ even if the topology on the latter is taken
as the induced one by $X$. (This implies that $H$ and $X$ induce the
same topology in $\Phi\left({\bf B}_\epsilon(H^0)\times  {\bf
B}_\epsilon(H^-)\right)$.)
\end{description}
\end{theorem}

\begin{corollary}[\cite{Lu1, Lu3}]\label{cor:A.2}{\rm (Shifting)}
Under the assumptions of Theorem~\ref{th:A.1},  for any Abelian
group ${\K}$ it holds that $C_q({\cal L}, 0;{\K})\cong
C_{q-m^-}({\cal L}^{\circ}, 0; {\K})$ for each $q=0, 1,\cdots$.
\end{corollary}

\begin{theorem}[\hbox{\cite[Th.6.4]{Lu3}}]\label{th:A.3}
Under the assumptions of Theorem~\ref{th:A.1}, suppose that
$\check{H}\subset H$ is a Hilbert subspace whose orthogonal
complementary in $H$ is  finite-dimensional and is contained in $X$.
Then $({\cal L}|_{\check{H}}, \check{H}, \check{X})$ with
$\check{X}:=X\cap\check{H}$ also satisfies the assumptions of
Theorem~\ref{th:A.1} around the critical point $0\in\check{H}$.
\end{theorem}

\begin{theorem}[\hbox{\cite[Th.6.3]{Lu3}}]\label{th:A.4}
 Under the assumptions of Theorem~\ref{th:A.1}, let
$(\widehat{H},\widehat{X})$ be another pair of Hilbert-Banach spaces
satisfying {\bf (S)}, and let $J:H\to\widehat{H}$ be a Hilbert space
isomorphism which can induce a Banach space isomorphism $J_X:X\to X$
(this means that $J(X)\subset \widehat{X}$ and $J|_X:X\to
\widehat{X}$ is a Banach space isomorphism). Set $\widehat{U}=J(U)$
(and hence $\widehat{U}_{\widehat
X}:=\widehat{U}\cap\widehat{X}=J(U_X)$) and $\widehat{\cal
L}:\widehat{U}\to\R$ by $\widehat{\cal L}={\cal L}\circ J^{-1}$.
Then $(\widehat{H},\widehat{X}, \widehat{U}, \widehat{\cal L})$
satisfies the assumptions of Theorem~\ref{th:A.1} too.
\end{theorem}

Let us give the following special case of the splitting lemma by
Ming Jiang \cite{Ji} (i.e. the case $X=Y$) with our  above
notations.

\begin{theorem}\label{th:A.5}{\rm (\cite[Th.2.5]{Ji})}
 Let the assumptions above {\bf (B1)} be satisfied when the neighborhood
$U_X$ therein is replaced by some ball ${\bf
B}_\varepsilon(X)\subset X$. (These imply that the conditions
(SP),\footnote{In the original Theorem~2.5 of \cite{Ji} the density
of $X$ in $H$ is not needed.} (FN1)-(FN3) in \cite[Th.2.5]{Ji} hold
in ${\bf B}_\varepsilon(X)$). Suppose also that the first condition
in {\bf (B1)} and the following condition are satisfied:
\begin{description}
\item[(CP1)] either $0$ is  not in the spectrum ${\sigma}(B(0))$ or it is an
isolated point of ${\sigma}(B(0))$. (This holds if  {\bf (B2)} is
true by Proposition~B.2 in \cite{Lu2}.)
\end{description}
Then  there exists a ball ${\bf B}_\delta(X)\subset {\bf
B}_\varepsilon(X)$, an origin-preserving local homeomorphism
$\varphi$ from ${\bf B}_\delta(X)$ to an open neighborhood of $0$ in
$U_X$ and a $C^1$ map $h:{\bf B}_\delta(H^0)\to X^++X^-$ such that
\begin{equation}\label{e:A.3}
{\cal L}|_X\circ\varphi(x)=\frac{1}{2}(B(0)x^\bot, x^\bot)_H+ {\cal
L}(h(z)+z)\;\forall x\in {\bf B}_\delta(X),
\end{equation}
where $z=P^0(x)$ and $x^\bot=(P^++P^-)(x)$. Moreover, the function
${\bf B}_\delta(X)\ni z\mapsto {\cal L}_X^{\circ}(z):={\cal
L}(h(z)+z)$ has the same properties as ${\cal L}^{\circ}$ in
Theorem~\ref{th:A.1}.
\end{theorem}

If $X^-:=X\cap H^-=H^-$ then $-B(0)|_{H^-}:H^-\to H^-$ is positive
definite and therefore $(-B(0)|_{H^-})^{\frac{1}{2}}:H^-\to H^-$ is
an isomorphism. It follows that
$$
\frac{1}{2}(B(0)P^+x,
P^+x)=\left(2^{-\frac{1}{2}}(-B(0)|_{H^-})^{\frac{1}{2}}P^-x,
2^{-\frac{1}{2}}(-B(0)|_{H^-})^{\frac{1}{2}}P^-x\right).
$$
Consider the Banach space isomorphism $\Gamma:X\to X$ given by
$$
X=X^-\dot{+} X^0\dot{+} X^+\ni x=x^-+ x^0+ x^+\to
2^{-\frac{1}{2}}(-B(0)|_{H^-})^{\frac{1}{2}}x^-+ x^0+ x^+.
$$
Replaceing $\varphi$ by $\varphi\circ\Gamma^{-1}$ and shrinking
$\delta>0$ suitably, (\ref{e:A.3}) becomes
\begin{equation}\label{e:A.4}
{\cal L}|_X\circ\varphi\circ\Gamma^{-1}(x)=\frac{1}{2}(B(0)x^+,
x^+)_H-\|x^-\|^2+ {\cal L}(h(z)+z)\;\forall x\in {\bf B}_\delta(X),
\end{equation}
where $z=P^0(x)$, $x^\star=P^\star x$, $\star=+,-$, and
$\|x^-\|^2=(x^-,x^-)_H$. In some applications such a form of
Theorem~\ref{th:A.5} is more convenient.

Clearly, the assumptions of Theorem~\ref{th:A.1} are stronger than
those of Theorem~\ref{th:A.5}. Under the conditions of
Theorem~\ref{th:A.1}  these two theorems take the same maps $h$ and
hence $\mathcal{L}^\circ_X$ and $\mathcal{L}^\circ$ are same near
$0\in H^0$.

 Following the methods of proof in
\cite[Th.5.1.17]{Ch05} we may obtain

\begin{corollary}\label{cor:A.6}{\rm (Shifting)}
Under the assumption of Theorem~\ref{th:A.5} suppose that
$H^-\subset X$ and $\dim (H^0+ H^-)<\infty$.  Then $C_q({\cal
L}|_{U_X}, 0;\K)=C_{q-m^-}({\cal L}_X^{\circ}, 0;\K)$ for each
$q=0,1,\cdots$, where $m^-:=\dim H^-$.
\end{corollary}

By simple arguments as in the proofs of \cite[Th.6.3, Th.6.4]{Lu3}
we may obtain the following corresponding results with
Theorems~\ref{th:A.3},~\ref{th:A.4}.

\begin{theorem}\label{th:A.7}
Under the assumptions of Theorem~\ref{th:A.5}, suppose that
$\check{H}\subset H$ is a Hilbert subspace whose orthogonal
complementary in $H$ is  finite-dimensional and is contained in $X$.
Then $({\cal L}|_{\check{H}}, \check{H}, \check{X},
\check{A}:=P_{\check{H}}\circ A|_{\check{X}},
\check{B}(\cdot):=P_{\check{H}}B(\cdot)|_{\check{H}})$, where
$\check{X}:=X\cap\check{H}$ and $P_{\check{H}}$ is the orthogonal
projection onto $\check{H}$, also satisfies the assumptions of
Theorem~\ref{th:A.5}.
\end{theorem}

\begin{theorem}\label{th:A.8}
Under the assumptions of Theorem~\ref{th:A.1}, let
$(\widehat{H},\widehat{X})$ be another pair of Hilbert-Banach spaces
satisfying {\bf (S)}, and let $J:H\to\widehat{H}$ be a Hilbert space
isomorphism which can induce a Banach space isomorphism $J_X:X\to
\widehat{X}$. Then $(\widehat{H},\widehat{X}, \widehat{\cal L}={\cal
L}\circ J^{-1})$ satisfies the assumptions of Theorem~\ref{th:A.5}
too.
\end{theorem}

Under Corollaries~\ref{cor:A.2},~\ref{cor:A.6} it holds that for any
$q=0,1,\cdots,$
\begin{equation}\label{e:A.5}
C_q({\cal L}, 0;\K)=C_{q-m^-}({\cal L}^\circ, 0;\K)=C_{q-m^-}({\cal
L}_X^\circ, 0;\K)=C_q({\cal L}_{U_X}, 0;\K).
\end{equation}
because ${\cal L}^\circ$ and ${\cal L}_X^\circ$ are same near $0\in
{\bf B}_\delta(X)$. Actually, a stronger result holds.

\begin{theorem}[\hbox{\cite[Cor.2.5]{Lu3}}]\label{th:A.9}
Under the assumptions of Theorem~\ref{th:A.1}  let $c={\cal L}(0)$.
For any open neighborhood $W$ of $0$ in $U$ and a field $\F$, write
$W_X=W\cap X$ as an open subset of $X$, then the inclusion
\begin{equation}\label{e:A.6}
\left({\cal L}_c\cap W_X, {\cal L}_c\cap W_X
\setminus\{0\}\right)\hookrightarrow\left({\cal L}_c\cap W, {\cal
L}_c\cap W\setminus\{0\}\right)
\end{equation}
induces isomorphisms
$$H_\ast\left({\cal L}_c\cap W_X,
{\cal L}_c\cap W_X\setminus\{0\};\F\right)\cong H_\ast\left({\cal
L}_c\cap W, {\cal L}_c\cap W\setminus\{0\};\F\right).
$$
\end{theorem}

It is not hard to prove that the corresponding conclusion also holds
true if we replace $\left({\cal L}_c\cap W_X, {\cal L}_c\cap W_X
\setminus\{0\}\right)$ and $\left({\cal L}_c\cap W, {\cal L}_c\cap
W\setminus\{0\}\right)$ in (\ref{e:A.6}) by $\bigl(\mathring{\cal
L}_c\cap W_X\cup\{0\}, \mathring{\cal L}_c\cap W_X \bigr)$ and
$\bigl(\mathring{\cal L}_c\cap W\cup\{0\}, \mathring{\cal L}_c\cap
W\bigr)$ respectively, where $\mathring{\cal L}_c=\{{\cal L}<c\}$.


\section{Appendix:\quad  Computations of gradients}
\label{app:B}\setcounter{equation}{0}

We shall make computations in general case. Let $E_{\gamma_0}={\rm
diag}(S_1,\cdots,S_{\sigma})\in\R^{n\times n}$ with ${\rm
ord}(S_1)\ge\cdots\ge{\rm ord}(S_{\sigma})$ be as in
Section~\ref{sec:5.3}. By the proof of Claim~\ref{cl:6.3} we have
$$
E_{\gamma_0^m}={\rm diag}(S_1^m,\cdots,S^m_{\sigma})={\rm
diag}(\hat{S}_1,\cdots,\hat{S}_{\tau})\quad\hbox{with}\quad {\rm
ord}(\hat{S}_1)\ge\cdots\ge{\rm ord}(\hat{S}_{\tau}),
$$
where each $\hat{S}_j$ is either $1$, or $-1$, or
${\scriptscriptstyle\left(\begin{array}{cc}
\cos\theta_j& \sin\theta_j\\
-\sin\theta_j& \cos\theta_j\end{array}\right)}$, $0<\theta_j<\pi$.

 For an integer $m\in\N$, assume that
 $L_{\gamma_0^m}\in C^1(\R\times\R^n\times\R^n, \R)$ and there
exist $a\in C(\R^+,\R^+)$, $b\in L^1([0,1],\R^+)$, $c\in
L^2([0,1],\R^+)$ such that $L_{\gamma_0^m}(t+1, x,
v)=L_{\gamma_0^m}(t, (E_{\gamma_0^m} x^T)^T, (E_{\gamma_0^m}
v^T)^T)\;\forall t$ and
\begin{eqnarray*}
&&|L_{\gamma_0^m}(t, x, v)|+|\partial_xL_{\gamma_0^m}(t, x, v)|\le a(|x|)(b(t)+ |v|^2),\\
&&|\partial_xL_{\gamma_0^m}(t, x, v)|\le a(|x|)(c(t)+
|v|)\quad\forall (t,x,v)\in [0,1]\times\R^n\times\R^n.
\end{eqnarray*}
Let $\hat{H}_{\gamma_0^m}$ be the Hilbert space above
(\ref{e:6.23}).
 Then the functional
$$
{\cal L}_{\gamma_0^m}:\hat{H}_{\gamma_0^m}\to\R,\;\xi\mapsto
\int^1_0L_{\gamma_0^m}(t,\xi(t),\dot{\xi}(t))dt
$$
is $C^1$ by Theorem~1.4 on the page 10 of \cite{MW}. Denote by
$\nabla^m{\cal L}_{\gamma_0^m}$ the gradient of ${\cal
L}_{\gamma_0^m}$. Note that $\hat{H}_{\gamma_0^1}=H_{\gamma_0}$ and
$\nabla^1{\cal L}_{\gamma_0^1}=\nabla{\cal L}_{\gamma_0}$.
 By the continuality of $\nabla^m
{\cal L}_{\gamma_0^m}$ we only need to compute $\nabla^m{\cal
L}_{\gamma_0^m}(\xi)$ for  $\xi\in{X}_{\gamma_0^m}$. Write
\begin{eqnarray}
&&\partial_v
L_{{\gamma_0}^m}\bigl(t,\xi(t),\dot{\xi}(t)\bigr)=(\mathfrak{V}^\xi_1(t),\cdots,
\mathfrak{V}^\xi_n(t)),\label{e:B.1}\\
&&\partial_x
L_{{\gamma_0}^m}\bigl(t,\xi(t),\dot{\xi}(t)\bigr)=(\mathfrak{X}^\xi_1(t),\cdots,
\mathfrak{X}^\xi_n(t)). \label{e:B.2}
\end{eqnarray}
Since $\xi(t+1)^T=E_{\gamma_0^m}\xi(t)^T$, by the assumptions on
$L_{\gamma_0^m}$ we have
\begin{eqnarray}
\partial_x L_{\gamma_0^m}(t+1,\xi(t+1),\dot\xi(t+1))&=&\partial_x
L_{\gamma_0^m}(t,\xi(t),\dot\xi(t))E_{\gamma_0^m}\quad\forall t\in\R,\label{e:B.3}\\
\partial_v L_{\gamma_0^m}(t+1,\xi(t+1),\dot\xi(t+1))&=&\partial_v
L_{\gamma_0^m}(t,\xi(t),\dot\xi(t))E_{\gamma_0^m}\quad\forall
t\in\R.\label{e:B.4}
\end{eqnarray}

\noindent{\bf Case 1}. ${\rm ord}(\hat{S}_1)=\cdots={\rm
ord}(\hat{S}_\tau)=1$ and so $\tau=n$.

By (\ref{e:B.1})-(\ref{e:B.2}) and (\ref{e:B.3})-(\ref{e:B.4}) we
have
\begin{eqnarray}
&&\mathfrak{V}^\xi_j(t+1)=-\mathfrak{V}^\xi_j(t),\quad
\mathfrak{X}^\xi_j(t+1)=-\mathfrak{X}^\xi_j(t)\;\forall
t\in\R\quad\hbox{if}\; \hat{S}_j=-1,\label{e:B.5}\\
 &&\mathfrak{V}^\xi_j(t+1)=\mathfrak{V}^\xi_j(t),\quad
\mathfrak{X}^\xi_j(t+1)=\mathfrak{X}^\xi_j(t)\;\forall
t\in\R\quad\hbox{if}\;\hat{S}_j=1.\label{e:B.6}
\end{eqnarray}
For the case of (\ref{e:B.5}), define $\mathfrak{G}^{\xi}_j:\R\to\R$
by
\begin{eqnarray}\label{e:B.7}
\mathfrak{G}^{\xi}_j(t)=\int^t_0\mathfrak{V}^\xi_j(s)ds-
\frac{1}{2}\int^1_0\mathfrak{V}^\xi_j(s)ds\quad\forall t\in\R.
\end{eqnarray}
It is the primitive function of the function $\mathfrak{V}^\xi_j(t)$,
and satisfies
\begin{eqnarray}\label{e:B.8}
\mathfrak{G}^{\xi}_j(t+1)=-\mathfrak{G}^{\xi}_j(t)\quad\forall
t\in\R.
\end{eqnarray}
For the case of (\ref{e:B.6}),   we define
\begin{equation}\label{e:B.9}
\mathfrak{G}^{\xi}_j(t)=\int^t_0\Bigl[\mathfrak{V}^\xi_j(s)-\int^{1}_0\mathfrak{V}^\xi_j(u)du\Bigr]ds,
\end{equation}
which is a $1$-periodic primitive function of the function
$s\mapsto\mathfrak{V}^\xi_j(s)-\int^{1}_0\mathfrak{V}^\xi_j(u)du$.
In these two cases it is easily proved that
\begin{eqnarray*}
\int^1_0\mathfrak{V}^\xi_j(t)\dot\eta_j(t)dt
=\int^1_0\dot{\mathfrak{G}}^{\xi}_j(t)\dot\eta_j(t)dt
=\langle\mathfrak{G}^{\xi}_j,\eta_j\rangle_{1,2,m}-m^2\int^1_0\mathfrak{G}^{\xi}_j(t)\eta_j(t)dt
\end{eqnarray*}
for any $\eta=(\eta_1,\cdots,\eta_n)\in \hat{H}_{\gamma_0^m}$ and
$j=1,\cdots,n$. By this and (\ref{e:B.1})-(\ref{e:B.2}) we obtain
\begin{eqnarray}\label{e:B.10}
d{\cal L}_{\gamma_0^m}(\xi) (\eta)
 \!\!\!&=&\sum^n_{j=1}\int^1_0\mathfrak{X}^\xi_j(t)\eta_j(t)dt+
\sum^n_{j=1}\int^1_0\mathfrak{V}^\xi_j(t)\dot\eta_j(t)dt\nonumber\\
&=&\sum^n_{j=1}\int^1_0[\mathfrak{X}^\xi_j(t)-m^2\mathfrak{G}^{\xi}_j(t)]\eta_j(t)dt+
\langle\mathfrak{G}^{\xi},\eta\rangle_{1,2,m}.
\end{eqnarray}
Since $\hat{H}_{\gamma_0^m}\ni\eta\mapsto
\sum^n_{j=1}\int^1_0[\mathfrak{X}^\xi_j(t)-m^2\mathfrak{G}^{\xi}_j(t)]\eta_j(t)dt$
is a bounded linear functional there exists a unique
$\mathfrak{F}^\xi=(\mathfrak{F}^\xi_1,\cdots,\mathfrak{F}^\xi_n)\in
\hat{H}_{\gamma_0^m}$ such that
\begin{eqnarray}\label{e:B.11}
\sum^n_{j=1}\int^1_0[\mathfrak{X}^\xi_j(t)-m^2\mathfrak{G}^{\xi}_j(t)]\eta_j(t)dt=
\langle\mathfrak{F}^\xi,\eta\rangle_{1,2,m}\quad\forall\eta\in
\hat{H}_{\gamma_0^m}.
\end{eqnarray}
This and (\ref{e:B.10}) lead to
\begin{equation}\label{e:B.12}
\nabla^m{\cal
L}_{\gamma_0^m}(\xi)=\mathfrak{G}^{\xi}+\mathfrak{F}^\xi.
\end{equation}
  Moreover, (\ref{e:B.11}) also implies
\begin{eqnarray}\label{e:B.13}
\int^1_0[\mathfrak{X}^{\xi}_j(t)-m^2\mathfrak{G}^{\xi}_j(t)]\eta_j(t)dt=
m^2\int^1_0\mathfrak{F}^\xi_j(t)\eta_j(t)dt+
\int^1_0\dot{\mathfrak{F}}^\xi_j(t)\dot{\eta}_j(t)dt
\end{eqnarray}
for any $W^{1,2}_{loc}$ map $\eta_j:\R\to\R^n$ satisfying
$\eta_j(t+1)=\pm\eta_j(t)\;\forall t\in\R$ if $\hat{S}_j=\pm 1$.

The following lemma may be proved directly.

\begin{lemma}\label{lem:B.1}
Let $m\in\N$, $\theta\in\R$, and let  $f\in L^1_{loc}(\R, \C^n)$ be
bounded. If  $f$ satisfies: $f(t+1)=e^{i\theta}f(t)\;\forall t$ then
 the equation $x''(t)- m^2x(t)=f(t)$
has a unique  solution
$$
x(t)=-\frac{1}{2m}\int^\infty_t e^{m(t-s)}f(s)\,ds -
\frac{1}{2m}\int_{-\infty}^t e^{m(s-t)}f(s)\,ds
$$
satisfying: $x(t+1)=e^{i\theta}x(t)\;\forall t$.
\end{lemma}

 From (\ref{e:B.1})-(\ref{e:B.2}) and (\ref{e:B.5})-(\ref{e:B.9}), Lemma~\ref{lem:B.1} and (\ref{e:B.13})
 we derive that
\begin{eqnarray}\label{e:B.14}
\mathfrak{F}_j(t)&=&\frac{1}{2m}\int^\infty_t
e^{m(t-s)}[\mathfrak{X}^\xi_j(s)-m^2\mathfrak{G}^{\xi}_j(s)]\,ds\nonumber\\
&+& \frac{1}{2m}\int_{-\infty}^t
e^{m(s-t)}[\mathfrak{X}^\xi_j(s)-m^2\mathfrak{G}^{\xi}_j(s)]\,ds
\end{eqnarray}
for $j=1,\cdots,n$. This and (\ref{e:B.12}) lead to
\begin{eqnarray}\label{e:B.15}
\nabla^m{\cal
L}_{\gamma_0^m}(\xi)(t)&=&\mathfrak{G}^{\xi}(t)+\frac{1}{2m}\int^\infty_t
e^{m(t-s)}\left[\partial_x
L_{\gamma_0^m}\bigl(s,\xi(s),\dot{\xi}(s)\bigr)-m^2\mathfrak{G}^{\xi}(s)\right]\,ds
\nonumber\\
&+&\frac{1}{2m}\int_{-\infty}^t e^{m(s-t)}\left[\partial_x
L_{\gamma_0^m}\bigl(s,\xi(s),\dot{\xi}(s)\bigr)-m^2\mathfrak{G}^{\xi}(s)\right]\,ds,
\end{eqnarray}
where by (\ref{e:B.7}) and (\ref{e:B.9})
$\mathfrak{G}^{\xi}=(\mathfrak{G}^{\xi}_1,\cdots,
\mathfrak{G}^{\xi}_n)$ are given by
\begin{eqnarray}\label{e:B.16}
\mathfrak{G}^{\xi}(t)&=&\int^t_0\partial_v
L_{\gamma_0^m}\bigl(s,\xi(s),\dot{\xi}(s)\bigr)ds\nonumber\\
&-& \left(\int^{1}_0\partial_v
 L_{\gamma_0^m}\bigl(s,\xi(s),\dot{\xi}(s)\bigr)ds\right) {\rm
diag}(a_1(t), \cdots, a_n(t))
\end{eqnarray}
with $a_j(t)=\frac{2tS_j+ 2t+1-S_j}{4}$, $j=1,\cdots,n$.

\noindent{\bf Case 2}. ${\rm ord}(\hat{S}_1)=\cdots={\rm
ord}(\hat{S}_\tau)=2$ and so $n$ is even and $2\tau=n$.

 By (\ref{e:B.3}) and (\ref{e:B.4}), for $l=1,\cdots,\tau=n/2$ we
 have
\begin{eqnarray*}
(\mathfrak{V}^\xi_{2l-1}(t+1),
\mathfrak{V}^\xi_{2l}(t+1))&=&(\mathfrak{V}^\xi_{2l-1}(t),
\mathfrak{V}^\xi_{2l}(t))S_{l}\quad\forall t\in\R,\\
(\mathfrak{X}^\xi_{2l-1}(t+1),
\mathfrak{X}^\xi_{2l}(t+1))&=&(\mathfrak{X}^\xi_{2l-1}(t),
\mathfrak{X}^\xi_{2l}(t))S_{l}\quad\forall t\in\R.
\end{eqnarray*}
Set $\mathscr{V}^\xi_l(t):=\mathfrak{V}^\xi_{2l-1}(t)+ i
\mathfrak{V}^\xi_{2l}(t)$ and
$\mathscr{X}^\xi_l(t):=\mathfrak{X}^\xi_{2l-1}(t)+
i\mathfrak{X}^\xi_{2l}(t)$ with $i=\sqrt{-1}$. Then
\begin{eqnarray} \label{e:B.17}
\mathscr{V}^\xi_l(t+1)=e^{i\theta_{l}}\mathscr{V}^\xi_l(t)\quad\forall
t\in\R,\quad
\mathscr{X}^\xi_l(t+1)=e^{i\theta_{l}}\mathscr{X}^\xi_l(t)\quad\forall
t\in\R.
\end{eqnarray}
 Since $e^{i\theta_{l}}\ne 1$ we may define
$V^\xi_{l}:\R\to\C$ by
\begin{equation}\label{e:B.18}
V^\xi_l(t)=\int^t_0\mathscr{V}^\xi_l(s)ds-\frac{1}{1-e^{i\theta_{l}}}\int^1_0\mathscr{V}^\xi_l(s)ds.
\end{equation}
It is the primitive function of the function $\mathscr{V}^\xi_l$, and
satisfies
\begin{equation}\label{e:B.19}
V^\xi_{l}(t+1)=e^{i\theta_{l}}V^\xi_{l}(t)\quad\forall
t\in\R,\;l=1,\cdots,\tau.
\end{equation}
Write  $\zeta_l(t):=\eta_{2l-1}(t)+ i \eta_{2l}(t)$,
$l=1,\cdots,\tau$. Then
\begin{eqnarray}\label{e:B.20}
d{\cal L}_{\gamma_0^m}(\xi) (\eta)
 \!\!\!&=&\sum^n_{j=1}\int^1_0\mathfrak{X}^\xi_j(t)\eta_j(t)dt+
\sum^n_{j=1}\int^1_0\mathfrak{V}^\xi_j(t)\dot\eta_j(t)dt\nonumber\\
\!\!\!&=&\sum^\tau_{l=1}{\rm
Re}\int^1_0\mathscr{X}^\xi_l(t)\overline{\zeta_l(t)}dt+
\sum^\tau_{l=1}{\rm Re}\int^1_0\mathscr{V}^\xi_l(t)\dot{\overline{\zeta_l(t)}}dt\nonumber\\
&=&\sum^\tau_{l=1}{\rm
Re}\int^1_0[\mathscr{X}^\xi_l(t)-m^2V^\xi_l(t)]\overline{\zeta_l(t)}dt+
{\rm Re}(V^{\xi},\zeta)_{1,2,m}
\end{eqnarray}
since
\begin{eqnarray}\label{e:B.21}
\int^1_0\mathscr{V}^\xi_l(t)\dot{\overline{\zeta_l(t)}}dt
=(V^{\xi}_l,\zeta_l)_{1,2,m}-m^2\int^1_0V^{\xi}_l(t)\overline{\zeta_l(t)}dt
\end{eqnarray}
for $l=1,\cdots,\tau$. Here $V^{\xi}=(V^{\xi}_1,\cdots,
V^{\xi}_\tau)$, $\zeta=(\zeta_1,\cdots,\zeta_\tau)$ and
\begin{eqnarray}
&&(V^{\xi}_l,\zeta_l)_{1,2,m}=m^2\int^1_0V^{\xi}_l(t)\overline{\zeta_l(t)}dt+
\int^1_0\dot{V}^\xi_l(t)\dot{\overline{\zeta_l(t)}}dt,\label{e:B.22}\\
&&(V^{\xi},\zeta)_{1,2,m}=\sum^\tau_{l=1}(V^{\xi}_l,\zeta_l)_{1,2,m}.\label{e:B.23}
\end{eqnarray}
Since
$\hat{H}_{\gamma_0^m}\ni\eta=(\eta_1,\cdots,\eta_n)\equiv(\eta_1+i\eta_2,\cdots,\eta_{2\tau-1}+i\eta_{2\tau})
\equiv(\zeta_1,\cdots,\zeta_\tau)\mapsto \sum^\tau_{l=1}{\rm
Re}\int^1_0[\mathscr{X}^\xi_l(t)-m^2V^\xi_l(t)]\overline{\zeta_l(t)}dt$
is a bounded linear functional there exists a unique
$\mathfrak{F}^\xi=(\mathfrak{F}^\xi_1,\cdots,\mathfrak{F}^\xi_n)\in
\hat{H}_{\gamma_0^m}$ such that for any $\eta\in
\hat{H}_{\gamma_0^m}$,
\begin{eqnarray}\label{e:B.24}
\sum^\tau_{l=1}{\rm
Re}\int^1_0[\mathscr{X}^\xi_l(t)-m^2V^\xi_l(t)]\overline{\zeta_l(t)}dt
= {\rm Re}\sum^\tau_{l=1}(\mathfrak{F}^\xi_{2l-1}+
i\mathfrak{F}^\xi_{2l},\zeta_l)_{1,2,m}.
\end{eqnarray}
This and (\ref{e:B.20}) lead to
\begin{equation}\label{e:B.25}
\nabla^m{\cal L}_{\gamma_0^m}(\xi)=\mathfrak{F}^\xi+ ({\rm
Re}V^\xi_1,{\rm Im}V^\xi_1,\cdots, {\rm Re}V^\xi_\tau,{\rm
Im}V^\xi_\tau).
\end{equation}
By Lemma~\ref{lem:B.1} and (\ref{e:B.24}) it is easily checked that
\begin{eqnarray*}
\mathfrak{F}^\xi_{2l-1}(t)+
i\mathfrak{F}^\xi_{2l}(t)&=&\frac{1}{2m}\int^\infty_t
e^{m(t-s)}[\mathscr{X}^\xi_l(s)-m^2V^\xi_l(s)]\,ds\\
& +& \frac{1}{2m}\int_{-\infty}^t
e^{m(s-t)}[\mathscr{X}^\xi_l(s)-m^2V^\xi_l(s)]\,ds
\end{eqnarray*}
for $l=1,\cdots,\tau$. From these, (\ref{e:B.25}) and (\ref{e:B.18})
we may derive
\begin{eqnarray}\label{e:B.26}
(\nabla^m{\cal
L}_{\gamma_0^m}(\xi))_{2l-1}&=&\mathfrak{F}^\xi_{2l-1}(t)+{\rm
Re}V^\xi_l\nonumber\\
&=&{\rm Re}V^\xi_l+\frac{1}{2m}\int^\infty_t
e^{m(t-s)}[\mathfrak{X}^\xi_{2l-1}(s)-m^2{\rm Re}V^\xi_l(s)]\,ds\nonumber\\
&& + \frac{1}{2m}\int_{-\infty}^t
e^{m(s-t)}[\mathfrak{X}^\xi_{2l-1}(s)-m^2{\rm
Re}V^\xi_l(s)]\,ds\\
(\nabla^m{\cal
L}_{\gamma_0^m}(\xi))_{2l}&=&\mathfrak{F}^\xi_{2l}(t)+{\rm
Im}V^\xi_l\nonumber\\
&=&{\rm Im}V^\xi_l+\frac{1}{2m}\int^\infty_t
e^{m(t-s)}[\mathfrak{X}^\xi_{2l}(s)-m^2{\rm Im}V^\xi_l(s)]\,ds\nonumber\\
&& + \frac{1}{2m}\int_{-\infty}^t
e^{m(s-t)}[\mathfrak{X}^\xi_{2l}(s)-m^2{\rm
Im}V^\xi_l(s)]\,ds\label{e:B.27}
\end{eqnarray}
with
\begin{eqnarray}\label{e:B.28}
\hspace{-10mm}{\rm
Re}V^\xi_l(t)&=&\int^t_0\mathfrak{V}^\xi_{2l-1}(s)ds-\frac{1}{2}
\int^1_0\mathfrak{V}^\xi_{2l-1}(s)ds+\frac{\sin\theta_l}{2-2\cos\theta_l}
\int^1_0\mathfrak{V}^\xi_{2l}(s)ds\\
\hspace{-10mm}{\rm
Im}V^\xi_l(t)&=&\int^t_0\mathfrak{V}^\xi_{2l}(s)ds-\frac{1}{2}
\int^1_0\mathfrak{V}^\xi_{2l}(s)ds-\frac{\sin\theta_l}{2-2\cos\theta_l}
\int^1_0\mathfrak{V}^\xi_{2l-1}(s)ds\label{e:B.29}
\end{eqnarray}
for $l=1,\cdots,\tau$. That is, we arrive at
\begin{eqnarray}\label{e:B.30}
\nabla^m{\cal
L}_{\gamma_0^m}(\xi)&=&\mathfrak{J}^\xi(t)+\frac{1}{2m}\int^\infty_t
e^{m(t-s)}[\partial_x
L_{\gamma_0^m}\bigl(s,\xi(s),\dot{\xi}(s)\bigr)-m^2\mathfrak{J}^\xi(s)]\,ds\nonumber\\
&+&\frac{1}{2m}\int^t_{-\infty} e^{m(s-t)}[\partial_x
L_{\gamma_0^m}\bigl(s,\xi(s),\dot{\xi}(s)\bigr)-m^2\mathfrak{J}^\xi(s)]\,ds,
\end{eqnarray}
where
\begin{eqnarray}\label{e:B.31}
\mathfrak{J}^\xi(t):&= &({\rm Re}V^\xi_1(t),{\rm
Im}V^\xi_1(t),\cdots, {\rm Re}V^\xi_\tau(t),{\rm Im}V^\xi_\tau(t))=
\int^t_0\partial_v
L_{\gamma_0^m}\bigl(s,\xi(s),\dot{\xi}(s)\bigr)ds\nonumber\\
&+&\int^1_0\partial_v
L_{\gamma_0^m}\bigl(s,\xi(s),\dot{\xi}(s)\bigr)ds
\left(\oplus^\tau_{l=1}\frac{\sin\theta_l}{2-2\cos\theta_l}{\scriptscriptstyle\left(\begin{array}{cc}
0& -1\\
1& 0\end{array}\right)}-\frac{1}{2}I_n\right).
\end{eqnarray}

\noindent{\bf Case 3}. ${\rm ord}(\hat{S}_1)=\cdots={\rm
ord}(\hat{S}_p)=2>{\rm ord}(\hat{S}_{p+1})=\cdots={\rm
ord}(\hat{S}_\tau)=1$ with $p<\tau$.

Since $1\le p<\tau=n-p$, combing the above two cases together we may
obtain
\begin{eqnarray}\label{e:B.32}
\nabla^m{\cal
L}_{\gamma_0^m}(\xi)&=&\mathfrak{K}^\xi(t)+\frac{1}{2m}\int^\infty_t
e^{m(t-s)}[\partial_x
L_{\gamma_0^m}\bigl(s,\xi(s),\dot{\xi}(s)\bigr)-m^2\mathfrak{K}^\xi(s)]\,ds\nonumber\\
&+&\frac{1}{2m}\int^t_{-\infty} e^{m(s-t)}[\partial_x
L_{\gamma_0^m}\bigl(s,\xi(s),\dot{\xi}(s)\bigr)-m^2\mathfrak{K}^\xi(s)]\,ds,
\end{eqnarray}
where
\begin{eqnarray}\label{e:B.33}
\mathfrak{K}^\xi(t)= \int^t_0\partial_v
L_{\gamma_0^m}\bigl(s,\xi(s),\dot{\xi}(s)\bigr)ds+
\int^1_0\partial_v L_{\gamma_0^m}\bigl(s,\xi(s),\dot{\xi}(s)\bigr)ds\times &&\nonumber\\
\times
\left(\oplus^p_{l=1}\frac{\sin\theta_l}{2-2\cos\theta_l}{\scriptscriptstyle\left(\begin{array}{cc}
0& -1\\
1& 0\end{array}\right)}-\frac{1}{2}I_{2p}\right)\oplus{\rm
diag}(a_{p+1}(t), \cdots, a_\tau(t))&&
\end{eqnarray}
with $a_j(t)=\frac{2t\hat{S}_j+ 2t+1-\hat{S}_j}{4}$,
$j=p+1,\cdots,\tau=n-p$.

These two formulas, (\ref{e:B.15})-(\ref{e:B.16}) and
(\ref{e:B.30})-(\ref{e:B.31}) can be written as a united way as done
for $m=1$ below.

\begin{corollary}\label{cor:B.2}
Let $E_{\gamma_0}={\rm diag}(S_1,\cdots,S_{\sigma})\in\R^{n\times
n}$ with ${\rm ord}(S_1)\ge\cdots\ge{\rm ord}(S_{\sigma})$, where
each ${S}_j$ is either $1$, or $-1$, or
${\scriptscriptstyle\left(\begin{array}{cc}
\cos\theta_j& \sin\theta_j\\
-\sin\theta_j& \cos\theta_j\end{array}\right)}$, $0<\theta_j<\pi$.
Then the gradient of ${\cal L}_{{\gamma_0}}$ on
$H_{{\gamma_0}}=\hat{H}_{{\gamma_0}}$ is given by
\begin{eqnarray}\label{e:B.34}
&&\nabla{\cal L}_{{\gamma_0}}(\xi)(t)=\nabla^1{\cal
L}_{{\gamma_0}}(\xi)(t)\nonumber\\
 &=&\frac{1}{2}\int^\infty_{t}
e^{t-s}\left[\partial_x
L_{{\gamma_0}}\bigl(s,\xi(s),\dot{\xi}(s)\bigr)-\mathfrak{R}^{\xi}(s)\right]\,ds
\nonumber\\
&+&\frac{1}{2}\int_{-\infty}^{t} e^{s-t}\left[\partial_x
L_{{\gamma_0}}\bigl(s,\xi(s),\dot{\xi}(s)\bigr)-\mathfrak{R}^{\xi}(s)\right]\,ds+\mathfrak{R}^{\xi}(t),
\end{eqnarray}
where as $2={\rm ord}(S_p)>{\rm ord}(S_{p+1})$ for some
$p\in\{0,\cdots,\sigma\}$,
\begin{eqnarray}\label{e:B.35}
\mathfrak{R}^\xi(t)= \int^t_0\partial_v
L_{\gamma_0}\bigl(s,\xi(s),\dot{\xi}(s)\bigr)ds+
\int^1_0\partial_v L_{\gamma_0}\bigl(s,\xi(s),\dot{\xi}(s)\bigr)ds\times &&\nonumber\\
\times \left(\oplus_{l\le
p}\frac{\sin\theta_l}{2-2\cos\theta_l}{\scriptscriptstyle\left(\begin{array}{cc}
0& -1\\
1& 0\end{array}\right)}-\frac{1}{2}I_{2p}\right)\oplus{\rm
diag}(a_{p+1}(t), \cdots, a_\sigma(t))&&
\end{eqnarray}
with $a_j(t)=\frac{2t{S}_j+ 2t+1-{S}_j}{4}$,
$j=p+1,\cdots,\sigma=n-p$. (As usual $p=0$ or $p=\sigma$ means there
is no  first or second term in (\ref{e:B.35}).)
\end{corollary}


\begin{thebibliography}{L3}

\bibitem{AbSc1} A. Abbondandolo, M. Schwarz,
 A smooth pseudo-gradient for the Lagrangian action
functional.  {\it Advanced Nonlinear Studies}, {\bf 9}(2009),
597-623, arXiv:0812.4364v1[math.DS], 2008.



\bibitem{BKl} V. Bangert, W. Klingenberg, Homology generated by iterated closed geodesics.
{\it Topology}, {\bf 22}(1983), 379-388.

\bibitem{BLo} V. Bangert, Y. Long, The existence of two closed geodesics
on every Finsler 2-sphere. {\it Math. Ann.}, {\bf 346}(2010),
335-366.

\bibitem{BaChSh} D. Bao, S. S. Chern, Z. Shen, An Introduction to
Riemann-Finsler Geometry. Springer, Berlin (2000).



\bibitem{Ber} M.S. Berger,  Nonlinearity and Functional Analysis.
Acad. Press, 1977.



\bibitem{Bo} R. Bott, Nondegenerate critical manifolds, {\it Ann.
Math.}, {\bf 60}(1954), 248-261.


\bibitem{Br} G.E. Bredon, Introduction to Compact Transformation Groups.
Acad. Press, 1972.

\bibitem{BiFrPi} L. Biliotti, F. Mercuri, P. Piccione, On a Gromoll-Meyer type
theorem in globally hyperbolic stationary spacetimes, {\it Comm. Anal. Geom.},
{\bf 16}(2008), 333-393.


\bibitem{Bre} H. Brezis,  Functional Analysis, Sobolev Spaces and
Partial Differential Equations. Springer 2011.





\bibitem{Ca} E. Caponio, The index of a geodesic in a randers space
and some remarks about the lack of regularity of the energy
functional of a Finsler metric. {\it Acta Mathematica Academiae
Paedagogicae Ny\'iregyh\'aziensis}, {\bf 26} (2010), 265-274.


\bibitem{CaJaMa1} E. Caponio, M.A. Javaloyes, A. Masiello, Morse theory of causal geodesics in a stationary spacetime
via Morse theory of geodesics of a Finsler metric. {\it Ann. Inst.
H. Poincar\'e Anal. Non Lin\'eaire}, {\bf 27}(2010), 857-876.


\bibitem{CaJaMa2} E. Caponio, M.A. Javaloyes, A. Masiello,
  Addendum to ``Morse theory of causal geodesics in a
stationary spacetime via Morse theory of geodesics of a Finsler
metric" [Ann. Inst. H. Poincar\'e Anal. Non Lin\'eaire, {\bf 27}(2010),
857--876]. {\it Ann. Inst. H. Poincar\'e Anal. Non Lin\'eaire}, {\bf 30}(2013),
961-968, arXiv:1211.3071.

\bibitem{CaJaMa3} E. Caponio, M.A. Javaloyes, A. Masiello, On the energy functional on Finsler manifolds and
applications to stationary spacetimes. {\it Math.Ann.}, {\bf
351}(2011), no.2, 365-392.

\bibitem{Ch83} K. C. Chang, A variant mountain pass lemma, {\bf
26}(1983), {\it Sci. Sinica Ser. A}, 1241-1255.

\bibitem{Ch05} K. C. Chang, {\it Methods in Nonlinear Analysis.}
Springer Monogaphs in Mathematics, Springer 2005.


\bibitem{Ch93} K. C. Chang, {\it Infinite Dimensional Morse Theory and Multiple Solution Problem.}
Birkh\"{a}user, 1993.

\bibitem{Ch94} K. C. Chang, $H^1$ versus $C^1$ isolated critical points, {\it C. R. Acad. Sci.
Paris S\'er. I Math.}, {\bf 319}(1994) 441-446.

\bibitem{ChGh} K. C. Chang and H. Ghoussoub, The Conley index
and the critical groups via an extension of Gromoll-Meyer theory.
{\it Topol. Methods in Nonlinear Analysis}, {\bf 7}(1996), 77-93.




\bibitem{CiDe} S. Cingolani and M. Degiovanni, On the Poincaré-Hopf theorem for
functionals defined on Banach spaces. {\it Adv. Nonlinear Stud.}
{\bf 9}(2009), no. 4, 679--699.

\bibitem{Cor} J. N. Corvellec, Morse theory for continuous
functionals, {\it J. Math. Anal. Appl}., {\bf 195}(1995), 1050-1072.

\bibitem{CorHa} J. N. Corvellec and A. Hantoute, Homotopical
stablity of isolated critical points of continuous functionals. {\it
Set-Valued Anal.}, {\bf 10}(2002), 143-164.

\bibitem{DeSo} A. A. De Moura, F. M. De Souza, A Morse lemma for degenerate critical points
with low differentiability. {\it Abstract and Applied Analysis}, {\bf
5}(2000), 113-118.

\bibitem{Di} T.T. Dieck, {\it Algebraic Topology}, European
Mathematical Society, Corrected 2nd printing 2010.

\bibitem{Fa} A. Fathi, {\it Weak Kam Theorem in Lagrangian Dynamics},
Pisa, Version 16 February, 2005.




\bibitem{GM1} D. Gromoll, W. Meyer, On differentiable functions with isolated
critical points, \emph{Topology}, \textbf{8} (1969), 361--369.

\bibitem{GM2} D. Gromoll and W. Meyer,  Periodic geodesics on compact
          Riemannian manifolds. {\it J. Diff. Geom}. 3. (1969). 493-510.

\bibitem{GroTa78} K. Grove; M, Tanaka,  On the number of invariant closed
geodesics. {\it Acta Math.} {\bf 140 }(1978), no. 1-2, 33-48.


\bibitem{Ji} M. Jiang, A generalization of Morse lemma and its applications. {\it Nonlinear Analysis},
{\bf 36}(1999), 943-960.

\bibitem{Kl} W. Klingenberg, {\it Lectures on Closed Geodesics},  Springer-Verlag
 Berlin, 1978.

\bibitem{KoKrVa} L.Kozma, A. Krist\'aly, C.Varga, Critical point
theorems on Finsler manifolds, {\it Beitr\"age zur Algebra und}
({\it Geometrie Contributions to Algebra and Geometry}), {\bf
45}(2004), no.1, 47-59.

\bibitem{Lo1} Y. Long, Multiple periodic points of the Poincar\'e map of Lagrangian
systems on tori. {\it Math. Z}. {\bf 233}(2000), 443-470.




\bibitem{LoLu} Y. Long and G. Lu,
 Infinitely many periodic solution
orbits of autonomous Lagrangian systems on tori. {\it J. Funct.
Anal.} {\bf 197}(2003), no. 2, 301--322.

\bibitem{Lu} G. Lu, Corrigendum: The Conley conjecture for Hamiltonian systems on the
cotangent bundle and its analogue for Lagrangian systems,
arXiv:0909.0609v1 [math.SG], 2009.



\bibitem{Lu0} G. Lu, The Conley conjecture for
Hamiltonian systems  on the cotangent bundle   and its analogue for
Lagrangian systems. {\it J. Funct. Anal.} {\bf 259}(2009),
2967--3034. arXiv:0806.0425v2 [math.SG], 2009.



\bibitem{Lu1} G. Lu, Corrigendum: The Conley conjecture for Hamiltonian
systems on the cotangent bundle and its analogue for Lagrangian
systems. {\it J. Funct. Anal.} {\bf 261}(2011), 542-589.


\bibitem{Lu3} G. Lu, The splitting lemmas for nonsmooth functionals
on Hilbert spaces I. {\it Discrete and continuous dynamical
systems}, {\bf 33}(2013), no.7, 2939-2990.
 arXiv:1211.2127.



\bibitem{Lu2} G. Lu, The splitting lemmas for nonsmooth functionals
on Hilbert spaces. arXiv:1102.2062v1.

\bibitem{Lu4} G. Lu, Methods of infinite dimensional Morse theory
for geodesics on Finsler manifolds, arXiv:1212.2078v2,v4.

\bibitem{LuW} G. Lu, Y.T. Wang, Stability of $G$-critical
groups, {\it submitted}.

\bibitem{Ma} H.H. Matthias, Zwei Verallgemeinerungen eines Satzes von Gromoll und Meyer,
Bonn Mathematical Publications, Universit\"at Bonn Mathematisches
Institut, 1980.


\bibitem{MW} J. Mawhin, M. Willem, {\it Critical point theory and
Hamiltonian systems}, Applied Mathematical Science, Springer-Verlag,
New York, 1989.




\bibitem{McSa} D. McDuff, D. Salamon, {\it $J$-holomorphic curves and
symplectic topology}, American Mathematical Society, Providence,
RI,2004.


\bibitem{Me} F. Mercuri, The critical points theory for the closed
geodesics problem, {\it Math. Z.}, {\bf 156}(1977), 231-245.


\bibitem{MePa} F. Mercuri, G. Palmieri, Morse theory with
low differentiablity.  {\it Boll. Unione Mat. Ital. B(7)}, {\bf 1}(1987), no.3, 621-631.


\bibitem{Mo} M. Morse, {\it The calculus of variations in the lagre},
Amer. Math. Soc. Colloqu. Publ. No 18, American Mathematical
Society, Providence, 1934.

\bibitem{Pa} R.S. Palais, Homotopy theory of infinite dimensional
manifolds. {\it Topology}, {\bf 5}(1966), 1-66.

\bibitem{Pa1} R. S. Palais, The principle of symmetric criticality, {\it Commun.
Math.Phys}., {\bf 69}(1979), 19-30.



\bibitem{Ra} H. B. Rademacher, Morse Theorie und geschlossene
Geod\"atische. {\it Bonner Math. Schriften} {\bf 229}(1992).




\bibitem{Shen} Z. Shen, Lectures on Finsler Geometry.
      World Scientific Publishing Co., New Jersey (2001).


\bibitem{Spa} E.H. Spanier, {\it Algebraic Topology},
Springer-Verlag 1966.


\bibitem{Tan82} M. Tanaka, On the existence of infinitely many isometry-invariant geodesics
 {\it J. Differential Geom.}  {\bf 17}(1982), no.2, 171-184.





\bibitem{Wa} Z. Q. Wang, Equivariant Morse theory for isolated critical
                 orbits and its applications to nonlinear problems.
                {\it Lect. Notes in Math.} No. 1306, Springer, (1988) 202-221.

\bibitem{Whi} G.W. Whitehead, {\it Elements of homotopy theory},
Graduate Texts in Mathematics 61, Springer 1978.


\end{thebibliography}
\end{document}